\documentclass[a4paper,11pt]{book}
\usepackage[a4paper]{geometry}
\usepackage[utf8]{inputenc}
\usepackage[T1]{fontenc}
\usepackage{tocloft}

\usepackage{lipsum}
\usepackage[english,francais]{babel}
\usepackage{amssymb,amsmath,amsthm,amscd}
\usepackage{mathrsfs}
\usepackage{mathabx}
\usepackage{subfig}
\usepackage{wrapfig}

\usepackage{tabularx} 
\usepackage{calc} 
\usepackage{graphicx} 
\usepackage{color} 

\usepackage[nottoc]{tocbibind}
\usepackage{longtable}
\usepackage[normalem]{ulem}

\usepackage{enumitem}

\usepackage{hyperref}

\usepackage{fancyhdr}
\makeatletter
\let\ps@plain=\ps@empty
\makeatother

\definecolor{grey}{gray}{0.4}
\definecolor{light-grey}{gray}{0.7}

\let\origdoublepage\cleardoublepage
\newcommand{\clearemptydoublepage}{%
  \clearpage
  {\pagestyle{empty}\origdoublepage}%
}
\let\cleardoublepage\clearemptydoublepage

\FrenchFootnotes 
\AddThinSpaceBeforeFootnotes 

\newcommand*\chapterstar[1]{%
  \chapter*{#1}%
  \addcontentsline{toc}{chapter}{#1}%
  \markboth{#1}{#1}}

\newenvironment{dedicace}{%
  \newpage\thispagestyle{empty}
  \hfill\begin{minipage}{100mm}\begin{flushright}\it}{%
  \end{flushright}\end{minipage}\vfill}

\newcommand{\executeiffilenewer}[3]{%
\ifnum\pdfstrcmp{\pdffilemoddate{#1}}%
{\pdffilemoddate{#2}}>0%
{\immediate\write18{#3}}\fi%
}

\theoremstyle{plain}
\newtheorem{theorem}{Theorem}[section]
\newtheorem{lemma}[theorem]{Lemma}
\newtheorem{corollary}[theorem]{Corollary}
\newtheorem{proposition}[theorem]{Proposition}
\newtheorem{fact}[theorem]{Fact}

\newtheorem*{theorem*}{Theorem}

\theoremstyle{definition}

\theoremstyle{remark}
\newtheorem{remark}[theorem]{Remark}

\newcommand{\tconst}[1]{\ensuremath{t_{\textnormal{\ref{#1}}}}}
\newcommand{\const}[1]{\ensuremath{c_{\textnormal{\ref{#1}}}}}
\newcommand{\Const}[1]{\ensuremath{C_{\textnormal{\ref{#1}}}}}

\def\C{\ensuremath{\mathbf{C}}}
\def\N{\ensuremath{\mathbf{N}}}
\def\Q{\ensuremath{\mathbf{Q}}}
\def\R{\ensuremath{\mathbf{R}}}
\def\Z{\ensuremath{\mathbf{Z}}}

\def\ep{\varepsilon}

\def\E{\ensuremath{\mathbf{E}}}
\def\P{\ensuremath{\mathbf{P}}}
\DeclareMathOperator{\Var}{Var}
\def\F{\ensuremath{\mathscr{F}}}

\def\Ind{\ensuremath{\mathbf{1}}}

\renewcommand\Re{\operatorname{Re}}
\renewcommand\Im{\operatorname{Im}}

\DeclareMathOperator{\supp}{supp}
\DeclareMathOperator{\med}{qu}

\def\to{\rightarrow}

\def\tand{\ensuremath{\text{ and }}}
\def\tif{\ensuremath{\text{ if }}}
\def\tas{\ensuremath{\text{ as }}}
\def\ton{\ensuremath{\text{ on }}}
\def\tor{\ensuremath{\text{ or }}}

\newcommand{\dd}{\mathrm{d}} 

\def\Bsh{\ensuremath{\Bh^*}}

\def\D{\ensuremath{\mathbb{D}}}
\def\Db{\ensuremath{\overline{\D}}}

\def\e{\ensuremath{e}}
\def\Et{\ensuremath{\widetilde{E}}}
\def\empt{\varnothing}
\def\Ft{\ensuremath{\widetilde{\F}}}

\def\G{\ensuremath{\mathcal{G}}}

\def\H{\ensuremath{\mathcal{H}}}

\def\lc{\lambda_c}
\def\lcb{\overline{\lambda_c}}

\def\Pt{\ensuremath{\widetilde{P}}}

\def\S{\ensuremath{\mathcal{S}}}
\def\T{\ensuremath{\mathcal{T}}}
\def\Tt{\ensuremath{\widetilde{\T}}}

\def\tree{\ensuremath{\mathfrak{t}}}

\def\bfr{\mathfrak{b}}
\def\sfr{\mathfrak{s}}
\def\Otilde{\widetilde{O}}

\def\G{\ensuremath{\mathscr{G}}}
\def\H{\ensuremath{\mathscr{H}}}
\def\p{\ensuremath{\mathbf{p}}}
\def\Ptilde{\widetilde{\P}}
\def\Etilde{\widetilde{\E}}

\def\Phat{\widehat{\P}}
\def\Ehat{\widehat{\E}}
\def\Qhat{\widehat{\Q}}
\def\Varhat{\widehat{\Var}}

\def\thbar{\overline{\theta}}
\def\Bfl{B$^\flat$}
\def\Bsh{B$^\sharp$}
\DeclareMathOperator{\Err}{Err}

\def\Mcal{\ensuremath{\mathcal{M}}}
\def\Ncal{\ensuremath{\mathcal{N}}}
\def\Mstar{\ensuremath{\Mcal^*}}
\def\Nstar{\ensuremath{\Ncal^*}}

\title{Mouvement brownien branchant avec s\'election}
\author{Pascal Maillard}

\begin{document}
\selectlanguage{french}


\pagestyle{empty}

\pdfbookmark[0]{Page de garde}{garde}

\thispagestyle{empty}

\begin{center}
  \begin{tabularx}{\textwidth}{m{3cm}Xm{3cm}X}
	\includegraphics[width = 2.5cm]{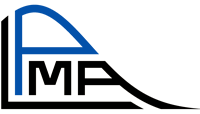}
	& \large \parbox{4cm}{Laboratoire de\\ Probabilités et\\ Modèles Aléatoires}
	& \includegraphics[width = 3cm]{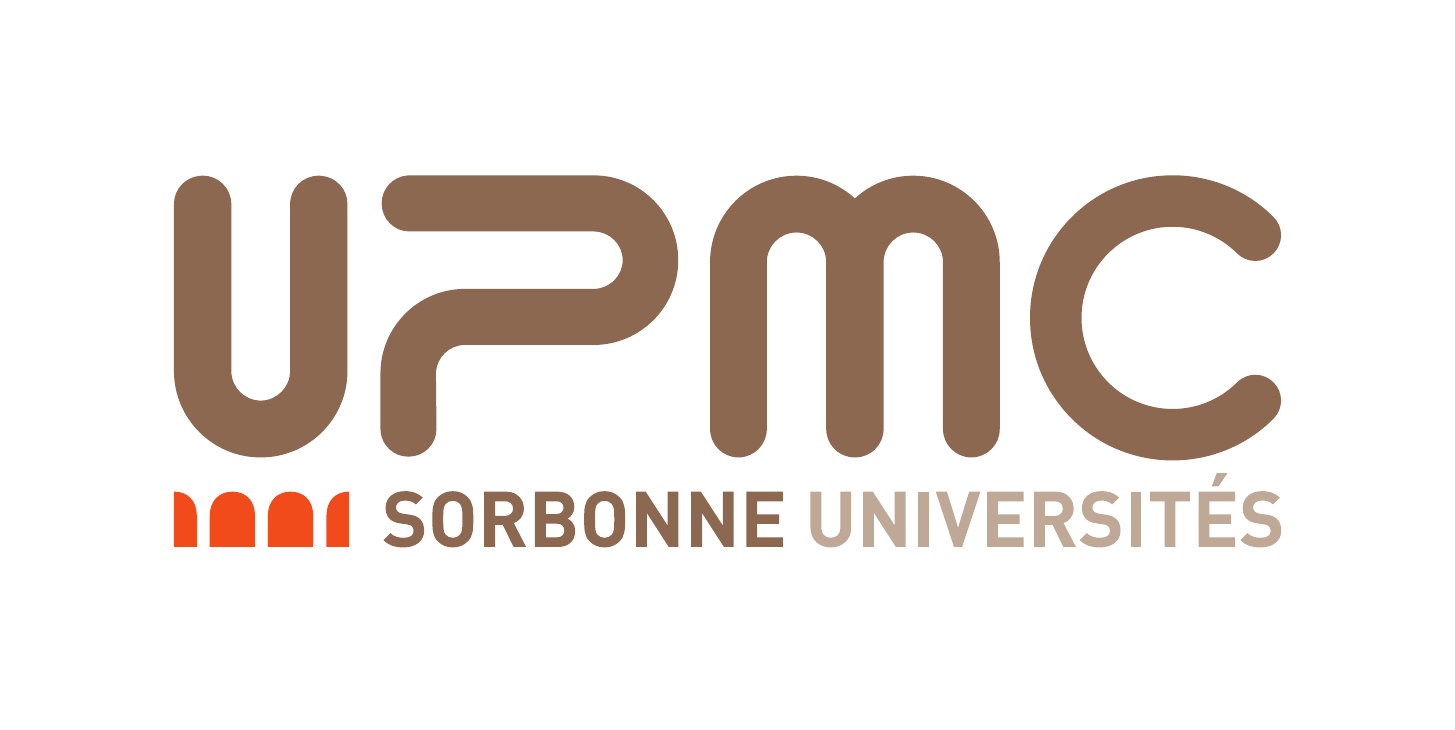}
	& \large \parbox{4cm}{Universit\'e Pierre\\ et Marie Curie}
  \end{tabularx}
\end{center}

\begin{center}
\vspace{\stretch{1}}
{\Large \textbf{École Doctorale Paris Centre}}

\vspace{\stretch{2}}

{\Huge \textsc{Thèse de doctorat}}

\vspace{\stretch{1}}

{\LARGE Discipline : Mathématiques}

\vspace{\stretch{3}}

{\large présentée par}
\vspace{\stretch{1}}

\textbf{{\LARGE Pascal \textsc{Maillard}}}

\vspace{\stretch{2}}
\hrule
\vspace{\stretch{1}}
{\LARGE \textbf{Mouvement brownien branchant avec s\'election}}

\vspace{\stretch{1}}

{\Large \textbf{(Branching Brownian motion with selection) }}
\vspace{\stretch{1}}
\hrule
\vspace{\stretch{1}}

{\Large 
dirigée par Zhan \textsc{Shi}}

\vspace{\stretch{3}}

 {\Large Rapporteurs :}

\vspace{\stretch{1}}

{\Large
\begin{tabular}{ll}
M. Andreas \textsc{Kyprianou} & University of Bath\\
M. Ofer \textsc{Zeitouni} & University of Minnesota \\
 & \& Weizmann Institute of Science
\end{tabular}
}

\vspace{\stretch{3}}

 {\Large Soutenue le 11 octobre 2012 devant le jury composé de :}

\vspace{\stretch{1}}
{\Large
\begin{tabular}{lll}
M\textsuperscript{me} Brigitte \textsc{Chauvin} & Université de Versailles &  examinatrice\\
M. Francis \textsc{Comets} & Université Paris 7 & examinateur\\
M. Bernard \textsc{Derrida} & Université Paris 6 et ENS & examinateur\\
M. Yueyun \textsc{Hu} & Université Paris 13 & examinateur\\
M. Andreas \textsc{Kyprianou} & University of Bath & rapporteur\\
M. Zhan \textsc{Shi} & Université Paris 6 & directeur\\
\end{tabular}
}

\end{center}


\newpage
\pagestyle{empty}

\vspace*{\fill}

\noindent\begin{center}
\begin{minipage}[t]{(\textwidth+2cm)/2}
Laboratoire de Probabilités et Modèles Aléatoires \\
Université Pierre et Marie Curie \\
4, place Jussieu \\
75 005 Paris
\end{minipage}
\hspace{1cm}
\begin{minipage}[t]{(\textwidth-5cm)/2}
\'Ecole Doctorale Paris Centre \\
Case courrier 188 \\ 
4, place Jussieu \\
75 252 Paris cedex 05
\end{minipage}
\end{center}

\newpage


\begin{dedicace}

  \vspace*{6cm}

	À mes parents, qui ont fait pousser la plante avec amour\\
	À Pauline, le soleil qui la fait vivre et lui est source de joie et de bonheur
\end{dedicace}


\cleardoublepage
\pdfbookmark[0]{Remerciements}{remerciements}
\chapter*{Remerciements}

\setlength{\parskip}{1em}

Je tiens à remercier chaleureusement mes deux rapporteurs Andreas Kyprianou et Ofer Zeitouni ainsi que les examinateurs Brigitte Chauvin, Francis Comets, Bernard Derrida et Yueyun Hu. J’admire les travaux de chacun et c’est un grand honneur qu’ils aient accepté d’évaluer mon travail.

Merci à mon directeur Zhan Shi qui m’a fait découvrir un sujet aussi riche et passionnant
que le mouvement brownien branchant, qui m’a toujours soutenu et prodigué ses conseils, tout en m’accordant une grande liberté dans mes recherches. Je remercie également les autres professeurs de l’année 2008/2009 du Master 2 Processus Stochastiques à l’UPMC pour la formation exceptionnelle en probabilités sans laquelle cette thèse n’aurait pas été possible. 
 
Certaines autres personnes ont eu une influence directe sur mes recherches mathématiques.
Je remercie vivement Julien Berestycki qui, par son enthousiasme infini, a su donner une
impulsion nouvelle à mes recherches. Merci à Olivier Hénard, mon ami, pour avoir partagé à
de nombreuses occasions sa fine compréhension des processus de branchement. Merci à Louigi Addario-Berry, Elie Aïdékon, Louis-Pierre Arguin, Jean Bérard, Nathanaël Berestycki, Jason Schweinsberg, Damien Simon, Olivier Zindy pour m’avoir soutenu, conseillé, éclairé.

Merci encore à Yacine~B, Pierre~B, Reda~C, Leif~D, Xan~D, Nic~F, Clément~F, Patrick~H, Mathieu~J, Cyril~L, Thomas~M, Bastien~M et Matthew~R pour ce ``aahh'' de compréhension survenu après une discussion spontanée ou après des heures de travail en commun. Merci aux anciens et actuels doctorants du LPMA pour un quotidien stimulant, agréable et drôle et pour avoir supporté mes mauvaises blagues au déjeuner et mes questions idiotes au GTT (ou l’inverse). Merci également à l’équipe administrative du LPMA pour leur travail efficace et leur dévouement.
 
Merci, Danke, Hvala, Thank you à mes amis, mathématiciens ou non, de Paris ou d’ailleurs, que j’ai eu la chance de rencontrer sur mon chemin. Merci aux superbes femmes de Pierrepont et à JP. Danke an meine Familie aus dem Norden.

Et finalement, les plus grands remerciements vont à Gérard, Elsbeth, Tobi, Karin et André
pour leur soutien inconditionnel, et surtout à Pauline, la seule personne qui sache aussi bien
me faire tourner la tête qu’y mettre de l’ordre quand elle tourne en rond.

\setlength{\parskip}{0em}

\cleardoublepage
\pdfbookmark[0]{Résumé/Abstract}{resume}
    

\begin{otherlanguage}{french}

  \vspace{1cm}

  \begin{center} \Large \bf Mouvement brownien branchant avec s\'election \end{center}

\section*{Résumé}

Dans cette thèse, le mouvement brownien branchant (MBB) est un système aléatoire de particules, où celles-ci diffusent sur la droite réelle selon des mouvements browniens et branchent à taux constant en un nombre aléatoire de particules d'espérance supérieure à~1. Nous étudions deux modèles de MBB avec sélection : le MBB avec absorption à une droite espace-temps et le $N$-MBB, où, dès que le nombre de particules dépasse un nombre donné $N$, seules les $N$ particules les plus à droite sont gardées tandis que les autres sont enlevées du système. Pour le premier modèle, nous étudions la loi du nombre de particules absorbées dans le cas où le processus s'éteint presque sûrement, en utilisant un lien entre les équations de Fisher--Kolmogorov--Petrovskii--Piskounov (FKPP) et de Briot--Bouquet. Pour le deuxième modèle, dont l'étude représente la plus grande partie de cette thèse, nous donnons des asymptotiques précises sur la position du nuage de particules quand $N$ est grand. 
Plus précisément, nous montrons qu'elle converge à l'échelle de temps $\log^3 N$ vers un processus de Lévy plus une dérive linéaire, tous les deux explicites, confirmant des prévisions de Brunet, Derrida, Mueller et Munier. Cette étude contribue à la compréhension de fronts du type FKPP sous l'influence de bruit. Enfin, une troisième partie montre le lien qui existe entre le MBB et des processus ponctuels stables.

\subsubsection*{Mots-clefs}

Mouvement brownien branchant, sélection, équation de Fisher--Kolmogorov--Petrovskii--Piskounov (FKPP) bruitée, équation de Briot--Bouquet, mesure aléatoire stable.

\end{otherlanguage}

\begin{otherlanguage}{english}

  \vspace{1cm}

  \begin{center} \rule{\textwidth/3}{1pt} \end{center}

  \vspace{1cm}

  \begin{center} \Large \bf Branching Brownian motion with selection \end{center}
  \section*{Abstract}
  
  In this thesis, branching Brownian motion (BBM) is a random particle system where the particles diffuse on the real line according to Brownian motions and branch at constant rate into a random number of particles with expectation greater than 1. We study two models of BBM with selection: BBM with absorption at a space-time line and the $N$-BBM, where, as soon as the number of particles exceeds a given number $N$, only the $N$ right-most particles are kept, the others being removed from the system. For the first model, we study the law of the number of absorbed particles in the case where the process gets extinct almost surely, using a relation between the Fisher--Kolmogorov--Petrovskii--Piskounov (FKPP) and the Briot--Bouquet equations. For the second model, the study of which represents the biggest part of the thesis, we give a precise asymptotic on the position of the cloud of particles when $N$ is large. 
More precisely, we show that it converges at the timescale $\log^3 N$ to a L\'evy process plus a linear drift, both of them explicit, which confirms a prediction by Brunet, Derrida, Mueller and Munier. This study contributes to the understanding of travelling waves of FKPP type under the influence of noise. Finally, in a third part we point at the relation between the BBM and stable point processes.

  \subsection*{Keywords} Branching Brownian motion, selection, Fisher--Kolmogorov--Petrovskii--Piskounov (FKPP) equation with noise, Briot--Bouquet equation, stable random measure.

\end{otherlanguage}   

\newpage

\selectlanguage{english}

\setcounter{tocdepth}{2} 
\pdfbookmark[0]{Contents}{tablematieres} 

\tableofcontents  



\newpage
\pagestyle{fancy}

\chapterstar{Introduction}

The ancestor of all branching processes is the \emph{Galton--Watson process}\footnote{See \cite{Kendall1966} for an entertaining historical overview. Note that it should be called the \emph{Bienaym\'e--Galton--Watson process}, but we will stick to the standard name.}. A Galton--Watson process $(Z_n)_{n\ge0}$ with reproduction law $q(k)_{k\ge0}$ is defined in the following way: $Z_0 = 1$ and $Z_{n+1} = Z_{n,1} + \cdots + Z_{n,Z_n}$, where the $Z_{n,i}$ are independent and identically distributed (iid) according to $q(k)_{k\ge0}$. The generating function $f(s) = \E[s^{Z_1}]$ plays an important role in the study of $(Z_n)_{n\ge0}$, because of the relation $\E[s^{Z_n}] = f^{(n)}(s) = f\circ\cdots\circ f(s)$, the $n$-fold composition of $f$ with itself. With this basic fact, one can see for example without difficulty that the probability of \emph{extinction} (i.e.\ the probability that $Z_n=0$ for some $n$) equals the smallest fixed point of $f$ in $[0,1]$. In particular, the extinction probability is one if and 
only if $\
E[L] = f'(1)$ is less than or equal to one. This motivates the classification of branching processes into supercritical, critical and subcritical, according to whether $\E[Z_1]$ is larger than, equal to or less than one. Furthermore, generating function techniques have been used extensively in the 1960's and 1970's in order to derive several limit theorems, one of the most famous being the \emph{Kesten--Stigum theorem}, which says that in the supercritical case, the martingale $Z_n/\E[Z_n]$ converges to a non-degenerate limit if and only if $\E[Z_1\log Z_1] < \infty$. 

The Galton--Watson process has two natural variants: First of all, one can define a branching process $(Z_t)_{t\ge 0}$ in continuous time, where each individual branches at rate $\beta>0$ into a random number of individuals, distributed according to $q(k)$. The generating function of $Z_t$ then satisfies two differential equations called Kolmogorov's forward and backward equations (see \eqref{eq_fwd} and \eqref{eq_bwd} in Chapter~\ref{ch:barrier}). In passing to continuous time one loses generality, because the discrete skeleton $(Z_{an})_{n\ge0}$ for any $a>0$ is a Galton--Watson process in the above sense. In fact, the question under which conditions a (discrete-time) Galton--Watson process can be embedded into a continuous-time process has been investigated in the literature (see \cite[Section~III.12]{Athreya1972}).

The second variant is to assign a \emph{type} to each individual and possibly let the reproduction of an individual depend on the type. In the simplest case, the case of a finite number of types, this yields to results which are similar to those of the single-type case. For example, the classification into supercritical, critical or subcritical processes now depends on the largest eigenvalue of a certain matrix and even the Kesten--Stigum theorem has an analogue (see \cite[Chapter V]{Athreya1972}).

\paragraph{Branching Brownian motion and FKPP equation.}
In this thesis, we study one-dimen\-sional \emph{branching Brownian motion (BBM)}, which is a fundamental example of a multitype branching processes in a non-compact state space, namely the real numbers\footnote{Strictly speaking, the type space of BBM is the space of continuous real-valued functions, but we ignore this fact here.}. Starting with an initial configuration of particles\footnote{We will often, but not always, use the terms ``particle'' and ``individual'' interchangeably.} located at the positions $x_1,\ldots,x_n\in\R$, the particles independently diffuse according to Brownian motions and branch at rate one into a random number of particles distributed according to the law $q(k)$. Starting at the position of their parent, the newly created particles then repeat this process independently of each other. It is the continuous counterpart of the \emph{branching random walk (BRW)}, a discrete-time multitype branching process where the offspring distribution of an individual at the position $x$ is 
given by a point process $\mathscr P$ translated by $x$. As in the single-type case, the BRW is a more general object than the BBM, because the discrete skeleton of a branching Brownian motion is itself a BRW. 

In losing generality, one gains in explicitness: When studying BBM one often has more tools at hand than for the BRW, because of the explicit calculations that are possible due to the Brownian motion. For example, let $u(x,t)$ denote the probability that there exists a particle to the right\footnote{For $a,b\in\R$, we say that $a$ is to the right of $b$ if $a > b$.} of $x$ at time $t$ in BBM started from a single particle at the origin. The function $u$ satisfies the so-called \emph{Fisher--Kolmogorov--Petrovskii--Piskounov (FKPP) equation}
\begin{equation}
\label{eq:FKPP}
\frac \partial {\partial t} u = \frac 1 2 \left(\frac \partial {\partial x}\right)^2 u + F(u),
\end{equation}
where the \emph{forcing} term is $F(u) = \beta (1-u - f(1-u))$.

Starting with McKean~\cite{McKean1975}, this fact has been exploited many times to give precise asymptotics on the law of the position of the right-most particle when the process is supercritical, i.e. $m = \sum_k(k-1)q(k) > 0$ \cite{Bramson1978,Bramson1983,Chauvin1988,Chauvin1990}. Recently, it has also been used for the study of the whole point process formed by the right-most particles \cite{Brunet2011a,Arguin2011,Arguin2010,Arguin2011a,Aidekon2011c}. Many of these results have later been proven for the branching random walk as well, either through the study of a functional equation which takes the role of \eqref{eq:FKPP} above \cite{Dekking1991,Rouault1993,Bachmann2000,Bramson2009,Webb2011} or using more probabilistic techniques \cite{McDiarmid1995,Hu2009,Addario-Berry2009,Aidekon2011b,Madaule2011}.

Let us state these results precisely. Fisher \cite{Fisher1937} and Kolmogorov, Petrovskii, Piskounov \cite{KPP}, who introduced the equation \eqref{eq:FKPP}, already noticed that it admits \emph{travelling wave solutions}, i.e.\ solutions of the form $u(x,t) = \varphi_c(x - ct)$ for every $c\ge c_0 = \sqrt{2\beta m}$. Furthermore, in \cite{KPP} it is proved that under the initial condition $u(x,0) = \Ind_{(x\le 0)}$, there exists a centring term $m(t)$, such that $u(x-m(t),t)\to\varphi_{c_0}$ and $m(t)\sim c_0t$ as $t\to\infty$. Together with tail estimates on the travelling wave, this implies a law of large numbers for the position of the right-most particle in BBM\footnote{An equivalent result for the BRW has been proven by Biggins \cite{Biggins1976,Big1977b}, after more restrictive versions by Hammersley \cite{Hammersley1974} and Kingman \cite{Kingman1975}.}. The next order of the centring term $m(t)$ was then studied first by McKean \cite{McKean1975}, who provided the estimate $m(t) \le c_0t - 1/(2c_0) \log t$ and then by Bramson \cite{Bramson1983}, who established almost fifty years after the discoverers of \eqref{eq:FKPP} that one could choose $m(t) = c_0t - 3/(2c_0) \log t$, a result which stimulated a wealth of research.

\paragraph{Travelling waves.}

The fact that (semi-linear) parabolic differential equations could describe wave-like phenomena has aroused great interest and spurred a lot of research, which is now a central pillar of the theory of parabolic differential equations (see for example \cite{Aronson1975,Volpert1994}) and has also been discussed to a great extent in the physics literature (see \cite{VanSaarloos2003} for an exhaustive account). The FKPP equation has a central place in this theory and is considered to be a basic prototype. 

Since the beginning of the 1990's, physicists have been especially interested in the effect of noise on wave propagation. The types of noise that one considers are mainly multiplicative white noise and discretisation of the wave profile (for an exhaustive list of references on this subject, see \cite{Panja2004}). The rationale behind the latter is that real-life systems consisting of a finite number of parts are only approximately described by differential equations such as \eqref{eq:FKPP}. For example, in the original work of Fisher \cite{Fisher1937}, the function $u(x,t)$ describes the proportion of an advantageous gene among a population in a one-dimensional habitat (such as a coast-line). In a population of size $N$, it can therefore ``in reality'' only take values which are multiples of $N^{-1}$. Discretisation therefore corresponds to ``internal'' noise of the system. A multiplicative white noise on the other hand models an external noise \cite{Brunet2001}. 

During the 90's, there were several studies which noticed that such a noise had a tremendous effect on the wave speed, causing a significant slowdown of the wave (see for example \cite{Breuer1994}). This was then brilliantly analysed by Brunet and Derrida \cite{Brunet1997}, who introduced the \emph{cutoff} equation, which is obtained by multiplying the forcing term $F(u)$ in \eqref{eq:FKPP} by $\Ind_{(u\ge N^{-1})}$. They found the solutions to this equation to have a wave speed slower than the original one by a difference of the order of $\log^{-2}N$ and verified this numerically \cite{Brunet1999} for an $N$-particle model, where each particle of generation $n+1$ chooses two parents uniformly from level $n$, takes the maximum of both positions and adds a noise term. But they did not stop there: In later works, with coauthors, they studied the fluctuations of such microscopic systems of ``FKPP type'' and developed an axiomatic phenomenological theory of fluctuating FKPP fronts which permits to describe the 
fluctuations of those systems. Among them, they proposed the particle system we call the $N$-BRW \cite{Brunet2006a}: At each time step, the particles reproduce as in the BRW, but only the $N$ right-most particles are kept, the others being removed from the system. This can be seen as a kind of \emph{selection mechanism}, which has an obvious biological interpretation: If one interprets the position of an individual as the value of its ``fitness'' \cite{Brunet2007}, i.e.\ a measure of how well the individual is adapted to an environment, then killing all but the $N$ right-most particles at each step is a toy model for \emph{natural selection}, a key concept of Darwinian evolution.

\paragraph{Further applications of BRW and BBM.}
Besides their role as prototypes of travelling waves, BRW and BBM have many other applications or interpretations, mainly because of their tree structure. For example, the BRW can be seen as a directed polymer on a disordered tree \cite{Derrida1988} and more generally as an infinite-dimensional version of the Generalised Random Energy Model (GREM) \cite{Derrida1985,Bovier2004a}. It also plays an important role in the study of the Gaussian Free Field on a 2D lattice box \cite{Bolthausen2001,Daviaud2006,Bolthausen2011,Bramson2012}. On a more basic level, BRW and BBM have been used as models for the ecological spread of a population or of a mutant allele inside a population \cite{Sawyer1976,Mollison1977,Sawyer1979}, especially in the multidimensional setting, which has been studied at least since \cite{Biggins1978}. Finally, a fascinating application of the BRW appears in the proof by Benjamini and Schramm that every graph with positive Cheeger constant contains a tree with positive Cheeger constant \cite{Benjamini1997}
.
\renewcommand{\descriptionlabel}[1]%
{\hspace{\labelsep}\uline{#1}}

\section*{Results}
We now come to the results obtained in this thesis on branching Brownian motion with \emph{selection}. As before, by selection we mean the process of killing particles, which can be interpreted as the effect of natural selection on a population, but should rather be viewed in the more global framework of fronts under the effect of noise. We will concentrate on two selection mechanisms:
\begin{description}
 \item[BBM with absorption.] Here, we absorb the particles at the space-time line $y=-x+ct$, i.e.\ as soon as a particle hits this line (a.k.a.\ the \emph{barrier}), it is killed immediately. This process is studied in Chapter~\ref{ch:barrier}.
 \item[$N$-BBM.] This is the continuous version of the $N$-BRW described above: Particles evolve according to branching Brownian motion and as soon as the number of particles exceeds $N$, we kill the \emph{left-most} particles, such that only the $N$ right-most remain. This process is studied in Chapter~\ref{ch:NBBM}.
\end{description}

In the next two paragraphs we describe the results we have obtained for those two models and place them into the context of the existing literature.

\paragraph{Notational convention.}

{\em We assume from now on that the branching rate satisfies $\beta = \beta_0 = 1/(2m)$, such that $c_0 = \sqrt{2\beta m} = 1$}. On can always reduce the situation to this case by rescaling time and/or space. This choice of parameters will also be made in Chapter~\ref{ch:NBBM}. However, in Chapter~\ref{ch:barrier} we will set $\beta$ to $1$, the reason being that it is based on the article \cite{Maillard2010} which was accepted for publication before the submission of this thesis and where this choice of parameters was made.

\paragraph{BBM with absorption.}

The study of branching diffusions with absorption goes back at least to Sevast'yanov \cite{Sevastyanov1958}, who studied the case of absorption at the border of a bounded domain. Watanabe \cite{Watanabe1967} considered branching diffusions in arbitrary domains under the condition that the probability of ultimate survival is positive. Kesten \cite{Kesten1978} was the first to consider the special case of one-dimensional BBM with absorption at a linear boundary: Starting with a single particle at $x>0$, he gives a constant drift $-c$ to the particles and kills them as soon as they hit the origin. He proves that the process gets extinct almost surely if and only $c\ge 1$ and provides detailed asymptotics for the number of particles in a given interval and for the probability that the system gets extinct before the time $t$ in the critical case $c=1$. 

The work of Neveu \cite{Neveu1988} in the case $c\ge 1$, where the process gets extinct almost surely, is of utmost importance to us. He made the simple but crucial observation that if one starts with one particle at the origin and absorbs particles at $-x$, then the process $(Z_x)_{x\ge 0}$, where $Z_x$ denotes the number of particles absorbed at $-x$, is a continuous-time Galton--Watson process (note that space becomes time for the process $(Z_x)_{x\ge0}$). This fact comes from the \emph{strong branching property}, which says that the particles absorbed at $-x$ each spawn independent branching Brownian motions, a fact that has been formalised and proven by Chauvin \cite{Chauvin1991} using the concept of \emph{stopping lines}. 

In the last two decades, there has been renewed interest in BBM and BRW with absorption and among the many articles on this subject we mention the article \cite{Biggins1991}, in which the criterion for almost sure extinction is established for the BRW, the articles by Biggins and Kyprianou \cite{BK1997,Biggins2004,Kyprianou2004}, who use it to study the system without absorption, the work \cite{HHK2006} which studies ``one-sided travelling waves'' of the FKPP equation, and the articles  \cite{Aldous1998,Pemantle2009,Derrida2007,Gantert2011,Berard2011,Aidekon2011e,Berestycki2010a} which study the survival probability at or near the critical drift. 

But let us get back to Neveu \cite{Neveu1988} and to the continuous-time Galton--Watson process $(Z_x)_{x\ge0}$. Neveu observed the maybe surprising fact that in the case of critical drift $c=1$, this process does not satisfy the conditions of the Kesten--Stigum theorem and that indeed $xe^{-x}Z_x$ converges as $x\to\infty$ to a non-degenerate limit $W$. 

This fact has recently aroused interest because of David Aldous' conjecture \cite{Aldous0000} that this result was true for the BRW as well. Specifically, he conjectured that $\E[Z\log Z]=\infty$, where $Z$ denotes the number of particles in branching random walk with critical drift that cross the origin for the first time\footnote{Aldous was actually interested in the total progeny of the process, but this question is equivalent, see \cite[Lemma~2]{Aidekon2011d}.}. Moreover, Aldous conjectured that in the case $c>1$, the variable $Z$ had a power-law tail. 

Aldous formulated its conjecture after Pemantle had already provided an incomplete proof \cite{Pem1999} in the Bernoulli case, based on singularity analysis of the generating function of $Z$. A complete proof in the case $c=1$ was then given by Addario-Berry and Broutin \cite{Addario-Berry2009} for general reproduction laws satisfying a mild integrability assumption. A\"{\i}d\'ekon \cite{Aidekon2010} further refined the results in the case of $b$-ary trees by showing that there are positive constants $\rho,C_1,C_2$, such that for every $x>0$, we have \[\frac{C_1x\e^{\rho x}}{n (\log n)^2} \le P(N_x > n) \le \frac{C_2x\e^{\rho x}}{n (\log n)^2}\quad\text{for large }n.\]
In Chapter~\ref{ch:barrier}, we prove a precise refinement of this result in the case of branching Brownian motion. Specifically, we prove the following:

\begin{theorem*}
Assume that the reproduction law admits exponential moments, i.e. that the radius of convergence of the power series $\sum_{k\ge 0}q(k)s^k$ is greater than 1.
\begin{itemize}[label=$-$]
\item In the critical speed area $(c = 1)$, as $n\to\infty$,
\[\quad P(Z_x = \delta n + 1) \sim  \frac{ x \e^{x}}{\delta  n^2 (\log n)^2} \quad\text{for each $x>0$}.\]
\item In the subcritical speed area $(c > 1)$ there exists a constant $K = K(c,f) > 0$, such that, as $n\to\infty$,
\[\quad P(Z_x = \delta n + 1) \sim \frac{\e^{\lcb x} - \e^{\lc x}}{\lcb-\lc} \frac{K}{n^{d+1}} \quad\text{for each $x>0$},\]
where $\lc<\lcb$ are the two roots of the quadratic equation $\lambda^2 - 2c\lambda + 1 = 0$, $d = \lcb/\lc$ and $\delta = \gcd\{k:q(k+1)>0\}$ if $q(0)=0$ and $\delta = 1$ otherwise.
\end{itemize}
\end{theorem*}

The proof of this theorem is inspired by Pemantle's incomplete proof mentioned above, in that it determines asymptotics on the generating function of $Z_x$ near its singularity $1$, which can be exploited to give the above asymptotics. The analysis of the generating function is made possible through a link between travelling waves of the FKPP equation and a classical differential equation in the complex domain, the \emph{Briot--Bouquet equation}. We will not go further into details here but instead refer to Chapter~\ref{ch:barrier}. 

If the reproduction law does not admit exponential moments, we can nevertheless apply Tauberian theorems to obtain the following result in the critical case:

\begin{theorem*}
Let $c=1$ and assume that $\sum q(k) k \log^2 k< \infty$. Then we have as $n\to\infty$,
\[P(Z_x > n) \sim \frac{ x \e^{ x}}{n(\log n)^2} \quad\text{for each $x>0$}.\]
\end{theorem*}

We further mention that A\"{\i}d\'ekon, Hu and Zindy \cite{Aidekon2011d} have recently proven these facts for the BRW without any (complex) analytical arguments and under even better moment conditions in the case $c>1$.

\paragraph{The $N$-BBM.} We now turn to the main part of this thesis: The study of the $N$-BBM. We have already outlined before the role that it takes as a prototype of a travelling wave with noise and will now present the heuristic picture obtained by Brunet, Derrida, Mueller and Munier \cite{Brunet2006,Brunet2006a}.

\begin{figure}[ht]
\begin{center}
\includegraphics[width=4in]{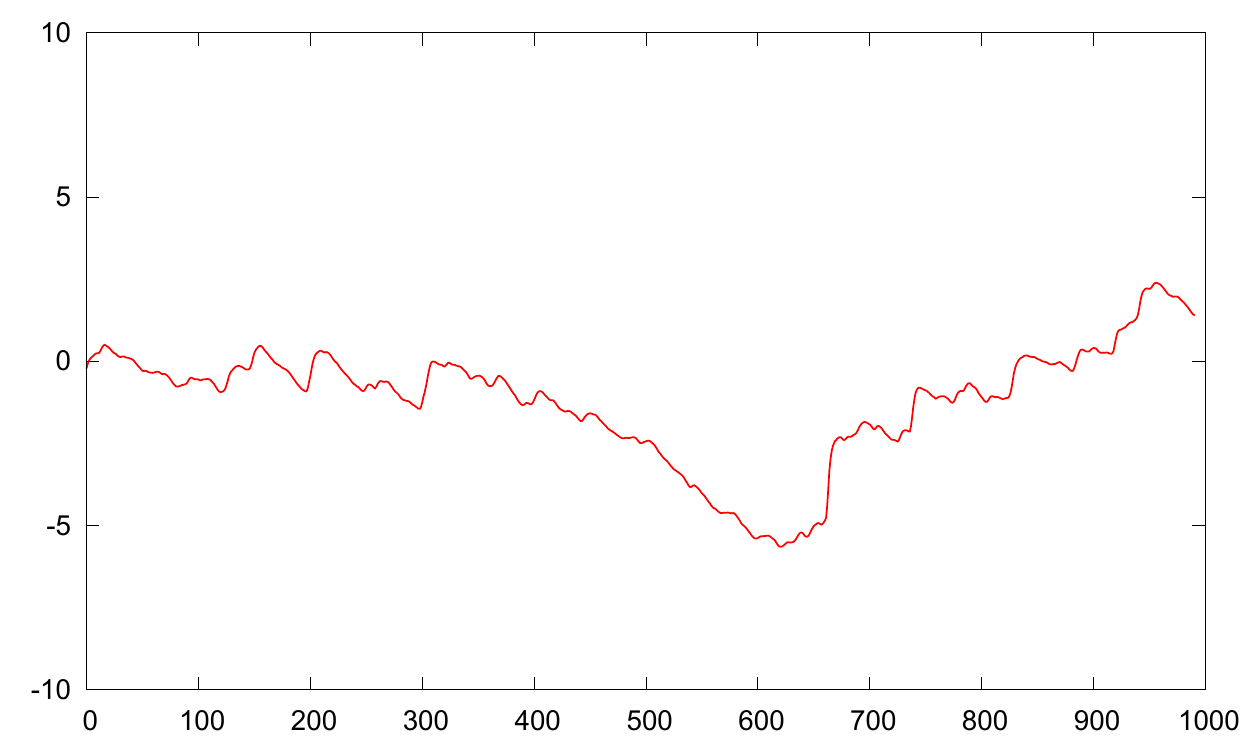} 
\end{center}
\begin{center}
\includegraphics[width=4in]{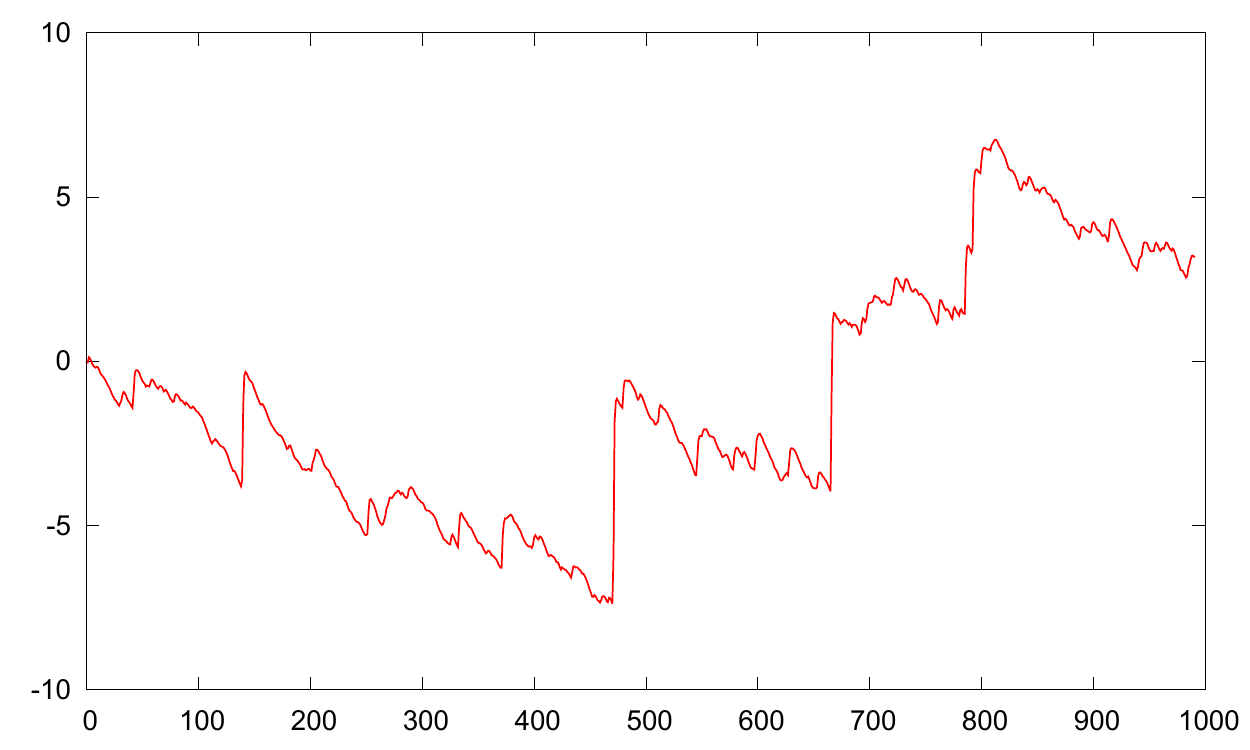} 
\end{center}
\begin{center}
\includegraphics[width=4in]{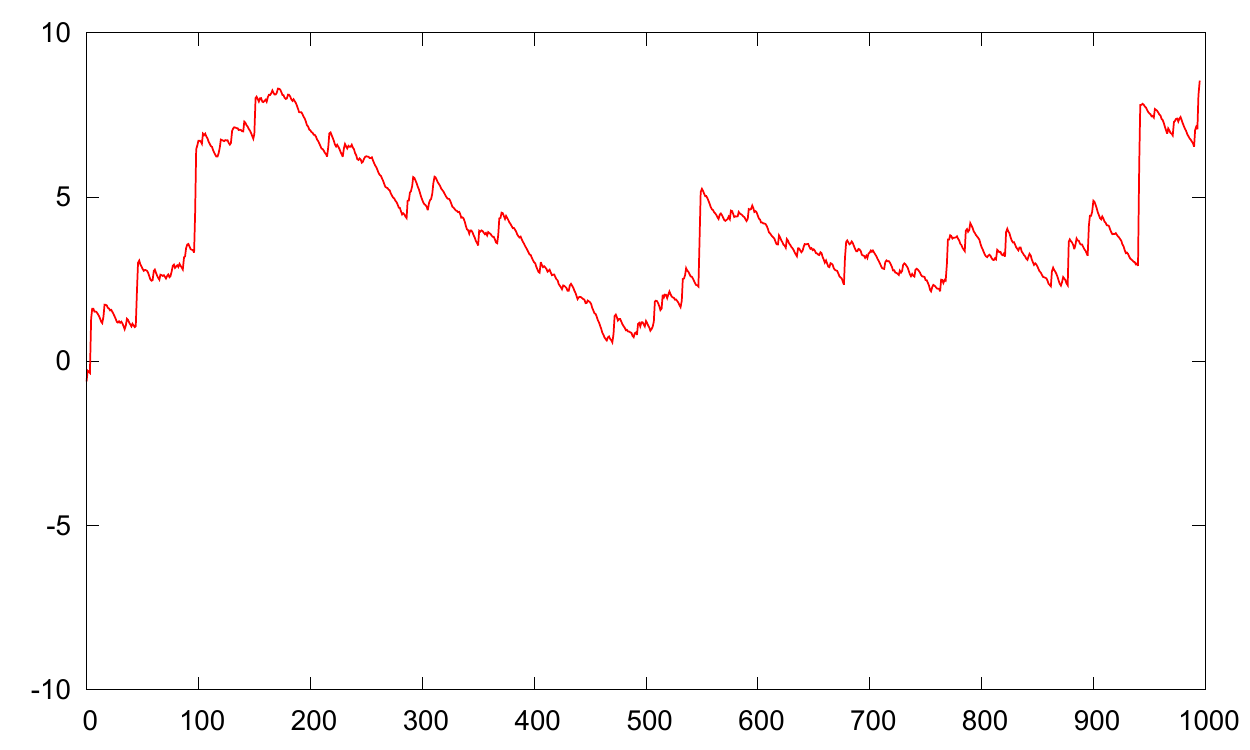} 
\end{center}
\caption{Simulation of the $N$-BRW with (from top to bottom) $N=10^{10}$, $N=10^{30}$ and $N=10^{60}$  particles (particles branch at each time step with probability $0.05$ into two and jump one step to the right with probability $0.25$). The graphs show the barycentre of the particles recentred roughly around $v_Nt$ (because of the slow convergence, small linear correction terms had to be added in order to fit the graphs into the rectangles). The horizontal axis shows time (at the scale $\log^3 N$), the vertical axis shows the position of the barycentre (without rescaling).}
\label{fig:nbbm}
\end{figure}

We recall the model: $N$ particles diffuse according to Brownian motions and branch at rate $\beta_0$ into a random number of particles given by the reproduction law $q(k)$. As soon as the number of particles exceeds $N$, only the $N$ right-most particles are kept and the others are immediately killed. This gives a cloud of particles, moving to the right with a certain speed $v_N <1$ and fluctuating around its mean. See Figure~\ref{fig:nbbm} at the end of the introduction for simulations. The authors of \cite{Brunet2006} provide detailed quantitative heuristics to describe this behaviour:

\begin{enumerate}
 \item Most of the time, the particles are in a meta-stable state. In this state, the diameter of the cloud of particles (also called the \emph{front}) is approximately $ \log N$, the empirical density of the particles proportional to $e^{- x} \sin(\pi  x/\log N)$, and the system moves at a linear speed $v_{\mathrm{cutoff}} = 1 - \pi^2/(2\log^2N)$. This is the description provided by the \emph{cutoff} approximation from \cite{Brunet1997} mentioned above.
 \item This meta-stable state is perturbed from time to time by particles moving far to the right and thus spawning a large number of descendants, causing a shift of the front to the right after a relaxation time which is of the order of $\log^2 N$. To make this precise, we fix a point in the bulk, for example the barycentre of the cloud of particles, and shift our coordinate system such that this point becomes its origin. Playing with the initial conditions of the FKPP equation with cutoff, the authors of \cite{Brunet2006} find that a particle moving up to the point $\log N + x$ causes a shift of the front by
\[
 \Delta = \log \Big(1+\frac{Ce^{x}}{\log^3 N}\Big),
\]
for some constant $C>0$. In particular, in order to have an effect on the position of the front, a particle has to reach a point near $\log N + 3 \log \log N$.
\item Assuming that such an event where a particle ``escapes'' to the point $\log N + x$ happens with rate $C e^{- x}$, one sees that the time it takes for a particle to come close to $\log N + 3 \log \log N$ (and thus causing shifts of the front) is of the order of $\log^3 N \gg \log^2 N$.
\item With this information, the full statistics of the position of the front (i.e.\ the speed $v_N$ and the cumulants of order $n\ge 2$) are found to be
\begin{equation}
 \label{eq:statistics}
\begin{split}
 v_N - v_{\mathrm{cutoff}} &\approx \pi^2 \frac{3\log\log N}{\log^3 N}\\
\frac{\text{[$n$-th cumulant]}}{t} &\approx \frac{\pi^2 n! \zeta(n)}{\log^3 N},\quad n\ge 2,
\end{split}
\end{equation}
where $\zeta$ denotes the Riemann zeta-function.
\end{enumerate}

In another paper \cite{Brunet2007}, the same authors introduce a related model, called the \emph{exponential model}. This is a BRW, where an individual at position $x$ has an infinite number of descendants distributed according to a Poisson process of intensity $e^{x-y}\,\dd y$. Again, at every step, only the $N$ right-most particles are kept. This model has some translational invariance properties which render it exactly solvable. Besides asymptotics on the speed of this system and the fluctuations, the authors of \cite{Brunet2007} then show that the genealogy of this system converges to the celebrated Bolthausen--Sznitman coalescent\footnote{The Bolthausen--Sznitman coalescent \cite{Bolthausen1998} is a process on the partitions of $\N$, in which a proportion $p$ of the blocks merge to a single one at rate $p^{-2}\,\dd p$. See \cite{Bertoin2006} and \cite{Berestycki2009} for an introduction to coalescent processes.}. This gives first insight into the fact that the genealogy of the $N$-BBM apparently 
converges to the same coalescent process, a fact that has been observed by numerical simulations \cite{Brunet2006a,Brunet2008}.

Despite (or because of) the simplicity of the $N$-BBM, it is very difficult to analyse it rigorously, because of the strong interaction between the particles, the impossibility to describe it exactly through differential equations and the fact that the shifts in the position of the system do not occur instantaneously but gradually over the fairly large timescale $\log^2 N$. For this reason, there have been few rigorous results on the $N$-BBM or the $N$-BRW: B\'erard and Gou\'er\'e \cite{Berard2010} prove the $\log^2 N$ correction of the linear speed of $N$-BRW in the binary branching case, thereby showing the validity of the approximation by a deterministic travelling wave with cutoff. Durrett and Remenik \cite{Durrett2009} study the empirical distribution of $N$-BRW and show that it converges as $N$ goes to infinity but time is fixed to a system of integro-differential equations with a moving boundary. Recently, Comets, Quastel and Ram\'{\i}rez \cite{Comets2012} studied a particle system expected to exhibit 
similar behaviour than the exponential model and show in particular that its recentred position converges to a totally asymmetric Cauchy process.

In \cite{Brunet2006}, the authors already had the idea of approximating the $N$-BBM by BBM with absorption at a linear barrier with the weakly subcritical slope $v_N$. This idea was then used with success in \cite{Berard2010} for the proof of the cutoff-correction to the speed of the $N$-BRW, relying on a result \cite{Gantert2011} about the survival probability of BRW with absorption at such a barrier and a result by Pemantle about the number of nearly optimal paths in a binary tree \cite{Pemantle2009}. In the same vein, Berestycki, Beres\-tycki and Schweinsberg \cite{Berestycki2010} studied the genealogy of BBM with absorption at this barrier and found that it converges at the $\log^3 N$ timescale to the Bolthausen--Sznitman coalescent, as predicted. In Chapter~\ref{ch:NBBM} of this thesis, we build upon their analysis in order to study the position of the $N$-BBM itself. The results that we prove are summarised in the following theorem, which confirms \eqref{eq:statistics}
 (see Theorem~\ref{th:1} of Chapter~\ref{ch:NBBM} for a complete statement).

\begin{theorem*}
Let $X_N(t)$ denote the position of the $\lceil N/2 \rceil$-th particle from the right in $N$-BBM. Then under ``good'' initial conditions, the finite-dimensional distributions of the process 
\[
\left(X_N\big(t\log^3 N\big) - v_N t\log^3 N\right)_{t\ge0}
\]
converge weakly as $N\to\infty$ to those of the L\'evy process $(L_t)_{t\ge0}$ with
\begin{equation*}
 \log E[e^{i\lambda (L_1-L_0)}] = i\lambda c + \pi^2\int_0^\infty e^{i\lambda x} - 1 - i\lambda x\Ind_{(x\le 1)}\,\Lambda(\dd x),
\end{equation*}
where $\Lambda$ is the image of the measure $(x^{-2}\Ind_{(x> 0)})\dd x$ by the map $x\mapsto \log(1+x)$ and $c\in\R$ is a constant depending only on the reproduction law $q(k)$.
\end{theorem*}

In order to prove this result, we approximate the $N$-BBM by BBM with absorption at a \emph{random} barrier instead of a linear one, a process which we call the B-BBM (``B'' stands for ``barrier''). This random barrier has the property that the number of individuals stays almost constant during a time of order $\log^3 N$. We then couple the $N$-BBM with two variants of the B-BBM, the \Bfl-BBM and the \Bsh-BBM, which in a certain sense bound the $N$-BBM from below and from above, respectively. For further details about the idea of the proof, we refer to Section~\ref{sec:NBBM_intro} in Chapter~\ref{ch:NBBM}.

\paragraph{Stable point processes occurring in branching Brownian motion.}
This paragraph describes the content of Chapter~\ref{ch:stable}, which is independent from the first two and has nothing to do with selection. It concerns the extremal particles in BBM and BRW without selection, i.e.\ the particles which are near the right-most. The study of these particles has been initiated again by Brunet and Derrida \cite{Brunet2011a}, who gave arguments for the following fact: The point process formed by the particles of BBM at time $t$, shifted to the left by $m(t) = t-3/2 \log t$, converges as $t\to\infty$ to a point process $Z$ which has the ``superposability'' property: the union of $Z$ translated by $e^\alpha$ and $Z'$ translated by $e^{\beta}$ has the same law as $Z$, where $Z'$ is an independent copy of $Z$ and $e^\alpha+e^\beta=1$. Moreover, they conjectured that this process, and possibly every process with the previous property, could be represented as a Poisson process with intensity $e^{-x}\,\dd x$, \emph{decorated} by an auxiliary point process $D$, i.e.\ each point $\xi_i$ of the Poisson process is replaced by an independent copy of $D$ translated by $\xi_i$.

In Chapter~\ref{ch:stable} we show that the ``superposability'' property has a classical interpretation in terms of \emph{stable point processes}, pointed out to us by Ilya Molchanov, and the above-mentioned representation is known in this field as the \emph{LePage series representation} of a stable point process. We furthermore give a short proof of this representation using only the theory of infinitely divisible random measures. For the BBM and BRW, the convergence of the extremal particles to such a process was proven by Arguin, Bovier, Kistler \cite{Arguin2011,Arguin2010,Arguin2011a}, A\"{\i}d\'ekon, Berestycki, Brunet, Shi \cite{Aidekon2011c} and Madaule \cite{Madaule2011}. Kabluchko \cite{Kabluchko2011} also has an interesting result for BRW started with an infinite number of particles distributed with density $e^{-cx}\,\dd x$, with $c>1$, which corresponds to travelling waves of speed larger than $1$.

\section*{Conclusion and open problems}

In this thesis, we have studied two models of branching Brownian motion with selection: the BBM with absorption at a linear barrier and the $N$-BBM. For the first model, we have given precise asymptotics on the number of absorbed particles in the case where the process gets extinct almost surely. For the second model, the study of which represents the major part of the thesis, we have shown that the recentred position of the particle system converges at the time-scale $\log^3 N$ to a L\'evy process which is given explicitly. Finally, in the last chapter, we have pointed out a relation between the extremal particles of BBM and BRW and stable point processes.

The study of the $N$-BBM constitutes an important step in the understanding of general fluctuating wave fronts, whose phenomenology is believed to be the same in many cases \cite{Brunet2006}. Furthermore, it is a natural and intricate example of a selection mechanism for BBM and BRW and an intuitive model of a population under natural selection.

We have not considered models of BBM with selection with density-dependent selection, i.e.\ where particles get killed with a rate depending on the number of particles in their neighbourhood. This kind of selection is indeed the most relevant for applications in ecology and has appeared in the literature mostly as branching-coalescing particle systems (see for example \cite{Shiga1986,Shiga1988,Ben-Avraham1998,Pechenik1999,Doering2003,Athreya2004}), but also as systems with a continuous self-regulating density \cite{Snyder2003,Hallatschek2011}. However, although we have not directly considered this type of selection mechanism, the $N$-BBM is closely related to a particular case: Suppose particles perform BBM and furthermore coalesce at rate $\ep$ when they meet (i.e.\ let two particles coalesce when their intersection local time is equal to an independent exponential variable of parameter $\ep$). Shiga \cite{Shiga1988} showed that this system is in \emph{duality} with the \emph{noisy FKPP equation}
\begin{equation}
\label{eq:noisy_FKPP}
\frac \partial {\partial t} u = \frac 1 2 \left(\frac \partial {\partial x}\right)^2 u + F(u) + \sqrt{\ep u(1-u)}\dot W, 
\end{equation}
where $\dot W$ is space-time white noise. This equation admits ``travelling wave'' solutions whose wave speed equals the speed of the right-most particle of the branching-coalescing Brownian motion (BCBM) \cite{Mueller2010}. Now, following the argumentation in \cite{Mueller2010}, the invariant measures of the BCBM are the Poisson process of intensity $\ep^{-1}$ and the configuration of no particles, the first being stable and the second unstable. Hence, if one looks at the right-most particles in the BCBM, they will ultimately form a wave-like profile with particles to the left of this wave distributed with a density $\ep^{-1}$. We have thus a similar picture than in the $N$-BBM, in that when particles get to the left of the front, they get killed quickly, in this case approximately with rate 1 instead of instantaneously. If $\ep$ is small, this does not make a difference because the typical diameter of the front goes to infinity as $\ep\to 0$. It is therefore plausible that our results about the $N$-BBM 
can be transferred to the BCBM and thus to the noisy FKPP equation \eqref{eq:noisy_FKPP}. Indeed, Mueller, Mytnik and Quastel \cite{Mueller2010} have showed that the wave speed of solutions to \eqref{eq:noisy_FKPP} is approximately $v_N$ (although they still have an error of $O(\log\log N/\log^3 N)$, where $N = \ep^{-1}$).

In what follows, we will outline several other open problems, mostly concerning the $N$-BBM.

\paragraph{Speed of the $N$-BBM.} If the reproduction law satisfies $q(0)=0$, then the \emph{speed} of the system is the constant $v_N^*$, such that $X_N(t)/t\to v_N^*$ as $t\to\infty$, where $X_N(t)$ is again the position of the $\lceil N/2 \rceil$-th particle from the right in $N$-BBM. It has been shown \cite{Berard2010} for the BRW with binary branching that this speed exists and it is not difficult to extend their proof to the BBM. Now, although our result about the position of the $N$-BBM describes quite precisely the fluctuations at the timescale $\log^3 N$ it tells us nothing about the behaviour of the system as time gets arbitrarily large. 
\emph{A priori}, it could be possible that funny things happen at larger timescales, which lead to stronger fluctuations or a different speed of the system. In their proof \cite{Berard2010} of the cutoff correction for the speed, B\'erard and Gou\'er\'e had to consider in fact timescales up to $\log^5 N$. It is therefore an open problem to prove that $v_N^*= v_N+c'+o(\log^{-3}N)$ for some constant $c'\in\R$ and that the cumulants scale as in \eqref{eq:statistics} as $t\to\infty$, which we conjecture to be true.

\paragraph{Empirical measure of the $N$-BBM.} Let $\nu^N_0(t)$ be the empirical measure of the particles at time $t$ in $N$-BBM seen from the left-most particle. Durrett and Remenik \cite{Durrett2009} show (for BRW) that the process $(N^{-1}\nu^N_0(t))_{t\in\R}$ is ergodic and therefore has an invariant probability $\pi^N$ and furthermore converges as $N\to\infty$ in law to a deterministic measure-valued process whose density with respect to Lebesgue measure solves a free boundary integro-differential equation. 
 Our result on the $N$-BBM and those prior to our work \cite{Brunet1997,Brunet2006,Berestycki2010} suggest that the $\pi^N$ should converge, as $N\to\infty$, to the Dirac measure concentrated on the measure $xe^{- x}\Ind_{x\ge 0}\,\dd x$. If one was to prove this, the methods of \cite{Durrett2009}, which are essentially based on Gronwall's inequality, would not directly apply, since they only work for timescales of at most $\log N$. On the other hand, our technique of the coupling with BBM with absorption has fairly restrictive requirements on the initial configuration of the particles. A method of proof would therefore be to show that these requirements are met after a certain time for any initial configuration and then couple the processes. 

\paragraph{The genealogy of the $N$-BBM.} As mentioned above, the genealogy of the $N$-BBM is expected to converge at the timescale $\log^3N$ to the Bolthausen--Sznitman coalescent as $N$ goes to infinity. A step towards this conjecture is the article \cite{Berestycki2010}, in which the authors show that this convergence holds for BBM with absorption at a weakly subcritical barrier; it should be easy to extend their arguments for the B-BBM defined in Section~\ref{sec:BBBM} of Chapter~\ref{ch:NBBM}. However, this would not be enough, since the monotone coupling we use to couple the $N$-BBM with the B-BBM (or rather its variants \Bfl-BBM and \Bsh-BBM) does not preserve the genealogical structure of the process. More work is therefore required to prove the conjecture for the $N$-BBM.

\paragraph{The $N$-BRW.} It should be possible to transfer the results obtained here for the $N$-BBM to the $N$-BRW, as long as the displacements have exponential right tails and the number of offspring of a particle has finite variance, say. Naturally, this would require a lot more work because of the lack of explicit expressions for the density of the particles; one can wonder whether this work would be worth it. However, it would be interesting to consider cases where one could get different behaviour than for the $N$-BBM, such as subexponentially decaying right tails of the displacements or models where particles have an infinite number of children, as in the exponential model. Another possible point of attack in the same direction would be to make rigorous for the $N$-BBM the findings in \cite{Brunet2012}. 
In that article, the authors weight the process by an exponential weight $e^{\rho R_N(t)}$, where $R_N(t)$ is the position of the right-most particle at time $t$, and obtain thus a one-parameter family of coalescent processes as genealogies.

\paragraph{Superprocesses.} If one considers $N$ independent BBMs with weakly supercritical reproduction (i.e.\ $m=a/N$ for some constant $a>0$), gives the mass $N^{-1}$ to each individual and rescales time by $N$ and space by $N^{-1}$, then one obtains in the limit a (supercritical) superprocess, such as the Dawson--Watanabe process in the case that the variance of the reproduction law is finite (see \cite{Etheridge2000} for a gentle introduction to superprocesses). One may wonder whether the results obtained in this thesis have analogues in that setting.  This is true for the results obtained in Chapter~\ref{ch:barrier}, as shown in \cite{Kyprianou2011}, in which the authors indeed transfer the results directly via the so-called \emph{backbone decomposition} of the superprocess \cite{Kyprianou2011a}. As for an analogue of the $N$-BBM: one could consider a similar model of a superprocess whose mass is kept below $N$ by stripping off some mass to the left as soon as the total mass exceeds $N$. 
It is possible that this process then exhibits similar fluctuations, which could again be related to those of the $N$-BBM by the backbone decomposition.

\numberwithin{equation}{section}

\chapter[Number of absorbed individuals in branching Brownian motion with a barrier]{The number of absorbed individuals in branching Brownian motion with a barrier}

\noindent This chapter is based on the article \cite{Maillard2010}.

\label{ch:barrier}

\section{Introduction}
We recall the definition of branching Brownian motion mentioned already in the introductory chapter: Starting with an initial individual sitting at the origin of the real line, this individual moves according to a 1-dimensional Brownian motion with drift $c$ until an independent exponentially distributed time with rate 1. At that moment it dies and produces $L$ (identical) offspring, $L$ being a random variable taking values in the non-negative integers with $P(L=1) = 0$. Starting from the position at which its parent has died, each child repeats this process, all independently of one another and of their parent. For a rigorous definition of this process, see for example \cite{INW} or \cite{Chauvin1991}.

We assume that $m = E[L]-1 \in (0,\infty)$, which means that the process is supercritical. At position $x>0$, we add an {\em absorbing barrier}, i.e.\ individuals hitting the barrier are instantly killed without producing offspring. Kesten proved \cite{Kesten1978} that this process becomes extinct almost surely if and only if the drift $c \ge c_0 = \sqrt{2 m}$ (he actually needed $E[L^2] < \infty$ for the ``only if'' part, but we are going to prove that the statement holds in general). Neveu \cite{Neveu1988} showed that the process $Z=(Z_x)_{x\ge0}$ is a continuous-time Galton--Watson process of finite expectation, but with $E[Z_x\log^+ Z_x]=\infty$ for every $x>0$, if $c=c_0$.

Let $\N = \{1,2,3,\ldots\}$ and $\N_0 = \{0\}\cup\N$. Define the infinitesimal transition rates (see \cite{Athreya1972}, p.\ 104, Equation (6) or \cite{Harris1963}, p.\ 95)\[q_n = \lim_{x\downarrow 0} \frac 1 x P(Z_x = n),\quad n \in \N_0\backslash\{1\}.\] We propose a refinement of Neveu's result:

\begin{theorem}
\label{th_tail}
Let $c=c_0$ and assume that $E[L(\log L)^{2}] < \infty$. Then we have as $n\to\infty$,
\[\sum_{k=n}^\infty q_k \sim \frac{c_0}{n (\log n)^2}\quad\tand\quad P(Z_x > n) \sim \frac{c_0 x \e^{c_0 x}}{n(\log n)^2} \quad\text{for each $x>0$}.\]
\end{theorem}

The heavy tail of $Z_x$ suggests that its generating function is amenable to singularity analysis in the sense of \cite{FO1990}. This is in fact the case in both the critical and subcritical cases if we impose a stronger condition upon the offspring distribution and leads to the next theorem.

Define $f(s)=E[s^L]$ the generating function of the offspring distribution. Denote by $\delta$ the span of $L-1$, i.e.\ the greatest positive integer, such that $L-1$ is concentrated on $\delta\Z$. Let $\lc \le \lcb$ be the two roots of the quadratic equation $\lambda^2 - 2c\lambda + c_0^2 = 0$ and denote by $d=\frac{\lcb}{\lc}$ the ratio of the two roots. Note that $c = c_0$ if and only if $\lc = \lcb$ if and only if $d=1$.

\begin{theorem}
\label{th_density}
Assume that the law of $L$ admits exponential moments, i.e. that the radius of convergence of the power series $E[s^L]$ is greater than 1.
\begin{itemize}[label=$-$]
\item In the critical speed area $(c = c_0)$, as $n\to\infty$,
\[q_{\delta n+1} \sim \frac{c_0}{\delta  n^2 (\log n)^2}\quad \tand \quad P(Z_x = \delta n + 1) \sim  \frac{c_0 x \e^{c_0x}}{\delta  n^2 (\log n)^2} \quad\text{for each $x>0$}.\]
\item In the subcritical speed area $(c > c_0)$ there exists a constant $K = K(c,f) > 0$, such that, as $n\to\infty$,
\[q_{\delta n+1} \sim \frac{K}{n^{d+1}}\quad\tand \quad P(Z_x = \delta n + 1) \sim \frac{\e^{\lcb x} - \e^{\lc x}}{\lcb-\lc} \frac{K}{n^{d+1}} \quad\text{for each $x>0$}.\]
\end{itemize}
Furthermore, $q_{\delta n + k} = P(Z_x = \delta n+ k) = 0$ for all $n\in\Z$ and $k\in\{2,\ldots,\delta\}$.
\end{theorem}

\begin{remark}
The idea of using singularity analysis for the study of $Z_x$ comes from Robin Pemantle's (unfinished) manuscript \cite{Pem1999} about branching random walks with Bernoulli reproduction.
\end{remark}
\begin{remark}
Since the coefficients of the power series $E[s^L]$ are real and non-negative, Pringsheim's theorem (see e.g.\ \cite{FS2009}, Theorem IV.6, p.\ 240) entails that the assumption in Theorem \ref{th_density} is verified if and only if $f(s)$ is analytic at 1.
\end{remark}
\begin{remark}
Let $\beta > 0$ and $\sigma > 0$. We consider a more general branching Brownian motion with branching rate given by $\beta$ and the drift and variance of the Brownian motion given by $c$ and $\sigma^2$, respectively. Call this process the $(\beta,c,\sigma)$-BBM (the reproduction is still governed by the law of $L$, which is fixed). In this terminology, the process described at the beginning of this section is the $(1,c,1)$-BBM. The $(\beta,c,\sigma)$-BBM can be obtained from $(1,c/(\sigma \sqrt{\beta}),1)$-BBM by rescaling time by a factor $\beta$ and space by a factor $\sigma/\sqrt{\beta}$. Therefore, if we add an absorbing barrier at the point $x>0$, the $(\beta,c,\sigma)$-BBM gets extinct a.s.\ if and only if $c\ge c_0 = \sigma\sqrt{2\beta m}$. Moreover, if we denote by $Z_x^{(\beta,c,\sigma)}$ the number of particles absorbed at $x$, we obtain that
\[(Z_x^{(\beta,c,\sigma)})_{x\ge0} \tand (Z_{x\sqrt{\beta}/\sigma}^{(1,c/(\sigma\sqrt{\beta}),1)})_{x\ge0} \text{ are equal in law.}\]
In particular, if we denote the infinitesimal transition rates of $(Z_x^{(\beta,c,\sigma)})_{x\ge0}$ by $q_n^{(\beta,c,\sigma)}$, for $n\in \N_0\backslash\{1\}$, then we have
\[q^{(\beta,c,\sigma)}_n = \lim_{x\downarrow0} \frac 1 x P\Big(Z_x^{(\beta,c,\sigma)} = n\Big) = \frac {\sqrt{\beta}} \sigma \lim_{x\downarrow0} \frac \sigma {x\sqrt{\beta}} P\Big(Z_{x\sqrt{\beta}/\sigma}^{(1,c/(\sigma\sqrt{\beta}),1)} = n\Big) = \frac {\sqrt{\beta}} \sigma q^{(1,c/(\sigma\sqrt{\beta}),1)}_n.\]
One therefore easily checks {\em that the statements of Theorems \ref{th_tail} and \ref{th_density} are still valid for arbitrary $\beta>0$ and $\sigma>0$, provided that one replaces the constants $c_0,\lc,\lcb,K$ by $c_0/\sigma^2$, $\lc/\sigma^2$, $\lcb/\sigma^2$, $\frac{\sqrt \beta}{\sigma} K(c/(\sigma\sqrt{\beta}),f)$, respectively.}
\end{remark}

The content of the paper is organised as follows: In Section \ref{sec_probability} we derive some preliminary results by probabilistic means. In Section \ref{sec_FKPP}, we recall a known relation between $Z_x$ and the so-called Fisher--Kolmo\-gorov--Petrovskii--Piskounov (FKPP) equation. Section \ref{sec_tail} is devoted to the proof of Theorem \ref{th_tail}, which draws on a Tauberian theorem and known asymptotics of travelling wave solutions to the FKPP equation. In Section \ref{sec_preliminaries} we review results about complex differential equations, singularity analysis of generating functions and continuous-time Galton--Watson processes. Those are needed for the proof of Theorem \ref{th_density}, which is done in Section \ref{sec_proofdensity}.

\section{First results by probabilistic methods}
\label{sec_probability}
The goal of this section is to prove
\begin{proposition}
\label{prop_heavytail}
Assume $c > c_0$ and $E[L^2] < \infty$. There exists a constant $C = C(x,c,L) > 0$, such that
\[P(Z_x > n) \ge \frac C {n^d}\quad \text{ for large }n.\]
\end{proposition}
This result is needed to assure that the constant $K$ in Theorem \ref{th_density} is non-zero. It is independent from Sections \ref{sec_FKPP} and \ref{sec_tail} and in particular from Theorem \ref{th_tail}. Its proof is entirely probabilistic and follows closely \cite{Aidekon2010}.

\subsection{Notation and preliminary remarks}
Our notation borrows from \cite{Kyprianou2004}. An {\em individual} is an element in the space of Ulam--Harris labels \[U=\bigcup_{n\in\N_0} \N^n,\]which is endowed with the ordering relations $\preceq$ and $\prec$ defined by
\[u \preceq v \iff \exists w\in U: v = uw\quad\tand\quad u\prec v \iff u\preceq v \tand u \ne v.\] The space of Galton--Watson trees is the space of subsets $\tree \subset U$, such that $\empt \in \tree$, $v\in \tree$ if $v\prec u$ and $u\in \tree$ and for every $u$ there is a number $L_u\in\N_0$, such that for all $j\in\N$, $uj\in \tree$ if and only if $j\le L_u$. Thus, $L_u$ is the number of children of the individual $u$.

Branching Brownian motion is defined on the filtered probability space $(\T,\F,(\F_t),P)$. Here, $\T$ is the space of Galton--Watson trees with each individual $u\in \tree$ having a mark $(\zeta_u,X_u)\in \R^+\times D(\R^+,\R\cup\{\Delta\})$, where $\Delta$ is a cemetery symbol and $D(\R^+,\R\cup\{\Delta\})$ denotes the Skorokhod space of cadlag functions from $\R^+$ to $\R\cup\{\Delta\}$. Here, $\zeta_u$ denotes the life length and $X_u(t)$ the position of $u$ at time $t$, or of its ancestor that was alive at time $t$. More precisely, for $v\in\tree$, let $d_v = \sum_{w\preceq v}\zeta_w$ denote the time of death and $b_v = d_v-\zeta_v$ the time of birth of $v$. Then $X_u(t) = \Delta$ for $t\ge d_u$ and if $v\preceq u$ is such that $t\in [b_v,d_v)$, then $X_u(t) = X_v(t)$.

The sigma-field $\F_t$ contains all the information up to time $t$, and $\F = \sigma\left(\bigcup_{t\ge 0} \F_t\right)$.

Let $y,c\in\R$ and $L$ be some random variable taking values in $\N_0\backslash\{1\}$. $P = P^{y,c,L}$ is the unique probability measure, such that, starting with a single individual at the point $y$,
\begin{itemize}[nolistsep, label=$-$]
\item each individual moves according to a Brownian motion with drift $c$ until an independent time $\zeta_u$ following an exponential distribution with parameter 1.
\item At the time $\zeta_u$ the individual dies and leaves $L_u$ offspring at the position where it has died, with $L_u$ being an independent copy of $L$.
\item Each child of $u$ repeats this process, all independently of one another and of the past of the process.
\end{itemize}
Note that often $c$ and $L$ are regarded as fixed and $y$ as variable. In this case, the notation $P^y$ is used. In the same way, expectation with respect to $P$ is denoted by $E$ or $E^y$.

A common technique in branching processes since \cite{LPP1995} is to enhance the space $\T$ by selecting an infinite genealogical line of descent from the ancestor $\empt$, called the {\em spine}. More precisely, if $T\in\T$ and $\tree$ its underlying Galton--Watson tree, then $\xi = (\xi_0,\xi_1,\xi_2,\ldots) \in U^{\N_0}$ is a \emph{spine} of $T$ if $\xi_0=\empt$ and for every $n\in\N_0$, $\xi_{n+1}$ is a child of $\xi_n$ in $\tree$. This gives the space
\[\Tt = \{(T,\xi)\in \T\times U^{\N_0}:\xi \text{ is a spine of }T\}\] of marked trees with spine and the sigma-fields $\Ft$ and $\Ft_t$. Note that if $(T,\xi)\in \Tt$, then $T$ is necessarily infinite.

Assume from now on that $m = E[L] - 1 \in (0,\infty)$. Let $N_t$ be the set of individuals alive at time $t$. Note that every $\Ft_t$-measurable function $f:\Tt\to\R$ admits a representation
\[f(T,\xi) = \sum_{u\in N_t} f_u(T)\Ind_{u\in\xi},\]
where $f_u$ is an $\F_t$-measurable function for every $u\in U$. We can therefore define a measure $\Pt$ on $(\Tt,\Ft,(\Ft_t))$ by
\begin{equation}
\label{eq_manyto1}
\int_{\Tt} f \,\dd\Pt = \e^{-m t} \int_{\T} \sum_{u\in N_t} f_u(T) P(\dd T).
\end{equation}
It is known \cite{Kyprianou2004} that this definition is sound and that $\Pt$ is actually a probability measure with the following properties:
\begin{itemize}[nolistsep, label=$-$]
\item Under $\Pt$, the individuals on the spine move according to Brownian motion with drift $c$ and die at an accelerated rate $m+1$, independent of the motion.
\item When an individual on the spine dies, it leaves a random number of offspring at the point where it has died, this number following the size-biased distribution of $L$. In other words, let $\widetilde{L}$ be a random variable with $E[f(\widetilde{L})] = E[f(L)L/(m+1)]$ for every positive measurable function $f$. Then the number of offspring is an independent copy of $\widetilde{L}$.
\item Amongst those offspring, the next individual on the spine is chosen uniformly. This individual repeats the behaviour of its parent.
\item The other offspring initiate branching Brownian motions according to the law $P$.
\end{itemize}
Seen as an equation rather than a definition, \eqref{eq_manyto1} also goes by the name of ``many-to-one lemma''.

\subsection{Branching Brownian motion with two barriers}
We recall the notation $P^y$ from the previous subsection for the law of branching Brownian motion started at $y\in\R$ and $E^y$ the expectation with respect to $P^y$. Recall the definition of $\Pt$ and define $\Pt^y$ and $\widetilde{E}^y$ analogously.

Let $a,b\in \R$ such that $y\in (a,b)$. Let $\tau = \tau_{a,b}$ be the (random) set of those individuals whose paths enter $(-\infty,a]\cup[b,\infty)$ and all of whose ancestors' paths have stayed inside $(a,b)$. For $u\in\tau$ we denote by $\tau(u)$ the first exit time from $(a,b)$ by $u$'s path, i.e.\
\[\tau(u) = \inf\{t\ge0: X_u(t) \notin (a,b)\} = \min\{t\ge0: X_u(t) \in \{a,b\}\},\] and set $\tau(u) = \infty$ for $u\notin\tau$. The random set $\tau$ is an (optional) stopping line in the sense of \cite{Chauvin1991}.

For $u\in\tau$, define $X_u(\tau) = X_u(\tau(u))$. Denote by $Z_{a,b}$ the number of individuals leaving the interval $(a,b)$ at the point $a$, i.e.\ \[Z_{a,b} = \sum_{u\in\tau} \Ind_{X_u(\tau) = a}.\]

\begin{lemma}
\label{lem_2barriers}
Assume $|c| > c_0$ and define $\rho = \sqrt{c^2 - c_0^2}$. Then
\[E^y[Z_{a,b}] = \e^{c(a-y)}\frac{\sinh((b-y)\rho)}{\sinh((b-a)\rho)}.\] If, furthermore, $V = E[L(L-1)] < \infty$, then
\[
\begin{split}
E^y[Z_{a,b}^2] = \frac{2V\e^{c(a-y)}}{\rho \sinh^3((b-a)\rho)}\Big[&\sinh((b-y)\rho)\int_a^y\e^{c(a-r)}\sinh^2((b-r)\rho)\sinh((r-a)\rho)\,\dd r\\
 + &\sinh((y-a)\rho)\int_y^b\e^{c(a-r)}\sinh^3((b-r)\rho)\,\dd r\Big] + E^y[Z_{a,b}].
\end{split}
\]
\end{lemma}
\begin{proof}
On the space $\Tt$ of marked trees with spine, define the random variable $I$ by $I = i$ if $\xi_i \in \tau$ and $I=\infty$ otherwise. For an event $A$ and a random variable $Y$ write $E[Y, A]$ instead of $E[Y \Ind_A]$. Then
\[E^y[Z_{a,b}] = E^y\Big[\sum_{u\in\tau} \Ind_{X_u(\tau)=a}\Big] = \Et^y[\e^{m \tau(\xi_I)}, I < \infty, X_{\xi_I}(\tau) = a]
\] by the many-to-one lemma extended to optional stopping lines (see \cite{Biggins2004}, Lemma 14.1 for a discrete version). But since the spine follows Brownian motion with drift $c$, we have $I < \infty$, $\Pt$-a.s.\ and the above quantity is therefore equal to
\[W^{y,c}[\e^{m T}, B_T = a],\]
where $W^{y,c}$ is the law of standard Brownian motion with drift $c$ started at $y$, $(B_t)_{t\ge 0}$ the canonical process and $T=T_{a,b}$ the first exit time from $(a,b)$ of $B_t$. By Girsanov's theorem, and recalling that $m=c_0^2/2$, this is equal to
\[
W^y[\e^{c(B_T-y)-\frac 1 2 (c^2 - c_0^2) T}, B_T = a],
\] where $W^y = W^{y,0}$. Evaluating this expression (\cite{BorodinSalminen}, p.\ 212, Formula 1.3.0.5) gives the first equality.

For $u\in U$, let $\Theta_u$ be the operator that maps a tree in $\T$ to its sub-tree rooted in $u$. Denote further by $C_u$ the set of $u$'s children, i.e.\ $C_u = \{uk: 1\le k \le L_u\}$. Then note that for each $u\in \tau$ we have
\[Z_{a,b} = 1 + \sum_{v\prec u} \mathop{\sum_{w\in C_v}}_{w\npreceq u} Z_{a,b}\circ \Theta_w,\]
hence
\begin{equation}
\label{eq_decomp}
\begin{split}
E^y[Z_{a,b}^2] &= E^y\Big[\sum_{u\in\tau}\Ind_{X_u(\tau)=a}Z_{a,b}\Big]\\
&= E^y[Z_{a,b}] + \Pt^y\Bigg[\e^{m \tau(\xi_I)} \sum_{v \prec \xi_I} \mathop{\sum_{w\in L_v}}_{w\npreceq \xi_I} Z_{a,b}\circ \Theta_w,\ X_{\xi_I}(\tau) = a\Bigg].
\end{split}
\end{equation}
Define the $\sigma$-algebras
\begin{align*}
\G &= \sigma(X_{\xi_I}(t);t\ge0),\\
\H &= \G \vee \sigma(\zeta_v;v\prec \xi_I),\\
\mathcal{I} &= \H \vee \sigma(\xi, I, (L_v;v\prec \xi_I)),
\end{align*}
such that $\G$ contains the information about the path of the spine up to the individual that quits $(a,b)$ first, $\H$ adds to $\G$ the information about the fission times on the spine and $\mathcal I$ adds to $\H$ the information about the individuals of the spine and the number of their children. Now, conditioning on $\mathcal I$ and using the strong branching property, the second term in the last line of \eqref{eq_decomp} is equal to
\[\Pt^y\Bigg[\e^{m \tau(\xi_I)}\sum_{v \prec \xi_I} (L_v-1) E^{X_v(d_v-)}[Z_{a,b}],\ X_{\xi_I}(\tau) = a\Bigg]\]
(recall that $d_v$ is the time of death of $v$). Conditioning on $\H$ and noting the fact that $L_v$ follows the size-biased law of $L$ for an individual $v$ on the spine, yields
\[\Pt^y\Bigg[\e^{m \tau(\xi_I)}\sum_{v \prec \xi_I} \frac{V}{m+1}E^{X_{\xi_I}(d_v)}[Z_{a,b}],\ X_{\xi_I}(\tau) = a\Bigg].\]
Finally, since under $\Pt$ the fission times on the spine form a Poisson process of intensity $m+1$, conditioning on $\G$ and applying Girsanov's theorem yields
\[
\begin{split}
&W^y\left[\e^{c(B_T-y)-\frac 1 2 \rho^2 T}\int_0^T V E^{B_t}[Z_{a,b}] \,\dd t,B_T = a\right]\\
&\hspace{2cm} = V \e^{c(a-y)} \int_a^b E^r[Z_{a,b}] W^y\left[\e^{-\frac 1 2 \rho^2 T} L_T^r, B_T = a\right] \,\dd r,
\end{split}
\]
where $L_T^r$ is the local time of $(B_t)$ at the time $T$ and the point $r$. The last expression can be evaluated explicitly (\cite{BorodinSalminen}, p.\ 215, Formula 1.3.3.8) and gives the desired equality.
\end{proof}

\begin{corollary}
\label{cor_2barriers}
Under the assumptions of Lemma \ref{lem_2barriers}, for each $b>0$ there are positive constants $C_b^{(1)}$, $C_b^{(2)}$, such that as $a\to-\infty$,
\begin{enumerate}[nolistsep, label=\alph*)]
\item $E^0[Z_{a,b}] \sim C_b^{(1)} \e^{(c+\rho)a}$,
\item if $c > c_0$, $E^0[Z_{a,b}^2] \sim C_b^{(2)}\e^{(c+\rho)a}$ and
\item if $c < -c_0$, $E^0[Z_{a,b}^2] \sim C_b^{(2)}\e^{2(c+\rho)a}$.
\end{enumerate}
\end{corollary}

The following result is well known and is only included for completeness. We emphasize that the only moment assumption here is $m = E[L]-1 \in (0,\infty)$. Recall that $Z_x$ denotes the number of particles absorbed at $x$ of a BBM started at the origin. For $|c|\ge c_0$, define $\lc$ to be the smaller root of $\lambda^2 - 2c\lambda + c_0^2$, thus $\lc = c - \sqrt{c^2-c_0^2}$.

\begin{lemma}
\label{lem_expectation_Z}
Let $x>0$.
\begin{itemize}[nolistsep, label=$-$]
\item If $|c|\ge c_0$, then $E[Z_x] = \e^{\lc x}$.
\item If $|c| < c_0$, then $E[Z_x] = +\infty$.
\end{itemize}
\end{lemma}
\begin{proof}
We proceed similarly to the first part of Lemma \ref{lem_2barriers}. Define the (optional) stopping line $\tau$ of the individuals whose paths enter $[x,\infty)$ and all of whose ancestors' paths have stayed inside $(-\infty,x)$. Define $I$ as in the proof of Lemma \ref{lem_2barriers}. By the stopping line version of the many-to-one lemma we have
\[E[Z_x] = E[\sum_{u\in\tau} 1] = \Et[\e^{m\tau(\xi_I)}, I<\infty].\] By Girsanov's theorem, this equals
\[W[\e^{cx - \frac 1 2 (c^2 - c_0^2) T_x}, T_x<\infty],\] where $W$ is the law of standard Brownian motion started at 0 and $T_x$ is the first hitting time of $x$. The result now follows from \cite{BorodinSalminen}, p.\ 198, Formula 1.2.0.1.
\end{proof}

\subsection{Proof of Proposition \ref{prop_heavytail}}
By hypothesis, $c > c_0$, $E[L^2] < \infty$ and the BBM starts at the origin. Let $x>0$ and let $\tau = \tau_x$ be the stopping line of those individuals hitting the point $x$ for the first time. Then $Z_x = |\tau_x|$.

Let $a < 0$ and $n\in\N$. By the strong branching property,
\[P^0(Z_x > n) \ge P^0(Z_x > n \mid Z_{a,x} \ge 1) P^0(Z_{a,x} \ge 1) \ge P^a(Z_x > n) P^0(Z_{a,x} \ge 1).\] If $P^0_-$ denotes the law of branching Brownian motion started at the point $0$ with drift $-c$, then
\[P^a(Z_x > n) = P^0_-(Z_{a-x} > n) \ge P^0_-(Z_{a-x,1} > n).\] In order to bound this quantity, we choose $a = a_n$ in such a way that $n = \frac 1 2 E^0_-[Z_{a_n-x,1}].$ By Corollary \ref{cor_2barriers} a), c) (applied with drift $-c$) and the Paley--Zygmund inequality, there is then a constant $C_1 > 0$, such that
\[P^0_-(Z_{a_n-x,1} > n) \ge \frac 1 4 \frac{E^0_-[Z_{a_n-x,1}]^2}{E^0_-[Z_{a_n-x,1}^2]} \ge C_1\quad\text{ for large }n.\]
Furthermore, by Corollary \ref{cor_2barriers} a) (applied with drift $-c$), we have
\[\frac 1 2 C^{(1)}_1 \e^{-\lc (a_n-x)} \sim n,\quad \tas n\to\infty,\] and therefore $a_n = -(1/\lc)\log n + O(1)$. Again by the Paley--Zygmund inequality and Corollary \ref{cor_2barriers} a), b) (applied with drift $c$), there exists $C_2 > 0$, such that for large $n$,
\[P^0(Z_{a_n,x} \ge 1) = P^0(Z_{a_n,x} > 0) \ge \frac{E^0[Z_{a_n,x}]^2}{E^0[Z_{a_n,x}^2]} \ge \frac{(C^{(1)}_x)^2}{2C^{(2)}_x} \e^{\lcb a_n} \ge \frac{C_2}{n^d}.\] This proves the proposition with $C=C_1C_2$.

\section{The FKPP equation}
\label{sec_FKPP}
As was already observed by Neveu \cite{Neveu1988}, the translational invariance of Brownian motion and the strong branching property immediately imply that  $Z=(Z_x)_{x\ge0}$ is a homogeneous continuous-time Galton--Watson process (for an overview to these processes, see \cite{Athreya1972}, Chapter\;III or \cite{Harris1963}, Chapter\;V). There is therefore an infinitesimal generating function
\begin{equation}
\label{eq_def_a}
a(s) = \alpha\left(\sum_{n=0}^\infty p_ns^n - s\right),\quad \alpha > 0,\ p_1 = 0,
\end{equation}
associated to it. It is a strictly convex function on $[0,1]$, with $a(0) \ge 0$ and $a(1)\le 0$. Its probabilistic interpretation is
\[\alpha = \lim_{x\to 0} \tfrac 1 x P(Z_x \ne 1)\quad\tand\quad p_n = \lim_{x\to0} P(Z_x = n | Z_x \ne 1),\] hence $q_n = \alpha p_n$ for $n\in\N_0\backslash\{1\}$. Note that with no further conditions on $c$ and $L$, the sum $\sum_{n\ge0} p_n$ need not necessarily be $1$, i.e.\ the rate $\alpha p_\infty$, where $p_\infty = 1-\sum_{n\ge0} p_n$, with which the process jumps to $+\infty$, may be positive.

We further define $F_x(s) = E[s^{Z_x}]$, which is linked to $a(s)$ by Kolmogorov's forward and backward equations (\cite{Athreya1972}, p.\;106 or \cite{Harris1963}, p.\;102):
\begin{align}
\label{eq_fwd}
\frac{\partial}{\partial x} F_x(s) &= a(s)\frac{\partial}{\partial s} F_x(s) & &\text{(forward equation)}\\
\label{eq_bwd}
\frac{\partial}{\partial x} F_x(s) &= a[F_x(s)] & &\text{(backward equation)}
\end{align}
The forward equation implies that if $a(1)=0$ and $\phi(x) = E[Z_x] = \frac{\partial}{\partial s} F_x(1-)$, then $\phi'(x) = a'(1)\phi(x)$, whence $E[Z_x] = \e^{a'(1)x}$. On the other hand, if $a(1) < 0$, then the process jumps to $\infty$ with positive rate, hence $E[Z_x] = \infty$ for all $x > 0$.

The next lemma is an extension of a result which is stated, but not proven, in \cite{Neveu1988}, Equation (1.1). According to Neveu, it is due to A.\ Joffe. To the knowledge of the author, no proof of this result exists in the current literature, which is why we prove it here.

\begin{lemma}
\label{lem_GW_psi}
Let $(Y_t)_{t\ge0}$ be a homogeneous Galton--Watson process started at 1, which may explode and may jump to $+\infty$ with positive rate. Let $u(s)$ be its infinitesimal generating function and $F_t(s) = E[s^{Y_t}]$. Let $q$ be the smallest zero of $u(s)$ in $[0,1]$.
\begin{enumerate}[nolistsep]
\item If $q < 1$, then there exists $t_-\in\R\cup\{-\infty\}$ and a strictly decreasing smooth function $\psi_-:(t_-,+\infty)\to (q,1)$ with $\lim_{t\to t_-}\psi_-(t) = 1$ and $\lim_{t\to \infty}\psi_-(t) = q$, such that on $(q,1)$ we have $u = \psi_-' \circ \psi_-^{-1}$, $F_t(s) = \psi_-(\psi_-^{-1}(s) + t)$.
\item If $q > 0$, then there exists $t_+\in\R\cup\{-\infty\}$ and a strictly increasing smooth function $\psi_+:(t_+,+\infty)\to (0,q)$ with $\lim_{t\to t_+}\psi_+(t) = 0$ and $\lim_{t\to \infty}\psi_+(t) = q$, such that on $(0,q)$ we have $u = \psi_+' \circ \psi_+^{-1}$, $F_t(s) = \psi_+(\psi_+^{-1}(s) + t)$.
\end{enumerate}
The functions $\psi_-$ and $\psi_+$ are unique up to translation.

Moreover, the following statements are equivalent:
\begin{itemize}[nolistsep, label=$-$]
\item For all $t > 0$, $Y_t < \infty$ a.s.
\item $q = 1$ or $t_- = -\infty$.
\end{itemize}
\end{lemma}

\begin{proof}
We first note that $u(s) > 0$ on $(0,q)$ and $u(s) < 0$ on $(q,1)$, since $u(s)$ is strictly convex, $u(0)\ge0$ and $u(1) \le 0$. Since $F_0(s) = s$, Kolmogorov's forward equation \eqref{eq_fwd} implies that $F_t(s)$ is strictly increasing in $t$ for $s\in(0,q)$ and strictly decreasing in $t$ for $s\in(q,1)$. The backward equation \eqref{eq_bwd} implies that $F_t(s)$ converges to $q$ as $t\to\infty$ for every $s\in[0,1)$. Repeated application of \eqref{eq_bwd} yields that $F_t(s)$ is a smooth function of $t$ for every $s\in[0,1]$.

Now assume that $q < 1$. For $n\in\N$ set $s_n = 1-2^{-n}(1-q)$, such that $q < s_1 < 1$, $s_n < s_{n+1}$ and $s_n \to 1$ as $n\to\infty$. Set $t_1 = 0$ and define $t_n$ recursively by \[t_{n+1} = t_n - t',\quad \text{where $t'>0$ is such that } F_{t'}(s_{n+1}) = s_n.\] Then $(t_n)_{n\in\N}$ is a decreasing sequence and thus has a limit $t_- \in \R\cup\{-\infty\}$. We now define for $t\in (t_-,+\infty)$,
\[\psi_-(t) = F_{t-t_n}(s_n),\quad \tif t \ge t_n.\] The function $\psi_-$ is well defined, since for every $n\in\N$ and $t\ge t_n$,
\[F_{t-t_n}(s_n) = F_{t-t_n}(F_{t_n-t_{n+1}}(s_{n+1})) = F_{t-t_{n+1}}(s_{n+1}),\]
by the branching property. The same argument shows us that if $s\in (q,1)$, $s_n > s$ and $t'>0$ such that $F_{t'}(s_n) = s$, then $F_t(s) = F_{t+t'}(s_n) = \psi_-(t+t'+t_n)$ for all $t\ge 0$. In particular, $\psi_-(t'+t_n) = s$, hence $F_t(s) = \psi_-(\psi_-^{-1}(s) + t)$. The backward equation \eqref{eq_bwd} now gives 
\[u(s) = \frac{\partial}{\partial t}F_t(s)_{\big|t=0} = \psi_-'(\psi_-^{-1}(s)).\] The second part concerning $\psi_+$ is proven completely analogously. Uniqueness up to translation of $\psi_-$ and $\psi_+$ is obvious from the requirement $\psi(\psi^{-1}(s) + t) = F_t(s)$, where $\psi$ is either $\psi_-$ or $\psi_+$.

For the last statement, note that $P(Y_t < \infty) = 1$ for all $t>0$ if and only if $F_t(1-) = 1$ for all $t>0$. But this is the case exactly if $q=1$ or $t_- = -\infty$.
\end{proof}

The following proposition shows that the functions $\psi_-$ and $\psi_+$ corresponding to $(Z_x)_{x\ge0}$ are so-called {\em travelling wave} solutions of a reaction-diffusion equation called the Fisher--Kolmo\-gorov--Petrovskii--Piskounov (FKPP) equation. This should not be regarded as a new result, since Neveu (\cite{Neveu1988}, Proposition 3) proved it already for the case $c\ge c_0$ and $L=2$ a.s.\ (dyadic branching). However, his proof relied on a path decomposition result for Brownian motion, whereas we show that it follows from simple renewal argument valid for branching diffusions in general.

Recall that $f(s) = E[s^L]$ denotes the generating function of $L$. Let $q'$ be the unique fixed point of $f$ in $[0,1)$ (which exists, since $f'(1) = m+1 > 1$), and let $q$ be the smallest zero of $a(s)$ in $[0,1]$.

\begin{proposition}
\label{prop_FKPP}
Assume $c\in\R$. The functions $\psi_-$ and $\psi_+$ from Lemma \ref{lem_GW_psi} corresponding to $(Z_x)_{x\ge0}$ are solutions to the following differential equation on $(t_-,+\infty)$ and $(t_+,+\infty)$, respectively.
\begin{equation}
\label{eq_FKPP}
\frac{1}{2} \psi'' - c \psi' = \psi - f\circ\psi.
\end{equation}
Moreover, we have the following three cases:
\begin{enumerate}[nolistsep]
\item If $c\ge c_0$, then $q = q'$, $t_- = -\infty$, $a(1) = 0$, $a'(1) = \lc$, $E[Z_x] = \e^{\lc x}$ for all $x>0$.
\item If $|c|<c_0$, then $q = q'$, $t_- \in\R$, $a(1) < 0$, $a'(1) = 2c$, $P(Z_x = \infty) > 0$ for all $x > 0$.
\item If $c\le -c_0$, then $q = 1$, $a(1) = 0$, $a'(1) = \lc$, $E[Z_x] = \e^{\lc x}$ for all $x>0$.
\end{enumerate}
\end{proposition}

\begin{proof}
Let $s\in (0,1)$ and define the function $\psi_s(x) = F_x(s) = E[s^{Z_x}]$ for $x\ge 0$.
By symmetry, $Z_x$ has the same law as the number of individuals $N$ absorbed at the origin in a branching Brownian motion started at $x$ and with drift $-c$. By a standard renewal argument (Lemma\;\ref{lem_difeq_branching_diffusion}), the function $\psi_s$ is therefore a solution of \eqref{eq_FKPP} on $(0,\infty)$ with $\psi_s(0+) = s$. This proves the first statement, in view of the representation of $F_x$ in terms of $\psi-$ and $\psi_+$ given by Lemma\;\ref{lem_GW_psi}.

Let $s\in(0,1)\backslash\{q\}$ and let $\psi(s) = \psi_-(s)$ if $s>q$ and $\psi(s) = \psi_+(s)$ otherwise. By \eqref{eq_FKPP},
\[a'(s) = \frac{\psi''\circ\psi^{-1}(s)}{\psi'\circ\psi^{-1}(s)} = 2c + 2\frac{\psi\circ\psi^{-1}(s) - f\circ\psi\circ\psi^{-1}(s)}{\psi'\circ\psi^{-1}(s)} = 2c + 2\frac{s - f(s)}{a(s)},\] whence, by convexity,
\begin{equation}
\label{eq_a}
a'(s)a(s) = 2ca(s) + 2(s-f(s)),\quad s\in[0,1].
\end{equation}

Assume $|c| \ge c_0$. By Lemma \ref{lem_expectation_Z}, $E[Z_x] = \e^{\lc x}$, hence $a(1) = 0$ and $a'(1) = \lc$, in particular, $a'(1) > 0$ for $c\ge c_0$ and $a'(1) < 0$ for $c \le -c_0$. By convexity, $q < 1$ for $c\ge c_0$ and $q = 1$ for $c \le -c_0$. The last statement of Lemma \ref{lem_GW_psi} now implies that $t_- = -\infty$ if $c\ge c_0$.

Now assume $|c| < c_0$. By Lemma \ref{lem_expectation_Z}, $E[Z_x] = +\infty$ for all $x>0$, hence either $a(1) < 0$ or $a(1) = 0$ and $a'(1) = +\infty$, in particular, $q < 1$ by convexity. However, if $a(1) = 0$, then by \eqref{eq_a}, $a'(1) = 2c - 2m/a'(1)$, whence the second case cannot occur. Thus, $a(1) < 0$ and $a'(1) = 2c$ by \eqref{eq_a}.

It remains to show that $q=q'$ if $q<1$. Assume $q\ne q'$. Then $a(q') \ne 0$ by the (strict) convexity of $a$ and $a'(q') = 2c$ by \eqref{eq_a}. In particular, $a'(q') \ge a'(1)$, which is a contradiction to $a$ being strictly convex.
\end{proof}

\section{Proof of Theorem \ref{th_tail}}
\label{sec_tail}
We have $c = c_0$ by hypothesis. Let $\psi_-$ be the travelling wave from Proposition \ref{prop_FKPP}, which is defined on $\R$, since $t_- = -\infty$. Let $\phi(x) = 1-\psi_-(-x)$, such that $\phi(-\infty) = 1-q$, $\phi(+\infty) =0$ and
\begin{equation}
\label{eq_phi}
\tfrac 1 2 \phi''(x) + c_0 \phi'(x) = f(1-\phi(x)) - (1-\phi(x)),
\end{equation}
by \eqref{eq_FKPP}. Furthermore, $a(1-s) = \phi'(\phi^{-1}(s))$ and $F_x(1-s) = 1-\phi(\phi^{-1}(s) - x)$.

Under the hypothesis $E[L(\log L)^{2}] < \infty$, it is known \cite{Yang2011} that there exists $K \in (0,\infty)$, such that $\phi(x) \sim Kx\e^{-c_0 x}$ as $x\to\infty$. Since $a(1) = 0$ and $a'(1) = c_0$ by Proposition \ref{prop_FKPP}, this entails that $\phi'(x) = a(1-\phi(x)) \sim -c_0Kx\e^{-c_0x}$, as $x\to\infty$.

Set $\varphi_1 = \phi'$ and $\varphi_2 = \phi$. By \eqref{eq_phi},
\[\frac{\dd}{\dd x} \begin{pmatrix} \varphi_1(x)\\\varphi_2(x) \end{pmatrix} = \begin{pmatrix} \phi''(x)\\\phi'(x) \end{pmatrix} = \begin{pmatrix} -2c_0\phi'(x) + 2[f(1-\phi(x)) - (1-\phi(x))] \\\phi'(x) \end{pmatrix}.\] Setting $g(s) = c_0^2s + 2[f(1-s)-(1-s)] = 2[f'(1)s + f(1-s)-1]$, this gives
\begin{equation}
\label{eq_system_critical}
\frac{\dd}{\dd x} \begin{pmatrix} \varphi_1(x)\\\varphi_2(x) \end{pmatrix} = M \begin{pmatrix} \varphi_1(x)\\\varphi_2(x) \end{pmatrix} + \begin{pmatrix} g(\varphi_2(x))\\0 \end{pmatrix},\quad \text{with }M = \begin{pmatrix} -2c_0 & -c_0^2\\ 1 & 0\end{pmatrix}.
\end{equation}
The Jordan decomposition of $M$ is given by
\begin{equation}
\label{eq_Jordan_critical}
J = A^{-1}MA = \begin{pmatrix}-c_0 & 1 \\ 0 & -c_0\end{pmatrix},\quad A= \begin{pmatrix}-c_0 & 1-c_0\\ 1&1\end{pmatrix}.
\end{equation}
Setting $\begin{pmatrix} \varphi_1\\ \varphi_2 \end{pmatrix}$ = $A\begin{pmatrix} \xi_1\\ \xi_2\end{pmatrix}$, we get with $\xi = \begin{pmatrix} \xi_1 \\ \xi_2 \end{pmatrix}$:
\[\xi'(x) = J\xi(x) + \begin{pmatrix} -g(\phi(x)) \\ g(\phi(x)) \end{pmatrix},\] which, in integrated form, becomes
\begin{equation}
\label{eq_xi_integrated}
\xi(x) = \e^{xJ}\xi(0) + \e^{xJ}\int_0^x \e^{-yJ} \begin{pmatrix} -g(\phi(y)) \\ g(\phi(y)) \end{pmatrix} \,\dd y.
\end{equation}
Note that
\begin{equation}
\label{eq_exJ}
\e^{xJ} = \begin{pmatrix} \e^{-c_0 x} & x\e^{-c_0 x} \\ 0 & \e^{-c_0 x} \end{pmatrix}.
\end{equation}

With the above asymptotic of $\phi$, the integral $\int_0^\infty \e^{c_0 y} g(\phi(y))\,\dd y$ is finite by Theorem B of \cite{Bingham1974} (see also Theorem 8.1.8 in \cite{Bingham1987}). Equations \eqref{eq_xi_integrated} and \eqref{eq_exJ} now imply that
\[
 \xi_2(x) \sim \e^{-c_0 x}\Big(\xi_2(0) + \int_0^\infty \e^{c_0 y} g(\phi(y))\,\dd y\Big),
\]
and
\[
 (\xi_1 + \xi_2)(x) \sim x \e^{-c_0 x}\Big(\xi_2(0) + \int_0^\infty \e^{c_0 y} g(\phi(y))\,\dd y\Big),
\]
and since $\phi = \xi_1 + \xi_2$ and $\xi_2 = \phi' + c_0 \phi$, this gives
\begin{equation}
\label{eq_phiprime_phi}
 (\phi' + c_0 \phi)(x) \sim \phi(x)/x \sim K \e^{-c_0 x}.
\end{equation}
With this information, one can now show by elementary calculus (see Section \ref{sec_appendix_tail}), that
\begin{align}
\label{eq_aprimeprime}
a''(1-s) &\sim \frac{c_0}{s (\log \frac 1 s)^2},\quad\tand\\
\label{eq_Fxprimeprime}
F_x''(1-s) &\sim \frac{c_0 x\e^{-c_0 x}}{s (\log \frac 1 s)^2},\quad\tas s\to 0.
\end{align}
By standard Tauberian theorems (\cite{Feller1971}, Section XIII.5, Theorem 5), \eqref{eq_aprimeprime} implies that \[U(n) = \sum_{k=1}^n k^2 q_k \sim c_0 \frac{n}{(\log n)^2}, \quad \tas n\to \infty.\] By integration by parts, this entails that \[\sum_{k=n}^\infty q_k = \int_{n-}^\infty x^{-2} U(\dd x) \sim c_0 \left(2 \int_n^\infty \frac{1}{x^2 (\log x)^2} \,\dd x - \frac{1}{n (\log n)^2}\right).\] But the last integral is equivalent to $1/(n (\log n)^2)$ (\cite{Feller1971}, Section VIII.9, Theorem 1), which proves the first part of the theorem. The second part is proven analogously, using \eqref{eq_Fxprimeprime} instead.

\section{Preliminaries for the proof of Theorem \ref{th_density}}
\label{sec_preliminaries}

In light of Proposition \ref{prop_heavytail}, one may suggest that under suitable conditions on $L$ one may extend the proof of Theorem \ref{th_tail} to the subcritical case $c > c_0$ and prove that as $n\to\infty$, $P(Z_x > n) \sim C'n^{-d}$ for some constant $C'$. In order to apply Tauberian theorems, one would then have to establish asymptotics for the $(\lfloor d \rfloor + 1)$-th derivatives of $a(s)$ and $F_x(s)$ as $s\to1$. In trying to do this, one quickly sees that the known asymptotics for the travelling wave ($1-\psi(x) \sim \text{const}\times \e^{-\lc x}$ as $x\to-\infty$, see \cite{Kyprianou2004}) are not precise enough for this method to work. However, instead of relying on Tauberian theorems, one can analyse the behaviour of the {\em holomorphic} function $a(s)$ near its singular point $1$. This method is widely used in combinatorics at least since the seminal paper by Flajolet and Odlyzko \cite{FO1990} and is the basis for our proof of Theorem \ref{th_density}. Not only does it work in both 
the critical and subcritical cases, it even yields asymptotics for the density instead of the tail only.
%

In the rest of this section, we will define our notation for the complex analytic part of the proof and review some necessary general complex analytic results.

\subsection{Notation}
\label{sec_notation}
In the course of the paper, we will work in the spaces $\C$ and $\C^2$, endowed with the Euclidean topology. An open connected set is called a {\em region}, a simply connected region containing a point $z_0$ is also called a {\em neighbourhood} of $z_0$. The closure of a set $D$ is denoted by $\overline{D}$, its border by $\partial D$. The disk of radius $r$ around $z_0$ is denoted by $\D(z_0,r) = \{z\in\C: |z-z_0| < r\}$, its closure and border by $\Db(z_0,r)$ and $\partial\D(z_0,r)$, respectively. We further use the abbreviation $\D = \D(0,1)$ for the unit disk. For $0 \le \varphi \le \pi$, $r > 0$ and $x\in\R$, we define
\begin{align*}
G(\varphi,r) &= \{z\in\D(1,r)\backslash\{1\}:|\arg(1-z)|<\pi-\varphi\}, & S_+(\varphi,x) &= [x,\infty) \times (-\varphi,\varphi),\\
\Delta(\varphi,r) &= \{z\in\D(0,1+r)\backslash\{1\}: |\arg(1-z)|<\pi-\varphi\}, & S_-(\varphi,x) &= (-\infty,x] \times (-\varphi,\varphi),\\
H(\varphi,r) &= \{z\in\D(0,r)\backslash\{0\}:|\arg z| < \varphi\}.
\end{align*}
Note that $H(\varphi,r) = 1-G(\pi-\varphi,r)$. Here and during the rest of the paper, $\arg(z)$ and $\log(z)$ are the principal values of argument and logarithm, respectively.

Let $G$ be a region in $\C$, $z_0\in \overline{G}$ and $f$ and $g$ analytic functions in $G$ with $g(z)\ne 0$ for all $z\in G$. We write
\begin{align*}
f(z) = o(g(z))\quad&\iff\quad \forall \ep>0\ \exists\delta>0\ \forall z\in G\cap \D(z_0,\delta):|f(z)|\le\ep |g(z)|,\\
f(z) = O(g(z))\quad&\iff\quad \exists C\ge 0\ \exists\delta>0\ \forall z\in G\cap \D(z_0,\delta):|f(z)|\le C |g(z)|,\\
f(z) = \Otilde(g(z))\quad&\iff\quad \exists K\in\C: f(z) = Kg(z) + o(g(z)),\\
f(z) \sim g(z)\quad&\iff\quad f(z) = g(z) + o(g(z)),
\end{align*}
specifying that the relations hold as $z\to z_0$.
\subsection{Complex differential equations}

In this section, we review some basics about complex differential equations. We start with the fundamental existence and uniqueness theorem (\cite{Bie1965}, p.\;1, \cite{Hil1976}, Theorem 2.2.1, p.\ 45 or \cite{Inc1944}, Section 12.1, p.\ 281).

\begin{fact}
\label{th_cde}
Let $G$ be a region in $\C^2$ and $(w_0,z_0)$ a point in $G$. Let $f:G\to \C$ be analytic in $G$, i.e.\ $f$ is continuous and both partial derivatives exist and are continuous. Then there exists a neighbourhood $U$ of $z_0$ and a unique analytic function $w:U\to\C$, such that
\begin{enumerate}[nolistsep]
\item $w(z_0) = w_0$,
\item $(w(z),z)\in G$ for all $z\in U$ and
\item $w'(z) = f(w(z),z)$ for all $z\in U$.
\end{enumerate}
In other words, the differential equation $w' = f(w,z)$ with initial condition $w(z_0) = w_0$ has exactly one solution $w(z)$ which is analytic at $z_0$.
\end{fact}

The following standard result is a special case of a theorem by Painlev\'e (\cite{Bie1965}, p.\ 11, \cite{Hil1976}, Theorem 3.2.1, p.\ 82 or \cite{Inc1944}, Section 12.3, p.\ 286f).

\begin{fact}
\label{fact_extend_border}
Let $H$ be a region in $\C$ and $w(z)$ analytic in $H$. Let $G$ be a region in $\C^2$, such that $(w(z),z)\in G$ for each $z\in H$ and suppose that there exists an analytic function $f:G\to \C$, such that $w'(z) = f(w(z),z)$ for each $z\in H$. Let $z_0\in\partial H$. Suppose that $w(z)$ is continuous at $z_0$ and that $(w(z_0),z_0)\in G$. Then $z_0$ is a regular point of $w(z)$, i.e.\ $w(z)$ admits an analytic extension at $z_0$.
\end{fact}

Let $[z_1,\ldots,z_k]_n$ denote a power series of the variables $z_1,\ldots,z_k$, converging in a neighbourhood of $(0,\ldots,0)$ and which contains only terms of order $n$ or higher. The complex differential equation
\begin{equation}
\label{eq_briot_bouquet}
zw' = \lambda w + pz + [w,z]_2,\quad\lambda,p\in\C,
\end{equation}
was introduced in 1856 by Briot and Bouquet \cite{BB1856} as an example of a complex differential equation admitting analytic solutions at a singular point of the equation. More precisely, they obtained (\cite{Hil1976}, Theorem 11.1.1, p.\;402):
\begin{fact}
\label{fact_bb_holo}
If $\lambda$ is not a positive integer, then there exists a unique function $w(z)$ which is analytic in a neighbourhood of $z=0$ and which satisfies \eqref{eq_briot_bouquet}. Furthermore, $w(0)=0$.
\end{fact}

The singular solutions to this equation were later investigated by Poincar\'e, Picard and others (for a full bibliography, see \cite{HKM1961}). We are going to need the following result (see \cite{HKM1961}, Paragraph III.9.$2^{\text{o}}$ or \cite{Hil1976}, Theorem 11.1.3, p.\;405, but note that the latter reference is without proof and the statement is slightly incomplete).
\begin{fact}
\label{fact_bb_general}
Assume $\lambda > 0$. There exists a function $\psi(z,u) = \sum_{jk\ge0} p_{jk} z^j u^k$, converging in a neighbourhood of $(0,0)$ and such that $p_{00} = 0$ and $p_{01} = 1$, such that the general solution of \eqref{eq_briot_bouquet} which vanishes at the origin is $w = \psi(z,u)$, with
\begin{itemize}[nolistsep, label=$-$]
\item $u = Cz^\lambda$, if $\lambda \notin \N$,
\item $u = z^\lambda(C+K\log z)$, if $\lambda \in \N$.
\end{itemize}
Here, $C\in\C$ is an arbitrary constant and $K\in\C$ is a fixed constant depending only on the right-hand side of \eqref{eq_briot_bouquet}.
\end{fact}
\begin{remark}
\label{rem_bb_general}
The above statement is slightly imprecise, in that the term \emph{solution} is not defined, i.e.\ what \emph{a priori} knowledge of $w(z)$ (regarding its domain of analyticity, smoothness, behaviour at $z=0$, \ldots) is required in order to guarantee that it admits the representation stated in Fact \ref{fact_bb_general}? Inspecting the proof (as in \cite{HKM1961}, for example) shows that it is actually enough to know that $w(z)$ satisfies \eqref{eq_briot_bouquet} on an interval $(0,\ep)$ of the real line and that $w(0+) = 0$. We briefly explain why:

In order to prove Fact \ref{fact_bb_general}, one shows that there exists a function $\psi$ of the form stated above, such that when changing variables by $w = \psi(z,u)$, the function $u(z)$ \emph{formally} satisfies one of the equations
\[zu' = \lambda u\quad\text{or}\quad zu' = \lambda u + K z^\lambda,\] according to whether $\lambda\notin\N$ or $\lambda \in \N$.

Now suppose that $w(z)$ satisfies the above conditions. By the implicit function theorem (\cite{Hor1973}, Theorem 2.1.2), we can invert $\psi$ to obtain a function $\varphi(w,z) = w + qz + [w,z]_2$, $q\in\C$, such that $\psi(z,\varphi(w,z)) = w$ in a neighbourhood of $(0,0)$. We may thus define $u(z) = \varphi(w(z),z)$ for all $z\in(0,\ep_1)$ for some $\ep_1>0$. Moreover, $u(z)$ now truly satisfies the above equations on $(0,\ep_1)$ and $u(0+) = 0$. Standard theory of ordinary differential equations on the real line now yields that $u$ is necessarily of the form stated in Fact \ref{fact_bb_general}.

We further remark that since $u(z)$ is analytic in the slit plane $\C\backslash(-\infty,0]$ and goes to $0$ as $z\to 0$ in $\C\backslash(-\infty,0]$, there exists an $r>0$, such that $(z,u(z))$ is in the domain of convergence of $\psi(z,u)$ for every $z\in H(\pi,r)$. Hence, every solution $w(z)$ can be analytically extended to $H(\pi,r)$.
\end{remark}

\subsection{Singularity analysis}
\label{sec_singanalysis}
We now summarise results about the singularity analysis of generating functions. The basic references are \cite{FO1990} and \cite{FS2009}, Chapter VI. The results are of two types: those that establish an asymptotic for the coefficients of functions that are explicitly  known, and those that estimate the coefficients of functions which are dominated by another function. We start with the results of the first type:

\begin{fact}
\label{th_singanalysis_explicit}
Let $d\in (1,\infty)\backslash\N$, $k\in\N$, $\gamma\in\Z\backslash\{0\}$, $\delta\in\Z$ and the functions $f_1$, $f_2$ defined by
\[
f_1(z) = (1-z)^d\quad\tand\quad f_2(z) = (1-z)^k \left(\log \frac 1 {1-z}\right)^\gamma \left(\log \log \frac 1 {1-z}\right)^\delta,
\]
for $z\in\C\backslash[1,+\infty)$. Let $(p^{(i)}_n)$ be the coefficients of the Taylor expansion of $f_i$ around the origin, $i=1,2$. Then $(p^{(i)}_n)$ satisfy the following asymptotics as $n\to\infty$:
\[
p^{(1)}_n \sim \frac{K_1}{n^{d+1}} \quad\tand\quad p^{(2)}_n \sim \frac{K_2 (\log n)^{\gamma - 1}(\log \log n)^\delta}{n^{k+1}},
\]
for some non-zero constants $K_1 = K_1(d),K_2 = K_2(k,\gamma,\delta)$. We have $K_2(1,-1,0) = 1$.
\end{fact}
\begin{proof}
For $f_1$, this is Proposition 1 from \cite{FO1990}. For $f_2$ this is Remark 3 at the end of Chapter 3 in the same paper. Note that the additional factors $\frac 1 z$ do not change the nature of the singularities, since $\frac 1 z$ is analytic at 1 (see the footnote on p.\ 385 in \cite{FS2009}). The last statement follows from Remark 3 as well.
\end{proof}

The results of the second type are contained in the next theorem. It is identical to Corollary 4 in \cite{FO1990}. Note that a potential difficulty here is that it requires analytical extension outside the unit disk.
\begin{fact}
\label{th_singanalysis_o}
Let $0<\varphi<\pi/2$, $r > 0$ and $f(z)$ be analytic in $\Delta(\varphi,r).$ Assume that as $z\to 1$ in $\Delta(\varphi,r)$,
\[f(z) = o\left((1-z)^\alpha L\left(\frac 1 {1-z}\right)\right),\quad\text{where }L(u) = (\log u)^\gamma (\log\log u)^\delta,\quad \alpha,\gamma,\delta\in\R.\]
Then the coefficients $(p_n)$ of the Taylor expansion of $f$ around 0 satisfy
\[p_n = o\left(\frac{L(n)}{n^{\alpha+1}}\right),\quad\tas n\to\infty.\]
\end{fact}

\subsection{An equation for continuous-time Galton--Watson processes}
\label{sec_GW}
In this section, let $(Y_t)_{t\ge0}$ be a homogeneous continuous-time Galton--Watson process starting at 1. Let $a(s)$ be its infinitesimal generating function and $F_t(s) = E[s^{Y_t}]$. Assume $a(1) = 0$ and $a'(1) = \lambda \in (0,\infty),$ such that $a(s) = 0$ has a unique root $q$ in $[0,1)$.

The following proposition establishes a relation between the infinitesimal generating function of a Galton--Watson process and its generating function at time $t$. For real $s$, the formulae stated in the proposition are well known, but we will need to use them for complex $s$, which is why we have to include some (complicated) hypotheses to be sure that the functions and integrals appearing in the formulae are well defined.

\begin{proposition}
\label{prop_exp_Ft_a}
Suppose that $a$ and $F_t$ have analytic extensions to some regions $D_a$ and $D_F$. Let $Z_a = \{s\in D_a: a(s) = 0\}$. Let there be simply connected regions $G\subset D_a\backslash Z_a$ and $D\subset G\cap D_F$ with $F_t(D) \subset G$ and $D\cap (0,1) \ne \empt$. Then the following equations hold for all $s\in D$:
\begin{equation}
\label{eq_int_Ft_a}
\int_s^{F_t(s)} \frac{1}{a(r)} \,\dd r = t,
\end{equation}
and
\begin{equation}
\label{eq_exp_Ft_a}
1-F_t(s) = \e^{\lambda t}(1-s)\exp\left(-\int_s^{F_t(s)} f^*(r) \,\dd r\right),
\end{equation}
where $f^*(s)$ is defined for all $s\in D_a\backslash Z_a$ as
\begin{equation}
\label{eq_fstar_def}
f^*(s) = \frac{\lambda}{a(s)} + \frac{1}{1-s},
\end{equation}
and the integrals may be evaluated along any path from $s$ to $F_t(s)$ in $G$.
\end{proposition}
\begin{proof}
For $s\in(0,1)\backslash\{q\}$, equation \eqref{eq_int_Ft_a} follows readily from Kolmogorov's backward equation \eqref{eq_bwd}, when the integral is interpreted as the usual Riemann integral (\cite{Athreya1972}, p.\;106). Now note that by definition of $G$, both $\frac{1}{a(s)}$ and $f^*$ are analytic in the simply connected region $G$ and therefore possess antiderivatives $g$ and $h$ in $G$. Thus, the functions
\[s\mapsto \int_s^{F_t(s)} \frac{1}{a(r)} \,\dd r = g(F_t(s))-g(s)\quad\tand\quad s\mapsto \int_s^{F_t(s)} f^*(r) \,\dd r = h(F_t(s))-h(s)\]
are analytic in $D$. By the analytic continuation principle, \eqref{eq_int_Ft_a} then holds for every $s\in D$, since $D\cap(0,1)\ne \empt$ by hypothesis. This proves the first equation. For the second equation, note that $-\log (1-s)$ is an antiderivative of $\frac{1}{1-s}$ in $G$, whence the right-hand side of \eqref{eq_exp_Ft_a} equals
\[\e^{\lambda t}(1-s)\exp\left(\log(1-F_t(s)) - \log(1-s) - \lambda \int_s^{F_t(s)} \frac{1}{a(r)}\,\dd r \right) = 1-F_t(s),\]
for all $s\in D$, by \eqref{eq_int_Ft_a}. This gives \eqref{eq_exp_Ft_a}.
\end{proof}

\begin{corollary}
\label{cor_1regular}
If 1 is a regular point of $a(s)$, then it is a regular point for $F_t(s)$ for every $t\ge0$.
\end{corollary}
\begin{proof}
Define $G = \{s\in\D: \Re s > q\}$. Then $G \cap Z_a = \empt$, since $q$ is the only zero of $a$ in $\D$ (every probability generating function $g$ with $g'(1) > 1$ has exactly one fixed point $q$ in $\D$; this can easily be seen by applying Schwarz's lemma to $\tau^{-1}\circ g \circ\tau$, where $\tau$ is the M\"obius transformation of the unit disk that maps $0$ to $q$). Let $s_1\in (q,1)$ be such that $F_t(s) \in G$ for every $s\in H = \{s\in\D:\Re s > s_1\}$. We can then apply Proposition \ref{prop_exp_Ft_a} to conclude that \eqref{eq_exp_Ft_a} holds for every $s\in H$.

Since $a(s)$ is analytic in a neighbourhood $U$ of 1 by hypothesis, it is easy to show that $f^*$ is analytic in $U$ as well. Thus, $f^*$ has an antiderivative $F^*$ in $H \cup U$. We define the function $g(s) = (1-s)\exp(F^*(s))$ on $H\cup U$. Since $g'(1) = -\exp(F^*(1)) \ne 0$, there exists an inverse $g^{-1}$ of $g$ in a neighbourhood $U_1$ of $g(1) = 0$. Let $U_2\subset U$ be a neighbourhood of 1, such that $\e^{\lambda t} g(s)\in U_1$ for every $s\in U_2$. Define the analytic function $\widetilde{F}_t(s) = g^{-1}(\e^{\lambda t} g(s))$ for $s\in U_2$. Then by \eqref{eq_exp_Ft_a}, we have $F_t(s) = \widetilde{F}_t(s)$ for every $s\in H\cap U_2$, hence $\widetilde{F}_t$ is an analytic extension of $F_t$ at 1.
\end{proof}

\begin{corollary}
\label{cor_analyticity_Ft}
Suppose that $a(s)$ has an analytic extension to $G(\varphi_0,r_0)$ for some $0 < \varphi_0 < \pi$ and $r_0 > 0$. Suppose further that there exist $c\in\R$, $\gamma > 1$, such that $a(1-s) = -\lambda s + \lambda c s/\log s + O(s/|\log s|^\gamma)$ as $s\to0$. Then for every $\varphi_0 < \varphi < \pi$ there exists $r>0$, such that $F_t(s)$ can be analytically extended to $G(\varphi,r)$, mapping $G(\varphi,r)$ into $G(\varphi_0,r_0)$.
\end{corollary}
\begin{proof}
Recall that $\lambda > 0$. By hypothesis, we can then assume that $a(s) \ne 0$ in $G(\varphi_0,r_0)$ by choosing $r_0$ small enough. Then $\lambda/a$ has an antiderivative $A$ on $G(\varphi_0,r_0)$. Define $B(s) = A(1-s)$ for $s\in H(\pi-\varphi_0,r_0)$, such that \[B'(s) = \frac{1}{s(1-c/\log s + O(|\log s|^{-\gamma}))} = \frac 1 s + \frac c {s\log s} + O\left(\frac 1 {s |\log s|^{-\min(\gamma,2)}}\right).\] We can therefore apply Lemma \ref{lem_inversion_A} to $B$ and deduce that there exist $\varphi_1\in(\varphi_0,\varphi)$ and $r_1,r \in (0,r_0)$, such that $A$ is injective on $G(\varphi_1,r_1)$ and such that $A(s) + \lambda t\in A(G(\varphi_1,r_1))$ for every $s\in G(\varphi,r)$. Hence, $\widetilde{F}_t(s) = A^{-1}(A(s) + \lambda t)$ is defined and analytic on $G(\varphi,r)$. By \eqref{eq_int_Ft_a}, $\widetilde{F}_t(s) = F_t(s)$ on $G(\varphi,r) \cap \D$, hence $\widetilde{F}_t$ is an analytic extension of $F_t$, mapping $G(\varphi,r)$ into $G(\varphi_1,r_1)\subset G(\varphi_0,r_0)$ by definition.
\end{proof}

\section{Proof of Theorem \ref{th_density}}
\label{sec_proofdensity}
We turn back to branching Brownian motion and to our Galton--Watson process $Z=(Z_x)_{x\ge 0}$ of the number of individuals absorbed at the point $x$. Throughout this section, we place ourselves under the hypotheses of Theorem \ref{th_density}, i.e.\ we assume that $c\ge c_0 = \sqrt{2 m}$ and that the radius of convergence of $f(s) = E[s^L]$ is greater than $1$. The equation $\lambda^2 - 2c\lambda + c_0^2 = 0$ then has the solutions $\lc = c - \sqrt{c^2-c_0^2}$ and $\lcb = c + \sqrt{c^2-c_0^2}$, hence $\lc = \lcb = c_0$ if $c = c_0$ and $\lc < c_0 < \lcb$ otherwise. The ratio $d = \lcb/\lc$ is therefore greater than or equal to one, according to whether $c > c_0$ or $c = c_0$, respectively. Recall further that $\delta\in\N$ denotes the span of $L-1$.

Let $a(s) = \alpha(\sum_{k\ge0} p_ks^k - s)$ be the infinitesimal generating function of $Z$ and let $F_x(s)=E[s^{Z_x}]$. We recall the equation \eqref{eq_a} from Section \ref{sec_FKPP}: For $s\in[0,1]$,
\begin{equation}
\label{eq_a_recalled}
a'(s)a(s) = 2ca(s) + 2(s-f(s)).
\end{equation}
By the analytic continuation principle, this equation is satisfied on the domain of analyticity of $a(s)$, in particular, on $\D$.

We now give a quick overview of the proof. Starting point is the equation \eqref{eq_a_recalled}. We are going to see that this equation is closely related to the Briot--Bouquet equation \eqref{eq_briot_bouquet} with $\lambda = d$. The representation of the solution to this equation given by Fact\;\ref{fact_bb_general} will therefore enable us to derive asymptotics for $a(s)$ near its singular point $s=1$ (Theorem\;\ref{th_asymptotics}). Via the results in Section\;\ref{sec_GW}, we will be able to transfer these to the functions $F_x(s)$ (Corollary\;\ref{cor_asymptotics_Ft}). Finally, the theorems of Flajolet and Odlyzko in Section \ref{sec_singanalysis} yield the asymptotics for $q_n$ and $P(Z_x = n)$.

More specifically, we will see that the main singular term in the expansion of $a(1-s)$ or $F_x(1-s)$ near $s=0$ is $s^d$, if $d\notin\N$ and $s^d \log s$, if $d\in\N$. At first sight, this dichotomy might seem strange, but it becomes evident if one remembers that we expect the coefficients of $F_x(s)$ (i.e.\ the probabilities $P(Z_x = n)$, assume $\delta = 1$) to behave like $1/n^{d+1}$, if $d > 1$ (see Proposition \ref{prop_heavytail}). In light of Fact \ref{th_singanalysis_explicit}, a logarithmic factor must therefore appear if $d$ is a natural number, otherwise $F_x(s)$ would be analytic at 1, in which case its coefficients would decrease at least exponentially.

We start by determining the singular points of $a(s)$ and $F_x(s)$ on the boundary of the unit disk, which is the content of the next three lemmas.

\begin{lemma}
\label{lem_delta}
Let $X$ be a random variable with law $(p_k)_{k\in\N_0}$ and let $x>0$. Then the spans of $X-1$ and of $Z_x-1$ are equal to $\delta$.
\end{lemma}
\begin{proof}
This follows from the fact that the BBM starts with one individual and the number of individuals increases by $l-1$ when an individual gives birth to $l$ children.
\end{proof}

\begin{lemma}
\label{lem_h}
If $\delta = 1$, then $a(s)$ and $(F_x(s))_{x>0}$ are analytic at every $s_0\in\partial\D\backslash\{1\}$. If $\delta \ge 2$, then there exist a function $h(s)$ and a family of functions $(h_x(s))_{x>0}$, all analytic on $\D$, such that
\[a(s) = sh(s^\delta)\quad\tand\quad F_x(s) = sh_x(s^\delta),\]
for every $s\in\D$. Furthermore, $h$ and $(h_x)_{x>0}$ are analytic at every $s_0\in\partial\D\backslash\{1\}$.
\end{lemma}
\begin{proof}
Assume first that $\delta \ge 2$. Define
\[h(s) = \alpha(\sum_n p_{1+\delta n} s^n - 1)\quad\tand\quad h_x(s) = \sum_n P(Z_x = 1+\delta n)s^n.\] By Lemma \ref{lem_delta}, $p_{k+\delta n} = P(Z_x = k+\delta n) = 0$ for every $k\in\{2,\ldots,\delta\}$ and $n\in\Z$, whence $a(s) = sh(s^\delta)$ and $F_x(s) = sh_x(s^\delta)$ for every $s\in\D$.

We now claim that $a$ and $F_x$ are analytic at every $s_0\in\partial\D$ with $s_0^\delta \ne 1$. Note that if $\delta \ge 2$, this implies that $h$ and $h_x$ are analytic at every $s_0\in\partial\D\backslash\{1\}$, since the function $s\mapsto s^\delta$ has an analytic inverse in a neighbourhood of any $s\ne0$.

First note that by \cite{Feller1971}, Lemma\;XV.2.3, p.\;475, we have $|\sum_n p_n s_0^n| < 1$ for every $s_0\in\partial\D$, such that $s_0^\delta \ne 1$, whence $a(s_0)\ne 0$. Now write the differential equation \eqref{eq_a_recalled} in the form
\begin{equation*}
\label{eq_a_with g}
a' = \frac{2ca + 2(s-f(s))}{a} =: g(a,s).
\end{equation*}
Since the radius of convergence of $f$ is greater than 1 by hypothesis, $g$ is analytic at $(a(s_0),s_0)$. Furthermore, $a$ is continuous at $s_0$, since $\sum_np_ns^n$ converges absolutely for every $s\in\Db$. Fact\;\ref{fact_extend_border} now shows that $a$ is analytic at $s_0$.

It remains to show that $F_x$ is analytic at $s_0$. Kolmogorov's forward and backward equations \eqref{eq_fwd} and \eqref{eq_bwd} imply that $a(s)F_x'(s) = a(F_x(s))$ on $[0,1]$, and the analytic continuation principle implies that this holds on $\D$. Now, let $s_0\in\partial\D$, such that $s_0^\delta \ne 1$. Then we have just shown that $a$ is analytic and non-zero at $s_0$. Furthermore, $|F_x(s_0)|<1$, by the above stated lemma in \cite{Feller1971} and Lemma \ref{lem_delta}. Thus, the function $f(w,s) = a(w)/a(s)$ is analytic at $(F_x(s_0),s_0)$, hence we can apply Fact\;\ref{fact_extend_border} again to conclude that $F_x$ is analytic at $s_0$ as well.
\end{proof}

The next lemma ensures that we can ignore certain degenerate cases appearing in the course of the analysis of \eqref{eq_a}. It is the analytic interpretation of the probabilistic results in Section \ref{sec_probability}.

\begin{lemma}
\label{lem_1_is_singularity}
1 is a singular point of $a(s)$. If $c = c_0$, then $a''(1) = +\infty$.
\end{lemma}
\begin{proof}
If $c=c_0$, the second assertion follows from Theorem \ref{th_tail} or from Neveu's result that $E[Z_x\log^+ Z_x] = \infty$ for $x > 0$ (see the remark before Theorem \ref{th_tail}). This implies that 1 is a singular point of $a(s)$. If $c > c_0$, Proposition \ref{prop_heavytail} implies that $E[s^{Z_x}] = \infty$ for every $s > 1$, whence $1$ is a singular point of the generating function $F_x(s)$ by Pringsheim's theorem (\cite{FS2009}, Theorem\;IV.6, p.\ 240). By Corollary \ref{cor_1regular}, it follows that 1 is a singular point of $a(s)$ as well.
\end{proof}

The next theorem is the core of the proof of Theorem \ref{th_density}.

\begin{theorem}
\label{th_asymptotics}
Under the assumptions of Theorem \ref{th_density}, for every $\varphi\in (0,\pi)$ there exists $r>0$, such that $a(s)$ possesses an analytical extension (denoted by $a(s)$ as well) to $G(\varphi,r)$. Moreover, as $1-s\to 1$ in $G(\varphi,r)$, the following holds.
\begin{itemize}[label=$-$]
\item If $d=1$, then
\begin{equation}
\label{eq_asymptotics_a_1}
a(1-s) = -c_0 s + c_0 \frac{s}{\log\frac{1}{s}} - c_0s \frac{\log\log\frac{1}{s}}{(\log\frac{1}{s})^2} + \Otilde\left(\frac{s}{(\log\frac{1}{s})^2}\right).
\end{equation}
\item If $d > 1$, then there is a $K = K(c,f)\in\C\backslash\{0\}$ and a polynomial $P(s)=\sum_{n=2}^{\lfloor d \rfloor} c_n s^n$, such that 
\begin{align}
\label{eq_asymptotics_a_not_N}
&\tif d\notin\N:\quad a(1-s) = -\lc s + P(s) + K s^d + o(s^d),\\
\label{eq_asymptotics_a_N}
&\tif d\in\N:\quad a(1-s) = -\lc s + P(s) + K s^d \log s + o(s^d).
\end{align}
\end{itemize}
\end{theorem}

\begin{proof}[Proof of Theorem \ref{th_asymptotics}]
We set $b(s) = a(1-s)$. By \eqref{eq_a_recalled},
\begin{equation}
\label{eq_b_first}
-b'(s)b(s) = 2cb(s) + 2 (1-s - f(1-s))\quad\ton\D(1,1).
\end{equation} Since $f$ is analytic at 1 by hypothesis, there exists $0 < \ep_1 < 1-q$ and a function $g$ analytic on $\D(0,\ep_1)$ with $g(0) = g'(0) = 0$, such that $f(1-s) = 1-(m+1)s+g(s)$ for $s\in \D(0,\ep_1)$.

As a first step, we analyse \eqref{eq_b_first} for real non-negative $s$. Since $\ep_1 < 1-q$, $b(s) < 0$ on $(0,\ep_1)$, whence we can divide both sides by $b(s)$ to obtain
\begin{equation}
\label{eq_b}
\frac{\dd b}{\dd s} = \frac{-2cb -c_0^2 s + 2 g(s)}{b}\quad\ton(0,\ep_1).
\end{equation}

Introduce the parameter $t(s) = \int_s^{\ep_1} \frac{\dd r}{-b(r)}$, $s\in(0,\ep_1]$, such that $t(\ep_1) = 0$, $t(0+) = +\infty$ and $t(s)$ is strictly decreasing on $(0,\ep_1]$. There exists then an inverse $s(t)$ on $[0,\infty)$, which satisfies $s'(t) = b(s(t))$. Hence, we have
\[\frac{\dd b}{\dd t} = \frac{\dd b}{\dd s}\frac{\dd s}{\dd t} = -2cb(t) - c_0^2 s(t) + 2 g(s(t))\quad\ton(0,\infty),\]
In matrix form, this becomes
\begin{equation}
\label{eq_bs}
\frac{\dd}{\dd t} \begin{pmatrix} b\\s \end{pmatrix} = M \begin{pmatrix} b\\s \end{pmatrix} + \begin{pmatrix} 2 g(s)\\0 \end{pmatrix},\quad M = \begin{pmatrix} -2c & -c_0^2\\ 1 & 0\end{pmatrix},
\end{equation} for $t\in(0,\infty)$. Note that this extends \eqref{eq_system_critical} to the subcritical case. This time, the Jordan decomposition of $M$ is given by
\begin{equation}
\label{eq_Jordan_subcritical}
A^{-1}MA = \begin{pmatrix}-\lcb & 0 \\ 0 & -\lc\end{pmatrix},\quad A = \begin{pmatrix}-\lcb & -\lc \\ 1 & 1\end{pmatrix},\quad\tif c > c_0,
\end{equation}
and by \eqref{eq_Jordan_critical}, if $c = c_0$. Setting 
\begin{equation}
\label{eq_bs_BS}
\begin{pmatrix} b \\ s \end{pmatrix} = A \begin{pmatrix}B\\S\end{pmatrix},
\end{equation}
transforms \eqref{eq_bs} into
\begin{alignat}{3}
\label{eq_BS_subcritical}
\frac{\dd B}{\dd t} &= -\lcb B + [B,S]_2, & \frac{\dd S}{\dd t} &= -\lc S + [B,S]_2, &\quad  \tif c &> c_0,\\
\label{eq_BS_critical}
\frac{\dd B}{\dd t} &= -c_0 B + S + [B,S]_2, &\quad \frac{\dd S}{\dd t} &= -c_0 S + [B,S]_2, & \tif c &= c_0,
\end{alignat}
for $t\in(0,\infty)$. Furthermore, by \eqref{eq_bs_BS}, we have 
\begin{align}
\label{eq_BSs}
s &= B + S,\\
\label{eq_S}
S &= \begin{cases}
(\lcb-\lc)^{-1} (b+\lcb s), & \tif d>1,\\
b + c_0 s, & \tif d=1,
\end{cases}\\
\label{eq_B}
B &= \begin{cases}
(\lc - \lcb)^{-1}(b+\lc s), & \tif d>1,\\
-b + (1 - c_0) s, & \tif d=1.
\end{cases}
\end{align}

From now on, let $\ep_2,\ep_3,\ldots$ be positive numbers that are as small as necessary. By the strict convexity of $b$ and the fact that $b'(0) = -\lc$ by Lemma \ref{lem_expectation_Z}, equation \eqref{eq_S} implies that $S$ is a strictly convex non-negative function of $s$ on $[0,\ep_2)$. This implies that the inverse $s = s(S)$ exists and is non-negative and strictly concave on $[0,\ep_3)$. It follows that $t(S)=t(s(S))$ exists on $[0,\ep_4)$. Equations \eqref{eq_BS_subcritical} and \eqref{eq_BS_critical} then yield for $S\in(0,\ep_4)$,
\begin{alignat}{2}
\label{eq_BS_subcritical_2}
\frac{\dd B}{\dd S} &= \frac{d B + [B,S]_2}{S + [B,S]_2}, &\quad & \tif c > c_0,\\
\label{eq_BS_critical_2}
\frac{\dd B}{\dd S} &= \frac{B - c_0^{-1}S + [B,S]_2}{S + [B,S]_2}, && \tif c = c_0.
\end{alignat}

By \eqref{eq_BSs} and the fact that $s(S)$ is strictly concave, $B$ is a strictly concave function of $S$ as well, hence strictly monotone on $(0,\ep_5)$. We claim that $B(S)^2=o(S)$ as $S\to 0$. For $d>1$, one checks by \eqref{eq_S} that $S(s) \sim s$, as $s\to0$, whence $B(S) = o(S)$, as $S\to 0$, by \eqref{eq_BSs}. If $d=1$, then $b'(0) = -c_0$ by Lemma \ref{lem_expectation_Z} and $b''(0) = +\infty$ by Lemma \ref{lem_1_is_singularity}. Equation \eqref{eq_S} then implies that $S(s)/s^2 \to +\infty$ as $s\to 0$, whence $s(S) = o(\sqrt S)$. The claim now follows by \eqref{eq_BSs}.

Proposition \ref{prop_reduction} now tells us that there exists a function $h(z) = [z]_2$, such that the function $\sfr(S) = S-h(B(S))$ has an inverse $S(\sfr)$ on $(0,\ep_6)$ and $\bfr(\sfr) = B(S(\sfr))$ satisfies the Briot--Bouquet equation
\begin{equation}
\label{eq_bfr_sfr}
\sfr \bfr' = \begin{cases}
d \bfr + [\bfr,\sfr]_2, &\tif d > 1,\\
\bfr - c_0^{-1} \sfr + [\bfr,\sfr]_2, &\tif d=1,
\end{cases}
\end{equation}
on $(0,\ep_6)$. By Fact \ref{fact_bb_general} and Remark \ref{rem_bb_general}, there exists then a function $\psi(z,u) = u+rz+[z,u]_2$, $r\in\C$, such that
$\bfr(\sfr) = \psi(\sfr,u(\sfr))$, where 
\[u(z) = Cz^d, \tif d\notin \N\quad\tand\quad u(z) = Cz^d\log z, \tif d\in\N,\] for some constant $C = C(c,f) \in\C$ (the form of $u$ in the case $d\in\N$ can be obtained from the one in Fact \ref{fact_bb_general} by changing $\psi$, $C$ and $K$). Moreover, comparing the coefficient of $\sfr$ on both sides of \eqref{eq_bfr_sfr}, we get, if $d>1$, $r = dr$, whence $r=0$ and if $d = 1$: $r + C = r - c_0^{-1}$, whence $C = -c_0^{-1}$.

Assume now $d>1$. Then $\bfr = u(\sfr) + [\sfr,u(\sfr)]_2$. Recall that $B = \bfr$ and $S = \sfr + h(\bfr)$. By \eqref{eq_BSs},
\[s = B+S = \bfr + \sfr + h(\bfr) = \sfr + u(\sfr) + [\sfr,u(\sfr)]_2,\] such that $s'(\sfr) = 1 + o(1)$ and $s(\sfr) = \sfr + [\sfr]_2 + o(\sfr^\gamma)$, as $\sfr\to0$, where $\gamma=(d+\lfloor d \rfloor)/2$, if $d\notin \N$ and $\gamma = d-1/2$, if $d\in\N$. By Lemmas \ref{lem_inversion} and \ref{lem_inversion_formula_d_ge_1}, for every $\varphi_0\in(0,\pi)$ there exists $r_0>0$, such that the inverse $\sfr(s)$ exists and is analytic on $H(\varphi_0,r_0)$ and satisfies
\[\sfr(s) = s + [s]_2 + o(s^\gamma),\quad\tas s\to0.\] This entails that 
\begin{alignat*}{2}
u(\sfr) &= C\sfr^d = C(s+o(s))^d = Cs^d + o(s^d),&&\quad\tif d\notin\N,\\
u(\sfr) &= C\sfr^d\log \sfr = C(s + o(s^{3/2}))^d \log(s+o(s)) = Cs^d\log s + o(s^d),&&\quad\tif d\in\N\backslash\{1\},\\
\sfr^n &= [s]_2 + o(s^{\gamma+1})+o(s^{\gamma^2}) = [s]_2 + o(s^d),&&\quad\text{ for all $n\ge2$}.
\end{alignat*}
It follows that
\[\bfr(s) = \bfr(\sfr(s)) = u(s) + [s]_2 + o(s^d),\quad\tas s\to0.\]We finally get by \eqref{eq_bs_BS}, \[b = -\lcb B - \lc S = -\lc s + (\lc - \lcb)\bfr = -\lc s + (\lc-\lcb) u(s) + [s]_2 + o(s^d),\] which proves \eqref{eq_asymptotics_a_not_N} and \eqref{eq_asymptotics_a_N}.

If $d=1$, recall that $u(z) = c_0^{-1} z\log \frac 1 z$ and $\bfr = u(\sfr) + r\sfr + [\sfr,u(\sfr)]_2$ for some $r\in\C$. By \eqref{eq_BSs}, \[s = B+S = \bfr + \sfr + h(\bfr) = u(\sfr) + (r+1)\sfr + [\sfr,u(\sfr)]_2,\] such that $s'(\sfr) = c_0^{-1}\log(\tfrac 1 \sfr) + O(1)$ and $s(\sfr) = c_0^{-1} \sfr \log \tfrac 1 \sfr + (r+1)\sfr + o(\sfr)$. Lemma \ref{lem_inversion} now implies that for every $\varphi_0\in(0,\pi)$ there exists $r_0>0$, such that the inverse $\sfr(s)$ exists and is analytic on $H(\varphi_0,r_0)$. Now, by \eqref{eq_bs_BS}, \[b = -c_0 s + S = -c_0 s + \sfr + h(\bfr) = -c_0 s + \sfr + O(\sfr^{3/2}).\] Lemma \ref{lem_inversion_formula_d_eq_1} now yields \eqref{eq_asymptotics_a_1}.
\end{proof}

\begin{remark}
The reason why we cannot explicitly determine the constant $K$ in Theorem \ref{th_asymptotics} is that we are analysing \eqref{eq_a} only {\em locally} around the point 1. Since the solution of \eqref{eq_a} with boundary conditions $a(q) = a(1) = 0$ is unique (this follows from the uniqueness of the travelling wave solutions to the FKPP equation), a {\em global} analysis of this equation should be able to exhibit the value of $K$. But it is probably easier to refine the probabilistic arguments of Section \ref{sec_probability}, which already give a lower bound that can be easily made explicit.
\end{remark}

The asymptotics established in Theorem \ref{th_asymptotics} for the infinitesimal generating function can now be readily transferred to the generating functions $F_x(s)$.

\begin{corollary}
\label{cor_asymptotics_Ft}
Under the assumptions of Theorem \ref{th_density}, for every $x>0$ and $\varphi \in (0,\pi)$ there exists $r>0$, such that $F_x(s) = E[s^{Z_x}]$ can be analytically extended to $G(\varphi,r)$. Furthermore, the following holds as $1-s\to1$ in $G(\varphi,r)$.
\begin{itemize}[label=$-$]
\item If $d=1$, then
\begin{equation}
\label{eq_asymptotics_Ft_1}
F_x(1-s) = 1-\e^{c_0x} s + c_0 x \e^{c_0 x} \left(\frac{s}{\log \frac 1 s} - \frac{s \log\log\frac 1 s}{(\log\frac 1 s)^2}\right) + \Otilde\left(\frac{s}{(\log\frac 1 s)^2}\right).
\end{equation}
\item If $d > 1$, then there is a polynomial $P_x(s)=\sum_{n=2}^{\lfloor d \rfloor} c_n s^n$, such that 
\begin{align}
\label{eq_asymptotics_Ft_not_N}
&\tif d \notin \N:\quad F_x(1-s) = 1-\e^{\lc x} s + P_x(s) + K_x d s^d + o(s^d),\\
\label{eq_asymptotics_Ft_N}
&\tif d\in\N:\quad F_x(1-s) = 1-\e^{\lc x} s + P_x(s) + K_x s^d \log s + o(s^d),
\end{align}
where $K_x = K(\e^{\lcb x} - \e^{\lc x})/(\lcb - \lc)$, with $K$ being the constant from Theorem \ref{th_asymptotics}.
\end{itemize}
\end{corollary}

\begin{proof}
Let $0 < \varphi_0 < \varphi$. By Theorem \ref{th_asymptotics}, there exists $r_0>0$, such that $a(s)$ can be analytically extended to $G(\varphi_0,r_0)$ and satisfies the hypothesis of Corollary \ref{cor_analyticity_Ft}. It follows that there exists $r > 0$, such that $F_x(s)$ can be analytically extended to $G(\varphi,r)$ and maps $G(\varphi,r)$ into $G(\varphi_0,r_0)$. Hence, the functions
\[w(s) = 1-F_x(1-s)\quad\tand\quad I(s) = \int_s^{w(s)} f^*(1-r)\,\dd r,\]
where $f^*(s)$ is defined as in \eqref{eq_fstar_def}, are analytic in $H(\pi-\varphi,r)$. In what follows, we always assume that $s\in H(\pi-\varphi,r)$. Appearance of the symbols $\sim,O,\Otilde,o$ means that we let $s$ go to $0$ in $H(\pi-\varphi,r)$.

First of all, we note that by Proposition \ref{prop_exp_Ft_a}, we have
\begin{equation}
\label{eq_w}
w(s) = s\e^{\lc x}\exp(I(s)) = s\e^{\lc x}\Big(1 + I(s) + \sum_{k=2}^\infty \frac{I(s)^k}{k!}\Big).
\end{equation}
Now assume $d > 1$. By Theorem \ref{th_asymptotics}, $a(1-s) = -\lc s + [s]_2 + u(s) + o(s^d)$, where $u(s) = K s^d$ or $u(s) = K s^d \log s$, according to whether $d\notin \N$ or $d\in\N$, respectively. It follows that
\[
\begin{split}
f^*(1-s) &= \frac{\lc}{a(1-s)} + \frac{1}{s} = -\frac{1}{s}\left(1-[s]_1-\tfrac{u(s)}{\lc s} + o(s^{d-1})\right)^{-1} + \frac 1 s\\
&= [s]_0 - \tfrac{u(s)}{\lc s^2} + o(s^{d-2}).
\end{split}
\]
Now, $\int_s^{w(s)} o(r^{d-2})\,\dd r = o(s^{d-1})$, since $w(s) \sim s\e^{\lc x}$ by Lemma \ref{lem_expectation_Z}. Thus,
\begin{equation}
\label{eq_I}
I(s) = \int_s^{w(s)} f^*(1-r)\,\dd r = [w(s),s]_1 - \int_s^{w(s)}\frac{u(r)}{\lc r^2} \,\dd r + o(s^{d-1}).
\end{equation}
Since $\int_s^{w(s)}r^{-2} u(r) \,\dd r = O(s^{d-1}\log s),$ equations \eqref{eq_w} and \eqref{eq_I} now give
\begin{equation}
\label{eq_w_dg1}
w(s) = s\e^{\lc x}\left(1 + [w(s),s]_1 - \int_s^{w(s)}\frac{u(r)}{\lc r^2} \,\dd r + o(s^{d-1})\right).
\end{equation}
If $d\ge2$, we deduce that $w(s) = s\e^{\lc x} + o(s^{3/2})$. Straightforward calculus now shows that
\begin{equation}
\label{eq_int_sing}
\int_s^{w(s)}\frac{u(r)}{\lc r^2} \,\dd r = \frac{K_x}{\e^{\lc x}}\frac{u(s)}{s} + [w(s),s]_1 + o(s^{d-1}),
\end{equation}
and \eqref{eq_w_dg1} and \eqref{eq_int_sing} now yield
\[w(s) = s\e^{\lc x} + [w(s),s]_2 - K_x u(s) + o(s^d).\]
Repeated application of this equation shows that $w(s) = s\e^{\lc x} + [s]_2 - K_x u(s) + o(s^d)$, which yields \eqref{eq_asymptotics_Ft_not_N} and \eqref{eq_asymptotics_Ft_N}.

In the critical case $d=1$, Theorem \ref{th_asymptotics} tells us that
\begin{equation}
\label{eq_fstar_asymptotics_d1}
f^*(1-s) = \frac{1}{s}\left(-\frac{1}{\log\frac 1 s} + \frac{\log\log\frac 1 s}{(\log\frac 1 s)^2} + \Otilde\left(\frac 1 {(\log\frac 1 s)^2}\right)\right).
\end{equation}
Write $\lambda = \lc = c_0$. For our first approximation of $w(s)$, we note that
\[
I(s) \sim - \int_s^{s\e^{\lambda x}} \frac 1 {r\log \frac 1 r} \,\dd r \sim - \frac{1}{\log \frac 1 s}\int_s^{s\e^{\lambda x}} \frac 1 r\,\dd r = -\frac{\lambda x}{\log \frac 1 s},
\]hence, by \eqref{eq_w},
\begin{equation}
\label{eq_w_first}
w(s) = s\e^{\lambda x}\left(1-\frac{\lambda x}{\log \frac 1 s} + o\left(\frac 1 {\log \frac 1 s}\right)\right).
\end{equation}

To obtain a finer approximation, we decompose $I(s)$ into
\[I(s) = \int_s^{s\e^{\lambda x}} f^*(1-r)\,\dd r + \int_{s\e^{\lambda x}}^{w(s)} f^*(1-r)\,\dd r =: I_1(s) + I_2(s).\]
We then have
\[
I_1(s) = -\frac{\lambda x}{\log \frac 1 s} + \lambda x\frac{\log\log \frac 1 s}{(\log \frac 1 s)^2} + \Otilde\left(\frac 1 {(\log \frac 1 s)^2}\right),\]
and, because of \eqref{eq_w_first},
\[
-I_2(s) \sim \int_{s\e^{\lambda x}}^{s\e^{\lambda x}(1+\frac{\lambda x}{\log s})} \frac 1 {r\log \frac 1 r} \,\dd r \sim \frac{\lambda x}{(\log\frac 1 s)^2}.
\]
Plugging this back into \eqref{eq_w} finishes the proof.
\end{proof}

\begin{proof}[Proof of Theorem \ref{th_density}]
Let $x>0$. We want to apply the methods from singularity analysis reviewed in Section \ref{sec_singanalysis} to the functions $a$ and $F_x$, if $\delta = 1$, or the functions $h$ and $h_x$ from Lemma \ref{lem_h}, if $\delta \ge 2$. Let $\varphi\in(0,\pi/2)$. By Theorem \ref{th_asymptotics} and Corollary \ref{cor_asymptotics_Ft}, there exists $r_0>0$, such that $a$ and $F_x$ can be analytically extended to $G(\varphi,r_0)$, which implies that for some $\varphi_1\in(\varphi,\pi/2)$ and $r_1\in(0,r)$, $h$ and $h_x$ can be extended to $G(\varphi_1,r_1)$, as well. Moreover, by Lemma \ref{lem_h}, each of these functions is analytic in a neighbourhood of every point of $C = \{s\in\partial\D: |1-s| \ge r_1/2\}$, which is a compact set. Hence, there exists a finite number of neighbourhoods which cover $C$. It is then easy to show that there exists $r>0$, such that the functions are analytic in $\Delta(\varphi_1,r)$.

If $\delta = 1$, we can then immediately apply Facts \ref{th_singanalysis_explicit} and \ref{th_singanalysis_o}, together with the asymptotics on $a$ and $F_x$ established in Theorem \ref{th_asymptotics} and Corollary \ref{cor_asymptotics_Ft}, to prove Theorem \ref{th_density}.

If $\delta \ge 2$, let $q(s)$ be the inverse of $s\mapsto s^\delta$ in a neighbourhood of 1, then $h(s) = a(q(s))/q(s)$ near 1, by Lemma \ref{lem_h}. But since $q'(1) = 1/\delta$, we have
\[
h(1-s) = a(1-(\frac 1 \delta s + c_2s^2 + c_3s^3 + \cdots))(1+c_1's+c_2's+\cdots),
\]
for some constants $c_n,c_n'$, and so equations \eqref{eq_asymptotics_a_1}, \eqref{eq_asymptotics_a_not_N} and \eqref{eq_asymptotics_a_N} transfer to $h$ with the coefficient of the main singular term divided by $\delta^d$. We can therefore use Facts \ref{th_singanalysis_explicit} and \ref{th_singanalysis_o} for the function $h$ to obtain the asymptotic for $(p_{\delta n+1})_{n\in\N}$ in Theorem \ref{th_density}. In the same way, equations \eqref{eq_asymptotics_Ft_1}, \eqref{eq_asymptotics_Ft_not_N} and \eqref{eq_asymptotics_Ft_N} yield asymptotics for $h_x$, such that we can use again Facts \ref{th_singanalysis_explicit} and \ref{th_singanalysis_o} to prove the second part of Theorem \ref{th_density}.
\end{proof}

\section{Appendix}

\subsection{A renewal argument for branching diffusions}
Let $W =(W_t)_{t\ge0}$ be a diffusion on an interval with endpoints $a'\le 0 < a$, such that $\lim_{x\downarrow 0} \P^x[T_0 < t] = 1$ for every $t>0$, where $T_0 = \inf\{t \ge 0: W_t = 0\}$ and $W_0 = x$, $\P^x$-almost surely. For $x\in (0,a)$, and only in the scope of this section, we define $P^x$ to be the law of the branching diffusion starting with a single particle at position $x$ where the particles move according to the diffusion $W$ and branch with rate $\beta$ according to the reproduction law with generating function $f(s)$. Moreover, particles hitting the point $0$ are absorbed at that point. Denote by $Z$ the number of particles absorbed during the lifetime of the process and define $u_s(x) = P^x[s^Z]$ for $s\in[0,1)$ and $x\in (0,a)$.

\begin{lemma}
\label{lem_difeq_branching_diffusion}
Let $s\in[0,1)$ and $\G$ be the generator of the diffusion $W$. Then
\[\G u_s = \beta(u_s-f\circ u_s)\quad\ton (0,a), \quad\text{with }u_s(0+) = s.\]
\end{lemma}
\begin{proof}
The proof proceeds by a renewal argument similar to the one in \cite{McKean1975}. As for the BBM, for an individual $u$, we denote by $\zeta_u$ its time of death, $X_u(t)$ its position at time $t$ and $L_u$ the number of $u$'s children. Define the event $A = \{\exists t\in [0,\zeta_\empt): X_\empt(t) = 0\}$. For $s\in[0,1)$ we have by the strong branching property
\[
\begin{split}
u_s(x) = E^x[s^Z] &= s P^x(A) + E^x\left[\left(E^{X_\empt(\zeta_\empt-)}[s^Z]\right)^{L_\empt}, A^c \right]\\
&= s \P^x(T_0 < \xi) + \E^x[f(u_s(W_\xi)), \xi \le T_0],
\end{split}
\] where $W=(W_t)_{t\ge0}$ is a diffusion with generator $\G$ starting at $x$ under $\P^x$, $T_0 = \inf\{t\ge0: W_t = 0\}$ and the random variable $\xi$ is exponentially distributed with rate $\beta$ and independent from $W$. Setting $v(x) = \P^x(T_0 < \xi)$ we get by integration by parts
\[v(x) = \int_0^\infty \beta \e^{-\beta t} \P^x(T_0 < t) \,\dd t = \int_0^\infty \e^{-\beta t} \P^x(T_0 \in \dd t) = \E^x\left[\e^{-\beta T_0}\right],\]
and therefore $\G v = \beta v$ on $(0,a)$ (\cite{BorodinSalminen}, Paragraph II.1.10, p.\ 18).

Denote the $\beta$-resolvent of the diffusion by $R_\beta$. By the strong Markov property,
\[\begin{split}
\E^x[f(u_s(W_\xi)), \xi \le T_0] &= \E^x\Big[\int_0^\infty \beta \e^{-\beta t} f(u_s(W_t)) \,\dd t\Big] - \E^x\Big[\int_{T_0}^\infty \beta \e^{-\beta t} f(u_s(W_t)) \,\dd t\Big]\\
 &= \beta R_\beta(f\circ u_s)(x) - \beta \E^x[\e^{-\beta T_0}]R_\beta(f\circ u_s)(0),
 \end{split}\]
hence $u_s = C_{s,\beta}v + \beta R_\beta(f\circ u_s)$, with $C_{s,\beta} = s - \beta R_\beta(f\circ u_s)(0)$. It follows that
\[\G u_s = \beta C_{s,\beta}v + \beta^2 R_\beta(f\circ u_s) - \beta(f\circ u_s) = \beta(u_s - f\circ u_s)\quad\ton (0,a).\]
By the above hypothesis on $W$, $\P^x(T_0<\xi) \to 1$ as $x\downarrow 0$, whence $u_s(0+) = s$.
\end{proof}

\subsection{Addendum to the proof of Theorem \ref{th_tail}}
\label{sec_appendix_tail}
With the notation used in the proof of Theorem \ref{th_tail}, recall that for some constant $K>0$ we have
\[
 (\phi' + c_0 \phi)(x) \sim \phi(x)/x \sim K \e^{-c_0 x},\quad \text{ as }x\to\infty.
\]
In what follows, formulae containing the symbols $\sim$ and $o()$ are meant to hold as $s\downarrow 0$. The above equation yields
\begin{equation}
\label{eq_a_asymptotic}
a(1-s) = \phi'(\phi^{-1}(s)) = -c_0s + (\phi' + c_0\phi)(\phi^{-1}(s)) = -c_0s + \frac{(c_0+o(1))s}{\log \frac 1 s}.
\end{equation}
Now, by \eqref{eq_a}, we have
\[
 a'(1-s)a(1-s) = 2c_0a(1-s) + c_0^2s - g(s),
\]
where we recall that $g(s)$ was defined as $g(s) = 2(f(1-s) - 1 + f'(1)s)$. From the above equation, one gets
\[
 a''(1-s) = -(a(1-s))^{-3}\Big(\big(c_0a(1-s) + c_0^2s-g(s)\big)^2 - g'(s) a(1-s)^2\Big).
\]
By definition, $g(s) = o(s)$ and $g'(s) = o(1)$. The previous equation together with \eqref{eq_a_asymptotic} then yields \eqref{eq_aprimeprime}.

Kolmogorov's forward and backward equations \eqref{eq_fwd} and \eqref{eq_bwd} give
\[F_x'(s) = \frac{a(F_x(s))}{a(s)},\]
and taking the derivative on both sides of this equation gives
\begin{equation}
\label{eq_Fxpp_asymp}
F_x''(s) = \frac{a(F_x(s))}{a(s)^2}\left(a'(F_x(s))-a'(s)\right).
\end{equation}
By \eqref{eq_aprimeprime} and $F_x'(1) = E[Z_x] = \e^{c_0x}$, we get
\[
a'(F_x(1-s))-a'(1-s) = -\int_s^{1-F_x(1-s)} a''(1-r)\,\dd r \sim -\frac{c_0^2x}{(\log s)^2}.
\]
This equation, together with \eqref{eq_Fxpp_asymp} now yields
\[F_x''(1-s) \sim \frac{-c_0\e^{c_0x}s}{c_0^2s^2} \left(-\frac{c_0^2x}{(\log s)^2}\right),\] which is \eqref{eq_Fxprimeprime}.

\subsection{Reduction to Briot--Bouquet equations}
\def\wfrak{\mathfrak{w}}
\def\zfrak{\mathfrak{z}}
In this section, we show how one can reduce differential equations as those obtained in the proof of Theorem \ref{th_density} to the canonical form \eqref{eq_briot_bouquet}. It is mostly based on pp.\;64 and 65 of \cite{Bie1965}.

\begin{lemma}
\label{lem_holo_solution}
Let $\lambda \in (0,1]$ and $p\in\C$. Then the equation
\begin{equation}
\label{eq_w_holo_solution}
w' = \frac{\lambda w + [w,z]_2}{z + pw + [w,z]_2}.
\end{equation}
has an analytic solution $w(z) = [z]_2$ in a neighbourhood of the origin.
\end{lemma}
\begin{proof}
We choose the {\em ansatz} $w = z\cdot w_1$. This transforms \eqref{eq_w_holo_solution} into
\[zw_1' + w_1 = \frac{\lambda z w_1+z^2[w_1,z]_0}{z+pzw_1+z^2[w_1,z]_0} = \frac{\lambda w_1+z[w_1,z]_0}{1+[w_1,z]_1}.\] Writing the inverse of the denominator as a power series in $w_1$ and $z$, this equals
\[(\lambda w_1+z[w_1,z]_0)(1+[w_1,z]_1) = \lambda w_1 + rz + [w_1,z]_2,\] for some $r\in\C$. This finally yields \[zw_1' = (\lambda-1) w_1 + rz + [w_1,z]_2.\]
Since $\lambda-1$ is not a positive integer, this equation now has an analytic solution $w_1(z) = [z]_1$ by Fact \ref{fact_bb_holo}, whence $w(z) = zw_1(z) = [z]_2$ solves \eqref{eq_w_holo_solution}.
\end{proof}
\begin{remark}
The important point in Lemma \ref{lem_holo_solution} is that the coefficient of $z$ in the numerator of \eqref{eq_w_holo_solution} is 0, which is why $w'(z) = 0$.
\end{remark}

\begin{proposition}
\label{prop_reduction}
Let $\lambda \ge 1$ and $p\in\R$. Suppose $w(z)$ is a strictly monotone real-valued function on $(0,\ep)$, $\ep>0$, with $w(z)^2 = o(z)$ as $z\to 0$ and satisfying
\begin{equation}
\label{eq_nonreduced}
w' = \frac{\lambda w + pz + [w,z]_2}{z + [w,z]_2}\quad \ton (0,\ep).
\end{equation}
Then there exists $h(z) = [z]_2$ and $\ep_1>0$, such that $\zfrak = z- h(w)$ has an inverse $z = z(\zfrak)$ on $(0,\ep_1)$ and such that
\begin{equation}
\label{eq_reduced}
\zfrak \frac{\dd w}{\dd\zfrak} = \lambda w + p\zfrak + [w,\zfrak]_2\quad\ton(0,\ep_1).
\end{equation}
\end{proposition}
\begin{proof}
By hypothesis, $w(z)$ is monotone on $(0,\ep)$ and therefore possesses an inverse $z = z(w)$ on $(0,\delta)$, $\delta>0$, which satisfies
\begin{equation}
\label{eq_inverse_z_w}
\frac{\dd z}{\dd w} = \frac{\lambda^{-1} z + [w,z]_2}{w + p\lambda^{-1} z + [w,z]_2}.
\end{equation}
By Lemma \ref{lem_holo_solution}, there exists then an analytic solution $z = g(w) = [w]_2$ to \eqref{eq_inverse_z_w}, since $\lambda^{-1}\in (0,1]$ by hypothesis.
Setting $\zfrak = z-g(w)$ transforms \eqref{eq_inverse_z_w} into a differential equation, which has $\zfrak = 0$ as a solution, hence it is of the form
\[\frac{\dd \zfrak}{\dd w} = \frac{\lambda^{-1} \zfrak + \zfrak[w,\zfrak]_1}{w + p\lambda^{-1} \zfrak + [w,\zfrak]_2}.\]
We have $\dd\zfrak/\dd z = 1 + g'(w(z))w'(z) = 1 + O(w(z)w'(z))$. By \eqref{eq_nonreduced},
\[w(z)w'(z) = O(w(z)^2/z + w(z)) = o(1),\] by hypothesis. Hence, there exists $\ep_1>0$, such that  $\zfrak(z)$ is strictly increasing on $(0,\ep_1)$ and therefore has an inverse. Thus, $w(\zfrak) = w(z(\zfrak))$ satisfies
\[\frac{\dd w}{\dd \zfrak} = \frac{w + p\lambda^{-1} \zfrak + [w,\zfrak]_2}{\lambda^{-1} \zfrak (1 + [w,\zfrak]_1)}\quad\ton (0,\ep_2),\] for some $\ep_2 > 0$. Expanding $(1+[w,\zfrak]_1)^{-1}$ as a power series at $(w,\zfrak) = (0,0)$ gives \eqref{eq_reduced}.
\end{proof}

\subsection{Inversion of some analytic functions}
The results in this section are needed in the proofs of Corollary \ref{cor_analyticity_Ft} and Theorem \ref{th_density}.

\begin{lemma}
\label{lem_injective_log}
Let $\varphi\in(0,\pi)$, $r>0$ and $h$ be an analytic function on $H(\varphi,r)$ with $h(z) = o(z)$ as $z\to 0$. Then there exists $r_1>0$, such that for all $z_1,z_2\in H(\varphi,r_1)$,
\[\log z_1 - \log z_2 + \int_{z_2}^{z_1} h(z)\,\dd z \ne 0.\]
\end{lemma}
\begin{proof}
Let $z_1,z_2\in H(\varphi,r)$. Write $z_i = a_i \e^{i\varphi_i}$, with $a_i>0$, $\varphi_i\in(-\varphi,\varphi)$, $i=1,2$. Define the paths
\[\gamma_1(t) = a_2 \e^{i(t\varphi_1 + (1-t)\varphi_2)}\quad\tand\quad \gamma_2(t) = (ta_1+(1-t)a_2)\e^{i\varphi_1},\quad t\in[0,1],\] such that their concatenation forms a path from $z_2$ to $z_1$ in $H(\varphi,r)$. Then
\[\int_{\gamma_1} h(s)\,\dd s = |\varphi_1 - \varphi_2|\cdot a_2o(1/a_2)\quad\tand\quad \int_{\gamma_2} h(s)\,\dd s = |\log a_1 - \log a_2| o(1).\] As a consequence,
\[\left|\int_{z_2}^{z_1} h(s)\,\dd s\right| = (|\varphi_1-\varphi_2|+|\log a_1 - \log a_2|)o(1) \le \sqrt 2 |\log z_1 - \log z_2| o(1).\]
This proves the statement.
\end{proof}

\begin{lemma}
\label{lem_inversion}
Let $r>0$ and $\varphi\in(0,\pi]$. Let $g$ and $h$ be analytic functions on $H(\varphi,r)$ with
$g'(z) = 1+o(1)$, $h'(z) = \log \frac 1 z +O(1)$, $g(z) \to 0$ and $h(z) \to 0$ as $z\to 0$ in $H(\varphi,r)$. Then for each $\varphi_0 \in (0,\varphi)$ and $\varphi_1\in(\varphi_0,\varphi)$ there exist $r_0,r_1>0$, such that $g$ and $h$ are injective on $H(\varphi_1,r_1)$ and the images of $H(\varphi_1,r_1)$ by $g$ and $h$ contain $H(\varphi_0,r_0)$.
\end{lemma}
\begin{proof}
By hypothesis, $g(z) = z+o(z)$ as $z\to0$ in $H(\varphi,r)$, whence $\arg g(z) = \arg z + o(1)$. Thus, there exists $r_1>0$, such that $g(H(\varphi_1,r_1)) \subset \C\backslash(-\infty,0]$.

Suppose that there exist $z_1,z_2\in H(\varphi_1,r_1)$, such that $g(z_1) = g(z_2)$. Let $\gamma$ be a path from $z_2$ to $z_1$ in $H(\varphi_1,r_1)$. Then $g\circ\gamma$ is a loop in $\C\backslash(-\infty,0]$, whence
\[0 = \int_{g\circ\gamma}\frac 1 z\,\dd z = \int_\gamma \frac{g'(z)}{g(z)}\,\dd z = \log z_1 - \log z_2 + \int_\gamma o(\tfrac 1 z)\,\dd z.\]
By Lemma \ref{lem_injective_log}, we can choose $r_1$ so small, that this equality cannot hold, whence $g$ is injective on $H(\varphi_1,r_1)$.

Since $g(z)\to0$ and $\arg g(z) = \arg z + o(1)$ as $z\to0$, there exists $r_0>0$, such that $g(\partial H(\varphi_1,r_1))$ encloses $H(\varphi_0,r_0)$. Now, since $g$ is injective on $H(\varphi_1,r_1)$, $H(\varphi_1,r_1)$ and $g(H(\varphi_1,r_1))$ are conformally equivalent,  whence $g(H(\varphi_1,r_1))$ is simply connected. It follows that $g(H(\varphi_1,r_1)) \supset H(\varphi_0,r_0)$.

Exactly the same arguments hold for $h$, since $h(z) = z(\log\tfrac 1 z + O(1))$ by hypothesis, whence $\arg h(z) = \arg z + o(1)$ and $h'(z)/h(z) = 1/z + o(1/z)$ as $z\to 0$.
\end{proof}

\begin{lemma}
\label{lem_inversion_A}
Let $r>0$, $\varphi\in(0,\pi]$ and $t \in \R$. Let $g$ be an analytic function on $H(\varphi,r)$ with
\[g'(z) = \frac 1 z + \frac{c}{z\log z}+O\left(\frac{1}{z|\log z|^\gamma}\right),\quad\tas z\to 0\text{ in }H(\varphi,r),\] for some $c\in\R$ and $\gamma > 1$. Then for each $0 < \varphi_0 < \varphi_1 < \varphi$ there exist $r_0,r_1>0$, such that $g$ is injective on $H(\varphi_1,r_1)$ and $g(z)+t\in g(H(\varphi_1,r_1))$ for every $z\in H(\varphi_0,r_0)$.
\end{lemma}
\begin{proof}
By the hypothesis on $g$, we have for $z_1,z_2\in H(\varphi,r)$, 
\[g(z_1)-g(z_2) = \log z_1 - \log z_2 + \int_{z_2}^{z_1} o(1/z)\,\dd z.\]
By Lemma \ref{lem_injective_log}, there exists therefore $r_1>0$, such that $g$ is injective on  $H(\varphi_1,r_1)$.

Since $1/(x|\log x|^\gamma)$ is integrable near 0, we have
\[g(z) = \log z + c \log(\log \tfrac 1 z) + o(1),\quad\tas z\to 0,\]
where we assume without loss of generalisation that the constant of integration is 0. It follows that $\Re g(z)\to -\infty$ and $\Im g(z) = \arg z + o(1)$ as $z\to 0$, since $c\in\R$. Hence, there exists an $R\in\R$, such that $g(\partial H(\varphi_1,r_1))$ encloses the strip $S = S_-(R,\varphi_1)$. As in the proof of Lemma \ref{lem_inversion}, it follows that $S\subset g(H(\varphi_1,r_1))$. Furthermore, again by the asymptotics of $\Re g$ and $\Im g$, there exists $r_0>0$, such that $g(s)+ t \in S$ for every $s\in G(\varphi_0,r_0)$. This concludes the proof.
\end{proof}

\begin{lemma}
\label{lem_inversion_formula_d_ge_1}
Let $w(z)$ be an analytic function on an open subset of $\C\backslash(-\infty,0]$, such that $w(z)\to 0$ as $z\to0$ and
\[z = w + a_2w^2 + \cdots + a_nw^n + o(w^\gamma),\quad\tas z\to0,\] for some $n\in\N$, $\gamma > n$ and $a_2,\ldots,a_n\in\C$. Then there exist $b_2,\ldots,b_n\in\C$, such that
\[w(z) = z + b_2z^2 + \cdots + b_nz^n + o(z^\gamma),\quad\tas z\to0.\]
\end{lemma}
\begin{proof}
For every $i\in\N$, we have by hypothesis 
\[z^i = w^i + a_{i,i+1}w^{i+1} + \cdots + a_{i,n}w^n + o(w^\gamma),\] for some $a_{i,i+1},\ldots,a_{i,n}\in\C$. For $2\le k\le n$, define recursively (with $b_1 = 1$) \[b_{k} = -(a_{1,k} + b_2 a_{2,k} + \cdots + b_{k-1} a_{k-1,k}).\] Then, $z + b_2 z^2 + \cdots + b^n z^n = w + o(w^\gamma).$ The statement now follows from the fact that $w(z)\sim z$ as $z\to0$ by hypothesis, whence $o(w^\gamma) = o(z^\gamma)$.
\end{proof}

\begin{lemma}
\label{lem_inversion_formula_d_eq_1}
Let $w(z)$ be an analytic function on an open subset of $\C\backslash(-\infty,0]$, such that $w(z)\to 0$ as $z\to0$ and
\[cz = w\log\tfrac 1 w + Cw + o(w),\quad\tas z\to 0,\] for some constants $c>0$, $C\in\C$. Then
\[w = \frac{cz}{\log \tfrac 1 z}\left(1-\frac{\log\log \tfrac 1 z + C-\log c + o(1)}{\log \tfrac 1 z}\right),\quad\tas z\to0.\]
\end{lemma}

\begin{proof}
Set $f(z) = w(z)/z$. By hypothesis, $\log z \sim \log w = \log z + \log f(z)$, whence $\log f(z) = o(\log z)$. Now define $g(z)$ by
\[w(z) = \frac{cz}{\log \tfrac 1 z} g(z),\]
such that $\log g(z) = \log f(z) - \log \log \tfrac 1 z = o(\log z).$
By hypothesis,
\[
cz \sim w\log\tfrac 1 w = \frac{cz}{\log \tfrac 1 z}g(z)\left(\log\log\tfrac 1 z+\log \tfrac 1 z-\log c-\log g(z)\right) \sim czg(z),
\]
whence $g(z)\sim 1$, which implies $\log g(z) = o(1)$. It now follows from the hypothesis that
\[
cz = \frac{cz}{\log \tfrac 1 z}g(z)\Big(\log\log \tfrac 1 z - \log cz + C + o(1)\Big),
\]
whence
\[g(z) = \left(1+\frac{\log\log \tfrac 1 z + C - \log c + o(1)}{\log \tfrac 1 z}\right)^{-1}.\] The statement now follows from the series representation of $(1+z)^{-1}$ at $z=0$.
\end{proof}

\section*{Acknowledgements}
\addcontentsline{toc}{section}{Acknowledgements}
The author is grateful to Elie A\"{\i}d\'ekon for stimulating discussions, to Louis Boutet de Monvel for several useful comments, to Jean Bertoin for useful comments on the proof of Lemma\;\ref{lem_difeq_branching_diffusion} and to Andreas Kyprianou and Yanxia Ren for having brought to my attention the reference \cite{Yang2011}. Furthermore, the author thanks two anonymous referees for their detailed comments, which led to a considerable improvement in the exposition of the paper.
%
%

\chapter{Branching Brownian motion with selection of the \texorpdfstring{$N$}{N} right-most particles}
\label{ch:NBBM}

\section{Introduction}
\label{sec:NBBM_intro}

In this chapter, we consider the $N$-BBM, whose definition we recall: Given a reproduction law $(q(k))_{k\ge 0}$ with $m=\sum (k-1)q(k) > 0$ and finite second moment, particles diffuse according to standard Brownian motion and branch at rate $\beta_0 = 1/(2m)$ into $k$ particles with probability $q(k)$. Furthermore, we fix a (large) parameter $N$ and as soon as the number of particles exceeds $N$, we keep the $N$ right-most particles and instantaneously kill the others.

For a finite counting measure $\nu$ on $\R$, define for $\alpha\in(0,1)$ and $N\in\N$, 
\[
\med^N_\alpha(\nu) = \inf\{x\in\R: \nu([x,\infty)) < \alpha N\}.
\]
Furthermore, set 
\[
x_\alpha = \inf\{x\ge 0: \int_x^\infty  y e^{- y}\,\dd y \le \alpha\}.
\]
For $N\in\N$ large enough, we then define $a_N = \log N + 3 \log\log N$ and
\[
\mu_N = \sqrt{1-\frac{\pi^2}{a_N^2}} = 1 - \frac{\pi^2}{2\log^2N} +  \frac{3\pi^2\log\log N}{\log^3 N} + o\Big(\frac 1 {\log^3 N}\Big).
\]
and set
\[
M^N_\alpha(t) = \med_\alpha^N(\nu^N_t) - \mu_N t,
\]
where $\nu^N_t$ is the counting measure formed by the positions of the particles of $N$-BBM at time~$t$. Our main theorem is then the following:

\begin{theorem}
  \label{th:1}
Suppose that at time $0$ there are $N$ particles distributed independently according to the density proportional to $\sin(\pi x/a_N)e^{- x}\Ind_{(0,a_N)}(x)$. Then for every $\alpha\in(0,1)$, the finite-dimensional distributions of the process 
\[
\left(M^N_\alpha\big(t\log^3 N\big)\right)_{t\ge0}
\]
converge weakly as $N\to\infty$ to those of the L\'evy process $(L_t+x_\alpha)_{t\ge0}$ with $L_0 = 0$ and 
\begin{equation}
\label{eq:laplace_levy}
 \log E[e^{i\lambda L_1}] = i\lambda c + \pi^2\int_0^\infty e^{i\lambda x} - 1 - i\lambda x\Ind_{(x\le 1)}\,\Lambda(\dd x).
\end{equation}
Here, $\Lambda$ is the image of the measure $(x^{-2}\Ind_{(x> 0)})\dd x$ by the map $x\mapsto \log(1+x)$ and $c\in\R$ is a constant depending on the reproduction law $q(k)$.
\end{theorem}

For a motivation of why this process is interesting, we refer to the introductory chapter. In the next subsection, we present a detailed sketch of the ideas in the proof of Theorem~\ref{th:1}.

\subsection{Heuristic ideas and overview of the results}
\label{sec:heuristic}

In the introductory chapter of this thesis, we have outlined the semi-deterministic description of the $N$-BBM from \cite{Brunet2006}: Most of the time the system is in a meta-stable state, where the particles are approximately distributed according to the density proportional to $e^{- x} \sin(\pi  x/\log N)\Ind_{(0,\log N)}(x)$. From time to time (with a rate of order $\log^{-3}N$), a particle goes far to the right and reaches a point near $a_N$ defined above. This particle then spawns a large number of descendants (of the order of $N$), which leads to a shift of the front. Our proof of Theorem~\ref{th:1} is inspired by this description and by the article by Berestycki, Berestycki and Schweinsberg \cite{Berestycki2010}, who made some of these ideas rigorous. We briefly recall their results and arguments.

They consider BBM with absorption at the origin and with drift $-\mu_N$. Their starting point is to introduce a second barrier at the point $a_{N,A} = a_N -  A$ for some large positive constant $A$ and divide the particles at time $t$ into two parts; on the one hand those that have stayed inside the interval $(0,a_{N,A})$, on the other hand those that have hit the point $a_{N,A}$ before hitting $0$. This corresponds roughly to the division of the process into a deterministic and a stochastic part. Indeed, killing the particles at $a_{N,A}$ prevents the number of particles from growing fast and thus permits to calculate expectations and variances of various quantities. For example, if at time $0$ we have $N$ particles distributed according to the meta-stable density, then the variance of the number of particles at the time $\log^3 N$ is of order of $e^{-A}N^2$.
For large $A$, the particles inside the interval $(0,a_{N,A})$ therefore behave almost deterministically at the timescale $\log^3 N$. Moreover, the leading term in the Fourier expansion of the transition density of Brownian motion (with drift $-\mu_N$ and killed at the border of the interval $(0,a_{N,A})$) is proportional to $e^{-\mu_N x} \sin(\pi x/a_{N,A})$, which explains the meta-stable density predicted by the physicists.

As for the particles that hit $a_{N,A}$, the authors of \cite{Berestycki2010} find that 1) the number of descendants at a later time of such a particle is of the order of $e^{-A}NW$, where $W$ is a random variable with tail $P(W > x) \sim 1/x$, as $x\to\infty$ and 2) the rate at which particles hit the right barrier is of the order of $e^A/\log^3N$. Putting the pieces together, they then show that the process which counts the number of particles of the system converges in the $\log^3 N$ timescale to Neveu's continuous-state branching process\footnote{A continuous-state branching process (CSBP) $(Z_t)_{t\ge0}$ is a time-changed L\'evy process without negative jumps: at time $t$, time is sped up by the factor $Z_{t-}$. CSBPs are scaling limits of Galton--Watson processes and thus have an inherent notion of genealogy. Neveu's CSBP is the CSBP with L\'evy measure $x^{-2}\Ind_{(x>0)}\,\dd x$, whose genealogy is given by the Bolthausen--Sznitman coalescent \cite{Bertoin2000}. As a wise reader, you have read the 
introductory chapter, such that you know 
what this coalescent is.} and its genealogy to the Bolthausen--Sznitman coalescent.

Our proof of Theorem~\ref{th:1} builds upon the ideas of \cite{Berestycki2010} presented above. The basic idea is to approximate the $N$-BBM by a BBM with absorption at a \emph{random} barrier, which is chosen in such a way that \emph{it keeps the number of particles almost constant}. We call the resulting system the ``B-BBM'' (B stands for ``barrier''). The B-BBM takes two (large) parameters $a$ and $A$, which have similar purposes than $a_{N,A}$ and $A$ above. We then set $\mu = \sqrt{1 - \pi^2/a^2}$ and start with BBM with drift $-\mu$ and an absorbing barrier at the origin. Now, at the beginning, this barrier stays at the origin and does not move. When and only when a particle hits $a$ \emph{and} spawns a lot of descendants do we increase the drift to the left. This increase is in order to kill particles and thus make the population size stay almost constant. Note that moving the barrier to the right is an equivalent 
operation, but increasing the drift is technically more convenient. After the system has relaxed (which takes a time of order $a^2$) the drift is set to $-\mu$ again and the process is repeated.

As in \cite{Berestycki2010}, an important quantity is 
%
\[
 Z_t = \sum w_Z(X_u(t)),\quad\text{ where }w_Z(x) = ae^{\mu(x-a)}\sin \frac{\pi x}{a}.
\]
Here, we sum over all the particles $u$ alive at time $t$ and $X_u(t)$ denotes the position of the particle $u$ at time $t$. The process $(Z_t)_{t\ge 0}$ is important because it is a martingale for BBM with absorption at $0$ and $a$ and drift $-\mu$. Furthermore, it gives the approximate number of particles at a time $t \gg a^2$. More precisely, set $N = \lceil 2\pi  e^{A}a^{-3}e^{\mu a}\rceil$ and suppose we kill particles at $0$ and $a$. The expected number of particles at a time $t$, with $a^2 \ll t \ll a^3$, is then approximately $Ne^{-A}Z_0$. Moreover, the variance is of the order of $N^2e^{-2A} Z_0$. Therefore, if $Z_0 \approx e^A$ then the number of particles is concentrated around its expectation for large $A$.

When a particle hits the right barrier at the time $\tau$, say, we absorb its descendants at the space-time line $\mathscr L_\tau: x = a-y + (1-\mu)(t-\tau)$, where $y$ is a large constant depending on $A$ only (this idea comes from \cite{Berestycki2010}). In doing so, the number of particles absorbed at the barrier has the same law as in BBM with drift $-1$ and absorption at $-y$, starting from a single particle at the origin. Moreover, the time it takes for all the particles to be absorbed depends only on $y$, not on $a$. Now, if $x_1,\ldots,x_n$ are the positions of the absorbed particles and $Z' = \sum_i w_Z(x_i)$, then the number of descendants of this particle at a later time is of the order of $e^{-A} N Z'$, provided that the drift stays constant. Consequently, we say that a \emph{breakout} occurs, whenever $Z' > \ep e^A$, where $\ep$ will be chosen such that $\ep \ll 1/A$. 

In order to define a breakout properly, we classify the particles into \emph{tiers}. Particles that have never hit the point $a$ form the particles of tier 0. As soon as a particle hits $a$ (at the time $\tau$, say) it advances to tier 1. Its descendants then belong to tier 1 as well, but whenever a descendant hits $a$ and has an ancestor which has hit the line $\mathscr L_\tau$ after $\tau$, it advances to tier 2 and so on. Whenever a particle advances to the next tier, it has a chance to break out. We can then define \emph{the time $T$ of the first breakout} and will indeed show that $T$ is approximately exponentially distributed with rate proportional to $\ep^{-1}a^{-3}$. Interestingly, we will see that with high probability breakouts only occur from particles which are of tier $0$ or $1$. In fact, the number of breakouts occuring from particles of tier $1$ between the times $0$ and $a^3$ is approximately proportional to $A$ (and the remaining $\approx \ep^{-1}$ breakouts occur from particles of tier $0$).

After the breakout, we will then increase the drift to the left slightly, in order to kill more particles than usual (remember that increasing the drift to the left corresponds to moving the barrier to the right). For this, we first choose a family of increasing smooth functions $(f_x)_{x\ge 0}$ with $f_x(0) = 0$ and $f_x(+\infty) = x$ for all $x\ge 0$. Such a function will be called a \emph{barrier function}. The drift after a breakout is then set to $\mu_t = \mu + (\dd/\dd t) f_\Delta(t/a^2)$, where $\Delta$ is the total amount by which we have to move the barrier. Thus, the only randomness in the choice of the barrier is in its total shift, not in its shape. Looking at the definition of $Z_t$ and the fact that $Z_0\approx e^A$, one easily guesses that we have to choose $\Delta = \log(1+e^{-A}Z')$ in order to get $Z_t$ ultimately back to its initial value. 
This already explains the convergence of the barrier to the L\'evy process given by \eqref{eq:laplace_levy}: On the one hand, we have $Z' \approx \pi W$, where $W$ is the random variable mentioned above. This implies that the law of $e^{-A}Z'$ conditioned on $Z'>\ep e^A$ is approximately $\ep x^{-2}\Ind_{(x\ge\ep)}\,\dd x$ for large $A$ and $a$.\footnote{The statement ``for large $A$ and $a$'' means that we let first $a$, then $A$ go to infinity, see Section \ref{sec:before_breakout_definitions}.} On the other hand, we will show that breakouts occur at a rate proportional to $\ep^{-1}a^{-3}$. Together with the definition of $\Delta$, this explains the L\'evy measure $\Lambda(\dd x)$ in \eqref{eq:laplace_levy}. One easily checks that the cumulants of this L\'evy process coincide with \eqref{eq:statistics}.

As for the shape of the barrier, it is determined by the fact that we want the number of particles at each time to be approximately $N$. By first-moment estimates, it will become clear in Section~\ref{sec:interval_number} that the correct barrier function to choose is
\[
 f_\Delta(t) =  \log\Big(1+(e^{\Delta}-1)\pi^{-2}e^{\pi^2 t/2}\frac{\dd}{\dd t}\theta(1,t)\Big),
\]
where the theta function $\theta(x,t)$ is defined in \eqref{eq:theta}.

We remark that in order to study the B-BBM up to the time $T$ of the first breakout, we need to study it \emph{conditioned to break out at time $t$} for every $t\ge 0$. We will see in Section~\ref{sec:fugitive} that this will lead to a decomposition of the particles into 1) a \emph{fugitive}, which is conditioned to break out at $t$ and which will effectively be a \emph{spine}, and 2) the other particles conditioned \emph{not to break out} before $t$. This is essentially a Doob transform of the process, which we will introduce in Section~\ref{sec:doob}. Furthermore, we will see that the tier $1$ particles will have an essential role: At the timescale $a^3$ they will lead to an additional shift of the barrier by an amount of the order of $A$. This term will play the role of the linear compensation that is necessary in order to obtain in the limit the L\'evy process of infinite variation stated in Theorem~\ref{th:1}.

We now describe how we use the results on the B-BBM in order to prove Theorem~\ref{th:1}. Initially, our plan was to \emph{couple} the $N$-BBM and the B-BBM, i.e.\ construct them on the same probability space. We would then assign a colour to each particle: \emph{blue} to the particles which appear in the $N$-BBM but not in the B-BBM, \emph{red} to those that appear in the B-BBM but not in the $N$-BBM, and \emph{white} to the particles that appear in both processes. Our aim was then to show that the number of blue and red particles was negligible after a time of order $a^3$. This, unfortunately, did not work out, because we were not able to handle the intricate dependence between the red and blue particles.

Instead, we couple the $N$-BBM with two different processes, the \Bfl- and the \Bsh-BBM, which are variants of the B-BBM and which bound the position of the $N$-BBM in a certain sense from below and above, respectively. The \Bfl-BBM is defined as follows: Initially, all particles are coloured white and evolve as in the B-BBM. A white particle is coloured red as soon as it has $N$ or more white particles to its right. Children inherit the colour of their parent. After a breakout and the subsequent relaxation, all the red particles are killed immediately and the process restarts with the remaining particles. It is intuitive that the collection of white particles then bounds the $N$-BBM from below (in some sense) because we kill ``more'' particles than in the $N$-BBM. Indeed, in Section~\ref{sec:coupling}, we show by a coupling method that the empirical measure of the white particles in \Bfl-BBM is stochastically dominated by the one of the $N$-BBM with respect to the usual stochastic ordering of measures.
It then remains to show that the number of red particles in \Bfl-BBM is negligible when $A$ and $a$ are large. We do this through precise estimates on the number of particles in the interval $[r,\infty)$ for every $r\ge0$. These allow us to estimate the expected number of particles which turn red at the point $r$. It turns out that this expectation is small enough, which permits to conclude.

The definition of the \Bsh-BBM, which is used to bound the $N$-BBM from above, is more intricate than the one of \Bfl-BBM. Again, we colour all initial particles white and particles evolve as in B-BBM with the following change: Whenever a white particle hits $0$ \emph{and} has less than $N$ particles to its right, instead of killing it immediately, we colour it blue and let it survive for a time of order $a^2$. More precisely, we cut time into intervals $I_n = [t_n,t_{n+1})$, with $t_n = Kna^2$ for some large constant $K$. A particle which gets coloured blue during $I_n$ then survives until the time $t_{n+2}$. At this time, all of its descendants to the left of the origin are killed and the others survive. It will turn out that this system bounds the $N$-BBM from above with high probability and that the number of blue particles will remain negligible during a time of order $a^3$, as long as $A$ and $a$ are large.

We note that although our technique of bounding the $N$-BBM from below and from above works well for the \emph{position} of the particles, it does not give us information about the genealogy; the reason being that the coupling deforms the genealogical tree of the process. Thus, although it should not be difficult to show that the genealogy of the B-BBM (and of the \Bfl- and \Bsh-BBM) converges to the Bolthausen--Sznitman coalescent we do not know at present how one could transfer this information to the $N$-BBM.



\subsection{Notation guide}

This chapter is quite long and therefore also uses a lot of different notation. Below is a list of recurrent symbols, roughly in the order of their first appearance. Following this list are some further remarks about notational conventions.

\begin{longtable}{lp{10.3cm}r}
\textbf{Symbol} & \textbf{Meaning} & \textbf{Sect.}\\[4pt]
$q(k)$ & Reproduction law & \ref{sec:NBBM_intro}, \ref{sec:critical_line}\\
$m$ & $m = \sum_k (k-1) q(k)$ & \ref{sec:NBBM_intro}, \ref{sec:critical_line}\\
$\beta_0$ & Branching rate, $\beta_0 = 1/(2m)$ & \ref{sec:NBBM_intro}, \ref{sec:critical_line}\\
$\theta,\thbar$ & Theta functions & \ref{sec:theta}\\
$W^x$ & Law of Brownian motion started at $x$ & \ref{sec:killed_bm}\\
$p_t^a(x,y)$ & Transition density of Brownian motion killed outside $[0,a]$ & \ref{sec:killed_bm}\\
$E_t$ & Error term & \ref{sec:killed_bm}\\
$I^a(x,S),J^a(x,S)$ & Integrals related to Brownian motion killed outside $[0,a]$ & \ref{sec:killed_bm}\\
$W^x_{\mathrm{taboo}},W^{x,t,y}_{\mathrm{taboo}}$ & Law of Brownian taboo process and its bridge & \ref{sec:taboo}\\
$\N$, $\N^*$ & The set of natural numbers including and excluding 0, resp. & \ref{sec:preliminaries_definition}\\
$U$ & The space of individuals & \ref{sec:preliminaries_definition}\\
$\Xi$ & A realisation of a branching Markov process & \ref{sec:preliminaries_definition}\\
$\zeta_u,b_u,d_u$ & Lifetime, birth and death times of an individual $u$ & \ref{sec:preliminaries_definition}\\
$\mathscr N(t)$ & Set of individuals alive at time $t$ & \ref{sec:preliminaries_definition}\\
$X_u(t)$ & Position of the individual $u$ at the time $t$ & \ref{sec:preliminaries_definition}\\
$\P^x,\E^x,\P^\nu,\E^\nu$ & Law of and expectation w.r.t.\ a branching Markov process (later: BBM) started at $x$ or with particles distributed according to a finite counting measure $\nu$ & \ref{sec:preliminaries_definition}\\
$\mathscr F_t$ & $\sigma$-algebra with information up to time $t$ & \ref{sec:stopping_lines}\\
$\mathscr F_{\mathscr L}$ & $\sigma$-algebra with information up to the stopping line $\mathscr L$ & \ref{sec:stopping_lines}\\
$\mathscr L_T$ & Stopping line generated by a stopping time $T$ & \ref{sec:stopping_lines}\\
$L$ & A random variable with law $q(k)$ & \ref{sec:critical_line}\\
$m_2$ & $m_2 = E[L(L-1)]$ & \ref{sec:critical_line}\\
$W$ & Seneta-Heyde martingale limit of BBM with absorption & \ref{sec:critical_line}\\
$\mu$ & $\mu = \sqrt{1-\pi^2/a^2}$ & \ref{sec:interval_notation}\\
$\p_t(x,y)$ & Density of BBM killed outside an interval & \ref{sec:interval_notation}\\
$\P_f,\E_f$ & Law of and expectation w.r.t.\ BBM with varying drift & \ref{sec:interval_notation}\\
$w_Z(x)$ & $w_Z(x) = a\sin(\pi x/a)e^{\mu(x-a)}\Ind_{(x\in[0,a])}$ & \ref{sec:ZY}\\
$w_Y(x)$ & $w_Y(x) = e^{\mu(x-a)}\Ind_{(x\ge 0)}$ & \ref{sec:ZY}\\
$Z_t,Y_t$ & Sums of $w_Z(x),w_Y(x)$ over the positions of time $t$  particles & \ref{sec:ZY}\\
$N_t$ & Number of particles at time $t$ & \ref{sec:interval_number}\\
$R_t$ & Number of particles hitting $a$ up to time $t$ & \ref{sec:right_border}\\
$\widetilde\P,\widetilde\E$ & Law/expectation of BBM weakly conditioned not to hit $a$ & \ref{sec:weakly_conditioned}\\
$A,a,\ep,\eta,y,\zeta$ & Parameters of BBM before a breakout and B-BBM & \ref{sec:before_breakout_definitions}\\
$\mathscr N_t^{(l)}$ & Stopping line of tier $l$ particles at time $t$ & \ref{sec:before_breakout_definitions}\\
$\mathscr S_t^{(l)}$ & Stopping line of tier $l$ particles hitting the critical line before $t$ & 
\ref{sec:before_breakout_definitions}\\
$\mathscr R_t^{(l)}$ & Stopping line of tier $l$ particles hitting $a$ before $t$ & 
\ref{sec:before_breakout_definitions}\\
$\mathscr S^{(u,t)}$ & Stopping line of descendants of $(u,t)$ hitting the critical line & 
\ref{sec:before_breakout_definitions}\\
$Z^{(l)}_t$ & Sum of $w_Z(x)$ over particles from $\mathscr N_t^{(l)}$ & \ref{sec:before_breakout_definitions}\\
$Z_t,Y_t,R_t,N_t$ & $Z_t = Z^{(0+)}_t$ (from Section~\ref{sec:before_breakout} on). Same for $Y_t$, $R_t$ and $N_t$ & \ref{sec:before_breakout_definitions}\\
$p_B$ & Probability of a breakout & \ref{sec:before_breakout_definitions}\\
$\Q^a$ & $\P^a$ conditioned not to break out & \ref{sec:before_breakout_definitions}\\
$T^{(l)}$ & Time of first breakout from tier $l$ & \ref{sec:before_breakout_definitions}\\
$\Phat,\Ehat$ (also $\Phat_\lambda$,$\Ehat_\lambda$) & Law/expectation of BBM conditioned not to break out before some fixed time & \ref{sec:before_breakout_particles}\\
$Z^{(l)}_\emptyset$ & Sum of $w_Z(x)$ over particles from $\mathscr S_t^{(l)}$ & \ref{sec:before_breakout_particles}\\
$\mathscr U$ & The fugitive (the particle which breaks out) & \ref{sec:fugitive}\\
$\widebar{\mathscr N}_t,\widecheck{\mathscr N}_t$ & The descendants of the fugitive at time $t$, whose most recent common ancestor with the fugitive is/is not a regular particle & \ref{sec:fugitive}\\
$\widebar Z_t$, $\widecheck Z_t$ & Sum of $w_Z(x)$ over particles from $\widebar{\mathscr N}_t$ and $\widecheck{\mathscr N}_t$, respectively & \ref{sec:fugitive}\\
$\widebar Z_t$, $\widecheck Z_t$ & Value of $Z_t$ restricted to particles from $\widebar{\mathscr N}_t$ and $\widecheck{\mathscr N}_t$, resp. & \ref{sec:fugitive}\\
$\widehat Z_t$ & Value of $Z_t$ restricted to particles not related to the fugitive & \ref{sec:BBBM}\\
$\Delta$ & The total shift of the barrier after a breakout & \ref{sec:BBBM_definition}\\
$X^{[n]}_t$ & The position of the first $n$ pieces of the barrier at time $t$ (barrier process) & \ref{sec:BBBM_definition}\\
$\Theta_n$ & Beginning of the $n+1$-th piece of B-BBM (i.e.\ time of the $n$-th breakout plus relaxation time) & \ref{sec:BBBM_definition}\\
$G_n$ & ``Good event'' related to the first $n$ pieces of B-BBM & \ref{sec:BBBM_definition}\\
$\gamma_0$ & $\gamma_0 = (\pi p_Be^A)^{-1}$.& \ref{sec:BBBM_definition}\\
\parbox{2cm}{$G_{\mathscr U},G_{\mathrm{fug}},\widehat G$,\\$\widecheck G,G_{\Delta},G_{\mathrm{nbab}}$} & Several good sets related to a piece of B-BBM & \ref{sec:piece_proof}
\end{longtable}

From Section~\ref{sec:interval} on, the symbol $C$ stands for a \emph{positive} constant, which may only depend on the reproduction law $q$ and the value of which may change from line to line. Furthermore, if $X$ is any mathematical expression, then the symbol $O(X)$ stands for a possibly random term whose absolute value is bounded by $C|X|$.

In Section~\ref{sec:before_breakout}, we introduce two parameters $A$ and $a$ and will first let $a$ then $A$ go to infinity. This will be expressed by the statements ``for large $A$ and $a$ we have\ldots'' or ``as $A$ and $a$ go to infinity\ldots'' (see Section~\ref{sec:before_breakout_definitions} for a precise definition). These phrases will become so common that in Sections~\ref{sec:BBBM} to~\ref{sec:Bsharp} they will often be used implicitly, although they will always be explicitly stated in the theorems, propositions, lemmas etc. Section~\ref{sec:before_breakout} furthermore introduces the notation $o(1)$, which stands for a (non-random) term that only depends on the reproduction law $q$ and the parameters $A$, $a$, $\ep$, $\eta$, $y$ and $\zeta$ and which goes to $0$ as $A$ and $a$ go to infinity.

Sections \ref{sec:Bflat} and \ref{sec:Bsharp} each use special notation which only appears in those sections. This notation is defined at the beginning of both sections. Moreover, in both sections, we sometimes denote quantities which refer to descendants of the fugitive \emph{after} a breakout by the superscript ``fug''.

\section{Brownian motion in an interval}
\label{sec:bm}

In this section, we recall some explicit formulae concerning real-valued Brownian motion killed upon exiting an interval. These formulae naturally involve Jacobi theta functions, since these are fundamental solutions of the heat equation with periodic boundary conditions. We will therefore first review their definition and some of their properties.

\subsection{A function of Jacobi theta-type}
\label{sec:theta}
We define for $x\in\R$ and $t>0$ the following function of Jacobi theta-type:
\begin{equation}
\label{eq:theta}
\theta(x,t) = \frac 1 2 + \sum_{n=1}^\infty e^{- \tfrac{\pi^2}{2} n^2 t} \cos(\pi n x).
\end{equation}

The definition \eqref{eq:theta} is a representation of $\theta$ as a Fourier series, which is particularly well suited for large $t$, but which does not reveal its behaviour as $t\to 0$. This is where the following representation comes in, which is related to \eqref{eq:theta} by the Poisson summation formula (see \cite{Bellman}, §9):
\begin{equation}
\label{eq:theta_gaussian}
\theta(x,t) = \sum_{n\in\Z} \frac{1}{\sqrt{2\pi t}} \exp\Big(-\frac{(x-2n)^2}{2t}\Big).
\end{equation}
One recognises immediately that for real $x$ and $t$, $\theta(x,t)$ is the probability density at time $t$ of Brownian motion on the circle $\R/2\Z$ started at $0$. In other words, $\theta(x,t)$ is the unique solution to the PDE
\begin{equation}
\label{eq:pde_theta}
\begin{cases}
\frac{\partial}{\partial t} u(x,t) = \frac 1 2 \left(\frac{\partial}{\partial x}\right)^2 u(x,t) & \text{(PDE)}\\
u(x,t) = u(x+2,t) & \text{(BC)}\\
u(x,0+) = \sum_{n\in\Z} \delta(x-2n) & \text{(IC)},
\end{cases}
\end{equation}
where $\delta(x)$ denotes the Dirac delta-function. This is the heat equation with periodic boundary condition and the Dirac comb as initial condition. Note that (PDE) and (BC) also follow directly from \eqref{eq:theta}.

We define
\[
 \thbar(t) = \frac{2}{\pi^2} e^{\frac{\pi^2}{2}t}\frac{\dd}{\dd t} \theta(1,t),
\]
which is a smooth function on $\R_+$. By \eqref{eq:theta} and \eqref{eq:theta_gaussian}, one can show that $\thbar$ is stricly increasing\footnote{More precisely, by elementary computations, \eqref{eq:theta} gives $t_0\in\R_+$, such that $\thbar$ is strictly increasing on $(t_0,\infty)$ and \eqref{eq:theta_gaussian} gives $t_1 > t_0$, such that $\thbar$ is strictly increasing on $[0,t_1)$.} with $\thbar(0) = 0$ and $\thbar(+\infty) = 1$.

\subsection{Brownian motion killed upon exiting an interval}
\label{sec:killed_bm}

Various quantities of Brownian motion killed upon exiting an interval can be expressed by theta functions. For $x\in\R$, let $W^x$ be the law of Brownian motion started at $x$, let $(X_t)_{t\ge0}$ be the canonical process and let $H_y = \inf\{t\ge0: X_t = y\}$. For $a>0$ and $x\in[0,a]$, denote by $W^x_{\text{killed},a}$ the law of Brownian motion started at $x$ and killed upon leaving the interval $(0,a)$. Let $p^a_t(x,y)$ be its transition density, i.e.
\begin{equation}
\label{eq:def_p}
p^a_t(x,y) \dd y = W^x_{\text{killed},a}(X_t\in\dd y) = W^x(X_t \in \dd y,\
H_0\wedge H_a > t), \quad x,y\in [0,a].
\end{equation}
Then $p^a_t(x,y)$ is the fundamental solution to the heat equation (PDE) with boundary condition
\[u(0,t) = u(a,t) = 0,\quad t\ge0.\]
Hence (see also \cite{ItoMcKean}, Problem 1.7.8 or \cite{BorodinSalminen}, formula 1.1.15.8),
\begin{equation}
\label{eq:p_theta}
p^a_t(x,y) = a^{-1}\left(\theta\Big(\frac{x-y}{a},\frac t {a^2}\Big) - \theta\Big(\frac{x+y}{a},\frac t {a^2}\Big)\right).
\end{equation}
Equation \eqref{eq:theta} then yields
\begin{equation}
\label{eq:p_sin}
p^a_t(x,y) = \frac{2}{a} \sum_{n=1}^\infty e^{-\frac{\pi^2}{2a^2} n^2 t} \sin(\pi n \tfrac x a) \sin(\pi n \tfrac y a).
\end{equation}
This representation is particularly useful for large $t$: Define 
\begin{equation}
\label{eq:def_E}
 E_t = \pi^2 \sum_{n=2}^{\infty} n^2e^{-\pi^2(n^2-1) t/2}.
\end{equation}
By \eqref{eq:p_sin} and the inequality $|\sin nx| \le n \sin x$, $x\in [0,\pi]$, one sees that
\begin{equation}
 \label{eq:p_estimate}
 \left|e^{\frac{\pi^2}{2a^2} t}p^a_t(x,y) - \frac 2 a \sin(\pi x/a) \sin(\pi y/a)\right| \le E_{t/a^2}\frac 2 a \sin(\pi x/a) \sin(\pi y/a).
\end{equation}
Note that the Green function (see e.g.\ \cite{Kallenberg1997}, Lemma 20.10, p379) is given by
\begin{equation}
\label{eq:p_potential}
 \int_0^\infty p^a_t(x,y)\,\dd t = W^x\Big(\int_0^{H_0\wedge H_a} \Ind_{(X_t\in\dd y)}\,\dd t\Big)/\dd y = 2 a^{-1} (x\wedge y)(a-x\vee y).
\end{equation}
Set $H=H_0\wedge H_a$ and define
\begin{equation}
\label{eq:def_r}
r^a_t(x) = W^x(H \in \dd t,\ X_H = a)/\dd t.
\end{equation}
Then (see \cite{BorodinSalminen}, formula 1.3.0.6),
\begin{equation}
\label{eq:r_theta}
r^a_t(x) = \frac{1}{a^2} \theta'\left(\frac{x}{a}-1,\frac t {a^2}\right),
\end{equation}
where $\theta'$ denotes the derivative of $\theta$ with respect to $x$.

The following two integrals are going to appear several times throughout the article, which is why we give some useful estimates here. For a measurable subset $S\subset \R$, define
\begin{equation}
 \label{eq:Ia} 
  I^a(x,S) = W^x\Big(e^{\frac{\pi^2}{2a^2}H_a} \Ind_{(H_0 > H_a\in S)}\Big) = \int_{S\cap (0,\infty)} e^{\frac{\pi^2}{2a^2} s} r_s^a(x)\dd s,
\end{equation}
and
\begin{equation}
 \label{eq:Ja}
  J^a(x,y,S) = \int_{S\cap (0,\infty)} e^{\frac{\pi^2}{2a^2} s} p_s^a(x,y)\dd s,
\end{equation}
which satisfy the scaling relations 
\begin{equation}
 \label{eq:IJ_scaling}
 I^a(x,S) = I\Big(\frac x a,\frac S {a^2}\Big),\quad J^a(x,y,S) = a J\Big(\frac x a,\frac y a,\frac S {a^2}\Big),
\end{equation}
with $I = I^1$ and $J=J^1$. The following lemma provides estimates on $I(x,S)$ and $J(x,y,S)$.

\begin{lemma}
 \label{lem:I_estimates}
There exists a universal constant $C$, such that for every $x\in[0,1]$ and every measurable $S\subset \R_+$, we have
\[
\begin{split}
 |I(x,S) - \pi \lambda(S)\sin(\pi x)| &\le C\Big(x \wedge E_{\inf S} (1\wedge \lambda(S))\sin(\pi x)\Big),\quad\tand\\
|J(x,y,S) - 2 \lambda(S)\sin(\pi x)\sin(\pi y)| &\le C\Big([(x\wedge y)(1-(x\vee y))] \wedge E_{\inf S} \sin(\pi x) \sin(\pi y)\Big),
\end{split}
\]
where $\lambda(S)$ denotes the Lebesgue measure of $S$ and $E_{\inf S}$ is defined in \eqref{eq:def_E}.
\end{lemma}
\begin{proof}
First note that $I(x,\cdot)$ is a positive measure on $\R_+$ for every $x\in[0,1]$, such that we have by \eqref{eq:Ia}, 
\[
 0\le I(x,S\cap[0,1]) \le I(x,[0,1]) \le W^x\Big(e^{\frac{\pi^2}{2}H_1} \Ind_{(H_0 > H_1 \in [0,1])}\Big) \le xe^{\frac{\pi^2}{2}},
\]
since the scale function of Brownian motion is $s(x) = W^x(H_1 < H_0) = x$. Furthermore, decomposing $I(x,S)$ into
\[
  I(x,S) = I(x,S\cap[0,1]) + I(x,S\cap (1,\infty)),
\]
it is enough to prove that $|I(x,S) - \pi \lambda(S)\sin(\pi x)| \le C (1\wedge \lambda(S))E_{\inf S} \sin(\pi x)$ for all $S$ with $\inf S > 0$. Now, by \eqref{eq:r_theta} and \eqref{eq:theta},
\[
\begin{split}
 I(x,S) &= \int_S e^{\frac{\pi^2}{2} s} \theta'(x,s) \dd s = \pi \int_S \sum_{n=1}^\infty e^{-\frac{\pi^2}{2}(n^2-1)s} (-1)^{n-1} n \sin(\pi n x)\dd s\\
  &= \pi \lambda(S) \sin(\pi x) + \pi \sum_{n=2}^\infty\Big(\int_S e^{-\frac{\pi^2}{2}(n^2-1)s} \dd s\Big) n(-1)^{n-1} \sin(\pi n x),
\end{split}
\]
where the exchange of integral and sum is justified by the uniform convergence of the sum for $s\ge \inf S$.
We now have for $n\ge 2$,
\[
  \int_S e^{-\frac{\pi^2}{2}(n^2-1)s}\dd s \le \int_{\inf S}^\infty e^{-\frac{\pi^2}{2}(n^2-1)s}\dd s = \frac{2}{\pi^2(n^2-1)} e^{-\frac{\pi^2}{2}(n^2-1)\inf S},
\]
as well as
\[
 \int_S e^{-\frac{\pi^2}{2}(n^2-1)s}\dd s \le \lambda(S)e^{-\frac{\pi^2}{2}(n^2-1)\inf S}
\]
Furthermore, we have for $n\ge 2$,
\[
 |n(-1)^{n-1} \sin(\pi n x)| \le n^2 \sin(\pi x) \le 2 (n^2 - 1) \sin(\pi x).
\]
It follows that
\[
 |I(x,S) - \pi \lambda(S)\sin(\pi x)| \le (\frac{4}{\pi} \wedge \pi \lambda(S)) E_{\inf S} \sin(\pi x).
\]
This proves the statement about $I$. The proof of the statement about $J$ is similar, drawing on \eqref{eq:p_sin} instead and on the following estimate:
\[
 J(x,y,[0,1]) = \int_0^1 e^{\frac{\pi^2}{2} t} p_t(x,y)\,\dd t \le e^{\frac{\pi^2}{2}} \int_0^\infty p_t(x,y)\,\dd t = e^{\frac{\pi^2}{2}} (x\wedge y)(1-(x\vee y)),
\]
by \eqref{eq:p_potential}.
\end{proof}

\subsection{The Brownian taboo process}
\label{sec:taboo}
The Markov process on $(0,a)$ with infinitesimal generator
\[\frac 1 2 \left(\frac \dd {\dd x}\right)^2 + \frac \pi a \cot \frac {\pi x} a \frac \dd {\dd x}\]
is called the \emph{Brownian taboo process} on $(0,a)$. It is a diffusion with scale function $s(x)$ and speed measure $m(\dd x)$, where
\begin{equation}
 \label{eq:scale_speed}
s(x) = \frac \pi a \cot \frac {\pi x} a \quad\quad\tand\quad\quad m(\dd x) = \frac {2 a^2} {\pi^2} \sin^2\frac {\pi x} a\ \dd x.
\end{equation}
The singular points $0$ and $a$ are therefore entrance-not-exit. For $x\in [0,a]$ we denote the law of the Brownian taboo process on $(0,a)$ started from $x$ by $W^x_{\text{taboo},a}$. Often we will drop the $a$ if its value is clear from the context.

The name of this process was coined by F.\ Knight \cite{Knight1969} who showed that it can be interpreted as Brownian motion conditioned to stay inside the interval $(0,a)$ (hence, $0$ and $a$ are \emph{taboo states}). The Brownian taboo process is also known as the \emph{three-dimensional Legendre process}, because of its relation to Brownian motion on the 3-sphere (see \cite{ItoMcKean}, p270). The 3-dimensional Bessel process is obtained by taking the limit in law as $a\to\infty$. Note that the normalisation of the scale function and speed measure in \eqref{eq:scale_speed} was chosen in such a way that they converge, respectively, to the scale function and speed measure of the 3-dimensional Bessel process, as $a\to\infty$.

Below we list some useful properties of the Brownian taboo process:
\begin{enumerate}[nolistsep]
 \item It satisfies the following \emph{scaling relation}: If $X_t$ is a Brownian taboo process on $(0,1)$, then $a X_{t/a^2}$ is a Brownian taboo process on $(0,a)$.
 \item It is the Doob transform of Brownian motion killed at $0$ and $a$, with respect to the space-time harmonic function $h(x,t) = \sin(\pi x/a) \exp(\pi^2t/(2a^2))$. In other words, for $x\in(0,a)$, $W^x_{\text{taboo}}$ is obtained from $W^x_{\text{killed}}$ by a Cameron--Martin--Girsanov change of measure with the martingale
\[
 Z_t = \left(\sin\frac {\pi x} a\right)^{-1} \sin\frac {\pi X_t} a \ \exp\frac {\pi^2}{2a^2} t.
\]
 \item As a consequence, its transition probabilities are given by
\begin{equation}
\label{eq:taboo_transition}
 p^{\text{taboo}(0,a)}_t(x,y) = W^x_{\text{taboo},a}\left(X_t \in \dd y\right)/\dd y = \frac {\sin (\pi y/a)} {\sin (\pi x/a)}e^{\frac {\pi^2}{2a^2} t}\, p^a_t(x,y).
\end{equation}
Equation \eqref{eq:p_estimate} now implies that
\begin{equation}
 \label{eq:ptaboo_estimate}
 p^{\text{taboo}(0,a)}_t(x,y) = \frac 2 a \sin^2(\pi y/a)(1+ O(1) E_{t/a^2}),\quad \text{for all }x,y\in[0,a],
\end{equation}
 \item As can be seen from above or directly, it admits the stationary probability measure
\[
(m(0,a))^{-1} m(\dd x) = 2/a \sin^2(\pi x/a)\ \dd x. 
\]
 \item If $W^{x,t,y}_{\text{taboo}}$ denotes the taboo bridge from $x$ to $y$ of length $t$, then $W^{x,t,y}_{\text{taboo}} = W^{x,t,y}_{\text{killed}}$ by the second property.
 \item As a consequence, the taboo process is self-dual in the sense that for a measurable functional $F$ and $t>0$, we have
\[
 W^{x,t,y}_{\text{taboo}}[F((X_s;0\le s\le t))] = W^{y,t,x}_{\text{taboo}}[F((X_{t-s};0\le s\le t))].
\]
\end{enumerate}

The following lemma will be needed in Sections \ref{sec:before_breakout} and \ref{sec:BBBM}.

\begin{lemma}
\label{lem:k_integral}
Let $c>0$ and define $k(x)=e^{-c x}$. There exists a constant $C$, depending only on $c$, such that we have for every $x,y\in[0,a]$,
\begin{equation}
 \label{eq:taboo_integral_error_1}
 W^{x}_{\mathrm{taboo}}\Big[\int_0^t k(X_s)\,\dd s\Big] \le C \Big(t/a^3 + \operatorname{err}(x)\Big),
\end{equation}
and for $t\ge a^2$, 
\begin{equation}
 \label{eq:taboo_integral_error_2}
 W^{x,t,y}_{\mathrm{taboo}}\Big[\int_0^t k(X_s)\,\dd s\Big] \le C \Big(t/a^3 + \operatorname{err}(x) + \operatorname{err}(y)\Big),
\end{equation}
with $\operatorname{err}(z) = 1\wedge z^{-1}$. If $t\le a^2$, we still have,
\begin{equation}
 \label{eq:taboo_integral_error_3}
 W^{x,t,y}_{\mathrm{taboo}}\Big[\int_0^t k(X_s)\,\dd s\Big] \le C.
\end{equation}
\end{lemma}

Note: one could probably show by induction that the law of the integral is dominated by an exponential law, by showing that its moments are bounded by $C^n n!$.
\begin{proof}
We first show that \eqref{eq:taboo_integral_error_1} implies \eqref{eq:taboo_integral_error_2}. By the self-duality of the taboo process, we have
\[
 W^{x,t,y}_{\mathrm{taboo}}\Big[\int_0^t k(X_s)\,\dd s\Big] = W^{x,t,y}_{\mathrm{taboo}}\Big[\int_0^{t/2} k(X_s)\,\dd s\Big] + W^{y,t,x}_{\mathrm{taboo}}\Big[\int_0^{t/2} k(X_s)\,\dd s\Big].
\]
It therefore remains to prove that
\[
 E(x,y) = W^{x,t,y}_{\mathrm{taboo}}\Big[\int_0^{t/2} k(X_s)\,\dd s\Big] \le C(t/a^3 + \operatorname{err}(x)).
\]
Conditioning on $\sigma(X_s; 0\le t\le t/2)$, this integral equals
\[
 E(x,y) = W^x_{\mathrm{taboo}}\Big[\frac{p^{\mathrm{taboo}}_{t/2}(X_{t/2},y)}{p^{\mathrm{taboo}}_t(x,y)}\int_0^{t/2} k(X_s)\,\dd s\Big].
\]
By \eqref{eq:ptaboo_estimate}, there exists a universal constant $C$, such that for $t\ge a^2$,
\[
 E(x,y) \le C\, W^x_{\mathrm{taboo}}\Big[\int_0^{t/2} k(X_s)\,\dd s\Big].
\]
Equation \eqref{eq:taboo_integral_error_1} therefore implies \eqref{eq:taboo_integral_error_2}.

Heuristically, one can estimate the left side of \eqref{eq:taboo_integral_error_1} in the following way: Since $k(x)$ is decreasing very fast, only the times at which $X_s$ is of order $1$ contribute to the integral. When started from the stationary distribution, the process takes a time of order $a^3$ to reach a point at distance $O(1)$ from $0$ \cite{Lambert2000} and it stays there for a time of order $1$, hence the integral is of order $t/a^3$. When started from the point $x$, an additional error is added, which is of order $1$, when $x$ is at distance of order $1$ away from $0$. Adding both terms gives the bound appearing in the statement of the lemma.

The exact calculations are most easily performed in the following way. Let $Y$ be a random variable with values in $(0,a)$ distributed according to $\widetilde{m}(\dd x) := 2/a \sin^2(\pi x/a)\,\dd x$, which is the stationary probability measure of the taboo process. Let $H_Y = \inf\{t>0: X_s = Y\}$. We then have
\[
\begin{split}
 W^x_{\mathrm{taboo}}\Big[\int_0^t k(X_s)\,\dd s\Big] &= W^x_{\mathrm{taboo}}\Big[\int_0^{H_Y} k(X_s)\,\dd s+\int_{H_Y}^t k(X_s)\,\dd s\Big]\\
&\le W^x_{\mathrm{taboo}}\Big[\int_0^{H_Y} k(X_s)\,\dd s\Big] + W^{\widetilde m}_{\mathrm{taboo}}\Big[\int_0^t k(X_s)\,\dd s\Big]\\
& =: I_1+I_2.
\end{split}
\]
By the inequality $\sin x \le x$ for $x\ge 0$,
\[
 I_2 = t \int_0^a \widetilde{m}(\dd y)\, k(y)\,\dd y \le 2 \pi^2 t/a^3 \int_0^\infty e^{-c y}(1+y) y^2\,\dd y \le CT/a^3,
\]
for some constant $C$ depending only on $c$.

Recall the definition of scale function and speed measure in \eqref{eq:scale_speed}. Define the Green functions
\[
 G_{y,a}(x,z) = s(x\wedge z) - s(y)\quad\tand\quad G_{0,y}(x,z) = s(y) - s(x\vee z).
\]
We then have (see e.g.\ \cite{Revuz1999}, Chapter 3, Corollary 3.8),
\begin{equation}
 \label{eq:proof_k_integral_1}
\begin{split}
 I_1 &= \int_0^x \widetilde{m}(\dd y) \int_y^a m(\dd z) G_{y,a}(x,z) k(z) + \int_x^a \widetilde{m}(\dd y) \int_0^y m(\dd z) G_{0,y}(x,z) k(z)\\
 &=: I_{11}+I_{12}.
\end{split}
\end{equation}
By Fubini's theorem, the first term in \eqref{eq:proof_k_integral_1} is easily bounded by
\[
 I_{11} \le \int_0^a m(\dd z) k(z) \int_0^z \widetilde{m}(\dd y) [s(z) - s(y)],
\]
and noticing that $\operatorname{sign}(s(z)) = -\Ind_{(z<a/2)} + \Ind_{(z>a/2)}$, we get
\[
\begin{split}
 I_{11} &\le \int_{a/2}^a m(\dd z) s(z) k(z)\int_0^z \widetilde{m}(\dd y) + \int_0^a m(\dd z) k(z) \int_0^z \widetilde{m}(\dd y) (-s(y))\\
 &\le C/a^3\Big(\int_{a/2}^a z^4 k(z)\,\dd z + \int_0^a z^4 k(z)\,\dd z\Big)\\
 &\le C/a^3,
\end{split}
\]
where again we made use of the inequality $\sin x \le x$ for $x\ge 0$.

For the term $I_{12}$ a little bit more care is needed. Using the fact that $\int_x^a \widetilde{m}(\dd y) \le 1$, we have
\[
\begin{split}
 I_{12} &\le \int_x^a \widetilde{m}(\dd y) s(y) \int_0^y m(\dd z) k(z) + (-s(x)\vee 0) \int_0^x m(\dd z) k(z)+ \int_x^{a/2} m(\dd z) |s(z)| k(z)\\
&=: I_{121} + I_{122} + I_{123}.
\end{split}
\]
To estimate the first two terms, note that
\[
 \int_0^y m(\dd z) k(z) \le C(1\wedge y^3),\quad\tand\quad\int_x^a \widetilde{m}(\dd y) s(y) \le C/a.
\]
such that
\[
 I_{121} + I_{122} \le C\Big(1/a + (1\wedge x^3)(-s(x)\vee 0)\Big) \le C\Big(1/a + (1\wedge x^{-1})\Big),
\]
because $-s(x) \le 1/x$ for $x\in[0,a]$. The third term is seen to be bounded by
\[
 I_{123} \le C \int_x^\infty z k(z)\,\dd z \le C (1+x)e^{-cx}.
\]
Altogether, we get
\[
 W^x_{\mathrm{taboo}}\Big[\int_0^t k(X_s)\,\dd s\Big] \le C\Big(t/a^3 + 1/a + \operatorname{err}(x)\Big).
\]
This proves \eqref{eq:taboo_integral_error_2} and therefore \eqref{eq:taboo_integral_error_1}.

In order to prove \eqref{eq:taboo_integral_error_3}, a different method is needed. We may assume that $x,y \le a/2$, otherwise we decompose the path at the first and/or last time it hits $a/2$ and bound the parts above $a/2$ trivially by $a^2e^{-ca/2} = O(1)$. The transition density of the taboo bridge can be written
\[
 W^{x,t,y}_{\mathrm{taboo}}\Big[X_s\in \dd z\Big] = \frac{p_s^a(x,z)p_{t-s}^a(z,y)}{p_t^a(x,y)}\,\dd z.
\]
If we denote by $p_t^0(x,y) = (2\pi t)^{-1/2} \exp(-(z^2+x^2)/2t)2\sinh(zx/t)$ the transition density of Brownian motion killed at $0$, then we have the trivial inequality $p_t^a(x,y) \le p_t^0(x,y)$ and furthermore $p_t^a(x,y) \ge C p_t^0(x,y)$ for $x,y\le a/2$ and $t\le a^2$ by Brownian scaling. It follows that
\[
 W^{x,t,y}_{\mathrm{taboo}}\Big[\int_0^t k(X_s)\,\dd s\Big] \le C R^{x,t,y}\Big[\int_0^t k(X_s)\,\dd s\Big],
\]
where $R^{x,t,y}$ denotes the law of the Bessel bridge of dimension 3. This Bessel bridge is the Doob transform of the Bessel process started at $x$ with respect to the space-time harmonic function $h_y(z,s) = p_{t-s}^0(z,y)/p_t^0(x,y)$. By the standard theory of Doob transforms, this is the Bessel process with additional drift
\[
 \frac{\dd}{\dd z}(\log h_y(z,s)) = -\frac{z^2}{t-s}+\frac{\dd}{\dd z}\log \sinh\frac{zy}{t-s} = -\frac{z^2}{t-s} + \frac{y}{t-s} \coth \frac{zy}{t-s}.
\]
This in an increasing function in $y$ and standard comparison theorems for diffusions (see e.g.\ \cite{Revuz1999}, Theorem IX.3.7) now yield that for $y_1\le y_2$, we have
\[
 R^{x,t,y_2}[k(X_s)] \le R^{x,t,y_1}[k(X_s)],
\]
since $k$ is a decreasing function. This is true in particular for $y_1 = 0$. Using the self-duality of the Bessel bridge, we can repeat the same reasoning with $x$. We thus have altogether
\[
 W^{x,t,y}_{\mathrm{taboo}}\Big[\int_0^t k(X_s)\,\dd s\Big] \le C R^{0,t,0}\Big[\int_0^t k(X_s)\,\dd s\Big],
\]
for any $x,y\le a/2$. This calculation can be done explicitly and yields \eqref{eq:taboo_integral_error_3}.
\end{proof}

\section{Preliminaries on branching Markov processes}
\label{sec:preliminaries}

In this section we recall some known results about branching Brownian motion and branching Markov processes in general.

\subsection{Definition and notation}
\label{sec:preliminaries_definition}

Branching Brownian motion can be formally defined using Neveu's \emph{marked trees} \cite{Neveu1986} as in \cite{Chauvin1988} and \cite{Chauvin1991}. We will follow this path here, but with slight differences, because we will need to consider more general branching Markov processes and the definition of branching Brownian motion in \cite{Chauvin1991} formally relied on the translational invariance of Brownian motion. 

We first define the space of Ulam--Harris labels, or \emph{individuals},
\[
 U=\{\emptyset\} \cup \bigcup_{n\ge 1} {\N^*}^n,
\]
where we use the notation $\N^* = \{1,2,3,\ldots\}$ and $\N = \{0\}\cup \N^*$. Hence, an element $u\in U$ is a \emph{word} over the alphabet $\N^*$, with $\emptyset$ being the empty word. For $u,w\in U$, we denote by $uw$ the \emph{concatenation} of $u$ and $w$. The space $U$ is endowed with the ordering relations $\preceq$ and $\prec$ defined by
\[
 u \preceq v \iff \exists w\in U: v = uw\quad\tand\quad u\prec v \iff u\preceq v \tand u \ne v.
\]
A \emph{tree} is by definition a subset $\mathfrak t\subset U$, such that 1) $\emptyset \in \mathfrak t$, 2) $u\in \mathfrak t$ and $v\prec u$ imply $v\in \mathfrak t$ and 3) for every $u\in\mathfrak t$ there is a number $k_u\in\N$, such that for all $j\in\N^*$, we have $uj\in \mathfrak t$ if and only if $j\le k_u$. Thus, $k_u$ is the number of children of the individual $u$. We denote the space of trees by $\mathscr T$ and endow it with the sigma-field $\mathscr A$ generated by the subsets $\mathscr T_u = \{\mathfrak t\in\mathscr T: u\in \mathfrak t\}$.

For a tree $\mathfrak t\in \mathscr T$ and $u\in \mathfrak t$, we define the \emph{subtree rooted at $u$} by
\[
 \mathfrak t^{(u)} = \{v\in U: uv \in \mathfrak t\}.
\]

Given a measurable space $\mathscr M$, a \emph{marked tree} (with space of marks $\mathscr M$) is a pair
\[
 \mathfrak t^{\mathscr M} = (\mathfrak t,(\eta_u;u\in\mathfrak t)),
\]
where $\mathfrak t\in \mathscr T$ and $\eta_u\in \mathscr M$ for all $u\in\mathfrak t$. The space of marked trees is denoted by $\mathscr T^\mathscr M$ and is endowed with the sigma-field $\mathscr A^{\mathscr M} = \pi^{-1}(\mathscr A)$, where $\pi: \mathscr T^\mathscr M\to \mathscr T$ is the canonical projection. Accordingly, we also define $\mathscr T^{\mathscr M}_u = \pi^{-1}(\mathscr T_u)$. The definition of a subtree extends as well to marked trees: For $u\in\mathfrak t$, we define
\[
 (\mathfrak t^{\mathscr M})^{(u)} = (\mathfrak t^{(u)},(\eta_{uv};v\in\mathfrak t^{(u)})).
\]

For our purposes, the space of marks $\mathscr M$ is always going to be a function space, namely, for a Polish space $\mathscr E$ and a cemetary symbol $\Delta \notin \mathscr E$, we define the Skorokhod space $D(\mathscr E)$ of functions $\Xi:[0,\infty)\to\mathscr E\cup \{\Delta\}$ which are right-continuous with left limits, with $\Xi(0)\ne \Delta$ and for which $\Xi(t) = \Delta$ implies $\Xi(s) = \Delta$ for all $s\ge t$. Then we define $\zeta(\Xi) = \inf\{t\ge 0: \Xi(t) = \Delta\}$. For an individual $u\in U$, its mark is denoted by $\Xi_u$ and we define $\zeta_u = \zeta(\Xi_u)$. The branching Markov process will then be defined on the space (we suppress the superscript $D(\mathscr E)$)
\[
 \Omega = \{\omega = (\mathfrak t,(\Xi_u;u\in\mathfrak t))\in\mathscr T^{D(\mathscr E)}: \forall u\in U\ \forall 1\le i\le k_u: \zeta_u < \infty \Rightarrow \Xi_u(\zeta_u-) = \Xi_{ui}(0)\},
\]
endowed with the sigma-field $\F = \Omega \cap \mathscr A^{D(\mathscr E)}$ generated by the sets $\Omega_u = \Omega \cap \mathscr T^{\mathscr M}_u$. We define for $u\in U$ the random variables
\[
 b_u = \sum_{v\prec u} \zeta_v,\quad d_u = b_u + \zeta_u = \sum_{v\preceq u} \zeta_v,
\]
which are the \emph{birth} and \emph{death} times of the individual $u$, respectively. We then define the set of individuals alive at time $t$ by 
\[
 \mathscr N(t) = \{u\in \mathfrak t: b_u\le t < d_u\}.
\]
The \emph{position} of $u$ at time $t$ is defined for $u\in \mathfrak t$ with $d_u > t$ by $X_u(t) = \Xi_v(t-b_v)$, where $v\in U$ is such that $v\in \mathscr N(t)$ and $v\preceq u$. If $d_u \le t$, then we set $X_u(t) = \Delta$.

Now suppose we are given a defective strong Markov process $\overline{X} = (\overline{X}_t)_{t\ge 0}$ on $\mathscr E$, with paths in $D(\mathscr E)$. The law of $\overline{X}$ started in $x\in \mathscr E$ will be denoted by $\overline{P}^x$. For simplicity, we will assume that for every $x\in \mathscr E$, we have $\zeta(\overline{X}) < \infty$, $\overline{P}^x$-almost surely. Furthermore, let $((q(x,k))_{k\in \N})_{x\in\mathscr E}$ be a family of probability measures on $\N$, measurable with respect to $x$. Then we define the branching Markov process with particle motion $\overline{X}$ and reproduction law $q$ as the (unique) family of probability measures $(\P^x)_{x\in\mathscr E}$ on $\Omega$ which satisfies
\begin{equation}
 \label{eq:bmp_renewal}
\P^x(\dd \omega) = \overline{P}^x(\dd X_\emptyset) q(X_\emptyset(\zeta_\emptyset-),k_\emptyset) \prod_{i=1}^{k_\emptyset} \P^{X_\emptyset(\zeta_\emptyset-)}(\dd \omega^{(i)}).
\end{equation}
Note that by looking at the space-time process $(X_t,t)_{t\ge 0}$, we can (and will) extend this definition to the time-inhomogeneous case.

\subsection{Stopping lines}
\label{sec:stopping_lines}

The analogue to stopping times for branching Markov processes are \emph{(optional) stopping lines}, for which several definitions exist. For branching Brownian motion, they have first been defined by Chauvin \cite{Chauvin1991}, the definition there is however too restricted for our purposes. Jagers \cite{Jagers1989} has given a definition of more general stopping lines for discrete-time branching processes; our definition will in fact mix both approaches. Note also that Biggins and Kyprianou \cite{Biggins2004} build up on Jagers' definition of stopping lines and define the subclasses of \emph{simple} and \emph{very simple} stopping lines (again for discrete-time processes). Chauvin's definition then corresponds to the class of very simple stopping lines.

We first define a \emph{(random) line} (called ``stopping line'' in \cite{Jagers1989}) to be a set $\ell = \ell(\omega) \subset U\times [0,\infty)$, such that
\begin{enumerate}[nolistsep]
 \item $u\in \mathscr N(t)$ for all $(u,t)\in \ell$ and
 \item $(u,t)\in\ell$ implies $(v,s)\notin \ell$ for all $v\preceq u$ and $s< t$.
\end{enumerate}
Note that a line is at most a countable set. For a pair $(u,t)\in U\times [0,\infty)$ and a line $\ell$, we write $\ell \preceq (u,t)$ if there exists $(v,s)\in \ell$, such that $v\preceq u$ and $s\le t$. For a subset $A\subset U\times [0,\infty)$, we write $\ell \preceq A$ if $\ell \preceq (u,t)$ for all $(u,t)\in A$. If $\ell_1$ and $\ell_2$ are two lines, we define the line $\ell_1 \wedge \ell_2$ to be the maximal line (with respect to $\preceq$), which is smaller than both lines.

We now define for each $u\in U$ two filtrations on $\Omega_u$ by
\begin{align*}
 \F_u(t) &= (\Omega_u \cap \sigma(\Xi_u(s);0\le s\le t-b_u)) \vee \bigvee_{v\prec u}(\Omega_v \cap \sigma(\Xi_v))\\ \F^{\mathrm{pre}}_u(t) &= (\Omega_u \cap \sigma(\Xi_u(s);0\le s\le t-b_u))\vee \bigvee_{v\nsucceq u} (\Omega_v \cap \sigma(\Xi_v)).
\end{align*}
Informally, $\F_u(t)$ contains the information on the path from $u$ to the root between the times $0$ and $t$, and $\F^{\mathrm{pre}}_u(t)$ contains this information and of all the other particles \emph{excluding} the descendants of $u$. In particular, we have $\F_u(t) \subset \F^{\mathrm{pre}}_u(t)$ The filtration $\F_u(t)$ is denoted by $\mathscr A_u(t)$ in Chauvin's paper \cite{Chauvin1991} and $\F^{\mathrm{pre}}_u(t)$ corresponds to the \emph{pre-$(u,t)$-sigma-algebra} as defined by Jagers \cite{Jagers1989}.

We can now define a \emph{stopping line} (``optional line'' in \cite{Jagers1989}) $\mathscr L$ to be a random line with the additional property
\begin{itemize}
 \item[3a.] $\forall (u,t)\in U\times [0,\infty):\{\omega\in \Omega_u:\mathscr L\preceq(u,t)\} \in \F^{\mathrm{pre}}_u(t).$
\end{itemize}
The sigma-algebra $\F_{\mathscr L}$ of the past of $\mathscr L$ is defined to be the set of events $E \in \F$, such that for all $(u,t)\in U\times [0,\infty)$,
\[
 E \cap \{\omega \in \Omega_u : \mathscr L\preceq(u,t)\} \in \F^{\mathrm{pre}}_u(t).
\]
For example, for any $t\ge0$, the set $\mathscr N(t)\times\{t\}$ is a stopping line. This permits us to define the filtration $(\F_t)_{t\ge0}$ by
\begin{equation*}
  \F_t = \F_{\mathscr N(t)\times\{t\}}.
\end{equation*}
Following Biggins and Kyprianou \cite{Biggins2004}, we now say that $\mathscr L$ is a \emph{simple stopping line} or a \emph{very simple stopping line}, if it satisfies the property 3b or 3c below, respectively:
\begin{itemize}
 \item[3b.] $\forall (u,t)\in U\times [0,\infty):\{\omega\in \Omega_u:\mathscr L\preceq(u,t)\} \in \Omega_u\cap\F_t.$
 \item[3c.] $\forall (u,t)\in U\times [0,\infty):\{\omega\in \Omega_u:\mathscr L\preceq(u,t)\} \in \F_u(t).$
\end{itemize}
Then $\mathscr N(t)\times\{t\}$ is obviously a very simple stopping line. Furthermore, if $T:D(\mathscr E)\to \R_+$ is a stopping time\footnote{In other words, $\{T(\Xi)\le t\} \in \sigma(\Xi(s);s\in[0,t])$ for every $t\ge 0$.}, then
\[
 \mathscr L_T = \{(u,t)\in U\times [0,\infty): u\in \mathscr N(t)\tand t = T(X_u)\}\}
\]
is a very simple stopping line as well. We recall that the definition of stopping lines in \cite{Chauvin1991} is equivalent to the definition of very simple stopping lines given here.

The first important property of stopping lines is the strong branching property. In order to state it, we define for $t\ge0$, $u\in \mathscr N(t)$,
\[
 \omega^{(u,t)} = (\mathfrak t^{(u)},(\Xi'_{uv};v\in\mathfrak t^{(u)})),
\]
with $\Xi'_u(\cdot) = \Xi_u(\cdot+t-b_u)$ and $\Xi'_{uv} = \Xi_{uv}$ for $v\in \mathfrak t^{(u)}\backslash \{\emptyset\}$. The \emph{strong branching property} (\cite[Proposition 2.1]{Chauvin1991}, \cite[Theorem 4.14]{Jagers1989}) then states that for every stopping line $\mathscr L$, conditioned on $\F_\mathscr L$, the subtrees $\omega^{(u,t)}$, for $(u,t)\in \mathscr L$, are independent with respective distributions $\P^{X_u(t)}$.

\subsection{Many-to-few lemmas and spines}
\label{sec:spine}

Another important tool in the theory of branching processes is the so-called Many-to-one lemma and its recently published extension, the Many-to-few lemma \cite{Harris2011} along with the \emph{spine decomposition} technique which comes along with it and has its origins in \cite{LPP1995}, although it appeared implicitly in the literature before that, see e.g.\ the references in the same paper. Here we state stopping line versions of these lemmas, which to the knowledge of the author have not yet been stated in this generality in the literature, although they belong to the common folklore. We will therefore only sketch how they can be derived from the existing literature.

We assume for simplicity that the strong Markov process $\overline{X}$ admits a representation as a conservative strong Markov process $X$ with paths in $D(\mathscr E)$, which is killed at a rate $R(x)$, where $R:\mathscr E\to [0,\infty)$ is measurable. The law of $X$ started at $x$ is denoted by $P^x$ and the time of killing by $\zeta$. Given a stopping time $T$ for $X$, we can then define a stopping time $\overline T$ for $\overline X$ by setting $\overline T = T$, if $T < \zeta$ and $\overline T = \infty$ otherwise. For simplicity, we write $\mathscr L_T$ for $\mathscr L_{\overline T}$. Finally, for every $x\in\mathscr E$, define $m(x) = \sum_{k\ge 0} (k-1) q(x,k)$, $m_1(x) = \sum_{k\ge 0} k q(x,k)$ and $m_2(x) = \sum_{k\ge 0} k(k-1) q(x,k)$.

We are now going to present the spine decomposition technique, following \cite{Hardy2006}. They assume that $q(x,0) \equiv 0$, but this restriction is actually not necessary, as noted in \cite{Harris2011}. Given a tree $\mathfrak t$, a \emph{spine} of $\mathfrak t$ is an element of the \emph{boundary} of $\mathfrak t$, i.e.\ it is a line of descent $\xi = (\xi_0 = \emptyset,\xi_1,\xi_2,\ldots)$ from the tree, which is finite if and only if the last element is a leaf of the tree. We augment our space $\Omega$ to the space $\Omega^*$ by
\[
 \Omega^* = \{(\omega,\xi): \omega\in\Omega,\ \xi \text{ is a spine of the tree underlying $\omega$}\}
\]
We are going to denote by $\xi_t$ the individual $u\in U$ that satisfies $u \in \mathscr N(t)$ and $u\in \xi$ if it exists, and $\xi_t = \emptyset$ otherwise. Instead of the redundant $X_{\xi_t}(t)$, we write $X_\xi(t)$. We also note that the definition of stopping lines can be extended to $\Omega^*$ by projection.

Now, for every $x\in \mathscr E$, one can define a probability measure $\P^{*,x}$ on $\Omega^*$ in the following way:
\begin{itemize}[nolistsep, label=$-$]
\item Initially, $X_\xi(0) = x$.
\item The individuals on the spine move according to the strong Markov process $X$ and die at the rate $m_1(y)R(y)$, when at the point $y\in\mathscr E$.
\item When an individual on the spine dies at the point $y\in\mathscr E$, it leaves $k$ offspring at the point where it has died, with probability $(m_1(x))^{-1} kq(x,\cdot)$ (this is also called the size-biased distribution of $q(x,\cdot)$\footnote{The size-biased distribution of the Dirac-mass at $0$ is again the Dirac-mass at $0$.}).
\item Amongst those offspring, the next individual on the spine is chosen uniformly. This individual repeats the behaviour of its parent (started at the point $y$).
\item The other offspring initiate independent branching Markov processes according to the law $\P^y$, independently of the spine.
\end{itemize}
This decomposition first appeared in \cite{Chauvin1988}. 

If we start with $n$ initial particles $1,\ldots,n$ at positions $x_1,\ldots,x_n$, we can extend this definition by defining a (non-probability) measure $\P^*$, which is the sum of $n$ probability measures $\P^*_1 + \cdots + \P^*_n$, where under $\P^*_i$, the particle $i$ follows the law $\P^{*,x_i}$ and the remaining particles $j\ne i$ follow the law $\P^{x_j}$.

We now have the important

\begin{lemma}[Many-to-one]
\label{lem:many_to_one}
 Let $\mathscr L$ be a simple stopping line. Define $T$ by $(\xi_T,T)\in \mathscr L$ if it exists, and $T=\infty$ otherwise. Let $Y$ be a random variable of the form
\[
 Y = \sum_{(u,t)\in \mathscr L} Y_u \Ind_{(u\in \xi)},
\]
where $Y_u$ an $\F_\mathscr L$-measurable random variable for every $u\in U$. Then
\begin{equation}
\E\Big[\sum_{(u,t)\in \mathscr L} Y_u\Big] = \E^*\Big[Ye^{\int_0^T R(X_\xi(t))m(X_\xi(t))\,\dd t}\Ind_{(T<\infty)}\Big] .
\end{equation}
\end{lemma}

Proofs of this result can be found for fixed time in \cite{Kyprianou2004}, \cite{Hardy2006} or \cite{Harris2011}. With simple stopping lines, it has been proven in the discrete setting \cite[Lemma~14.1]{Biggins2004} and their arguments can be used to adapt the proofs in \cite{Hardy2006} and \cite{Harris2011} to yield the result stated here.

Often, we will use a simpler version of the Many-to-one lemma, which is the following
\begin{lemma}[Simple Many-to-one]
\label{lem:many_to_one_simple}
 Let $T = T(X)$ be a stopping time for the strong Markov process $X$ which satisfies $P^x(T<\infty) = 1$ for every $x\in\mathscr E$. Let $f:\mathscr E\to [0,\infty)$ be measurable. Then we have
\[
 \E^x\Big[\sum_{(u,t)\in \mathscr L_T} f(X_u(t))\Big] = E^x\Big[e^{\int_0^T R(X_t) m(X_t)\,\dd t} f(X_T)\Big]
\]
\end{lemma}

The next lemma tells us about second moments of sums of the previous type. To state it, we define for a stopping time $T$ for $X$ the density of the branching Markov process before $\mathscr L_T$, by
\begin{equation}
\label{eq:def_density}
 p_T(x,\dd y,t) = \E^x\Big[\sum_{u\in \mathscr N(t)} \Ind_{(X_u(t)\in \dd y,\ t < T(X_u))}\Big].
\end{equation} 
\begin{lemma}
\label{lem:many_to_two}
 Let $H$ be the hitting time functional of a closed set $F\subset \mathscr E$ on $D(\mathscr E)$ which satisfies $P^x(H<\infty) = 1$ for every $x\in\mathscr E$. Let $f:\mathscr E\to [0,\infty)$ be measurable. Then we have
\begin{multline}
\label{eq:many_to_two}
 \E^x\Big[\Big(\sum_{(u,t)\in \mathscr L_H} f(X_u(t))\Big)^2\Big] = \E^x\Big[\sum_{(u,t)\in \mathscr L_H} (f(X_u(t)))^2\Big]\\
 + \int_0^\infty \int_{\mathscr E} p_H(x,\dd y,s) R(y)m_2(y) \Big(\E^y\Big[\sum_{(u,t)\in \mathscr L_H} f(X_u(t))\Big]\Big)^2\,\dd s
\end{multline}
\end{lemma}
\begin{remark}
  \label{rem:many_to_two}
This lemma can be proven using the Many-to-few lemma from \cite{Harris2011} (which is valid for stopping lines as well by the same argument as the one above) or with Lemma~\ref{lem:many_to_one}, by noting that
\[
 \Big(\sum_{(u,t)\in \mathscr L_H} f(X_u(t))\Big)^2 = \sum_{(u,t)\in \mathscr L_H} (f(X_u(t)))^2 + \sum_{(u,t)\in \mathscr L_H} \Big(f(X_u(t)) \sum_{(v,s)\in \mathscr L_H,\ v\ne u} f(X_v(s))\Big).
\]
It also has an intuitive explanation (see, e.g.\  the proof of Proposition 18 in \cite{Berestycki2010}): Write the above sum slightly differently as
\begin{equation}
  \label{eq:2nd_moment_decomposition}
 \Big(\sum_{(u,t)\in \mathscr L_H} f(X_u(t))\Big)^2 = \sum_{(u,t)\in \mathscr L_H} (f(X_u(t)))^2 + \sum_{(u,t),(v,s)\in \mathscr L_H,\,u\ne v} f(X_u(t))f(X_v(s)),
\end{equation}
Now, a particle at the point $y$ spawns an expected number of $R(y)m_2(y)\,\dd s$ of ordered pairs of particles during the time interval $[s,s+\dd s]$, and by the strong branching property, the two particles in a pair evolve independently. This yields \eqref{eq:many_to_two}. Note that this heuristic argument can be seen as a decomposition of the ordered pairs of particles in the second sum in \eqref{eq:2nd_moment_decomposition} according to their most recent common ancestor, an approach which can be made rigorous in a discrete setting.
\end{remark}

Taking for $X$ the space-time process $(Y_t,t)_{t\ge 0}$ of a possibly non-homogeneous strong Markov process $(Y_t)_{t\ge 0}$ with paths in $D(\mathscr E)$ and the closed set $F = \mathscr E\times \{t\}$, for some $t\ge 0$, we obtain the following useful corollary, which appeared already in \cite{Watanabe1967} and \cite{Sawyer1976} in the homogeneous case.
\begin{lemma}
\label{lem:many_to_two_fixedtime}
Let $f:\mathscr E\times \R_+\to [0,\infty)$ be measurable and let $t\ge 0$. Then we have
\begin{multline}
\label{eq:many_to_two_fixedtime}
 \E^{(x,0)}\Big[\Big(\sum_{u\in \mathscr N(t)} f(Y_u(t),t)\Big)^2\Big] = \E^{(x,0)}\Big[\sum_{u\in \mathscr N(t)} (f(Y_u(t),t))^2\Big]\\
 + \int_0^t \int_{\mathscr E} p(x,\dd y,s) R(y,s)m_2(y,s) \Big(\E^{(y,s)}\Big[\sum_{u\in \mathscr N(t)} f(Y_u(t),t)\Big]\Big)^2\,\dd s
\end{multline}
\end{lemma}

\subsection{Doob transforms}
\label{sec:doob}
As in the previous subsection, we assume for simplicity that the strong Markov process $\overline{X}$ admits a representation as a conservative strong Markov process $X$ with paths in $D(\mathscr E)$, which is killed at a rate $R(x)$, where $R:\mathscr E\to [0,\infty)$ is measurable. Let $H$ be the hitting time functional of a closed set $F\subset \mathscr E$ on $D(\mathscr E)$. Furthermore, let $h:F\to [0,1]$ be a measurable function. We extend the function $h(x)$ to $\mathscr E$ by setting
\[
 h(x) = \E^x\Big[\prod_{(u,t)\in \mathscr L_H} h(X_u(t))\Big],
\]
%
We are going to assume that $h(x) > 0$ for all $x\in \mathscr E\backslash F$. Then for all such $x$ we can define a law $\P_h^x$ on $\Omega$ by
\[
 \P_h^x(\dd \omega) = (h(x))^{-1} \prod_{(u,t)\in \mathscr L_H} h(X_u(t))\times \P^x(\dd \omega),
\]
where the multiplication is in the sense of a Radon--Nikodym derivative. Now define
\[
 Q(x) = \sum_{k\ge 0} q(x)h(x)^{k-1},\quad\tand\quad q_h(x,k) = \frac{q(x)h(x)^{k-1}}{Q(x)}.
\]
By \eqref{eq:bmp_renewal}, we now have (dropping the symbol $\emptyset$ for better reading and setting $H = H(X_\emptyset)$)
\begin{multline*}
 h(x) \P_h^x(\dd \omega) = \overline{P}^x(\dd X)
\Big(
    \Ind_{(H < \zeta)}
	h(X(H))q(X(\zeta-),k)
	\prod_{i=1}^{k} \P^{X(\zeta-)}(\dd \omega^{(i)})\\
    +\Ind_{(\zeta \le H)}
	h(X(\zeta-))^{k}q(X(\zeta-),k)
	\prod_{i=1}^{k} \P^{X(\zeta-)}(\dd \omega^{(i)})
	\prod_{(u,t)\in \mathscr L_H(\omega^{(i)})}h(X_u(t))
\Big).
\end{multline*}
If we denote by $X^H$ the process $X$ stopped at $H$, and the law of $X^H$ under $\overline{P}^x$ by $(\overline{P}^x)^H$, then the last equation and the strong Markov property give
\begin{equation}
\label{eq:bmp_conditioned_renewal}
 \begin{split}
h(x) \P_h^x(\dd \omega) = (\overline{P}^x)^H(\dd X^H)
  \Big(
    \Ind_{(H < \zeta)} h(X(H))
  + \Ind_{(\zeta \le H)} h(X(\zeta-))Q(X(\zeta-))
  \Big)\\
\times
  \Big(
    \Ind_{(H < \zeta)} \P^{X(H)}(\dd \omega^{(\emptyset,H)})
    +\Ind_{(\zeta \le H)}
	q_h(X(\zeta-),k)
	\prod_{i=1}^{k} \P_h^{X(\zeta-)}(\dd \omega^{(i)})
  \Big).
 \end{split}
\end{equation}
In particular, integrating over $k$, $\omega^{(i)}$, $i=1,2,\ldots$, and $X(t)$ for $t\in [H,\zeta)$, we get that
\[
 h(x) = (\overline{E}^x)^H\Big(\Ind_{(H < \zeta)} h(X(H)) + \Ind_{(\zeta \le H)} h(X(\zeta-))Q(X(\zeta-))\Big).
\]
We can therefore define a law $\overline{P}_h^x$ on the paths in $D(\mathscr E)$ stopped at $H$ by
\[
 \overline{P}_h^x(\dd X) = (h(x))^{-1} \Big(
    \Ind_{(H(X) < \zeta)} h(X(H))
  + \Ind_{(\zeta \le H(X))} h(X(\zeta-))Q(X(\zeta-))
  \Big)\times (\overline{P}^x)^H(\dd X),
\]
where the multiplication is again in the sense of a Radon--Nikodym derivative. Then \eqref{eq:bmp_conditioned_renewal} yields the following decomposition of the law $\P_h^x$:
\begin{itemize}[nolistsep,label=$-$]
 \item As long as a particle has not hit the set $F$ yet, it moves according to the law $\overline{P}_h^x$. If it gets killed at the point $y$, it spawns $k$ offspring according to the law $q_h(y,\cdot)$, which initiate independent branching Markov processes according to the law $\P_h^y$.
 \item When a particle hits the set $F$ at the point $y$, it continues as a branching Markov process according to the law $\P^y$.
\end{itemize}
If $R(x) \equiv R$, one gets a simpler characterisation of the law $\overline{P}_h^x$: In this case, $h(x)$ is a harmonic function for the law of the stopped process $X^H$ under $P^x$, whence we can define the Doob transform
\[
 P_h^x(\dd X) = (h(x))^{-1} \Big(\Ind_{(H = \infty)} + \Ind_{(H < \infty)} h(X(H))\Big)P^x(\dd X^H).
\]
Then the law $\overline{P}_h^x$ is obtained from the law $P_h^x$ by killing the process at the time-dependent rate $R Q(x)\Ind_{(t < H)}$.

\section{BBM with absorption at a critical line}
\label{sec:critical_line}
From this section on, $q(k)$ will denote a law on $\{0,1,2,\ldots\}$ and $L$ a random variable with law $q(k)$. We define $m = \E[L-1]$ and $m_2 = \E[L(L-1)]$ and suppose that $m>0$ and $m_2 < \infty$. We study the branching Markov process where, starting with a single particle at the origin, particles move according to standard Brownian motion with drift $-1$ and branch at rate $\beta_0 = 1/(2m)$ into $k$ particles according to the reproduction law $q(k)$. At the point $-y$, we add an \emph{absorbing barrier} to the process, i.e.\ particles hitting this barrier are instantly killed. Formally, we are considering the process up to the stopping line $\mathscr L_{H_{-y}}$, where $H_{-y}$ is the hitting time functional of the point $-y$. It is well-known since Kesten \cite{Kesten1978} that this process gets extinct almost surely. As a consequence, the number of particles absorbed at the barrier, i.e.\ the random variable
\[
N_y = \#\mathscr L_{H_{-y}},
\]
is almost surely finite. By the strong branching property and the translational invariance of Brownian motion, one sees that the process $(N_y)_{y\ge 0}$ is a continuous-time Galton--Watson process, a fact which was first noticed by Neveu \cite{Neveu1988} (see \cite{Athreya1972}, Chapter III or \cite{Harris1963}, Chapter V for an introduction to continuous-time Galton--Watson processes). Let $u(s)$ be its infinitesimal generating function. Neveu stated that $u = \psi'\circ \psi^{-1}$, where $\psi$ is a so-called travelling wave of the FKPP (Fisher--Kolmogorov--Petrovskii--Piskounov) equation: Write $f(s) = \sum_k s^k q(k)$. Then $\psi$ is a solution of the equation
\begin{equation}
 \label{eq:travelling_wave}
 \frac 1 2 \psi'' - \psi' = \beta_0(\psi - f\circ \psi),
\end{equation}
with $\psi(-\infty) = 1$ and $\psi(+\infty)$ is the extinction probability of the process, i.e.\ the smaller root of $f(s) = s$. For a proof of these results, see Section~\ref{sec_FKPP} of Chapter~\ref{ch:barrier}.

In the same paper \cite{Neveu1988}, Neveu introduced his \emph{multiplicative martingales}, which he used to derive the Seneta-Heyde norming for the martingale $e^{- y}N_y$. He proved that in the case of binary branching, one has
\begin{equation}
 \label{eq:def_Wy}
 W_y :=  y e^{- y} N_y \to W\quad\text{almost surely as $y\to\infty$},
\end{equation}
where $W > 0$ almost surely. His proof relied on a known asymptotic for the travelling wave $\psi$, namely that
\begin{equation}
 \label{eq:tw_asymptotics}
 1-\psi(-x) \sim K x e^{- x},\quad\tas x \to \infty,
\end{equation}
for some constant $K> 0$.
It was recently shown \cite{Yang2011} that this asymptotic is true if and only if $E[L\log^2 L] <\infty$ and the proof of \eqref{eq:def_Wy} works in this case as well. We also still have in this case, for every $x\in \R$,
\begin{equation}
 \label{eq:W_laplace}
 E[e^{-e^{ x}W}] = \psi(x),
\end{equation}
a fact which was already proven by Neveu \cite{Neveu1988} for dyadic branching.

In \cite{Berestycki2010}, further properties of the limit $W$ have been established under the hypothesis of dyadic branching, namely
\begin{equation}
\label{eq:W_tail}
 \P(W > x) \sim \frac{1}{x},\quad \tas x\to\infty,
\end{equation}
and
\begin{equation}
 \label{eq:W_expec}
 \E[W\Ind_{(W\le x)}] - \log x \to \const{eq:W_expec},\quad \tas x\to\infty,
\end{equation}
for some constant $\const{eq:W_expec}\in\R$. Equation \eqref{eq:W_tail} has been proven in Propositions 27 and 40 of \cite{Berestycki2010}, and \eqref{eq:W_expec} appears in the proof of Proposition 39 of the same paper. Their arguments were very ingenious but indirect and although they could be extended to general reproduction laws with finite variance, we will reprove them here directly under (probably) minimal assumptions, based on methods of \cite{Maillard2010}. The main result in this section is
\begin{proposition}
 \label{prop:W}
If $E[L\log^2L] < \infty$, then \eqref{eq:W_tail} holds. If $E[L\log^3L]<\infty$, then \eqref{eq:W_expec} holds.
\end{proposition}
See also \cite{Buraczewski2009} for a proof of \eqref{eq:W_tail} in the case of branching random walk.
Before proving this result in the next subsection, we state a lemma which is immediate from \eqref{eq:def_Wy} and the fact that $N_y$ is almost surely finite (see also Corollary 25 in \cite{Berestycki2010}):
\begin{lemma}
\label{lem:y_zeta}
 Suppose $E[L\log^2 L]<\infty$. For any $\eta > 0$, there exist $y$ and $\zeta$, such that $y \ge \eta^{-1}$ and
\begin{multline*}
\P(|W_y - W| > \eta) + \P(N_y > \zeta) + \P(\mathscr L_{H_{-y}} \nsubseteq U\times [1,\zeta])\\
 + \P(\sup_{0\le t\le \zeta} \#\{u\in\mathscr N(t): (u,t)\preceq \mathscr L_{H_{-y}}\} > \zeta) \le \eta.
\end{multline*}
\end{lemma}

\subsection{Proof of Proposition \ref{prop:W}}

Define $\chi(\lambda) = E[e^{-\lambda W}]$ for $\lambda \ge 0$. Our first result is:

\begin{lemma}
\label{lem:chi}
Suppose $\E[L\log^2L]<\infty$. Then $\chi''(\lambda) \sim \lambda^{-1}$ as $\lambda \to 0+$. Furthermore, $\E[L\log^3 L] < \infty$ if and only if $r(\lambda) = \lambda^{-1} - \chi''(\lambda) \ge 0$ for $\lambda \ge 0$, with $\int_0^1 r(\lambda)\,\dd\lambda < \infty$.
\end{lemma}

\begin{proof}
Define $\phi(x) = 1-\psi(-x)$, such that $u(s) = \phi'(\phi^{-1}(s))$. By \eqref{eq:tw_asymptotics} and the hypothesis  $\E[L\log^2L]<\infty$,
\begin{equation}
 \label{eq:phi_asymptotic}
 \phi(x) \sim K x e^{-x},\quad\tas x\to \infty.
\end{equation}
Furthermore, by \eqref{eq:travelling_wave}, we have
\begin{equation}
\label{eq:phi} 
\frac 1 2 \phi''(x) + \phi'(x) = \beta_0(f(1-\phi(x)) - (1-\phi(x))).
\end{equation}
Setting $g(s) = 2\beta_0[f(1-s) - 1 + f'(1)s] \ge 0$ and $\rho = \phi + \phi'$, we get from \eqref{eq:phi}, and the fact that $\beta_0 = 1/(2m)$ and $f'(1) = m+1$ by definition,
\begin{equation}
 \label{eq:rho}
 \rho'(x) = -\rho(x) + g(\phi(x)).
\end{equation}
As in the proof of Theorem~\ref{th_tail} from Chapter~\ref{ch:barrier}, we will study the function $\rho$ through the integral equation corresponding to \eqref{eq:rho}, namely
\begin{equation}
 \label{eq:rho_integral}
 \rho(x) = e^{-x}\Big(\rho(0) + \int_0^x e^yg(\phi(y))\,\dd y\Big) = e^{-x}\Big(\rho(0) + \int_{\phi(x)}^{\phi(0)} \frac{e^{\phi^{-1}(s)}g(s)}{-u(s)}\ \dd s\Big).
\end{equation}
Now, by Theorem B of \cite{Bingham1974} (see also Theorem 8.1.8 in \cite{Bingham1987}) we have for every $d\ge 0$,
\begin{equation}
 \label{eq:bingham}
 \int_0^1 \frac{\log^d \frac 1 s}{s^2}g(s)\,\dd s < \infty\quad\Longleftrightarrow\quad\int_0^1 \frac{\log^d \frac 1 s}{s}g'(s)\,\dd s < \infty\quad\Longleftrightarrow\quad E[L\log^{1+d}L]<\infty.
\end{equation}
Furthermore, by Proposition~\ref{prop_FKPP} from Chapter~\ref{ch:barrier}, we have $-u(s) \sim s$ as $s\to 0$, and by \eqref{eq:phi_asymptotic}, we have $e^{\phi^{-1}(s)} \sim (\log 1/s)/s$ as $s\to 0$. By the hypothesis $E[L\log^2 L]<\infty$, this gives
\[
 \int_0^{\phi(0)} \frac{e^{\phi^{-1}(s)}g(s)}{-u(s)}\ \dd s < \infty,
\]
whence, by \eqref{eq:rho_integral},
\begin{equation}
\label{eq:rho_asymptotic}
\rho(x) \sim K e^{-x},\quad\tas x\to\infty,
\end{equation}
where the constant $K$ is actually the same as the one in \eqref{eq:phi_asymptotic}, see again the proof of Theorem~\ref{th_tail} from Chapter~\ref{ch:barrier}. Now, from \eqref{eq:W_laplace}, we get $\chi(\lambda) = 1-\phi(-\log \lambda)$, whence, by \eqref{eq:rho} and \eqref{eq:rho_integral},
\begin{equation}
 \label{eq:chi_rho}
\begin{split}
\chi''(\lambda) = -\frac 1 {\lambda^2}\rho'(-\log \lambda) &= \frac K {\lambda} + \frac 1 {\lambda^2}\Big(-\lambda\int_{-\log\lambda}^\infty e^yg(\phi(y))\,\dd y - g(\phi(-\log\lambda))\Big)\\
& = \frac K {\lambda} + \frac 1 \lambda\Big(\int_{-\log\lambda}^\infty e^y\phi'(y)g'(\phi(y))\,\dd y\Big),
\end{split}
\end{equation}
where the last equation follows from integration by parts. This proves the first statement, with the constant $K$ instead of $1$, since the last integral vanishes as $\lambda \to 0$. Now, setting
\[
 r(\lambda) = - \frac 1 \lambda\Big(\int_{-\log\lambda}^\infty e^y\phi'(y)g'(\phi(y))\,\dd y\Big),
\]
we first remark that $r(\lambda) \ge 0$, since the integrand is negative for $y\in\R$. By the Fubini--Tonelli theorem, we then have
\begin{align*}
 \int_0^1r(\lambda)\dd \lambda & = \int_0^1 -\frac 1 \lambda \int_{-\log \lambda}^\infty e^y\phi'(y)g'(\phi(y))\,\dd y\,\dd \lambda\\
& = \int_0^\infty - ye^y\phi'(y)g'(\phi(y))\,\dd y\\
& = \int_0^{\phi(0)}e^{\phi^{-1}(y)} \phi^{-1}(y) g'(y)\,\dd y,
\end{align*}
which is finite if and only if $E[L\log^3 L]<\infty$, by \eqref{eq:bingham} and the fact that $e^{\phi^{-1}(y)} \phi^{-1}(y) \sim (\log^2 1/s)/s.$ This proves the second statement, again with the constant $K$ instead of $1$.

The previous arguments worked for every travelling wave $\psi$. In order to show that that the constant $K$ is equal to 1 in our case, we use Neveu's multiplicative martingale (this idea was also used in \cite{Liu2000}, Theorem~2.5). It was observed by Neveu \cite{Neveu1988} (see also \cite{Chauvin1991} for a rigorous proof), that $((1-\phi(x+y))^{N_y})_{y\ge 0}$ is a martingale for every $x\in\R$ with values in $[0,1]$. By \eqref{eq:def_Wy} and \eqref{eq:phi_asymptotic}, we then get by dominated convergence, for every $x\in\R$,
\[
 \chi(Ke^x) = \lim_{y\to \infty} E[e^{-Kye^{x-y}N_y}] = \lim_{y\to\infty}E[(1-\phi(y-x))^{N_y}] = 1-\phi(-x) = \chi(e^x).
\]
This yields $K=1$.
\end{proof}

\begin{remark}
 Choosing arbitrary initial points $x_0,x_1\in\R$ instead of $0$ in \eqref{eq:rho_integral}, one sees that
\[
 e^{x_0}\rho(x_0) + \int_{x_0}^\infty e^yg(\phi(y))\,\dd y = e^{x_1}\rho(x_1) + \int_{x_1}^\infty e^yg(\phi(y))\,\dd y.
\]
In particular, since $\rho$ is bounded, letting $x_0\to -\infty$ and $x_1\to +\infty$ yields
\[
 \int_{-\infty}^\infty e^yg(\phi(y))\,\dd y = 1.
\]
One could hope (see the proof of Proposition \ref{prop:W} below) that this helps in determining the constant \const{eq:W_expec}, but apparently this does not seem to be the case.
\end{remark}

\begin{proof}[Proof of Proposition \ref{prop:W}]
For $n\in\N,$ we define the function
\[
 V_n(x) = \int_0^x y^n P(W\in\dd y) = E[W^n\Ind_{(W\le x)}],
\]
such that with $\chi^{(n)}$ denoting the $n$-th derivative of $\chi$, we have for $\lambda > 0$,
\[
 \chi^{(n)}(\lambda) = (-1)^n \int_0^\infty e^{-\lambda x} \dd V_n(x).
\]
If $E[L\log^2L]<\infty$, Proposition \ref{prop:W} and Karamata's Tauberian theorem (\cite{Feller1971}, Theorem XIII.5.2 or \cite{Bingham1987}, Theorem 1.7.1) now yields
\begin{equation}
 \label{eq:V2}
 V_2(x) \sim x,\quad\tas x\to\infty.
\end{equation}
By an integration by parts argument (see also \cite{Feller1971}, Theorem VIII.9.2 or \cite{Bingham1987}, Theorem 8.1.2), we get \eqref{eq:W_tail}. Now suppose that $E[L\log^3L] < \infty$. By Lemma~\ref{lem:chi}, we have $\chi'(\lambda) - \log \lambda \to c\in \R$, as $\lambda \to 0$. By Theorem 3.9.1 from \cite{Bingham1987} (with $\ell(x) \equiv 1$), this yields
\[
 V_1(x) - \log x \to \gamma - c,\quad\tas x\to \infty,
\]
where $\gamma$ is the Euler--Mascheroni constant. This is exactly \eqref{eq:W_expec}.


\end{proof}

\section{BBM in an interval}
\label{sec:interval}

In this section we study branching Brownian motion killed upon exiting an interval. Many ideas in this section (except for Section~\ref{sec:weakly_conditioned} and parts of Section~\ref{sec:interval_number}) stem from Sections~2 and~3 of \cite{Berestycki2010} and for completeness, we will reprove some of their results with streamlined proofs. However, we will also extend their results to the case of Brownian motion with variable drift.

\subsection{Notation}
\label{sec:interval_notation}
During the rest of the paper, the symbol $C$ stands for a \emph{positive} constant, which may only depend on the reproduction law $q$. Its value may change from line to line. If a subscript is present, then this subscript is the number of the equation where this constant appears for the first time (example: $\Const{eq:Rt_nu_expec}$). In this case, this constant is fixed after its value has been chosen in the corresponding equation. If $X$ is any mathematical expression, then the symbol $O(X)$ stands for a possibly random  term whose absolute value is bounded by $C|X|$.

Recall the definition of $q(k)$, $m$, $m_2$ and $\beta_0$ from Section \ref{sec:critical_line} and the hypotheses on $m$ and $m_2$. In this section, we let $a \ge \pi$ and set
\begin{equation}
 \label{eq:mu}
\mu = \sqrt{1 - \frac{\pi^2}{a^2}}.
\end{equation}
From \eqref{eq:mu}, one easily gets the basic estimate
\begin{equation}
 \label{eq:c0_mu}
 0\le 1 - \mu \le \frac{\pi^2}{2\mu a^2}.
\end{equation}

We then denote by $\P^x$ the law of the branching Markov process where, starting with a single particle at the point $x\in\R$, particles move according to Brownian motion with variance $1$ and drift $-\mu$ and branch at rate 1 into $k$ particles according to the reproduction law $q(k)$. Expectation with respect to $\P^x$ is denoted by $\E^x$. On the space of continuous functions from $\R_+$ to $\R$, we define $H_0$ and $H_a$ to be the hitting time functionals of $0$ and $a$. We further set $H = H_0\wedge H_a$. Then note that the density of the branching Brownian motion before $\mathscr L_H$, as defined in \eqref{eq:def_density}, has a density with respect to Lebesgue measure given for $t>0$ and $x,y \in (0,a)$ by
\begin{equation}
 \label{eq:density_bbm}
\p_t(x,y) = e^{\mu(x-y) + \frac{\pi^2}{2 a^2} t} p_t^a(x,y),
\end{equation}
where $p_t^a$ was defined in \eqref{eq:def_p}.

Now, let $f:\R_+\to\R_+$ be non-decreasing, with $f(0)=0$, continuous and such
that the left-derivative $f'$ exists everywhere and is of bounded variation.
Such a function will be called a \emph{barrier function}. We define
\begin{equation}
\label{eq:def_f_norm}
 \|f\| = \max\Big\{\|f\|_\infty, \|f'\|_\infty,\|f'\|_\infty^2,\int_0^\infty|\dd f'(s)|\Big\},
\end{equation}
where $\|\cdot\|_\infty$ is the usual supremum norm. Furthermore, we set $\Err(f,t) = \|f\|\Big(\frac{1}{a} + \frac{t}{a^3}\Big)$.
Now define
\begin{equation}
 \label{eq:mu_t}
\mu_t = \mu + \frac{\dd}{\dd t} f(t/a^2) = \mu + \frac{1}{a^2} f'(t/a^2),
\end{equation}
such that $\mu_0 = \mu$ and $\mu_t \ge \mu$ for all $t\ge0$. We denote by $\P^x_f$ the law of the branching Brownian motion described above, but with infinitesimal drift $-\mu_t$. Expectation with respect to $\P^x_f$ is denoted by $\E^x_f$ and the density of the process is denoted by $\p^f_t(x,y)$.

The above definitions can be extended to arbitrary initial configurations of particles distributed according to a counting measure $\nu$ on $(0,a)$. In this case the superscript $x$ is replaced by $\nu$ or simply omitted if $\nu$ is known from the context.

\subsection{The processes \texorpdfstring{$Z_t$ and $Y_t$}{Z\_t and Y\_t}}
\label{sec:ZY}

Recall from Section \ref{sec:preliminaries} that the set of particles alive at time $t$ is denoted by $\mathscr{N}(t)$. We define
\[
 \mathscr N_{0,a}(t) = \{u\in \mathscr N(t): H(X_u) > t\},
\]
where $H$ was defined in the previous subsection. Now set $w_Z(x) = ae^{\mu (x-a)}\sin(\pi x/a)\Ind_{(x\in[0,a])}$ and $w_Y(x) = e^{\mu(x-a)}\Ind_{(x\ge 0)}$ and define
\[
 Z_t = \sum_{u\in \mathscr N_{0,a}(t)} w_Z(X_u(t))\quad\tand\quad Y_t = \sum_{u\in \mathscr N_{0,a}(t)} w_Y(X_u(t)).
\]
Then $Z_t$ is a martingale under $\P^x$, the proof of which is standard and relies on the branching property, the Many-to-one lemma (Lemma~\ref{lem:many_to_one_simple}) and the fact that $e^{t/2} w_Z(B_t)$ is a martingale for a Brownian motion with drift $-\mu$ killed at $0$ and $a$, which is easily seen by It\={o}'s formula, for example. Furthermore, it is easy to see as well that $Z_t$ is a supermartingale under $\P^x_f$.

The following lemma relates the density of BBM with variable drift to BBM with fixed drift. 
\begin{lemma}
 \label{lem:mu_t_density} For all $x,y\in[0,a]$, $t\ge0$,
\[
 \p^f_t(x,y) = \p_t(x,y) e^{- f(t/a^2) + O(\Err(f,t))}.
\]
\end{lemma}
\begin{proof}
 By the Many-to-one lemma and Girsanov's theorem, we have
\begin{equation}
\label{eq:005}
\begin{split}
 \p^f_t(x,y) &= e^{\beta_0 mt}W_{-\mu_t}^x\left(B_t\in \dd y,\, H > t\right)\\
&= \exp\left(t/2 - \int_0^t\frac{\mu_s^2-\mu^2}{2} \dd s\right)W_{-\mu}^x\left(\exp\Big(-\int_0^t \mu_s-\mu\,\dd B_s\Big),\,B_t\in\dd y, H > t\right).
\end{split}
\end{equation}
By integration by parts, we have
\begin{equation}
 \label{eq:ipp_mu_B}
 \int_0^t\mu_s\,\dd B_s = \mu_t B_t - \mu B_0 - \int_0^tB_s\,\dd \mu_s.
\end{equation}
Since $B_t \in (0,a)$ for all $t\ge 0$, we have
\begin{equation}
 \label{eq:010}
 \Big|\int_0^t B_s\,\dd \mu_s\Big| \le a \int_0^t |\dd \mu_s| \le \frac{\int_0^\infty|\dd
   f'(s)|}{a}
\end{equation}
Furthermore,
\[
 \frac{\mu_t^2}{2} = \frac{\mu^2}{2} + \frac{\mu}{a^2}f'(t/a^2)+\frac{f'(t/a^2)^2}{2a^4},
\]
such that
\begin{equation}
 \label{eq:020}
 \Big|\int_0^t \frac{\mu_s^2}{2}\,\dd s - \frac{\mu^2}{2}t - \mu f(t/a^2)\Big| \le \frac{\|f'\|_\infty^2 t}{2a^4}.
\end{equation}
Finally,
\begin{equation}
 \label{eq:030}
 \Big|\mu_tB_t - \mu B_t\Big| = \Big|\frac{1}{a^2}f'(t/a^2) B_t\Big| \le \frac{\|f'\|_\infty}{a}.
\end{equation}
Equations \eqref{eq:density_bbm}, \eqref{eq:005}, \eqref{eq:ipp_mu_B}, \eqref{eq:010}, \eqref{eq:020} and \eqref{eq:030} now give
\[
  \p^f_t(x,y) = \p_t(x,y) e^{-\mu f(t/a^2) + O(\Err(f,t))},
\]
and the lemma now follows from \eqref{eq:c0_mu}.
\end{proof}

\begin{proposition}
 \label{prop:quantities}
Under any initial configuration of particles, for every $t\ge 0$, we have
\begin{equation}
 \label{eq:Zt_expectation}
\E_f[Z_t] = Z_0e^{- f(t/a^2) + O(\Err(f,t))},
\end{equation}
and if in addition $\mu \ge 1/2$, then
\begin{equation}
 \label{eq:Zt_variance} 
\Var_f(Z_t) \le C e^{- f(t/a^2)+O(\Err(f,t))} \Big(\frac{t}{a^3} Z_0 + Y_0 \Big).
\end{equation}
Furthermore, we have for every $t\ge 0$ (without hypothesis on $\mu$),
\begin{equation}
 \label{eq:Yt_weak}
\E_f[Y_t] \le Ce^{- f(t/a^2)+O(\Err(f,t))}Y_0.
\end{equation}
and for $t\ge a^2$,
\begin{equation}
  \label{eq:Yt}
\E_f[Y_t] \le Ce^{- f(t/a^2)+O(\Err(f,t))} \frac{Z_0}{a}.
\end{equation}
Moreover, for every $a^2\le t\le a^3$, we have
\begin{equation}
  \label{eq:Yt_variance}
\Var_f(Y_t) \le Ce^{- f(t/a^2)+O(\Err(f,t))} \frac{Y_0}a.
\end{equation}
\end{proposition}

\begin{proof} 
Equation \eqref{eq:Zt_expectation} follows from Lemma~\ref{lem:mu_t_density} and the fact that $Z_t$ is a martingale under $\P^x$. In order to show \eqref{eq:Yt_weak} and \eqref{eq:Yt}, it suffices by Lemma~\ref{lem:mu_t_density} to consider the case without variable drift. We first suppose that $t\ge a^2$. By \eqref{eq:density_bbm} and \eqref{eq:p_estimate}, we get
\[
 \E^x[Y_t] \le e^{\mu (x-a)}\int_0^a e^{\frac{\pi^2}{2a^2}t} p_t(x,y)\,\dd y
\le C e^{\mu (x-a)} \sin(\pi x/a)\int_0^a \frac 2 a \sin(\pi y/a)\,\dd y.
\]
The last integral is independent of $a$. Summing over $x$ yields \eqref{eq:Yt} as well as \eqref{eq:Yt_weak} in the case $t\ge a^2$. Now, if $t< a^2$, by the Many-to-one lemma and Girsanov's theorem, we have
\[
 \E^x[Y_t] = e^{\beta_0 mt} W^x_{-\mu}\Big[e^{\mu (X_t-a)},\ H_0\wedge H_a < t\Big] = e^{\pi^2t/(2a^2)} W^x[H_0\wedge H_a < t] e^{\mu (x-a)}.
\]
Summing over $x$ yields \eqref{eq:Yt_weak}.

In order to prove \eqref{eq:Zt_variance}, we have by Lemma~\ref{lem:many_to_two_fixedtime},
\begin{equation}
 \label{eq:bbm_2ndmoment}
 \E^x_f[Z_t^2] = \E_f^x\Big[\sum_{u\in \mathscr N_{0,a}(t)}w_Z(X_u(t))^2\Big] + \beta_0 m_2\int_0^a \int_0^t \p_s^f(x,y) (\E_f^{(y,s)}[Z_t])^2\,\dd s\,\dd y.
\end{equation}
By Lemma~\ref{lem:mu_t_density} and the fact that $Z_t$ is a martingale with respect to the law $\P^x$, this yields
\begin{equation}
 \label{eq:050}
 \E_f^x[Z_t^2] \le C e^{- f(t/a^2)+O(\Err(f,t))} \left(\E^x\Big[\sum_{u\in \mathscr N_{0,a}(t)}w_Z(X_u(t))^2\Big] + \int_0^a \int_0^t \p_s(x,y) w_Z(y)^2\,\dd s\,\dd y\right).
\end{equation}
Now we have for $x\in (0,a)$,
\[
 w_Z(x)^2 = (a\sin(\pi x/a) e^{-\mu(a-x)})^2 \le \pi^2 (a-x)^2e^{-2\mu(a-x)} \le C w_Y(x),
\]
because $\mu \ge 1/2$ by hypothesis. This yields
\begin{equation}
 \label{eq:S1}
 S_1 := \E^x\Big[\sum_{u\in \mathscr N_{0,a}(t)}w_Z(X_u(t))^2\Big] \le C \E^x[Y_t] \le C w_Y(x),
\end{equation}
by \eqref{eq:Yt_weak}. Now, by \eqref{eq:density_bbm} and \eqref{eq:Ja}, we have
\[
 S_2 := \int_0^a \int_0^t \p_s(x,y) w_Z(y)^2\,\dd s\,\dd y = ae^{\mu (x-a)}\int_0^aae^{\mu (y-a)}\sin^2(\pi y/a) J^a(x,y,t)\,\dd y.
\]
Lemma~\ref{lem:I_estimates} now gives
\begin{equation}
\label{eq:S2}
\begin{split}
 S_2 &\le C ae^{\mu (x-a)}\int_0^a e^{-\mu y}\sin^2(\pi y/a)\Big(t\sin(\pi x/a)\sin(\pi y/a) + ay\Big)\,\dd y\\
&\le C ae^{\mu (x-a)} \Big(\sin(\pi x/a)\frac t{a^3} + \frac 1 a\Big)\int_0^\infty e^{-\mu y} y^3\,\dd y,
\end{split}
\end{equation}
the last line following again from the change of variables $y\mapsto a-y$ and the inequality $\sin x\le x$. Using again the fact that $\mu \ge 1/2$, equations \eqref{eq:050}, \eqref{eq:S1} and \eqref{eq:S2} now imply
\begin{equation}
 \label{eq:Zt_2ndmoment}
\E^x_f[Z_t^2] \le C e^{- f(t/a^2)+O(\Err(f,t))} \Big(\frac{t}{a^3} w_Z(x) + w_Y(x) \Big).
\end{equation}
If we write the positions of the initial particles as $x_1,\ldots,x_n$, then by the independence of their contributions to $Z_t$,
\begin{equation}
\label{eq:41}
 \Var_f(Z_t) = \sum_i \Var^{x_i}_f(Z_t) \le \sum_i \E^{x_i}_f[Z_t^2].
\end{equation}
Equations \eqref{eq:Zt_2ndmoment} and \eqref{eq:41} now prove \eqref{eq:Zt_variance}. Equation \eqref{eq:Yt_variance} is proven similarly.
\end{proof}

\subsection{The number of particles}
\label{sec:interval_number}
In this subsection, we establish precise first and second moment estimates for the number of particles alive at a time $t$. These estimates extend those of \cite{Berestycki2010}, which are effective only when $t \gg a^2$. For $r\in[0,a]$ and $t\ge0$, we denote by $N_t(r)$ the number of particles in $[r,a]$ at time $t$.

\begin{proposition}
\label{prop:N_expec}
Suppose $\mu \ge 1/2$. Let $t\ge0$, $x,r\in [0,a)$ and suppose that $x\ge (r+a/{20})\Ind_{(t\le a^2)}$. Then
\begin{equation}
\label{eq:N_expec}
 \E^x[N_{t}(r)] = e^{\mu (x-r)}\left[-\frac{2(1+\mu r)e^{\pi^2 t/(2a^2)}}{a^2} \theta'\Big(\frac x a,\frac t {a^2}\Big) + O\left((1+r^3)\frac{\sin(\pi x/a) \vee \Ind_{(t\le a^2)}}{a^4}\right)\right].
\end{equation}
In particular, if $a/{20}+r\Ind_{(t\le a^2)}\le x < a$, then
\begin{multline}
 \label{eq:N_expec_2}
\E^x[N_{t}(r)] = \frac{2\pi(1+\mu r)e^{\mu (a-r)}}{ a^3}\Big[w_Z(x)\thbar\left(\frac t {a^2}\right)\\
+O\left(\frac{(a-x)^2+(1+r^2)}{a}\left(\frac{w_Z(x)}{a} \vee w_Y(x)\Ind_{(t\le a^2)}\right)\right)\Big].
\end{multline}
Moreover, for every $x\in(0,a)$ and $t > 0$, we have
\begin{equation}
 \label{eq:N_expec_large_t}
\E^x[N_{t}(r)] = \frac{2\pi (1+\mu r)e^{\mu (a-r)}}{ a^3} w_Z(x) \Big(1+\mathrm{error}\Big),
\end{equation}
with $|\mathrm{error}| \le E_{t/a^2}(1 + O(1/a^2))$.
\end{proposition}

\begin{proof}
 Let $t\ge0$ and $0\le x < 1$. By \eqref{eq:density_bbm} and Brownian scaling, we have
\begin{equation}
 \label{eq:250}
 \E^{ax}[N_{ta^2}(r)] = e^{\mu ax + \pi^2 t/2} \int_{r/a}^1 e^{-\mu a z} p_t^1(x,z)\,\dd z.
\end{equation}
Note that $(\partial/\partial z) p_t^1(x,z)|_{z=0} = -2\theta'(x,t)$, by \eqref{eq:p_theta}. Taylor's formula then implies that there exists $\xi\in [0,z]$, such that
\begin{equation}
 \label{eq:252}
 p_t^1(x,z) = - 2z \theta'(x,t) + z^2 \frac{\partial^2}{\partial z^2} p_t^1(x,z)\big|_{z=\xi}.
\end{equation}
We first note that for all $x\in(0,1)$ and $t\ge1$, or for all $x\in[\frac 1 {40},1]$ and $t\le 1$, one has by \eqref{eq:theta} and \eqref{eq:theta_gaussian},
\begin{equation}
  \label{eq:254}
  |\theta'(x,t)| \le C (\sin(\pi x) \vee \Ind_{(t\le 1)}).
\end{equation}
Furthermore, if $t\ge 1$, we have
\begin{equation}
 \label{eq:255}
\begin{split}
 \Big|\frac{\partial^2}{\partial z^2} p_t^1(x,z)\big|_{z=\xi}\Big| &= \pi^2\Big|\sum_{n=1}^\infty e^{-\pi^2 n^2 t/2} n^2 \sin(\pi n x)\sin(\pi n \xi)\Big|\\
&\le \pi^4 \xi \sin(\pi x) \sum_{n=1}^\infty n^4 e^{-\pi^2 n^2 t/2} \le C z \sin(\pi x) e^{-\pi^2 t/2},
\end{split}
\end{equation}
by the inequality $|\sin nx| \le n \sin x$, $x\in[0,\pi]$. Equations \eqref{eq:252}, \eqref{eq:254} and \eqref{eq:255} now give
\begin{equation}
 \label{eq:256}
 \int_{r/a}^1 e^{-\mu a z} p_t^1(x,z)\,\dd z = - \frac{2\theta'(x,t)}{\mu^2 a^2} (1+\mu r) e^{-\mu r} + O\Big((1+r^3) e^{-\mu r}\frac{\sin(\pi x)}{a^4} e^{-\pi^2 t/2}\Big).
\end{equation}
Now suppose that $t\le1$, $x,r\in [0,a)$ and $x\ge r+a/{20}$.  By \eqref{eq:p_theta} and the mean value theorem, we have $(\partial^2/\partial z^2) p_t^1(x,z) = 2 z \theta'''(x+\xi',t)$ for some $\xi'\in[0,z]$. With \eqref{eq:theta_gaussian}, one then easily sees that for $z\le x-a/40$, one has
\begin{equation}
 \label{eq:258}
 \Big|\frac{\partial^2}{\partial z^2} p_t^1(x,z)\Big| \le C z e^{-
\Const{eq:258}/t}.
\end{equation}
Equations \eqref{eq:252}, \eqref{eq:254} and \eqref{eq:258} now give, for $x' = x-a/40$,
\begin{equation}
\label{eq:260}
 \int_{r/a}^{x'/a} e^{-\mu a z} p_t^1(x,z)\,\dd z = - \frac{2\theta'(x,t)}{\mu^2 a^2} (1+\mu r)e^{-\mu r} + O\Big(\frac{(1+r^3)e^{-\mu r}}{a^4}\Big).
\end{equation}
Furthermore, we have
\begin{equation}
 \label{eq:262}
\int_{x'/a}^1 e^{-\mu a z} p_t^1(x,z)\,\dd z \le e^{-\mu x'} \int_{x'/a}^1 p_t^1(x,z)\,\dd z \le e^{-\mu x'}.
\end{equation}
Equations \eqref{eq:250}, \eqref{eq:256}, \eqref{eq:260}, \eqref{eq:262} and the hypothesis on $\mu$ now imply \eqref{eq:N_expec}.

In order to prove \eqref{eq:N_expec_2}, let $ x\in[1/{20},1]$. By \eqref{eq:pde_theta} and Taylor's formula, there exists $\xi\in[x,1]$, such that
\begin{equation}
 \label{eq:267}
-\theta'(x,t) = 2(1-x)\frac{\partial}{\partial t}\theta(1,t) - (1-x)^2\theta'''(\xi,t).
\end{equation}
With \eqref{eq:theta} and \eqref{eq:theta_gaussian}, one now easily sees that $|\theta'''(\xi,t)| \le C(e^{-\pi^2t/2}\sin(\pi x) \vee \Ind_{(t\le 1)})$.  Equations \eqref{eq:c0_mu}, \eqref{eq:N_expec} and \eqref{eq:267} now readily imply \eqref{eq:N_expec_2}. For the last equation, by \eqref{eq:density_bbm} and \eqref{eq:p_estimate}, we have for every $x\in(0,a)$ and $t>0$,
\begin{equation}
 \label{eq:271}
\E^x[N(ta^2)] = (1+\mathrm{err}) w_Z(x) \frac {2 e^{\mu a}} {a^2} \int_0^a e^{-\mu y} \sin(\pi y/a)\,\dd y,
\end{equation}
with $|\mathrm{err}| \le E_t$. Evaluating this integral and using \eqref{eq:c0_mu} and the hypothesis on $\mu$ yields \eqref{eq:N_expec_large_t}.
\end{proof}

We will often use the following handier upper bounds on $\E^x[N_{t}(r)]$:

\begin{lemma}
  \label{lem:N_expec_upper_bound}Suppose $\mu \ge 1/2$. Let $t\ge0$, $x,r\in
  [0,a]$ and suppose that $x\ge (r+a/{20})\Ind_{(t\le a^2)}$.
  \begin{equation}
    \label{eq:N_expec_upper_bound_1}
\E^x[N_{t}(r)] \le C (1+r)\frac{e^{\mu (a-r)}}{a^3}\Big(w_Z(x)+\frac{1+r^2} a w_Y(x)\Ind_{(t\le a^2)}\Big).
  \end{equation}
Furthermore, we have for all $x,r\in[0,a]$ and $t\ge 0$,
\begin{equation}
  \label{eq:N_expec_upper_bound_2}
  \E^x[N_{t}(r)] \le C e^{\mu(x-r)}.
\end{equation}
\end{lemma}
\begin{proof}
  One sees from \eqref{eq:theta} and \eqref{eq:theta_gaussian} that $\theta''(x,t) \le C$ for
  $(x,t)\in[0,1]\times[1,\infty)\cup [1/{20},1]\times[0,1]$. Equation \eqref{eq:N_expec_upper_bound_1} now
  follows from \eqref{eq:N_expec} and the mean value theorem. For \eqref{eq:N_expec_upper_bound_2}, we have by \eqref{eq:Yt_weak},
\[
\E^x[N_{t}(r)] \le e^{-\mu r}\E^x\Big[\sum_{u\in\mathscr N(t)}e^{\mu X_u(t)}\Big] \le C e^{\mu(x-r)}.
\]
\end{proof}

\begin{lemma}
\label{lem:N_2ndmoment}
Suppose $\mu\ge 1/2$ and $r\le 9a/10$. For every $t\ge 0$ and $x\in[0,a]$, we have for large $a$,
\begin{equation}
 \label{eq:N_2ndmoment}
 \E^x_f[N_{t}(r)^2] \le Ce^{- f(t/a^2)+O(\Err(f,t))} \left(\frac{e^{\mu a}}{a^3}\right)^2 (1+r^2)\Big(1+\frac{1+r^4}{a^2}\Big)e^{-2 \mu r}\Big(w_Z(x)  \frac t {a^3} + w_Y(x)\Big)
\end{equation}
\end{lemma}
\begin{proof}
As in the proof of \eqref{eq:Zt_variance}, we have by Lemmas~\ref{lem:many_to_two_fixedtime} and \ref{lem:mu_t_density},
\begin{equation}
\label{eq:272}
 \E^x[N_{t}(r)^2] \le e^{- f(t/a^2)+O(\Err(f,t))} \Big(\E^x[N_{t}(r)] + \beta_0 m_2 \int_0^a \int_0^t \p_s(x,z) (\E^z[N_{t-s}(r)])^2\,\dd s\,\dd z\Big).
\end{equation}
By \eqref{eq:density_bbm}, Lemma~\ref{lem:N_expec_upper_bound} and the hypotheses on $x$ and $r$, we have after a change of variables $z\to a-z$ in the integral,
 \begin{multline}
\label{eq:278}
  \int_0^a \p_s(x,z) (\E^z[N_{t-s}(r)])^2\,\dd z \le C \left(\frac{e^{\mu a}}{a^3}\right)^2 e^{\mu(x-a-2r)+\pi^2s/(2a^2)}\\
\times \Big( (1+r^2) \int_0^a p_s^a(a-x,z) e^{-\mu z}\Big(z^2+\frac{1+r^4}{a^2}\Big) \,\dd z +  a^6  \int_{a/{20}}^a p_s^a(a-x,z) e^{-\mu z} \,\dd z\Big).
 \end{multline}
Integrating \eqref{eq:278} over $s$ from $0$ to $t$ and splitting the interval
at $t\wedge a^2$, we have by \eqref{eq:p_potential} (for the first piece) and \eqref{eq:p_estimate} (for the second),
\begin{multline}
\label{eq:280}
\int_0^a \int_0^t \p_s(x,z) (\E^z[N_{t-s}(r)])^2\,\dd s\,\dd z \le C\left(\frac{e^{\mu a}}{a^3}\right)^2\\
\times (1+r^2)\Big(1+\frac{1+r^4}{a^2}\Big)e^{-2\mu r} (1+a^7e^{-\mu a/{20}})\Big( w_Z(x) \frac t {a^3}+w_Y(x) \Big)
\end{multline}
Equations \eqref{eq:N_expec_upper_bound_2}, \eqref{eq:272} and \eqref{eq:280} and the hypotheses on $\mu$ and $r$ now imply the lemma.
\end{proof}

\subsection{The particles hitting the right border}
\label{sec:right_border}
In this section we recall some formulae from \cite{Berestycki2010} about the number of particles hitting the right border of the interval. We reprove these formulae here for completeness and because Lemma~\ref{lem:I_estimates} makes their proofs straightforward. For most formulae we will assume that $f\equiv 0$, i.e.\ that we are working under the measure $\P$. Only Lemma~\ref{lem:R_f} contains an upper bound on the expected number of particles for general $f$, which will be useful in Section \ref{sec:BBBM}.

For a measurable subset $S\subset \R$, define $R_S$ to be the number of particles killed at the right border during the (time) interval $S$, i.e.
\[
 R_S = \#\{(u,t):u\in \mathscr N(t)\tand H_0(X_u)>H_a(X_u) = t\in S\}.
\]
The following lemma gives exact formulae of the expectation and the second moment of $R_S$.

\begin{lemma}
 \label{lem:R_moments}
For every $x\in (0,a)$, we have
 \begin{align}
  \label{eq:RS_x_expec}
  \E^x[R_S] &= e^{\mu(x-a)} I^a(x,S),\\
  \label{eq:RS_x_2ndmoment}
  \E^x[R_S^2] &= \E^x[R_S] + m_2 e^{\mu(x-a)} \int_0^a \dd y\, e^{\mu(y-a)} \int_0^\infty \dd t\,e^{\frac{\pi^2}{2a^2}t}p_t^a(x,y) I^a(y,S-t)^2
 \end{align}
\end{lemma}
We will first prove a more general result, which will be needed in Section \ref{sec:fugitive}.
\begin{lemma}
\label{lem:R_many_to_one}
 For every $x\in (0,a)$ and any measurable function $f:\R_+\to \R_+$, we have
 \[
  \E^x\Big[\sum_{(u,t)\in \mathscr L_H} f(t) \Ind_{(X_u(t) = a)}\Big] = e^{\mu(x-a)}\int_0^t f(s)I^a(x,\dd s).
 \]
\end{lemma}

\begin{proof}
Recall that $H_0$ and $H_a$ denote the hitting time functionals of $0$ and $a$ and $H = H_0\wedge H_a$. Then note that $W^x_{-\mu}(H<\infty) = 1$ for all $x\in [0,a]$. We then have
\begin{align*}
 \E^x\Big[\sum_{(u,t)\in \mathscr L_H} f(t) \Ind_{(X_u(t) = a)}\Big] &= W^x_{-\mu}\Big[e^{\beta_0 mH_a} f(H_a)\Ind_{(H_0 > H_a\le t)}\Big] && \text{by Lemma~\ref{lem:many_to_one_simple}}\\
&= e^{\mu(x-a)} W^x_0\Big[e^{\frac{\pi^2}{2a^2}H_a}f(H_a)\Ind_{(H_0 > H_a\le t)}] && \text{by Girsanov's transform}\\
&= e^{\mu(x-a)} \int_0^t f(s)I^a(x,\dd s) && \text{by \eqref{eq:Ia}.}
\end{align*}
Note that in the second line we used the fact that $\beta_0 m = 1/2$ by definition.
\end{proof}
\begin{proof}[Proof of Lemma~\ref{lem:R_moments}]
Equation \eqref{eq:RS_x_expec} follows from Lemma~\ref{lem:R_moments} and \eqref{eq:IJ_scaling} by taking $f=\Ind_S$. Equation \eqref{eq:RS_x_2ndmoment} follows from Lemma~\ref{lem:many_to_two} and \eqref{eq:RS_x_expec}.
\end{proof}

\begin{lemma}
\label{lem:Rt}
 For any initial configuration $\nu$ and any $0\le s\le t$, we have
\begin{equation}
 \label{eq:Rt_nu_expec}
 |\E[R_{[s,t]}] - \frac{\pi (t-s)}{a^3} Z_0| \le \Const{eq:Rt_nu_expec} \Big(Y_0 \wedge E_{s/a^2} (1\wedge (t-s)/a^3) Z_0\Big),
\end{equation}
where $E_s$ is defined in \eqref{eq:def_E}. Furthermore, if $\mu \ge 1/2$ and $0\le t\le a^3$, then for each $x\in (0,a)$,
\begin{equation}
 \label{eq:Rt_nu_variance}
 \E^x[R_t^2]\le \Const{eq:Rt_nu_variance} \Big(\frac{t}{a^3} w_Z(x) + w_Y(x)\Big),
\end{equation}
\end{lemma}

\begin{proof}
We have $\E[R_t] = \int \nu(\dd x) \E^x[R_t]$, such that \eqref{eq:Rt_nu_expec} follows from \eqref{eq:RS_x_expec} and Lemma~\ref{lem:I_estimates}. For the second moment, we have by \eqref{eq:RS_x_2ndmoment},
\[
\begin{split}
 \E^x [R_t^2 - R_t] &= \beta m_2e^{\mu(x-a)} \int_0^a \dd y\, e^{\mu(y-a)} \int_0^t \dd s\,e^{\frac{\pi^2}{2a^2}t}p_s^a(x,y)I^a(y,t-s)^2\\
 &\le C e^{\mu(x-a)} \int_0^a \dd y\, e^{\mu(y-a)} I^a(y,t)^2 J^a(x,y,t)\\
 &\le C e^{\mu(x-a)} \int_0^a \dd y\, e^{\mu(y-a)} (t/a^2 \sin(\pi y/a) + 1)^2\\
&\quad\quad\quad\quad\quad\times (a t/a^2 \sin(\pi x/a)\sin(\pi y/a) + a^{-1}(x\wedge y)(a-(x\vee y)))
\end{split}
\]
Performing the change of variables $y\mapsto a-y$ in the integral and making use of the inequalities $a^{-1}(x\wedge y)(a-(x\vee y)) \le a-y$ and $\sin x \le x$, we get
\[
\begin{split}
 \E^x [R_t^2 - R_t] &\le C e^{\mu(x-a)} (\sin(\pi x/a) t/a^2 + 1) \int_0^\infty \dd y\,e^{-\mu y} (y + y^2 t/a^3 + y^3 t^2/a^6)\\
 &\le C e^{\mu(x-a)} (\sin(\pi x/a) t/a^2 + 1) (1+ t^2/a^6),
\end{split}
\]
where we used the hypothesis $\mu \ge 1/2$. The last inequality, together with \eqref{eq:Rt_nu_expec} and the hypothesis $t\le a^3$  yields \eqref{eq:Rt_nu_variance}.
\end{proof}

\begin{lemma}
 \label{lem:R_f}
Let $f$ be a function as in Section \ref{sec:interval_notation}. Then for every $x\in (0,a)$, we have
\[
 \E^x_f[R_S] \le \E^x[R_S].
\]
\end{lemma}
\begin{proof}
 As in the proof of Lemma~\ref{lem:R_moments}, we have
\[
 \E^x_f[R_S] = W^x_{-\mu_t}\Big[e^{\beta_0 mH_a} \Ind_{(H_0 > H_a\in S)}\Big] = W^x_0\Big[e^{-\int_0^{H_a} \mu_t\,\dd B_t - \int_0^{H_a} \mu_t^2/2\,\dd t + H_a/2} \Ind_{(H_0 > H_a\in S)}\Big],
\]
by Girsanov's theorem and the definition of $\beta_0$. Now, we have by \eqref{eq:ipp_mu_B}, on the event $\{H_0 > H_a\}$,
\[
 \int_0^{H_a} \mu_t\,\dd B_t = \mu(a-x) + a (\mu_{H_a} - \mu) - \int_0^{H_a} B_t\,\dd \mu_t \ge \mu(a-x),
\]
since $B_t \in [0,a]$ for $t\in [0,H_a]$. This gives
\[
 \E^x_f[R_S] \le e^{\mu(x-a)}W^x_0\Big[e^{(1-\mu^2) H_a/2} \Ind_{(H_0 > H_a\in S)}\Big] = \E^x[R_S],
\]
by the proof of Lemma~\ref{lem:R_moments}.
\end{proof}

We finish this section with a lemma which links BBM with absorption at a critical line to our BBM with selection model.

\begin{lemma}
\label{lem:critical_line}
 Let $\zeta \ge 1$, $y\ge 1$, $\mu \ge 1/2$ and $f$ be a barrier function (defined in Section \ref{sec:interval_notation}). Suppose that $\sqrt a \ge y + \zeta$ and $\|f\|\le \sqrt a$. Let $(x_i,t_i)_{i=1}^N$ be a collection of space-time points with
\[
 x_i = a-y+(1-\mu)t_i-f(s/a^2),\quad i=1,\ldots,N,
\]
and $t_i \le \zeta$ for all $i$. Define $Z = \sum_i w_Z(x_i)$, $Y = \sum_i w_Y(x_i)$ and $W_y = ye^{- y} N$. Then,
\[
 Z =\pi W_y\Big(1+O\Big(\frac 1 a \Big)\Big)\quad\tand\quad  Y = \frac 1 { y} W_y \Big(1+O\Big(\frac{1}{a}\Big)\Big).
\]
In particular, for large $a$, we have
\[
 Y \le Z/y.
\]

\end{lemma}
\begin{proof}
By \eqref{eq:c0_mu} and the hypotheses $\mu \ge 1/2$ and $\zeta \ge 1$, we have for all $i$,
\[
 x_i = a-y + O\Big(\frac{\zeta(1+\|f\|)}{a^2}\Big).
\]
Hence, by \eqref{eq:c0_mu} and the hypotheses $\mu \ge 1/2$, $\|f\|\le \sqrt a$ and $\sqrt a \ge y + \zeta$,
\begin{equation}
 \label{eq:240}
 e^{\mu x_i} = e^{\mu a-y}\Big(1+O\Big(\frac{1}{a}\Big)\Big).
\end{equation}
Furthermore, since $x-x^2/3\le \sin x \le x$ for $x\ge 0$ and by the hypotheses $y\ge 1$ and $a \ge y + \zeta$,
\begin{equation}
 \label{eq:245}
 \sin(\pi x_i/a) =  \sin(\pi (a-x_i)/a) =\frac \pi a y\Big(1 + O\Big(\frac y a\Big)\Big).
\end{equation}
The lemma now follows by summing over \eqref{eq:240} and \eqref{eq:245}.
\end{proof}

\subsection{Penalizing the particles hitting the right border}
\label{sec:weakly_conditioned}
In this section, let $(U_u)_{u\in U}$ be iid random variables, uniformly distributed on $(0,1)$, independent of the branching Brownian motion. Furthermore, let $p:\R_+\in (0,1]$ be measurable and such that $p(t) = 0$ for large enough $t$. Recall that $H = H_0 \wedge H_a$. We define the event
\[
 E = \{\nexists (u,t)\in \mathscr L_H: X_u(t) = a \tand U_u \le p(t)\}.
\]
Our goal in this section is to describe the law $\Ptilde_f^x = \P_f^x(\cdot|E)$. We first note that
\begin{equation}
\label{eq:180}
 \P_f^x(\dd \omega, E) = \P_f^x(\dd \omega) \prod_{(u,t)\in L_t}\Big(\Ind_{(X_u(t) \ne a)} + p(t)\Ind_{(X_u(t) = a)}\Big).
\end{equation}
In order to apply the results from Section \ref{sec:doob}, we define
\begin{align}
\label{eq:def_gamma}
 h(x,t) &= \P_f^{(x,t)}(E)\\
\label{eq:def_beta}
 Q(x,t)      &= \sum_{k=0}^\infty q(k)h(x,t)^{k-1}\\
\label{eq:def_btilde}
 \widetilde q(x,t,k) &= q(k) h(x,t)^{k-1} / Q(x,t)
\end{align}
By the results from Section \ref{sec:doob}, under the law $\Ptilde_f^x$, the BBM stopped at $\mathscr L_H$ is the branching Markov process where
\begin{itemize}[nolistsep,label={--}]
 \item particles move according to the Doob transform of Brownian motion with drift $-\mu_t$ (stopped at $0$ and $a$) by the space-time harmonic function $h(x,t)$ and
 \item a particle located at the point $x\in(0,a)$ at time $t$ branches at rate $\beta_0 Q(x,t)\Ind_{x\in (0,a)}$, throwing $k$ offspring with probability $\widetilde q(x,t,k)$.
\end{itemize}

We have the following useful Many-to-one lemma for the conditioned process stopped at the stopping line $\mathscr L_t = \mathscr L_{H \wedge t}$: Define the function
\begin{equation}
 \label{eq:def_e}
 e(x,t) = \beta_0 \sum_{k\ge 0} k (1-h(x,t)^{k-1}) q(k) \le \beta_0 m_2(1-h(x,t)).
\end{equation}
\begin{lemma}
\label{lem:many_to_one_conditioned}
 For any measurable function $g:[0,a]\to \R_+$, $x\in (0,a)$ and $t\ge 0$, we have
\begin{equation}
\label{eq:many_to_one_conditioned}
 \Etilde_f^x\Big[\sum_{u\in \mathscr L_t} g(X_u(t))\Big] = W_{-\mu_t}^x\Big[g(X_{H\wedge t})\frac{h(X_{H\wedge t},H\wedge t)}{h(x,0)} e^{(H\wedge t)/2 -\int_0^{H\wedge t} e(X_s,s) \,\dd s}\Big].
\end{equation}
In particular, if we denote by $\widetilde{\p}^f_t(x,y)$ the density of the $\Ptilde_f^x$-BBM, then
\begin{equation}
 \label{eq:ptilde_estimate}
\widetilde{\p}_t(x,y) =  \frac{h(y,t)}{h(x,0)} \p_t(x,y) W^{x,t,y}_{\mathrm{taboo}}\Big[e^{-\int_0^t e(X_s,s) \,\dd s}\Big] \le (h(x,0))^{-1} \p_t(x,y),
\end{equation}
and for general $f$,
\begin{equation}
  \label{eq:ptilde_estimate_f}
  \widetilde{\p}^f_t(x,y) \le (h(x,0))^{-1} \p^f_t(x,y).
\end{equation}
\end{lemma}

\begin{proof}
By Lemma~\ref{lem:many_to_one_simple} and the description of the law $\Ptilde_f^x$ given above, we have
\[
 \Etilde_f^x\Big[\sum_{u\in \mathscr L_t} g(X_u(t))\Big] = W^x_{-\mu_t}\Big[g(X_{H\wedge t})\frac{h(X_{H\wedge t},H\wedge t)}{h(x,0)} e^{\int_0^{H\wedge t} \beta_0 \widetilde  m(X_s,s)Q(X_s,s) \,\dd s}\Big],
\]
where $\widetilde m(x,t) = \sum_k(k-1)\widetilde q(x,t,k)$, which yields \eqref{eq:many_to_one_conditioned}. Equation \eqref{eq:ptilde_estimate} follows from \eqref{eq:many_to_one_conditioned} applied to the Dirac Delta-function $g = \delta_y$, $y\in (0,a)$, together with \eqref{eq:density_bbm}.
\end{proof}

The previous lemma immediately gives an upper bound for the quantities we are interested in:
\begin{corollary}
\label{cor:conditioned_upperbound}
 Let  $x\in (0,a)$, $t\ge 0$ and  $g:[0,a]\to \R_+$ be measurable with $g(0) = g(a) = 0$. Define $S_t = \sum_{u\in \mathscr L_t} g(X_u(t))$. Then,
 \begin{align}
  \label{eq:conditioned_upperbound_1}
  \Etilde_f^x[S_t] &\le (h(x,0))^{-1} \E_f^x[S_t],\\
  \label{eq:conditioned_upperbound_2}
  \Etilde_f^x[S_t^2] &\le (h(x,0))^{-1} \E_f^x[S_t^2].
 \end{align}
\end{corollary}
\begin{proof}
Equation \eqref{eq:conditioned_upperbound_1} immediately follow from \eqref{eq:ptilde_estimate}. In order to prove the second-moment estimates, we note that by Lemma~\ref{lem:many_to_two_fixedtime} and the description of the conditioned process,
\[
 \begin{split}
  \Etilde_f^x[S^2_t] &= \Etilde_f^x\Big[\sum_{u\in \mathscr L_t} g(X_u(t))^2\Big] + \int_0^t\int_0^a\widetilde{\p}_s(x,y) \widetilde{m_2}(y,s)\beta_0 Q(y,s)\left(\Etilde_f^{(y,s)}[S_t]\right)^2\,\dd y\,\dd s,
 \end{split}
\]
where $\widetilde{m_2}(x,t) = \sum_{k\ge 0} k(k-1) \widetilde{q}(x,t,k).$ By \eqref{eq:def_btilde}, we have $\widetilde{m_2}(x,t)Q(x,t)\le h(x,t) m_2$. Equation \eqref{eq:conditioned_upperbound_2} then follows from \eqref{eq:bbm_2ndmoment}, \eqref{eq:ptilde_estimate} and \eqref{eq:conditioned_upperbound_1}.
\end{proof}

The following lemma gives a good lower bound on the first-moment estimates in the case where $f\equiv 0$.
\begin{lemma}
\label{lem:conditioned_Z_lowerbound}
Suppose $\mu \ge 1/2$, $t\le a^3$ and $p(s) = 0$ for all $s\ge a^3$. Let $S_t$ be as in Corollary \ref{cor:conditioned_upperbound}. We have
\begin{equation}
\label{eq:conditioned_Z_lowerbound}
 \Etilde^x[S_t] \ge \E^x[S_t](1-\Const{eq:conditioned_Z_lowerbound}\|p\|_\infty ).
\end{equation}
\end{lemma}

This follows from the following estimate on $h(x,0)$, which will be sharpened in Lemma~\ref{lem:T0}.

\begin{lemma} Suppose $p(s) = 0$ for all $s\ge a^3$. Then for all $x\in(0,a)$, we have
 \label{lem:h_estimate}
\[
1-h(a-x,0) \le C \|p\|_\infty   (x + 1) e^{-\mu x}.
\]
\end{lemma}
\begin{proof}
By Markov's inequality, we have
\[
 \begin{split}
  1-h(x,0) &= \P_f^x(\#\{(u,s)\in \mathscr L_{a^3}: X_u(s)=a,\ U_u \le p(s)\} \ge 1) \\
    &\le \E_f^x(\#\{(u,s)\in \mathscr L_{a^3}: X_u(s)=a,\ U_u \le \|p\|_\infty\})\\
    &\le \|p\|_\infty  \E_f^x(R_{a^3}),
 \end{split}
\]
The lemma now follows from Lemmas~\ref{lem:Rt} and \ref{lem:R_f} and the inequality $\sin x\le x$, $x\in[0,\pi]$.
\end{proof}



\begin{proof}[Proof of Lemma~\ref{lem:conditioned_Z_lowerbound}]
By \eqref{eq:ptilde_estimate},
\begin{equation}
 \label{eq:210}
 \Etilde^x[S_t] \ge \E^x[S_t]\inf_{y\in(0,a)}\Big(h(y,t)W_{\mathrm{taboo}}^{x,t,y}\Big[e^{-\int_0^t e(X_s,s) \,\dd s}\Big]\Big).
\end{equation}
By Lemma~\ref{lem:h_estimate} and the hypothesis on $\mu$, we have for every $y\in (0,a)$,
\begin{equation}
 \label{eq:218}
h(a-y,t) \ge h(a-y,0) \ge 1-C\|p\|_\infty e^{-y/3}.
\end{equation}
This gives
\begin{equation}
 \label{eq:219}
 W_{\mathrm{taboo}}^{x,t,y}\Big[e^{-\int_0^t e(X_s,s) \,\dd s}\Big] \ge 1 - W_{\mathrm{taboo}}^{x,t,y}\Big[\int_0^t e(X_s,s) \,\dd s\Big] \ge 1-C\|p\|_\infty,
\end{equation}
by \eqref{eq:def_e}, \eqref{eq:218}, Lemma~\ref{lem:k_integral}, the inequality $e^{-x} \ge 1-x$ for $x\ge 0$ and the fact that the law of the Brownian taboo process is preserved under the map $y\mapsto a-y$. The lemma now follows from \eqref{eq:210}, \eqref{eq:218} and \eqref{eq:219}.
\end{proof}

Finally, we study the law of $R_t$ under the new probability.
\begin{lemma}
 \label{lem:conditioned_Rt}
We have for every $x\in [0,a]$,
\begin{equation}
\label{eq:conditioned_Rt}
 \E_f^x[R_t] - \|p\|_\infty \E_f^x[R_t^2] \le \Etilde_f^x[R_t] \le (h(x,0))^{-1}\E_f^x[R_t],
\end{equation}
and if there is a $p\in[0,1]$, such that $p(s) \equiv p$ for $s\le t$, then we even have
\begin{equation}
\label{eq:conditioned_Rt_constant}
 \Etilde_f^x[R_t] \le \E_f[R_t].
\end{equation}
\end{lemma}

\begin{proof}
Let $\mathscr R_t$ be the stopping line
\[
 \mathscr R_t = \{(u,s)\in \mathscr L_{H_a}: s\le t\}.
\]
We have by definition of the law $\Ptilde$,
\begin{equation}
 \label{eq:230}
 \Etilde_f^x[R_t] = \frac{\E_f^x\Big[R_t \prod_{(u,s)\in\mathscr R_t} (1-p(X_u(s)))\Big]}{\E_f^x\Big[\prod_{(u,s)\in\mathscr R_t} (1-p(X_u(s)))\Big]}.
\end{equation}
Now the denominator is $h(x,0)$ by \eqref{eq:def_gamma}, which yields the right-hand side of \eqref{eq:conditioned_Rt}. The left-hand side follows by noticing that
\[
 \E_f^x\Big[R_t \prod_{(u,s)\in\mathscr R_t} (1-p(X_u(s)))\Big] \ge \E_f^x[R_t(1-\|p\|_\infty )^{R_t}] \ge \E_f^x[R_t] - \|p\|_\infty \E_f^x[R_t^2].
\]
For \eqref{eq:conditioned_Rt_constant}, we note that if $p(s) \equiv p$ for $s\le t$, then by \eqref{eq:230},
\[
 \Etilde_f^x[R_t] = \frac{\E_f^x[R_t(1-p)^{R_t}]}{\E_f^x[(1-p)^{R_t}]}.
\]
Since $(1-p)^k$ is decreasing in $k$, this yields \eqref{eq:conditioned_Rt_constant}.
\end{proof}

\section{BBM with absorption before a breakout}
\label{sec:before_breakout}

In this section, we are studying branching Brownian motion with drift $-\mu$ and absorption at $0$ until a \emph{breakout} occurs, an event which will be defined in Section \ref{sec:before_breakout_definitions} and which corresponds to a particle going far to the right and spawning a large number of descendants. In \eqref{eq:decomposition}, we decompose the system into a particle conditioned to break out at a specific time $T$ (this particle will be called the \emph{fugitive}) and the remaining particles, which are conditioned not to break out before time $T$. These two parts will be studied separately, the former in Section \ref{sec:fugitive} and the latter in Section \ref{sec:before_breakout_particles}. Before that, in Section \ref{sec:time_breakout}, we study the law of the time of the first breakout, showing that it is approximately exponentially distributed. First of all, however, we start with the necessary definitions:

\subsection{Definitions}
\label{sec:before_breakout_definitions}
We will introduce several parameters which will be used during the rest of the paper. The two most important parameters are $a$ and $A$, which are both large positive constants. The meaning of $a$ is as in the previous sections: It is the right border of an interval in which the particles are staying most of the time. The parameter $A$ has a more subtle meaning and controls the number of particles of the system and with it the intensity at which particles hit the point $a$. In Section \ref{sec:BBBM}, we will indeed choose the initial conditions such that $Z_0 \approx  e^A$.


When we study the system for large $A$ and $a$, we first let $a$ go to infinity, then $A$. Thus, the statement ``For large $A$ and $a$ we have\ldots'' means: ``There exist $A_0$ and a function $a_0(A)$, both depending on the reproduction law $q$ only, such that for $A \ge A_0$ and $a \ge a_0(A)$ we have\ldots''. Likewise, the statement ``As $A$ and $a$ go to infinity\ldots'' means ``For all $A$ there exists $a_0(A)$, depending on the reproduction law $q$ only, such that as $A$ goes to infinity and $a\ge a_0(A)$\ldots''. These phrases will become so common that in Sections~\ref{sec:BBBM} to~\ref{sec:Bsharp} they will often be used implicitly, although they will always be explicitly stated in the theorems, propositions, lemmas etc. We further introduce the notation $o(1)$, which stands for a (non-random) term that only depends on the reproduction law $q$ and the parameters $A$, $a$, $\ep$, $\eta$, $y$ and $\zeta$  and which goes to $0$ as $A$ and $a$ go to infinity.

The remaining parameters we introduce are all going to depend on $A$, but not on $a$. First of all, there is the small parameter $\ep$, which controls the intensity of the breakouts. Indeed, when $Z_0 \approx  e^A$, the mean time one has to wait for a breakout will be approximately proportional in $\ep$. Morally, one could choose $\ep$ such that $e^{-A/2} \ll \ep \ll A^{-1}$, but for technical reasons we will require that
\begin{align}
 \label{eq:ep_upper}
 \ep &\le \Const{eq:ep_upper} A^{-17},\quad \tand\\
 \label{eq:ep_lower}
 \ep &\ge \Const{eq:ep_lower} e^{-A/6}.
\end{align}
Another protagonist is $\eta$, which we will choose as small as we need and which will be used to bound the probability of very improbable events, as well as the contribution of the variable $Y$. It will be enough to require that
\begin{equation}
 \label{eq:eta}
 \eta \le e^{-2A},
\end{equation}
which, by \eqref{eq:ep_lower}, implies
\begin{equation}
 \label{eq:eta_ep}
\eta \le C \ep^{12}.
\end{equation}
The last parameters are $y$ and $\zeta$, which are defined as in Lemma~\ref{lem:y_zeta}, with $\eta$ there being the $\eta$ defined above. Note that the parameters $\eta$, $y$ and $\zeta$ appeared already in \cite{Berestycki2010} and had the same meaning there.

We can now proceed to the definition of the process. Recall the definition of $\mu$ in \eqref{eq:mu}.
As in Section \ref{sec:interval_notation}, we denote by $\P^x$ and $\E^x$ the law and expectation of branching Brownian motion with drift $-\mu$ starting from a particle at the point $x\in\R$; we extend this definition to general initial distributions of particles according to a counting measure $\nu$. Recall from Section \ref{sec:preliminaries_definition} that $\mathscr N(t)$ denotes the set of individuals alive at time $t$. We want to absorb the particles at 0 and do this formally by setting
\[
 \mathscr N_0(t) = \{u\in \mathscr N(t): H_0(X_u) > t\},
\]
where $H_0$ is again the hitting time functional of 0.

\begin{figure}
 \centering
\begin{tabular}{cc}
\executeiffilenewer{sec6_1.svg}{sec6_1.pdf}%
{inkscape -z -D --file=sec6_1.svg %
--export-pdf=sec6_1.pdf --export-latex --export-area-drawing}%
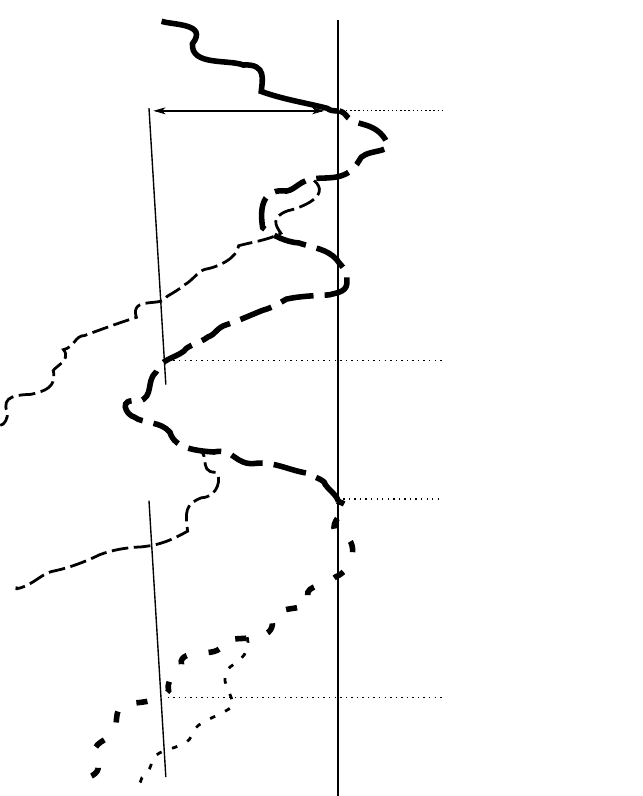%
 & \def\svgwidth{8cm}%
\executeiffilenewer{sec6_1b.svg}{sec6_1b.pdf}%
{inkscape -z -D --file=sec6_1b.svg %
--export-pdf=sec6_1b.pdf --export-latex --export-area-drawing}%
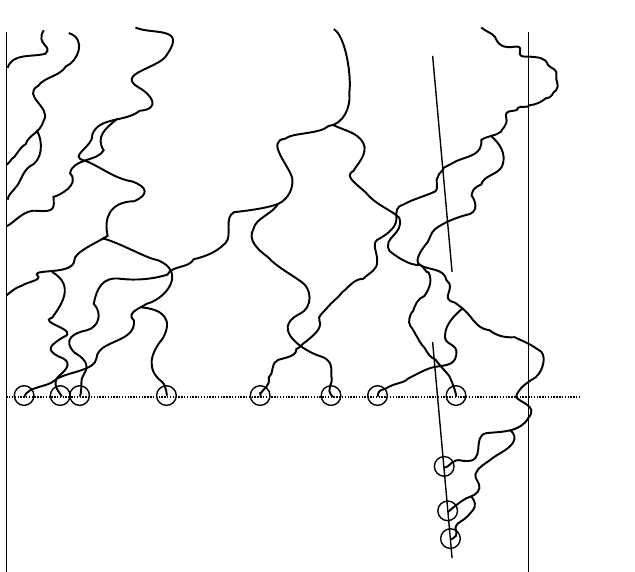%

\end{tabular}
 
\caption{Left picture: A graphical view of the \emph{tiers}. The trajectory of one particle $u$ is singled out (thick line) and the times $\tau_l(u)$ and $\sigma_l(u)$ are shown. The tier 0 particles are drawn with straight lines, the tier 1 particles with dashed and the tier 2 particles with dotted lines. Right picture: The stopping line $\mathscr N_t$ (encircled particles).}
\label{fig:tiers}
\end{figure}

Instead of absorbing particles at $a$, we are now going to classify them into \emph{tiers} as described in the introduction. Let $u\in U$, $t\ge 0$. We define two sequences of random times $(\tau_n(u))_{n\ge-1}$ and $(\sigma_n(u))_{n\ge0}$ by $\tau_{-1}(u) = 0$, $\sigma_0(u) = 0$ and for $n\ge 0$:
\begin{equation}
 \label{eq:def_tau_sigma}
\begin{split}
 \tau_n(u) &= \inf\{s\ge \sigma_n(u):X_u(s) = a\},\\
 \sigma_{n+1}(u) &= \inf\{s\ge \tau_n(u):X_u(s)=a-y+(1-\mu)(s-\tau_n(u))\},
\end{split}
\end{equation}
where we set $\inf \emptyset = \infty$. See Figure~\ref{fig:tiers} for a graphical description. We now define for $t\ge 0$ the stopping lines
\begin{align}
 \label{eq:def_Rline}
\mathscr{R}^{(l)}_t &= \{(u,s)\in U\times[0,t]: s=\tau_l(u)\tand u\in\mathscr{N}_0(s)\},\ l \ge -1,\tand\\
 \label{eq:def_Lline}
\mathscr{S}^{(l)}_t &= \{(u,s)\in U\times[0,t]: s=\sigma_l(u)\tand u\in\mathscr{N}_0(s)\},\ l\ge 0.
\end{align}

That means, $\mathscr{R}^{(l)}_t$ contains the particles of tier $l$ at the moment at which they touch the right barrier and $\mathscr{S}^{(l)}_t$ contains the particles of tier $l$ at the moment at which they come back to the critical line. Note that the sets $\mathscr{R}^{(l)}_t$ and $\mathscr{S}^{(l)}_t$ are increasing in $t$ and $\mathscr{R}^{(l-1)}_t \preceq \mathscr{S}^{(l)}_t \preceq \mathscr{R}^{(l)}_t$ for every $l\ge0$. We also set
\[
R^{(l)}_t = \#\mathscr{R}^{(l)}_t.
\]
In order to extend the definitions of the variables $Z_t$ and $Y_t$ to the current setup, we could simply replace $\mathscr N_{0,a}(t)$ by $\mathscr N_0(t)$ in their definition (see Section \ref{sec:interval_notation}). However, it will be more useful to take special care of the individuals $u$ for which $\tau_l(u) \le t < \sigma_{l+1}(u)$ for some $l\ge 0$, since these are in some kind of ``intermediary'' state which is difficult to analyse. We therefore define the stopping lines
\begin{multline}
 \label{eq:def_Nline_l}
 \mathscr{N}^{(l)}_t = \{(u,s)\in U\times \R_+: u\in\mathscr{N}_0(s),\ \tau_{l-1}(u) \le s < \tau_{l}(u),\\
\text{ and either }t<\sigma_l(u) = s\tor \sigma_l(u) \le t = s\},\ l\ge 0,
\end{multline}
and
\begin{equation}
\label{eq:def_Nline}
 \mathscr{N}_t = \bigcup_{l\ge 0}  \mathscr{N}^{(l)}_t.
\end{equation}
In other words, the stopping line $\mathscr{N}^{(l)}_t$ contains the particles of tier $l$ that have already come back to the critical line at time $t$, as well as the descendants of those that haven't, at the moment at which they hit the critical line. We then define for $l\ge 0$ (recall the definitions of $w_Z$ and $w_Y$ from Section~\ref{sec:interval}),
\[
 Z_t^{(l)} = \sum_{u\in\mathscr{N}^{(l)}_t} w_Z(X_u(t)),\quad Y_t^{(l)} = \sum_{u\in\mathscr{N}^{(l)}_t} w_Y(X_u(t)),
\]
Furthermore, for any symbol $S$ and $0\le k\le l\le \infty$, we write,
\[
 S^{(k;l)} = \sum_{i=k}^l S^{(i)},\quad S^{(l+)} = S^{(l;\infty)},\quad S = S^{(0+)}.
\]
For a particle $(u,s) \in \mathscr{R}^{(l)}_t$, we now define the stopping line
\begin{equation*}
  \mathscr S^{(u,s)} = \{(v,r)\in U\times \R_+: v\in \mathscr N_0(r),\ (u,t)\preceq (v,r) \tand r = \sigma_{l+1}(v)\}.
\end{equation*}
This stopping line yields a collection $(X_v(r),r-s)_{(v,r)\in \mathscr S^{(u,s)}}$ of space-time points and we denote by $Z^{(u,s)}$, $Y^{(u,s)}$ and $W_y^{(u,s)}$ the quantities from Lemma~\ref{lem:critical_line} corresponding to this collection of points (in particular, $W_y^{(u,s)} = ye^{- y}\#\mathscr S^{(u,s)}$). Of course, we have chosen the stopping line in such a way that the variable $W_y^{(u,s)}$ follows the same law as the variable $W_y$ defined in \eqref{eq:def_Wy}. We also define $\tau_{\mathrm{max}}^{(u,s)} = \max_{(v,r)\in \mathscr S^{(u,s)}} (r-s)$. We then define the ``good'' event
\begin{multline}
  \label{eq:116}
  G^{(u,s)} = \{\#\mathscr S^{(u,s)} \le \zeta\} \cap \{\mathscr S^{(u,s)} \subseteq U\times [1,\zeta]\}\\ \cap\{\sup_{s\le r\le \zeta} \#\{v\in\mathscr N_0(r): (u,s)\preceq (v,r)\preceq \mathscr S^{(u,s)}\} \le \zeta\}
\end{multline}
and the event of a \emph{breakout},
\begin{equation}
 \label{eq:def_B}
B^{(u,s)} = \{Z^{(u,s)} > \ep e^A\}\cup (G^{(u,s)})^c,
\end{equation}
(the inclusion of the ``bad'' event $(G^{(u,s)})^c$ is for technical reasons). If the event $B^{(u,s)}$ occurs, the particle $u$ is then also called the \emph{fugitive}. We set
\begin{equation}
 \label{eq:def_pB}
p_B = \P^a(B^{(\emptyset,0)}),
\end{equation}
and define the law of BBM started at $a$ with the first particle conditioned not to break out:
\[
 \Q^a(\cdot) = \P^a(\cdot\,|\, B^c) = \frac{\P^a(\cdot,\ B^c)}{1-p_B},
\]
where we set $B = B^{(\emptyset,0)}$. We further set $Z = Z^{(\emptyset,0)}$ and $W_y = W_y^{(\emptyset,0)}$ and note that by Lemmas~\ref{lem:critical_line} and \ref{lem:y_zeta}, we have for large $a$,
\begin{equation}
 \label{eq:Z_W}
\P^a(|Z-\pi W| > 2\eta) + \P^a((G^{(\emptyset,0)})^c) < \eta,
\end{equation}
where $W$ is defined as in \eqref{eq:def_Wy}. Hence, by \eqref{eq:eta}, \eqref{eq:W_tail} and Lemma~\ref{lem:y_zeta}, we get
\begin{equation}
 \label{eq:pB}
p_B = \Big(\pi + o(1)\Big)\frac 1 {\ep e^A},
\end{equation}
which goes to 0 as $A$ and $a$ go to infinity, by \eqref{eq:ep_lower}. Furthermore, \eqref{eq:Z_W} yields for large $A$ and $a$,
\begin{equation}
\begin{split}
\label{eq:QaZ}
 \Q^a[Z] &= (\E[\pi W\Ind_{(\pi W\le \ep e^A(1+o(1)) + O(\eta))}] + O(\eta \ep e^A))(1+O(p_B)) \\
&= \pi(A+\log\ep+\const{eq:QaZ}+ o(1)),
\end{split}
\end{equation}
by \eqref{eq:W_expec}, \eqref{eq:ep_upper}, \eqref{eq:eta} and \eqref{eq:pB}. In particular, we have for $A\ge 1$ and large $a$,
\begin{equation}
 \label{eq:QaZC}
\Q^a[Z] \le CA.
\end{equation}
Moreover, by \eqref{eq:W_tail}, \eqref{eq:Z_W} and \eqref{eq:eta}, we have for $A\ge 1$ and large $a$,
\begin{equation}
 \label{eq:QaZsquared}
\Q^a[Z^2] \le C\ep e^A.
\end{equation}
Finally, note that by Lemma~\ref{lem:critical_line}, we have $Y\le \eta Z$, $\Q^a$-almost surely, a fact that will often be used without further reference.

We now define for every $l\in\N$ the time of the first breakout of a particle of tier $l$,
\begin{equation}
 \label{eq:def_Tl}
T^{(l)}(\omega) = \inf\{t\ge0: \omega \in \bigcup_{(u,s)\in\mathscr{R}^{(l)}_t} B^{(u,s)}\},
\end{equation}
and set
\begin{equation}
\label{eq:def_T}
 T^{(0;l)} = \min_{0\le j\le l} T^{(l)},\quad\text{with }T = T^{(0;\infty)} = \min_{j\ge 0} T^{(l)}.
\end{equation}

Now fix $t > 0$ and $l\in\N\cup\{\infty\}$. We want to describe the system conditioned on $T^{(0;l)}=t$. For this, suppose that at time $0$ the particles are distributed according to a counting measure $\nu = \sum_{i=1}^n \delta_{x_i}$. We denote by $\mathscr U$ the fugitive of the breakout that happened at time $T^{(0;l)}$ and define $p_i = \P^\nu(i \preceq \mathscr U\,|\,T^{(0;l)} = t)$. This yields a law $(p_i)_{i=1}^n$ on the initial particles, depending on $\nu$ and $t$. Since the variable $T^{(0;l)}$, the time of the first breakout, is the minimum of the variables $T^{(0;l)}_i$, $i=1,\ldots,n$, the times of the first breakout of the BBM descending from the particle $i$, we can decompose the process into
\begin{equation}
\label{eq:decomposition}
 \P^\nu\Big(\prod_{i=1}^n \dd \omega^{(i)}\,\Big|\,T^{(0;l)}=t\Big) = \sum_{i=1}^n p_i \times\P^{x_i}(\dd \omega^{(i)}\,|\,T^{(0;l)}=t) \times \prod_{j\ne i} \P^{x_j}(\dd \omega^{(j)}\,|\,T^{(0;l)}>t).
\end{equation}
That is, we first choose according to the law $(p_i)_{i=1}^n$ the initial particle that is going to cause the breakout. This particle spawns a BBM conditioned to break out at time $t$. The remaining particles spawn independent BBM conditioned not to break out before time $t$.

\begin{remark}
\label{rem:hat_f}
  Note that many results in this section can be done for BBM with varying drift given by a barrier function $f$. For example, with Lemma~\ref{lem:R_f} one gets immediately that $\P_f(T^{(0)}>t) \ge \P(T^{(0)} > t)$ for all $t\ge0$ and it will be clear from the next section that in fact $\P_f(T>t) \ge \P(T>t)$ as well. However, in order to simplify notation and because we will not need the results often in this generality, we state them here only for BBM with fixed drift.
\end{remark}

\subsection{The time of the first breakout}
\label{sec:time_breakout}

For large $A$ and $a$, define
\begin{equation}
 \label{eq:tcrit}
 \tconst{eq:tcrit}= \tfrac{1}{20 A} a^3.
\end{equation}

We want to prove that the random variable $T$ defined in the previous section is approximately exponentially distributed with parameter $p_B\pi Z_0/a^3$, which is the statement of the following proposition:

\begin{proposition}
\label{prop:T}
 Let $0\le t\le  \tconst{eq:tcrit}$. For $A$ and $a$ large enough, we have
\begin{equation}
\label{eq:T}
\P(T> t) = \exp\Big(-\pi p_B Z_0\frac t {a^3}\Big(1 + O(At/a^3 + p_B)\Big)+O(p_B Y_0)\Big).
\end{equation}
\end{proposition}

The proof proceeds by a sequence of lemmas. Lemma~\ref{lem:T0} gives a estimate on $\P(T^{(0)} > t)$. This is used in Lemma~\ref{lem:T_a}, in order to obtain an estimate on $\Q^a(T > t)$, using a recursive argument. Finally, Proposition \ref{prop:T} is proven by combining Lemmas~\ref{lem:T0} and \ref{lem:T_a}.

\begin{lemma}
 \label{lem:T0}
Let $0\le t\le a^3$.  Suppose that $p_B \le 1/2$. Then,
\begin{equation}
\label{eq:T0}
\P(T^{(0)}> t) = \exp\Big( -\pi p_B Z_0 \frac t {a^3} \big(1+O(p_B)\big) + O(p_BY_0)\Big).
\end{equation}
\end{lemma}
\begin{proof}
 Let $x_1,\ldots,x_n$ be the positions of the initial particles. Since the initial particles spawn independent branching Brownian motions, we have
\begin{equation}
\label{eq:340}
 \P(T^{(0)} > t) = \prod_i \P^{x_i}(T^{(0)}>t).
\end{equation}
We have for every $x\in(0,a)$,
\begin{equation}
 \label{eq:308}
  \P^x(T^{(0)} > t) = \E^x\Big[\prod_{(u,s)\in\mathscr{R}^{(0)}_t} \Ind_{B^{(u,s)}}\Big] = \E^x\Big[(1-p_B)^{R^{(0)}_t}\Big],
\end{equation}
since by the strong branching property, the random variables $Z^{(u,s)}$ are independent conditioned on $\mathscr{R}^{(0)}_t$. By Lemma~\ref{lem:Rt} and the assumption $t\le a^3$, we have
\begin{align}
 \label{eq:310}
&|\E^x[R^{(0)}_t] - \pi w_Z(x) t/a^3| \le C w_Y(x),\\
 \label{eq:315}
&\E^x[(R^{(0)}_t)^2] \le C (w_Z(x) t/a^3+w_Y(x)).
\end{align}
By Jensen's inequality, \eqref{eq:310} and the inequality $|\log(1-z)| \le z+z^2$ for $z\le 1/2$,
\begin{equation}
 \label{eq:317}
 \E^x\Big[(1-p_B)^{R^{(0)}_t}\Big] \ge \exp\Big( -\pi p_B w_Z(x) \frac t {a^3} (1+ Cp_B) - Cp_Bw_Y(x)\Big).
\end{equation}
Furthermore, the inequality $(1-p)^n \le 1-np+n(n-1)p^2/2$ and Equations \eqref{eq:310} and \eqref{eq:315} give
\begin{equation}
\label{eq:325}
 \E^x\Big[(1-p_B)^{R^{(0)}_t}\Big] \le 1-\pi p_B w_Z(x) \frac t {a^3}  (1- Cp_B) + Cp_Bw_Y(x).
\end{equation}
The lemma now follows from \eqref{eq:308}, \eqref{eq:317} and \eqref{eq:325} together with the inequality $1+z \le e^z$ for $z\in\R$.
\end{proof}

In the following lemma, note that according to the definition of the tiers, a particle starting at $a$ starts immediately in tier 1.
\begin{lemma}
 \label{lem:T_a}
Let $0\le t\le \tconst{eq:tcrit}$. Then, for large $A$ and $a$,
\begin{equation}
 \label{eq:T_a}
\Q^a(T > t) \ge \exp\Big(-C p_B A (\tfrac t {a^3} + \eta)\Big).
\end{equation}
\end{lemma}
\begin{proof}
We have
\begin{equation}
 \label{eq:360}
 \Q^a(T > t) = \Q^a\Big[\prod_{(u,s)\in\mathscr{S}^{(1)}_t}\P^{X_u(s)}(T > t-s)\Big] \ge \Q^a[\P^\nu(T > t)],
\end{equation}
where $\nu = \sum_{(u,s)\in\mathscr{S}^{(1)}_t} \delta_{X_u(s)}$. Since $T>t$ implies $T^{(0)}>t$, we have
\begin{equation}
 \label{eq:362}
\P^\nu(T > t) = \P^\nu(T > t|T^{(0)}>t) \P^\nu(T^{(0)} > t).
\end{equation}
Let $Z=Z^{(\emptyset,0)}$ and $Y=Y^{(\emptyset,0)}$, such that $Y \le \eta Z$, $\Q^a$-almost surely, by Lemma~\ref{lem:critical_line} and the definition of the ``bad'' event $B^{(\emptyset,0)}$. By Lemma~\ref{lem:T0}, we have for large $A$,
 \begin{equation}
\label{eq:363}
\P^\nu(T^{(0)} > t) \ge \exp\Big(- C p_B Z\left(\tfrac t {a^3}+\eta\right)\Big).
\end{equation} 
Furthermore, with the notation from Section \ref{sec:weakly_conditioned}, with $p(s) \equiv p_B$,
\begin{align*}
 \P^\nu(T > t|T^{(0)}>t) &= \widetilde{\P}^\nu(T > t) = \widetilde{\P}^\nu\Big[\prod_{(u,s)\in\mathscr{R}^{(0)}_t}\Q^a(T>t-s)\Big] \ge \widetilde{\P}^\nu\Big[\Q^a(T>t)^{R^{(0)}_t}\Big].
\end{align*}
By Jensen's inequality and Lemmas~\ref{lem:conditioned_Rt} and \ref{lem:Rt}, this implies
\begin{equation}
\label{eq:364}
 \P^\nu(T > t|T^{(0)}>t) \ge \Q^a(T>t)^{\widetilde{\E}^\nu[R^{(0)}_t]} \ge \Q^a(T>t)^{\pi Z (t/a^3 + O(\eta))}.
\end{equation}
Equations \eqref{eq:360}, \eqref{eq:362}, \eqref{eq:363} and \eqref{eq:364}, together with Jensen's inequality and \eqref{eq:QaZ}, now yield for large $A$ and $a$,
\begin{equation}
\label{eq:370}
 \Q^a(T > t) \ge \Q^a(T>t)^{\pi^2A (t/a^3 + O(\eta))}\times\exp\Big(- C p_B A\left(\tfrac t {a^3}+\eta\right)\Big).
\end{equation}
By the hypothesis on $t$ and \eqref{eq:eta}, the exponent of $\Q^a(T>t)$ in \eqref{eq:370} is smaller than $1/2$ for large $A$ and $a$. This yields the statement.
\end{proof}

\begin{proof}[Proof of Proposition \ref{prop:T}]
The upper bound follows from Lemma~\ref{lem:T0} and the trivial inequality $T\le T^{(0)}$. For the lower bound, we note that as in the proof of Lemma~\ref{lem:T_a}, we have by Jensen's inequality and Lemma \ref{lem:conditioned_Rt},
\begin{equation}
\label{eq:380}
 \P(T > t) = \P(T > t\,|\,T^{(0)} > t)\P(T^{(0)} > t) \ge \Q^a(T > t)^{\E[R_t^{(0)}]}\P(T^{(0)} > t).
\end{equation}
By Lemmas~\ref{lem:T_a} and \ref{lem:Rt}, we have
\begin{equation}
\label{eq:382}
 \Q^a(T > t)^{\E[R_t^{(0)}]} \ge \exp\Big(-C p_B A (\tfrac t {a^3} + \eta)(\tfrac t {a^3} Z_0 + Y_0)\Big),
\end{equation}
The lower bound in \eqref{eq:T} now follows from \eqref{eq:380}, \eqref{eq:382} and Lemma~\ref{lem:T0}, together with the hypothesis on $t$, \eqref{eq:eta} and \eqref{eq:pB}.
\end{proof}

\begin{lemma}
\label{lem:moments_T}
Define $\gamma = (\pi p_B Z_0)^{-1}$. Suppose that $Y_0\le C$ and let $\alpha \ge 0$ and $n\in\N$. Then, for large $A$, for every $l\in\N\cup\{\infty\}$,
\begin{equation}
\label{eq:T_moments}
 \E[(T^{(0;l)}/a^3 + \alpha)^n\Ind_{(T^{(0;l)}\le a^3)}] \le C\sum_{k=0}^n\frac{n!\alpha^k(2\gamma)^{n-k}}{k!}
\end{equation}
Furthermore, if $0\ge \delta = A^{-1}o(1)$, then for large $A$ and $a$,
\begin{equation}
 \label{eq:T_expec}
 \E[(T/a^3)\Ind_{(T \le \delta a^3)}] = \gamma(1+O(A\delta + p_B)) + O((\delta+\gamma)e^{-O(\delta/\gamma)})
\end{equation}
\end{lemma}
\begin{proof}
We first note that we have, for $n\ge 1$ and $\rho > 0$,
\begin{equation}
 \label{eq:398}
\int_0^\infty (t+\alpha)^{n-1} e^{-\rho t}\,\dd t = \sum_{k=0}^{n-1} \frac{(n-1)!\alpha^k}{k!\rho^{n-k}}.
\end{equation}
Now, we have
\[
\begin{split}
\E[(T^{(0;l)}/a^3 + \alpha)^n\Ind_{(T^{(0)}\le a^3)}] &= \int_0^1 (t+\alpha)^n\,\P(T^{(0;l)}/a^3 \in \dd t) \\
&\le n\int_0^1 (t+\alpha)^{n-1}\,\P(T^{(0;l)}>ta^3)\,\dd t + \alpha^n.
\end{split}
\]
The inequality \eqref{eq:T_moments} now follows from Lemma~\ref{lem:T0} and \eqref{eq:398} and the hypothesis on $Y_0$. For the second part, we note that
\[
 \E[(T/a^3)\Ind_{(T \le \delta a^3)}] = \int_0^{\delta} \P(T > ta^3)\,\dd t - \delta \P(T > \delta a^3),
\]
and by Proposition \ref{prop:T} and the hypothesis on $Y_0$, we have for $t\le \delta$ and large $A$ and $a$,
\[
 \P(T > ta^3) = (1+O(p_B))\exp(-\gamma^{-1} t(1+O(A\delta + p_B))).
\]
Equation \eqref{eq:T_expec} now follows from the last two equations.
\end{proof}

We now show how we can couple the variable $T$ with an exponentially distributed variable:
\begin{lemma}
\label{lem:coupling_T}
 Suppose that $e^{-A}Z_0 = 1 + O(\ep^{3/2})$ and that $Y_0 \le \eta Z_0$. Then there exists a coupling $(T,V)$, such that $V$ is a random variable which is exponentially distributed with parameter $\pi p_Be^A$, $T$ is $\sigma(V)$-measurable and $\P(|T/a^3 - V| > \ep^{3/2}) \le C\ep^2$ for large $A$ and  $a$.
\end{lemma}

\begin{proof}
For brevity, set $\rho := \pi p_Be^A$. Let $F$ be the tail distribution function of $T$, i.e. $F(t) := \P(T \ge t)$. The number of individuals in a BBM tree being at most countable, $T$ has no atoms except $\infty$. We can therefore define a random variable $U$ which is uniformly distributed on $(0,1)$ by setting
\[
  U = F(T)\Ind_{(T<\infty)} + U' F(\infty)\Ind_{(T=\infty)},
\]
where $U'$ is a uniformly distributed random variable on $(0,1)$, independent of $T$. Now we define $V = -\rho^{-1} \log U$. Then $V$ is exponentially distributed with parameter $\rho$ and $T = F^{-1}(e^{-\rho V})$, where $F^{-1}$ denotes the generalised right-continuous inverse of $F$. Hence, $T$ is $\sigma(V)$-measurable. On $\{T<\infty\}$, we have by Proposition \ref{prop:T}, for $a$ large enough,
\begin{equation}
\label{eq:390}
\begin{split}
 V &= -\rho^{-1} (\pi p_Be^A e^{-A} Z_0 T/a^3(1+ O(AT/a^3 + p_B))) + O(p_BY_0)\\
&= T/a^3(1 + O(\ep^{3/2} + AT/a^3 + p_B)) + O(p_B e^A\eta),
\end{split}
\end{equation}
by the hypotheses on $Z_0$ and $Y_0$. Hence, by \eqref{eq:ep_lower}, \eqref{eq:eta_ep} and \eqref{eq:pB}, we have for $a$ large enough,
\[
 |T/a^3 - V| = O(\ep^{3/2}T/a^3 + A(T/a^3)^2) + O(\ep^{2}).
\]
But now we have by Lemma~\ref{lem:T0}, for large $A$ and $a$,
\begin{equation}
\label{eq:392}
 \P(T/a^3 > \ep^{7/8}) \le \P(T^{(0)}/a^3 > \ep^{7/8}) \le Ce^{-O(\ep^{-1/8})} \le \ep^2/2.
\end{equation}
The statement now follows from \eqref{eq:390} and \eqref{eq:392} together with \eqref{eq:ep_upper}
\end{proof}

\subsection{The particles that do not participate in the breakout}
\label{sec:before_breakout_particles}

In this section, we fix $\tau\le  \tconst{eq:tcrit}$. We define the law of BBM conditioned not to break out before $\tau$ in the tiers $0,\ldots,l$ by
\[
 \Phat_l(\cdot) = \P(\cdot\,|\,T^{(0;l)} > \tau),\quad\text{with }\Phat = \Phat_\infty.
\]
Expectation w.r.t. $\Phat_l$ is denoted by $\Ehat_l$. Under the law $\Phat_l$, the process stopped at $\mathscr L_H$ then follows the law $\Ptilde$ from Section \ref{sec:weakly_conditioned}, with
\begin{equation}
 \label{eq:def_ps}
p(t) = p_B\Ind_{(t\le \tau)} + (1-p_B)\Q^a(T^{(0;l)} \le \tau-t),
\end{equation}
such that by Lemma~\ref{lem:T_a}, \eqref{eq:eta}, \eqref{eq:QaZC} and the trivial inequality $T^{(0;l)} \ge T$, we have for large $A$ and $a$,
\begin{equation}
 \label{eq:pbar_estimate}
\|p\|_\infty  \le C p_B.
\end{equation}
In particular, by Lemma~\ref{lem:h_estimate}, \eqref{eq:c0_mu} and \eqref{eq:pbar_estimate}, we have for large $A$ and $a$,
\begin{equation}
  \label{eq:hbar}
 (h(x,0))^{-1} \le 1+Cp_B.
\end{equation}
Lemmas \ref{lem:many_to_one_conditioned} and ~\ref{lem:conditioned_Z_lowerbound}, together with \eqref{eq:pbar_estimate} and \eqref{eq:hbar} now immediately imply the following:

\begin{corollary}
  \label{cor:hat_many_to_one}
 Let  $x\in (0,a)$, $0\le t\le \tau\le \tconst{eq:tcrit}$ and  $g:[0,a]\to \R_+$ be measurable with $g(0) = g(a) = 0$. Define $S^{(0)}_t = \sum_{u\in \mathscr N^{(0)}_t} g(X_u(t))$ or $S^{(0)}_t = \R^{(0)}_t$. Then, with $\Ehat_l$ as above,
\begin{equation*}
 \Ehat_l^x[S^{(0)}_t] = (1+O(p_B)) \E^x[S^{(0)}_t] \quad\tand\quad \Ehat_l^x[(S^{(0)}_t)^2] \le (1+O(p_B)) \E^x[(S^{(0)}_t)^2].
\end{equation*}
\end{corollary}
Moreover, as in the proof of Lemma~\ref{lem:conditioned_Rt}, one can show that
\begin{equation}
 \label{eq:QhatZ}
\Qhat^a_l[Z] = (1+O(p_B)) \Q^a[Z]\quad\tand\quad \Qhat^a_l[Z^2] \le (1+O(p_B)) \Q^a[Z^2].
\end{equation}

We define two filtrations $(\G_l)_{l\ge 0}$ and $(\H_l)_{l\ge 0}$ by 
\[
 \G_l = \F_{\mathscr{S}^{(l)}_\tau},\quad \H_l = \F_{\mathscr{R}^{(l)}_\tau},
\]
such that $\G_l \subset \H_l \subset \G_{l+1}$ for every $l$. Now define for measurable $I\subset [0,\tau]$,
\[
 Z_{\emptyset,I}^{(l)} = \sum_{(u,t)\in\mathscr{S}^{(l)}_\tau\cap (U\times I)} Z^{(u,t)},\quad Y_{\emptyset,I}^{(l)} = \sum_{(u,t)\in\mathscr{S}^{(l)}_\tau\cap (U\times I)} Y^{(u,t)},
\]
with $Z_\emptyset^{(l)} = Z_{\emptyset,[0,\tau]}^{(l)}$ and $Y_\emptyset^{(l)} = Y_{\emptyset,[0,\tau]}^{(l)}$. Furthermore, recall from Section \ref{sec:before_breakout_definitions} the definition of $Z_\emptyset^{(k;l)}$, $Z_\emptyset^{(l+)}$ and $Z_\emptyset$ and the corresponding quantities for $Y$.
\begin{lemma}
 \label{lem:Zl_exp_upperbound}
We have for all $0\le k\le l$, $x\in (0,a)$ and large $A$ and $a$,
\begin{equation}
 \label{eq:Zl_exp_upperbound}
 \Ehat_l[Z_\emptyset^{(k+1)}\,|\,\G_{k}] \le (\pi+Cp_B)\Q^a[Z](\tfrac \tau {a^3} Z_\emptyset^{(k)} + C Y_\emptyset^{(k)}).
\end{equation}
In particular, for large $A$ and $a$,
\begin{equation}
 \label{eq:Zl_exp_upperbound_2}
 \Ehat_l[Z_\emptyset^{(k+1)}] \le C A(\tfrac \tau {a^3} Z_0 + C Y_0)\Big(\pi^2A(\tfrac \tau {a^3} + \Const{eq:Zl_exp_upperbound_2} \eta)\Big)^{k}.
\end{equation}
In the case $k=0$, we also have for large $A$ and $a$,
\begin{equation}
 \label{eq:Z1_exp_lowerbound}
 \Ehat_l[Z_\emptyset^{(1)}] \ge (\pi-C p_B)\Q^a[Z](\tfrac \tau {a^3}Z_0 - CY_0) .
\end{equation}
\end{lemma}

\begin{proof}
We have for $0\le k\le l$,
\begin{equation}
\label{eq:Zl_Hl}
 \Ehat_l[Z_\emptyset^{(k+1)}\,|\,\H_{k}] = \Ehat_l\Big[\sum_{(u,t)\in\mathscr R^{(k)}_\tau} Z^{(u,t)}\,\Big|\,\H_{k}\Big] = \Qhat_{l-k}^a[Z] R_\tau^{(k)},
\end{equation}
since conditioned on $\H_{k}$, the random variables $Z^{(u,t)}$ are iid under $\Phat_l$ of the same law as $Z$ under $\Qhat_{l-k}^a$ and independent of $\H_{k}$, by the strong branching property. Now, we have
\[
 \Ehat_l[R_\tau^{(k)}\,|\,\G_{k}] = \sum_{(u,t)\in\mathscr S^{(k)}_\tau} \E^{X_u(t)}[R_{\tau-t}^{(0)}\,|\,T > \tau-t],
\]
such that by Lemma~\ref{lem:Rt} and Corollary~\ref{cor:hat_many_to_one},
\begin{equation}
\label{eq:400}
\Ehat_l[R_\tau^{(k)}\,|\,\G_{k}] \le (1+O(p_B))\Big(\pi \frac \tau {a^3} Z_\emptyset^{(k)} + \Const{eq:Rt_nu_expec} Y_\emptyset^{(k)}\Big).
\end{equation}
Equations \eqref{eq:Zl_Hl} and \eqref{eq:400} together with \eqref{eq:QhatZ} give \eqref{eq:Zl_exp_upperbound}. Equation \eqref{eq:Zl_exp_upperbound_2} follows easily from \eqref{eq:Zl_exp_upperbound} by \eqref{eq:QaZC} and the fact that $Y_\emptyset^{(k)} \le \eta Z_\emptyset^{(k)}$, $\Phat_l$-almost surely for $1\le k\le l$. Now, in the case $k=0$, we have $\G_0 = \F_0$ by definition. Denote the positions of the initial particles by $x_1,\ldots,x_n$. By Lemmas~\ref{lem:Rt}, \ref{lem:conditioned_Rt} and \eqref{eq:pbar_estimate},
\[
\begin{split}
 \Ehat_l[R_\tau^{(0)}] = \sum_{i=1}^n \Ehat_l^{x_i}[R_\tau^{(0)}] &\ge \sum_{i=1}^n \E^{x_i}[R_\tau^{(0)}] - \|p\|_\infty  \E^{x_i}[(R_\tau^{(0)})^2] \\
&\ge \pi \frac \tau {a^3} Z_0 - \Const{eq:Rt_nu_expec} Y_0 - \|p\|_\infty  \Const{eq:Rt_nu_variance}\left(\tfrac \tau {a^3} Z_0 +Y_0\right),
\end{split}
\]
which yields \eqref{eq:Z1_exp_lowerbound}.
\end{proof}

%
%
%
%
%

\begin{lemma}
\label{lem:hat_quantities}
Suppose that $t\le \tau$. Then for large $A$ and $a$, we have
\begin{equation*}
 \Ehat[Z_t] = Z_0\left(1+\pi \Q^a[Z] \tfrac t {a^3} + O\left(p_B + \left(A\tfrac t {a^3}\right)^2\right)\right)+O(AY_0),\quad \Ehat[Y_t] \le C \left(Y_0 + \eta A\tfrac t {a^3}Z_0\right).
\end{equation*} 
If moreover $t\ge 2a^2$, then
\begin{equation*}
\Ehat[Y_t\Ind_{(R_{[t-a^2,t]} = 0)}] \le \frac C a (Z_0 + AY_0)\quad\tand\quad \Phat(R_{[t-a^2,t]} \ne 0) \le C\eta(Y_0+ Z_0).
\end{equation*} 
\end{lemma}
\begin{proof}
  First note that we have $\inf_{x\in[0,a]}h(x,0) \ge 1/2$ for large $A$, by \eqref{eq:hbar}. Furthermore, by the hypothesis $\tau \le \tconst{eq:tcrit}$ and \eqref{eq:eta}, we have $\pi^2A(\tfrac \tau {a^3} + C \eta)\le 1/2$ for large $A$ and $a$. The first two inequalities now follow from Proposition~\ref{prop:quantities}, Corollary~\ref{cor:hat_many_to_one},  Lemma~\ref{lem:Zl_exp_upperbound}, \eqref{eq:QhatZ} and \eqref{eq:eta} by summing over $k$. In particular, if $t\ge 2a^2$, then for large $A$ and $a$, $\Ehat[Z_{t-a^2}] \le C(Z_0 + AY_0)$ and $\Ehat[Y_{t-a^2}] \le C\eta(Z_0+Y_0)$. Together with Proposition~\ref{prop:quantities}, Lemma~\ref{lem:Rt} and Corollary~\ref{cor:hat_many_to_one}, this proves the other two inequalities.
\end{proof}

In order to estimate second moments, we will make use of the following extension to the Many-to-two lemma (Lemma~\ref{lem:many_to_two_fixedtime}). For $x,z\in(0,a)$ and $0\le t\le\tau$, we define $\widehat{m_2}(x,t)$ to be the quantity $\widetilde{m_2}(x,t)$ from Section~\ref{sec:weakly_conditioned} corresponding to the penalisation from \eqref{eq:def_ps} and  $\widehat{\p}^{(l)}_t(x,z)$ to be the density of tier $l$ particles at position $z$ and time $t$ under the law $\Phat^x$, not counting the particles $u$ with $t < \sigma_l(u)$. Then set $\widehat{\p}_t = \widehat{\p}^{(0+)}_t$.

\begin{lemma}
  \label{lem:many_to_two_with_tiers}
Let $w:[0,a]\to\R_+$ be measurable and define $S_t = \sum_{(u,t)\in \mathscr N_t} w(X_u(t))$. Then
\begin{multline}
  \label{eq:many_to_two_with_tiers}
  \Ehat^x[S_t^2] = \Ehat^x[\sum_{(u,s)\in\mathscr N_t} w(X_u(s))^2] + \beta_0 \int_0^t\int_0^a
   \widehat{m_2}(z,s) \widehat{\p}_s(x,z) \Ehat^{(z,s)}[S_t]^2\,\dd z\,\dd s\\
+ \Ehat^x\Big[\sum_{(u,s)\in\mathscr R_t} \sum_{(v_1,s_1),(v_2,s_2)\in\mathscr
  S^{(u,s)},\,v_1\ne v_2} \Ehat^{(X_{v_1}(s_1),s_1)}[S_t]\Ehat^{(X_{v_2}(s_2),s_2)}[S_t]\Big].
\end{multline}
\end{lemma}
\begin{proof}
For $(v_1,s_1),(v_2,s_2)\in \mathscr N_t$ we write $(v_1,s_1)\wedge (v_2,s_2)$ for their most recent common ancestor. Then define for $l\ge 1$,
\begin{multline*}
\mathscr A_1^{(l)} = \{((v_1,s_1),(v_2,s_2))\in \mathscr N_t^2: v_1\ne v_2 \text{ and if } (v_1,s_1)\wedge (v_2,s_2) = (v_0,s_0),\\
 \text{ then } \tau_{l-1}(v_0) \le s_0 < \sigma_l(v_0)\},
\end{multline*}
and for $l\ge 0$,
\begin{multline*}
\mathscr A_2^{(l)} = \{((v_1,s_1),(v_2,s_2))\in \mathscr N_t^2: v_1\ne v_2 \text{ and if } (v_1,s_1)\wedge (v_2,s_2) = (v_0,s_0),\\
 \text{ then } \sigma_l(v_0)\le s_0 < \tau_l(v_0)\},
\end{multline*}
We then have
\begin{multline}
  \label{eq:42}
  \Ehat^x\Big[\sum_{((v_1,s_1),(v_2,s_2))\in \mathscr A^{(l)}_1} w(X_{v_1}(s_1))w(X_{v_2}(s_2))\,\Big|\,\mathscr G_l\Big] \\
= \sum_{(u,s)\in\mathscr R^{(l)}_t} \sum_{(v_1,s_1),(v_2,s_2)\in\mathscr S^{(u,s)},\,v_1\ne v_2} \Ehat^{(X_{v_1}(s_1),s_1)}[S_t]\Ehat^{(X_{v_2}(s_2),s_2)}[S_t],
\end{multline}
and by Lemma~\ref{lem:many_to_two_fixedtime} (see also Remark~\ref{rem:many_to_two}),
\begin{multline}
  \label{eq:54}
  \Ehat^x\Big[\sum_{((v_1,s_1),(v_2,s_2))\in \mathscr A^{(l)}_2} w(X_{v_1}(s_1))w(X_{v_2}(s_2))\Big] \\
= \beta_0 \int_0^t\int_0^a m_2(z,s) \widehat{\p}^{(l)}_s(x,z) \Ehat^{(z,s)}[S_t]^2\,\dd z\,\dd s.
\end{multline}
Equations \eqref{eq:42} and \eqref{eq:54}, together with \eqref{eq:2nd_moment_decomposition} and the tower property of conditional expectation yield the lemma.
\end{proof}
\begin{remark}
  \label{rem:many_to_two_with_tiers}
An analogous result holds for $R_t$.
\end{remark}

\begin{lemma}
  \label{lem:Zt_variance}
We have for every $t\le \tau$,
\[
\Varhat(Z_t) \le C\ep e^A\Big(\frac{t}{a^3} Z_0 + Y_0\Big)\quad\tand\quad \Varhat(R_t) \le C\ep e^A\Big(\frac{t}{a^3} Z_0 + Y_0\Big).
\]
\end{lemma}
\begin{proof}
  By Lemma~\ref{lem:hat_quantities} and the hypothesis  $t\le \tau\le \tconst{eq:tcrit}$, we have for every $x\in(0,a)$ and $s\le t$,
  \begin{equation}
    \label{eq:56}
    \Ehat^{(x,s)}[Z_t] \le CA \Big(w_Z(x) \frac t {a^3} + w_Y(x)\Big) \le CAe^{-\frac 3 4 (a-x)}.
  \end{equation}
Now, as in the proof of \eqref{eq:Zt_2ndmoment}, we have by \eqref{eq:56} and Corollary~\ref{cor:hat_many_to_one},
\begin{equation}
  \label{eq:58}
 \int_0^t\int_0^a \widehat{\p}^{(0)}_s(x,z) \Ehat^{(z,s)}[Z_t]^2\,\dd z\,\dd s \le CA^2 \Big(\frac{t}{a^3} w_Z(x) + w_Y(x) \Big),
\end{equation}
which yields
\begin{equation}
\label{eq:61}
\begin{split}
  \int_0^t\int_0^a &\widehat{m_2}(z,s) \widehat{\p}_s(x,z) \Ehat^{(z,s)}[Z_t]^2\,\dd z\,\dd s \\
&\le C \Ehat^x\Big[\sum_{(u,s_0)\in\mathscr S_t} \int_{s_0}^t\int_0^a \widehat{\p}^{(0)}_{s-s_0}(X_u(s_0),z) \Ehat^{(z,s)}[Z_t]^2\,\dd z\,\dd s\Big]\\
 &\le CA^2 \Big(\frac{t}{a^3} \Ehat^x[Z_{\emptyset,[0,t]}] + \Ehat^x[Y_{\emptyset,[0,t]}]\Big)\\
 &\le CA^2 \Big(\frac{t}{a^3} w_Z(x) +  w_Y(x)\Big), 
\end{split}
\end{equation}
by Lemma~\ref{lem:Zl_exp_upperbound} and the hypothesis $t\le \tau\le \tconst{eq:tcrit}$. Furthermore, by \eqref{eq:56},  \eqref{eq:eta} and Lemma~\ref{lem:critical_line}, we have
\begin{multline}
  \label{eq:60}
\Ehat^x\Big[\sum_{(u,s)\in\mathscr R_t} \sum_{(v_1,s_1),(v_2,s_2)\in\mathscr S^{(u,s)},\,v_1\ne v_2} \Ehat^{(X_{v_1}(s_1),s_1)}[Z_t]\Ehat^{(X_{v_2}(s_2),s_2)}[Z_t]\Big] \\
\le C \Q^a[Z^2] \Ehat^x[R_t] \le C\ep e^A \Big(\frac{t}{a^3} w_Z(x) +  w_Y(x)\Big),
\end{multline}
by \eqref{eq:QaZsquared} and Lemmas~\ref{lem:Rt} and Corollary~\ref{cor:hat_many_to_one}. Lemma~\ref{lem:many_to_two_with_tiers}, together with \eqref{eq:61} and \eqref{eq:60} as well as Lemma~\ref{lem:hat_quantities} and the inequality $w_Z^2 \le Cw_Y$ gives
\begin{equation}
  \label{eq:62}
  \Varhat^x(Z_t) \le \Ehat^x[(Z_t)^2] \le C \ep e^A\Big(\frac{t}{a^3} w_Z(x) +  w_Y(x)\Big).
\end{equation}
Summing over the initial particles yields the inequality for $\Varhat(Z_t)$. The proof of the second inequality is analogous, relying on Lemma~\ref{lem:Rt} and Remark~\ref{rem:many_to_two_with_tiers}.
\end{proof}

We finish the section by a corollary which will be useful in the next section.
\begin{corollary}
\label{cor:hat_quantities}
 Suppose that $t\le \tau$ and $x\in(0,a)$. Then for large $A$ and $a$, we have $\Ehat^x[Z_t] \le CAe^{-(a-x)/2}$, $\Ehat^x[Z_t^2] \le C\ep e^Ae^{-(a-x)/2}$ and $\Ehat^x[Y_t] \le Ce^{-(a-x)/2}$.
\end{corollary}
\begin{proof}
  Immediate from Lemmas~\ref{lem:hat_quantities} and \ref{lem:Zt_variance} and \eqref{eq:c0_mu}.
\end{proof}

\subsection{The fugitive and its family}
\label{sec:fugitive}
We now describe the BBM conditioned to break out at a given time. Recall that $\mathscr U$ denotes the fugitive. For simplicity, we write $\tau_l$ and $\sigma_l$ for $\tau_l(\mathscr U)$ and $\sigma_l(\mathscr U)$, respectively, $l\in\N$. On the event $T=T^{(l)}$, we define  $\mathscr U_j$ to be the ancestor of $\mathscr U$ alive at the time $\sigma_j$, $j=0,\ldots,l$. By the strong branching property, we have the following decomposition:

\begin{lemma}
  \label{lem:decomposition_Tl}
Let $k\in\N$, $l\in \N\cup\{\infty\}$ with $l\ge k$ and $t\ge 0$. Conditioned on $\F_0$, $T^{(0;l)}=T^{(k)}=t$, $\tau_0$ and $\mathscr U_0$, the BBM admits the following recursive decomposition:
\begin{enumerate}[nolistsep]
\item The initial particles $u\ne \mathscr U_0$ spawn independent BBM conditioned on $T^{(0;l)}>t$,
\item independently, the particle $\mathscr U_0$ spawns BBM conditioned on a particle (call it $\mathscr U_0'$) hitting $a$ for the first time at the time $\tau_0$, all the children of which born before $\tau_0$ being conditioned on $T^{(0;l)}>t$.
\begin{enumerate}
\item If $k=0$, this particle $\mathscr U_0'$ spawns BBM conditioned on the event $B$ of a breakout.
\item If $k>0$, it spawns BBM starting at the space-time point $(a,\tau_0)$ conditioned on $T^{(0;l)}=T^{(k)}=t$. In particular, if we write $\mathscr S = \mathscr S^{(\mathscr U_0',\tau_0)}$, then conditioned on $\F_{\mathscr S}$, the particles in $\mathscr S$ spawn BBM starting from the collection of space-time points $\mathscr S$, conditioned on $T^{(0;l-1)}=T^{(k-1)} = t$.
\end{enumerate}
\end{enumerate}
\end{lemma}
Note that in the case 2b above, the subtree spawned by $(\mathscr U_0',\tau_0)$ follows the law $\Q$ \emph{conditioned} on $T^{(0;l)}=T^{(k)} = t$, hence the law of $\mathscr S^{(\mathscr U_0',\tau_0)}$ is not the same as under $\Q$. In particular,  $\E[Z^{(\mathscr U_0',\tau_0)}] \not\le CA$. Indeed, conditioning on one of its descendants breaking out at a later time corresponds to a kind of size-bias on the number of particles. However, it is still true that $Z^{(\mathscr U_0',\tau_0)} < \ep e^A$, by the definition of a breakout.

Lemma~\ref{lem:decomposition_Tl} gives a decomposition of the BBM conditioned on $T^{(k)}$ into $k+1$ pieces. In order to describe what happens in a single piece, define $\F_{\mathscr U}$ to be the $\sigma$-field generated by the family of events $(\{T^{(0;l)}=T^{(k)}\})_{k\in\N}$ and by $\mathscr U$, $(X_{\mathscr U}(s))_{\sigma_j\le s < \tau_j}$, $(d_u,k_u)_{u\preceq\mathscr U}$ and  $\mathscr S^{(\mathscr U,\tau_j)}$, for $j=0,\ldots,k$, on the event $\{T^{(0;l)}=T^{(k)}\}$. Note that in particular, on the event $\{T^{(0;l)}=T^{(k)}\}$, $\tau_j$, $(\mathscr U_j,\sigma_j)$ are $\F_{\mathscr U}$-measurable, $j=0,\ldots,k$. It is plain that conditioned on $\F_{\mathscr U}$, the subtrees spawned by the children of the ancestors of $\mathscr U$ during the intervals $[\sigma_j,\tau_j)$ are independent BBM conditioned not to break out before $t$.

It remains to describe the trajectory of the fugitive and its reproduction. By a decomposition at the first time of branching as in Section~\ref{sec:doob}, we could describe the law of BBM starting from a single particle at position $x$ at time $0$ and conditioned on a particle hitting $a$ for the first time at a time $s \le t$, all the children of which born before $s$ being conditioned on $T>t$. However, it is faster to use the Many-to-one lemma instead, which is the method of proof of the following lemma, which we state for general penalisations $p(s)$. This result is essentially \cite[Theorem~1]{Chauvin1991a}.
\begin{lemma}
 \label{lem:many_to_one_fugitive}
Let $t \ge 0$, $x\in(0,a)$ and $p:\R_+\to[0,1]$ be measurable with $p(s) = 0$ for $s$ large enough. Denote by $\Ptilde$ the law associated to $p(s)$ as in Section \ref{sec:weakly_conditioned}. Recall the definition of $e(x,t)$ from \eqref{eq:def_e}. Then, given a family of $\F_{\mathscr R_t^{(0)}}$-measurable non-negative random variables $(Y_u)_{u\in U}$, we have
\begin{equation*}
  \Etilde^x\Big[Y_u\,\Big|\,\exists u \in U : (u,t) \in \mathscr R_t^{(0)}\Big] = \frac{\Etilde^{*,x}\Big[Y_{\xi_t}e^{-\int_0^t e(\xi_s,s)\,\dd s} \,\Big|\,H_0(\xi) > H_a(\xi)=t\Big]}{\Etilde^{*,x}\Big[e^{-\int_0^t e(\xi_s,s)\,\dd s} \,\Big|\,H_0(\xi) > H_a(\xi)=t\Big]},
\end{equation*}
where under $\Ptilde^{*,x}$, the spine follows standard Brownian motion and spawns particles with rate $\widetilde{m_1}(x,t) \beta_0 Q(x,s)$ according to the reproduction law $((\widetilde{m_1}(x,t))^{-1}k\widetilde q(x,t,k))_{k\ge 0}$, which start independent $\Ptilde$-BBM. In particular, conditioned on $(X_u(s))_{0\le s \le t}$ and $(d_v,k_v)_{v\preceq u}$, the children of the ancestors of $u$ follow independent $\Ptilde$-BBM.
\end{lemma}
\begin{proof}
We have
\begin{equation}
 \label{eq:600}
 \Etilde^x\Big[Y_u\,\Big|\,\exists u\in U : (u,t) \in \mathscr R_t^{(0)}\Big]
  = \frac{\Etilde^x\Big[\sum_{(u,s)\in \mathscr R_t^{(0)}}\Ind_{(H_a(X_u) \in \dd t)} Y_u\Big]}{\Etilde^x\Big[\sum_{(u,s)\in \mathscr R_t^{(0)}}\Ind_{(H_a(X_u) \in \dd t)}\Big]},
\end{equation}
and by Lemma~\ref{lem:many_to_one},
\begin{equation}
  \label{eq:47}
 \Etilde^x\Big[\sum_{(u,s)\in \mathscr R_t^{(0)}}\Ind_{(H_a(X_u) \in \dd t)} Y_u\Big] = \Etilde^{*,x}_{-\mu}\Big[Y_{\xi_t}e^{\int_0^t \widetilde m(\xi_s,s)\beta_0 Q(\xi_s,s)\,\dd s} \Ind_{(H_0(\xi) > H_a(\xi)\in\dd t)}\Big].
\end{equation}
According to the description of the conditioned process in Section \ref{sec:weakly_conditioned} and the description of the spine in Section \ref{sec:spine}, the particles on the spine follow the Doob transform of Brownian motion with drift $-\mu$ by the harmonic function $h(x,t)$ and spawn particles with rate $\widetilde{m_1}(x,t) \beta_0 Q(x,s)$ according to the reproduction law $((\widetilde{m_1}(x,t))^{-1}k\widetilde q(x,t,k))_{k\ge 0}$. With Girsanov's transform, \eqref{eq:47} yields,
  \begin{equation}
    \label{eq:57}
  \Etilde^x\Big[\sum_{(u,s)\in \mathscr R_t^{(0)}}\Ind_{(H_a(X_u) \in \dd t)} Y_u\Big] = e^{-\mu(a-x)}\frac{h(a,t)}{h(x,0)}\Etilde^{*,x}\Big[Y_{\xi_t}e^{-\int_0^t e(\xi_s,s)\,\dd s} \Ind_{(H_0(\xi) > H_a(\xi)\in\dd t)}\Big].
  \end{equation}
Equations \eqref{eq:600} and \eqref{eq:57} yield the first statement. The second statement follows by taking appropriate test functionals $(Y_u)_{u\in U}$.
\end{proof}

\begin{corollary}
\label{cor:many_to_one_fugitive}
In addition to the assumptions in Lemma~\ref{lem:many_to_one_fugitive}, suppose $t\in [0,a^3]$ and $\|p\|_\infty  \le c_p$, where $c_p$ is some universal constant implicitly defined below. Then,
\begin{equation}
\label{eq:many_to_one_fugitive}
 \Etilde^x\Big[Y_u\,\Big|\,\exists u\in U : (u,t) \in \mathscr R_t^{(0)}\Big]
 \le C\Etilde^{*,x}\Big[Y_{\xi_t}\,\Big|\, H_0(\xi) > H_a(\xi)\in\dd t\Big].
\end{equation}
\end{corollary}
\begin{proof}
By Lemma~\ref{lem:k_integral}, the inequality $e^{-x} \ge 1-x$ and the hypothesis on $t$, we have
\begin{equation}
  \label{eq:610}
 \Etilde^{*,x}\Big[e^{-\int_0^t e(\xi_s,s)\,\dd s} \,\Big|\,H_0(\xi) > H_a(\xi)=t\Big] = W^{x,t,a}_{\mathrm{taboo}}[e^{-\int_0^t e(X_s,s)\,\dd s}] \ge (1-C\|p\|_\infty).
\end{equation}
The statement now follows from \eqref{eq:610} and Lemma~\ref{lem:many_to_one_fugitive}.
\end{proof}

\begin{figure}
\centering
\def\svgwidth{3cm}
\executeiffilenewer{sec6_4.svg}{sec6_4.pdf}%
{inkscape -z -D --file=sec6_4.svg %
--export-pdf=sec6_4.pdf --export-latex --export-area-drawing}%
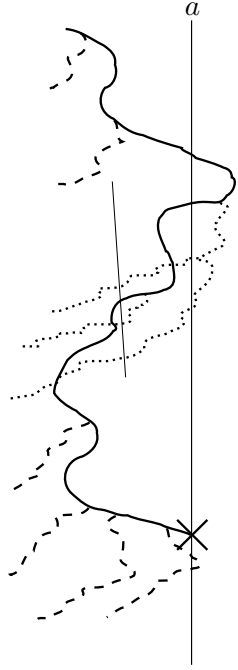%

\caption{A visualization of the bar- and check-particles: The path of the fugitive until the breakout (depicted by a cross) is drawn with a solid line, the bar-particles (spawned between the times $\sigma_l$ and $\tau_l$) with dashed lines and the check-particles (spawned between the times $\tau_l$ and $\sigma_{l+1}$) with dotted lines.}
\label{fig:sec6_4}
\end{figure}

We come back to the BBM conditioned to break out at a given time and set up the important definitions. Recall that $\mathscr U$ denotes the fugitive. We will denote by a \emph{bar} the quantities referring to the particles spawned between the times $\sigma_l(\mathscr U)$ and $\tau_l(\mathscr U)$ for some $l\ge 0$, and by a \emph{check} those referring to the particles spawned between the times $\tau_l(\mathscr U)$ and $\sigma_{l+1}(\mathscr U)$. See Figure~\ref{fig:sec6_4} for a visualization. Formally, we set
\[
 \mathscr{\widebar N}_t = \{u\in \mathscr N_t:(u,t)\wedge (\mathscr U,T) \in U\times \bigcup_{l\ge 0}[\sigma_l(\mathscr U),\tau_l(\mathscr U))\},
\]
\[
  \mathscr{\widebar R}_t = \{(u,s)\in \bigcup_{l\ge 0}\mathscr R^{(l)}_t:(u,s)\wedge (\mathscr U,T) \in U\times \bigcup_{l\ge 0}[\sigma_l(\mathscr U),\tau_l(\mathscr U)))\}
\]
and
\[
 \mathscr{\widecheck N}_t = \{u\in \mathscr N_t:(u,t)\wedge (\mathscr U,T) \in U\times \bigcup_{l\ge 1}[\tau_{l-1}(\mathscr U),\sigma_l(\mathscr U))\}.
\]
We then define
\[
 \widebar Z_t = \sum_{u\in \mathscr{\widebar N}_t}w_Z(X_u(t)),\quad \widebar Y_t = \sum_{u\in \mathscr{\widebar N}_t}w_Y(X_u(t)),\quad \widebar R_t = \#\mathscr{\widebar R}_t,
\]
and
\[
 \widecheck Z_t = \sum_{u\in \mathscr{\widebar N}_t}w_Z(X_u(t)),\quad \widecheck Y_t = \sum_{u\in \mathscr{\widebar N}_t}w_Y(X_u(t)).
\]
Note that on the event $T=T^{(0)}$, we have $\mathscr{\widecheck N}(T) = \emptyset$ by definition. 

By Corollary~\ref{cor:hat_quantities}, the following $\F_{\mathscr U}$-measurable functional will be of use in the study of the bar-quantities.
\begin{equation}
  \label{eq:E_U}
\mathscr E_{\mathscr U} = \sum_{i=0}^l \sum_{u\preceq \mathscr U} \Ind_{(\sigma_i\le d_u < \tau_i)}(k_u-1)e^{-(a-X_u(d_u))/2} ,\quad\ton \{T=T^{(l)}\}.
\end{equation}
where we recall that $d_u$ and $k_u$ are respectively the time of death and the number of children of the individual $u$.

\begin{lemma}
  \label{lem:E_U}
For large $A$ and $a$, we have for every $l\in\N$,
\[
\E[\mathscr E_{\mathscr U}\,|\,T=T^{(l)} \le \tconst{eq:tcrit}] \le C(l+1).
\]
\end{lemma}
\begin{proof}
By Lemma~\ref{lem:decomposition_Tl} and Corollary~\ref{cor:many_to_one_fugitive} (which can be applied because of \eqref{eq:pbar_estimate}), we have for every $l\in\N$,
\[
\E[\mathscr E_{\mathscr U}\,|\,T=T^{(l)} \le \tconst{eq:tcrit}]\le C  \sum_{i=0}^l \E\Big[W_i\Big[\int_0^{\tau_i-\sigma_i} e^{-(a-B_t)/2}\,\dd t\Big]\,\Big|\,T=T^{(l)}\le \tconst{eq:tcrit}\Big],
\]
where $W_i$ is the law of a Brownian bridge of length $\tau_i-\sigma_i$ from $X_{\mathscr U}(\sigma_i)$ to $X_{\mathscr U}(\tau_i)$. The statement now follows readily from Lemma~\ref{lem:k_integral}.
\end{proof}

We can now study the probability that the fugitive stems from a given tier.
\begin{lemma}
\label{lem:T1_T2}
Suppose that $C_1 e^A \le Z_0 \le C_2e^A$ and $Y_0 \le \eta Z_0$. Then for large $A$ and $a$, we have for $l\in\{1,2\}$,
 \[
  \P(T^{(l+)} < T^{(0;l-1)}) \le C(\ep A)^l,
 \]
\end{lemma}
\begin{proof}
By \eqref{eq:ep_lower}, \eqref{eq:pB} and \eqref{eq:QhatZ}, we have for large $A$ and $a$,
\begin{equation}
 \label{eq:777}
\P(T^{(0)} > \tconst{eq:tcrit}) \le \exp(-C A/\ep) \le \exp(-C/\sqrt\ep).
\end{equation}
Now, for the rest of the proof, let $t\le \tconst{eq:tcrit}$ and let $\nu = \sum_{i=1}^n \delta_{x_i}$. We have by the decomposition \eqref{eq:decomposition} of the process conditioned on $T^{(0;l-1)}=t$,
\begin{multline}
  \label{eq:1}
 \P^\nu(T^{(l+)} > t\,|\,T^{(0;l-1)} = t) = \sum_{i=1}^n p_i \Phat^{\nu-\delta_{x_i}}(T^{(l+)} > t\,|\,T^{(0;l-1)} > t)\\
\times\P^{x_i}(T^{(l+)} > t\,|\,T^{(0;l-1)} = t)
\end{multline}
Define $c_t = \pi^2A( t /{a^3} + \Const{eq:Zl_exp_upperbound_2} \eta)$, which is less than $1/2$ for large $A$ and $a$ by the hypothesis on $t$. By Lemma~\ref{lem:Zl_exp_upperbound} and the hypothesis, we have for every $j\le l$,
\begin{equation}
  \label{eq:11}
\Ehat_{l-j-1}[Z_\emptyset^{(l-j)}] \le CZ_0 c_t^{l-j},
\end{equation}
and moreover for every $(x,s)\in (0,a)\times[0,t]$,
 \begin{equation}
   \label{eq:12}
\Ehat_{l-j-1}^{(x,s)}[Z_\emptyset^{(l-j)}] \le CA e^{-(a-x)/2}c_t^{l-j-1},
 \end{equation}
For every $k\le l-1$, \eqref{eq:11} applied to the particles in $\widecheck{\mathscr S}^{(j)}_t$, $j=1,\ldots,k$ yields 
\begin{equation}
  \label{eq:13}
  \E[\widecheck Z^{(l)}_\emptyset\,|\,T^{(k)} = T^{(0;l-1)} = t] \le C\ep e^A \sum_{j=1}^lc_t^{l-j} \le C\ep e^A c_t^{l-k}.
\end{equation}
Moreover, by \eqref{eq:12}, Corollary~\ref{cor:hat_quantities} and Lemma~\ref{lem:E_U}, we have
\begin{equation}
  \label{eq:14}
  \E[\widebar Z^{(l)}_\emptyset\,|\,T^{(k)} = T^{(0;l-1)} = t] \le CA c_t^{l - k}.
\end{equation}
In total, we get by \eqref{eq:11}, \eqref{eq:12} and \eqref{eq:13}, the hypothesis on $Z_0$ and \eqref{eq:ep_lower},
\begin{equation}
  \label{eq:15}
 \E[\widehat Z^{(l)}_\emptyset+\widecheck Z^{(l)}_\emptyset+\widebar Z^{(l)}_\emptyset\,|\,T^{(k)} = T^{(0;l-1)} = t] \le C e^A\Big(c_t^l + \ep c_t^{l - k}\Big).
\end{equation}
Jensen's inequality and Proposition~\ref{prop:T} together with \eqref{eq:15} and the inequality $1-e^{-x} \le x$ give
\begin{equation}
  \label{eq:115}
 \P(T^{(l+)} < t\,|\,T^{(k)} = T^{(0;l-1)} = t) \le C c_t(\ep^{-1}c_t^{l} + c_t^{l-k}).
\end{equation}
Summing \eqref{eq:115} over $k$, the law of total expectation gives
\begin{equation}
  \label{eq:114}
 \P(T^{(l+)} < t\,|\, T^{(0;l-1)} = t) \le C c_t(\ep^{-1}c_t^{l} + c_t).
\end{equation}
The lemma now follows by integrating \eqref{eq:114} over $t$ from $0$ to $\tconst{eq:tcrit}$ and using Lemma~\ref{lem:moments_T} and \eqref{eq:777}.
\end{proof}

\begin{remark}
  One may wonder whether one can simply calculate $\P(T^{(l)}\in\dd t\,|\,T^{(0;l-1)}> t)$ for every $l\ge 1$ and $t\le\tconst{eq:tcrit}$, using only the tools from Section~\ref{sec:before_breakout_particles}. This would require fine estimates on the density of the point process formed by the particles from tier $l-1$ hitting $a$ just before $t$. These estimates can be most easily obtained if one stops descendants of the particles hitting $a$ at a (large) fixed time $\zeta$ instead of the line from Lemma~\ref{lem:critical_line}, with which the results in this paper would hold as well. However, in order not to lose generality, we stick to Lemma~\ref{lem:T1_T2}, which is enough for our purposes.
\end{remark}

\section{The B-BBM}
\label{sec:BBBM}

\def\numcst{\epsilon}

We will now define properly the BBM with the moving barrier, also called the B-BBM (the ``B'' stands for ``barrier''), which will be used in the subsequent sections to approximate the $N$-BBM. We will still use
all the definitions from Section \ref{sec:before_breakout_definitions}, with
one notational change: Recall that by \eqref{eq:decomposition}, we can
decompose the process into two parts; the first part consisting of the
particles spawned by the ancestor of the fugitive and the second part
consisting of the remaining particles. As in Section \ref{sec:fugitive}, the
quantities which refer to the particles of the first part will be denoted by a
bar (e.g.\ $\widebar Z$) or check (e.g.\ $\widecheck Z$). The quantities of the
second part will be denoted with a hat in this section (e.g.\ $\widehat Z$),
in reference to the law $\Phat$ from Section
\ref{sec:before_breakout_particles}. 

\subsection{Definition of the model}
\label{sec:BBBM_definition}
Suppose that we are given a family $(f_x)_{x\ge 0}$ of non-decreasing functions $f_x \in \mathscr C^2(\R,\R_+)$, such that for each $x\ge 0$, $f_x(t) = 0$ for $t\le 0$, $f_x(+\infty) = x$ and for each $\delta > 0$ small enough there exist $x_0 = x_0(\delta)$, $t_0 = t_0(\delta)$, such that
\begin{enumerate}[nolistsep]
 \item $x_0(\delta) \to \infty$ as $\delta \to 0$,
 \item $\|f_x\| \le \delta^{-1}$ for all $x\in [0,x_0]$,
 \item $f_x(t) \ge (1-\delta)x$ for all $t\ge t_0$ and
 \item $t_0 \le \delta^{-1}$.
\end{enumerate}
where $\|f\|$ is defined in \eqref{eq:def_f_norm}. A family of functions with the above properties exists, for example the following, which we will choose in Sections~\ref{sec:Bflat} and \ref{sec:Bsharp}.
 \begin{equation}
   \label{eq:wall_fn}
 f_x(t) = \log\left(1+(e^{ x}-1)\thbar(t)\right).
 \end{equation}

\begin{figure}[h]
 \centering
\def\svgwidth{12cm}
\executeiffilenewer{sec7.svg}{sec7.pdf}%
{inkscape -z -D --file=sec7.svg %
--export-pdf=sec7.pdf --export-latex --export-area-drawing}%
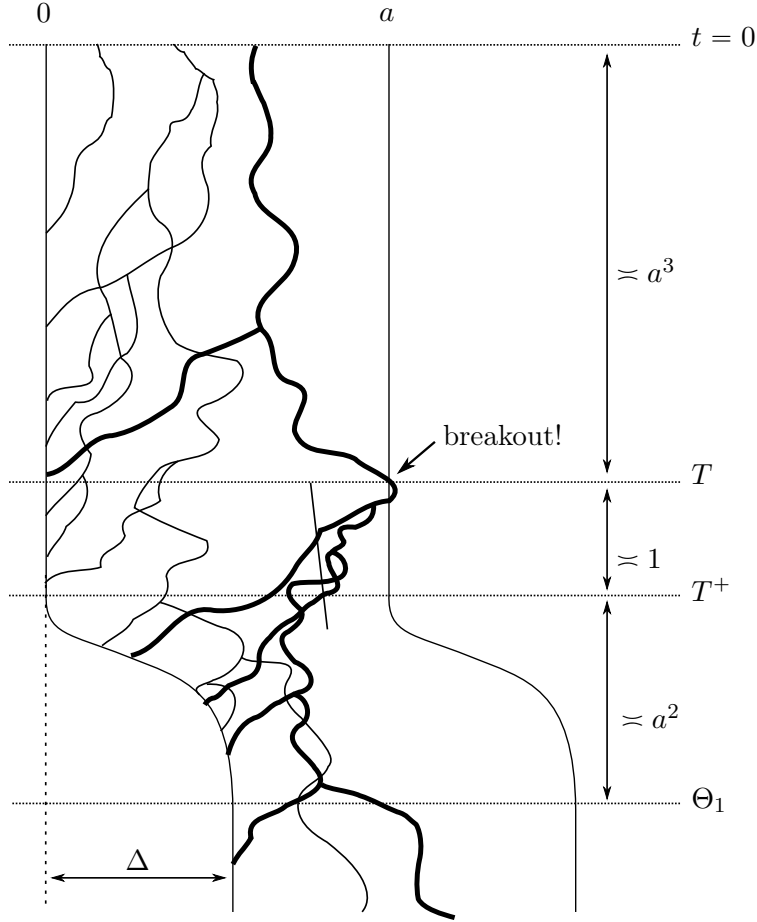%

\caption{A caricatural graphical description of the B-BBM until the time $\Theta_1$. The fugitive and its descendants are drawn with thick lines, the other particles with thin lines. A breakout happens at time $T$ and the barrier is moved from the time $T^+$ on. Note that technically we increase the drift to the left instead of moving the barrier. The three important timescales ($1$, $a^2$ and $a^3$) are shown as well.}
\label{fig:sec7}
\end{figure}

For every $A$ and $a$, let $\nu_0 = \nu_0^{A,a}$ be a (possibly random) finite counting measure. We now define the B-BBM with  initial configuration of particles $\nu_0$. Starting from BBM with constant drift $-\mu$ with initial configuration $\nu_0$, we define for each $n\in\N$ a stopping time $\Theta_n$ and for each $n\in\N^*$ the barrier process $(X^{[n]}_t)_{t\in [\Theta_{n-1},\Theta_n]}$ as follows:
\begin{enumerate}
 \item We set $\Theta_0 = 0$.
 \item Denote by $T=T_1$ the time of the first breakout of the BBM absorbed at 0 and by $\mathscr U$ the fugitive, as in Section \ref{sec:before_breakout_definitions}. We set $X^{[1]}_t = 0$ for $t\in [0,T]$.
\item
Define $T^- = (T-e^Aa^2)\wedge 0$ and 
\begin{equation}
\label{eq:Delta_def}
\Delta =  \log\Big(e^{-A}(\widehat Z_{T^-} + \widecheck Z_{\emptyset\vee T^-} +Z^{(\mathscr U,T)})\Big)\vee 0.
\end{equation}
where $\widecheck Z_{\emptyset\vee T^-} = \widecheck Z^{(1)}_{T^-}$ if $T^- \ge \tau_0(\mathscr U)$ and $=\widecheck Z^{(1)}_\emptyset$ otherwise.  The above quantities are defined in Section \ref{sec:fugitive}, after Corollary~\ref{cor:many_to_one_fugitive}.
 \item Define $T^+ = T^+_1 = T+\tau^{(\mathscr U,T)}_{\mathrm{max}}$ and $\Theta_1 = (T +
   e^Aa^2)\vee T^+$. Note that $T^+$ and therefore also $\Theta_1$ is a stopping time for the BBM. Now define
\[
X^{[1]}_t = f_{\Delta}\Big(\frac{t-T^+}{a^2}\Big), t\in[T^+,\Theta_1].
\]
We then give to the particles an additional drift $-(\dd/\dd t)X^{[1]}_t$ for $t\in [T^+,\Theta_1]$, in the meaning of Section \ref{sec:interval_notation}.
 \item We have now defined $T_1$, $T^+_1$, $\Theta_1$ and $X^{[1]}$. We further define $\nu_1$ to be the measure formed by the particles at time $\Theta_1$, which have never hit $0$. To define $T_2$, $T^+_2$, $\Theta_2$ and $X^{[2]}$, we repeat the above steps with the process formed by the BBM started from those particles, with the definitions changed such that the barrier process starts at $X^{[2]}_{\Theta_1} = X^{[1]}_{\Theta_1}$, time starts at $\Theta_1$ etc.
\item We now construct the barrier process $X^{[\infty]}_t$ from the pieces by $X^{[\infty]}_t = X^{[n]}_t$, if $t\in [\Theta_{n-1},\Theta_n]$.
\end{enumerate}

Define the event $G_0 = \{\supp \nu_0 \subset (0,a),\ |Z_0 - e^A| \le \ep^{3/2}e^A,\ Y_0 \le \eta\}$.
Recall the definition of the phrase ``As $A$ and $a$ go to infinity'' from Section \ref{sec:before_breakout_definitions}. Define the predicate
\begin{itemize}
 \item[(HB)] (The law of) $\nu_0$ is such that $\P(G_0)\to 1$ as $A$ and $a$ go to infinity.
 \item[(HB$_0$)] $\nu_0$ is deterministic and such that $G_0$ is verified.
\end{itemize}
Furthermore, if for some $n>0$, $(t^{A,a}_j)_{j=1}^n\in\R_+^n$ for all $A$ and $a$, then define the predicate
\begin{itemize}
\item[(Ht)] There exists $0\le t_1 < \dots < t_n$, such that for all $j\in\{1,\ldots,n\}$, $t^{A,a}_j\to t_j$ as $A$ and $a$ go to infinity, with $t_1^{A,a}\equiv 0$ if $t_1=0$.
\end{itemize}

\begin{theorem}
 \label{th:barrier}
Suppose (HB). Let $n\ge 1$ and $(t^{A,a}_j)_{j=1}^n\in\R_+^n$ such that (Ht) is verified. Define the process
\[
(X_t)_{t\ge 0} = (X^{[\infty]}_{a^3t} -  \pi^2 A t)_{t\ge 0}.
\]
Then, as $A$ and $a$ go to infinity, the law of the vector $(X_{t^{A,a}_j})_{j=1}^n$ converges to the law of $(L_{t_j})_{j=1}^n$, where $(L_t)_{t\ge0}$ is the L\'evy process from Theorem~\ref{th:1}.
\end{theorem}

A stronger convergence than convergence in the sense of finite-dimensional distributions is convergence in law with respect to Skorokhod's ($J_1$-)topology (see \cite[Chapter~3]{Ethier1986}). Obviously, the convergence in Theorem~\ref{th:barrier} does not hold in this stronger sense, because the barrier is continuous but the L\'evy process is not and the set of continuous functions is closed in Skorokhod's $J_1$-topology\footnote{One could prove convergence in the $M_1$-topology (cf.\ \cite[Section~3.3]{Whitt2002}), but the use of this less standard topology would lead us away from our real goal: the proof of Theorem~\ref{th:1}.}. However, if we create artificial jumps, we can rectify this:

\begin{theorem}
 \label{th:barrier2}
Suppose (HB). Define $J_t = X^{[\infty]}_{\Theta_n}$, if $t\in [\Theta_n,\Theta_{n+1})$, for $n\in\N$. Then as $A$ and $a$ go to infinity, the process $(X'_t)_{t\ge 0} = (J_{a^3t}-\pi^2 At)_{t\ge 0}$ converges in law with respect to Skorokhod's topology to the L\'evy process defined in the statement of Theorem~\ref{th:barrier}.
\end{theorem}

For each $n\ge 1$, we define the event $G_n$ to be the intersection of $G_{n-1}$ with the following events (see Section~\ref{sec:before_breakout_definitions} for the definition of $\mathscr N_t$):
\begin{itemize}
 \item $\supp\nu_n\subset (0,a)$,
 \item $\mathscr N_{\Theta_n} \subset U\times \{\Theta_n\}$ and $\Theta_n > T_n^+$ (for $n>0$),
 \item $|e^{-A}Z_{\Theta_n} -1| \le \ep^{3/2}$ and $Y_{\Theta_n} \le \eta$.
\end{itemize}
Note that $G_n\in\F_{\Theta_n}$ for every $n\ge 0$.

The core of the proof of Theorems \ref{th:barrier} and \ref{th:barrier2} will be the following proposition:
\begin{proposition}
\label{prop:piece}
Fix $\lambda \in \R$ and define $\gamma_0 = (\pi p_Be^A)^{-1}$.
There exists a numerical constant $\numcst > 0$, such that for large $A$ and $a$, we have $\P(G_n^c\,|\,\F_0)\Ind_{G_0} \le n\ep^{1+\numcst}$ for every $n\ge 0$ and
\begin{equation}
\label{eq:piece_fourier}
 \log \E\Big[e^{i\lambda X^{[\infty]}_{\Theta_n}}\,\Big|\,\F_0\Big]\Ind_{G_0} = \Big(n\gamma_0(\kappa(\lambda) + i\lambda \pi^2A + o(1) + O(\ep^\numcst))\Big)\Ind_{G_0}.
\end{equation} 
where $\kappa(\lambda)$ is the cumulant from \eqref{eq:laplace_levy} and where $o(1)$ and $O(\ep^\numcst)$ may depend on $\lambda$.
\end{proposition}

\subsection{Proof of Proposition \ref{prop:piece}}
\label{sec:piece_proof}
First note that conditioning on $\F_0$, we can and will assume without loss of generality that the hypothesis (HB$_0$) holds and will only state it again in the lemmas. In the same vein, \emph{we will always assume that $A$ and $a$ are large and will omit stating this fact, except again in the statements of the lemmas}.

Furthermore, it is only necessary to treat the case $n=1$, since the rest can be obtained by induction: Suppose that we have shown the result for $n=1$. Since the process starts afresh at the stopping time $\Theta_n$, we then have $\P(G_n\backslash G_{n+1}\,|\,\F_{\Theta_n}) \le  \ep^{1+\numcst}\Ind_{G_n} $ for every $n\ge0$, which implies the statement on the probability of $G_n$. Similarly, if we set $\kappa'(\lambda) = \gamma_0\left(\kappa(\lambda) + i\lambda \pi^2 A + o(1)\right)$, where $o(1)$ is the term from \eqref{eq:piece_fourier}, then
\begin{align*}
  \E\Big[e^{i\lambda (X^{[\infty]}_{\Theta_{n+1}}-X^{[\infty]}_{\Theta_{n}})}\Ind_{G_{n+1}}\,|\,\F_{\Theta_n}\Big] &= \Big(e^{\kappa'(\lambda) + O(\ep^\numcst)} + O\left(\P(G_n\backslash G_{n+1}\,|\,\F_{\Theta_n})\right)\Big)\Ind_{G_n}\\
&= \left(1 + O(\ep^{1+\numcst})\right)e^{\kappa'(\lambda) + O(\ep^\numcst)}\Ind_{G_n},
\end{align*}
because $\gamma_0 = \ep(1+o(1))$ by \eqref{eq:pB}. It follows that
\begin{align*}
   \E\Big[e^{i\lambda X^{[\infty]}_{\Theta_{n+1}}}\Ind_{G_{n+1}}\Big] &= e^{\kappa'(\lambda) + O(\ep^\numcst)}\E\Big[e^{i\lambda X^{[\infty]}_{\Theta_n}} \left(1 + O(\ep^{1+\numcst})\right)\Ind_{G_n}\Big] \\
&= e^{\kappa'(\lambda) + O(\ep^\numcst)}\E\Big[e^{i\lambda X^{[\infty]}_{\Theta_n}} \Ind_{G_n}\Big]
\end{align*}
By induction, and the fact that $\P(G_{n+1}^c)\le n\ep^{1+\numcst} \le e^{O(\ep^{1+\numcst})}$, we then obtain \eqref{eq:piece_fourier}.

Throughout the proof, we will work on several ``good sets'', which we all define here for easy reference, since they will be reused in the following sections. For this reason, some of them (for example $\widehat G$) are actually more restrictive than what would be necessary for the proof of Proposition~\ref{prop:piece}.

\begin{itemize}
\item $G_{\mathscr U} =\{e^Aa^2\le T\le \sqrt \ep a^3,\,T=T^{(0;1)},\,\mathscr E_{\mathscr U} \le e^{A/3}\}$.
\item $G_{\mathrm{fug}} = \{Z^{(\mathscr U,T)} \le e^A/\ep\} \cap G^{(\mathscr U,T)}$ (the event $G^{(\mathscr U,T)}$ was defined in \eqref{eq:116}).
\item $\widehat G = \{|\widehat Z_{T^-} - e^A| \le \ep^{1/4} e^A,\ \widehat Y_{T^-} \le a^{-1}e^{3A/2}\}$.
\item $\widecheck G = \{\widecheck Z_{\emptyset\vee T^-} \le \ep^{1/4}e^A,\ \widecheck Y_{\emptyset\vee T^-} \le e^{-A/2}\}$ (we define $\widecheck Y_{\emptyset\vee T^-}$ analogously to $\widecheck Z_{\emptyset\vee T^-}$).
\item $G_\Delta = G_{\mathscr U}\cap G_{\mathrm{fug}}\cap \widehat G\cap \widecheck G$.
\item $G_{\mathrm{nbab}}$ (``nbab'' stands for ``no breakout after breakout'') is the event that no bar-particle breaks out between $T$ and $\Theta_1+e^Aa^2$ and that no hat-, fug- or check-particle (fug-particles = descendants of $\mathscr S^{(\mathscr U,T)}$) hits $a$ between $T^-$ and $\Theta_1+e^Aa^2$.
\end{itemize}

Recall the definition of  $\F_{\mathscr U}$ from Section~\ref{sec:fugitive} and set $\F_\Delta := \F_{\mathscr U} \vee \widehat \F_{T^-}\vee \widecheck \F_{T^-}$, such that $G_\Delta\in\F_\Delta$ and the random variable $\Delta$ is measurable with respect to $\F_\Delta$. Define further $Z_\Delta = \widehat Z_{T^-}+Z^{(\mathscr U,T)}+\widecheck Z_{\emptyset\vee T^-}$ and $Y_\Delta$ analogously. Finally, define $\P_{\mathrm{nbab}} = \P(\cdot\,|\,G_{\mathrm{nbab}})$.
 
\paragraph{The probability of $G_\Delta$.}

\begin{lemma} 
\label{lem:GDelta}
Suppose (HB$_0$). For large $A$ and $a$, we have $\P(G_\Delta) \ge 1 - C\ep^{5/4}$.
\end{lemma}
\begin{proof}
By Proposition~\ref{prop:T} and Lemma~\ref{lem:T1_T2}, we have
\begin{equation}
  \label{eq:21}
  \P(e^Aa^2 \le T \le \sqrt{\ep}a^3,\,T=T^{(0;1)}) \ge 1 - CA^2\ep^2.
\end{equation}
Since $\sqrt\ep \le \tconst{eq:tcrit}$ for large $A$ by \eqref{eq:ep_upper}, Equation \eqref{eq:21} and Markov's inequality applied to Lemma~\ref{lem:E_U} yield
\begin{equation}
  \label{eq:32}
  \P(G_{\mathscr U}) \ge 1 - CA^2\ep^2.
\end{equation}
Furthermore, by \eqref{eq:Z_W}, we have,
\begin{equation}
 \label{eq:Gfug}
\P(G_{\mathrm{fug}}^c) \le p_B^{-1}\left(\P^a\left(\pi W > e^A/\ep-2\eta \right) + \eta\right) \le C\ep^2,
\end{equation}
by \eqref{eq:W_tail}, \eqref{eq:ep_lower}, \eqref{eq:eta} and \eqref{eq:pB}.

In order to estimate the probability of $\widehat G$, we will calculate first- and second moments of the quantities in the definition of $\widehat G$. These estimates will be needed again later on, when we calculate the Fourier transform of the variable $\Delta$.

By Lemma~\ref{lem:hat_quantities} and the hypotheses, together with \eqref{eq:pB}, \eqref{eq:ep_lower} and \eqref{eq:eta},
\begin{equation}
  \label{eq:710}
 \E[e^{-A}\widehat Z_{T^-}\,|\,\F_{\mathscr U}]\Ind_{G_{\mathscr U}} = 1+\pi \Q^a[Z] \tfrac T {a^3} + O\left(\ep^{3/2} + \left(A\tfrac T {a^3}\right)^2\right),
\end{equation}
which together with Lemma~\ref{lem:Zt_variance} gives
\begin{equation}
\label{eq:712}
\begin{split}
\E[(e^{-A}\widehat Z_{T^-}-1)^2\,|\,\F_{\mathscr U}]\Ind_{G_{\mathscr U}} &\le \Var(e^{-A}\widehat Z_{T^-}\,|\,\F_{\mathscr U})\Ind_{G_{\mathscr U}}+ C\left((A\tfrac T {a^3})^2+ \ep^3\right) \\
&\le C\left((A\tfrac T {a^3})^2+\ep \tfrac T {a^3}+\ep^3\right),
\end{split}
\end{equation}
by the hypotheses on the initial configuration. Lemma~\ref{lem:moments_T} and \eqref{eq:710}, together with  \eqref{eq:ep_upper}, \eqref{eq:eta_ep} and \eqref{eq:QaZ}, now give
\begin{equation}
 \label{eq:714}
 \E[e^{-A}\widehat Z_{T^-}\,|\,G_{\mathscr U}] = 1+\pi^2\gamma_0 (A+\log \ep+\const{eq:QaZ}+o(1)).
\end{equation}
Similarly, \eqref{eq:710} and \eqref{eq:712} and Lemma~\ref{lem:moments_T} give
\begin{equation}
 \label{eq:716}
 \E[(e^{-A}\widehat Z_{T^-}-1)^2\,|\,G_{\mathscr U}] = O(A^2\ep^2).
\end{equation}
Furthermore, by Lemma~\ref{lem:hat_quantities} and Markov's inequality, together with \eqref{eq:eta} and the hypotheses, we have
\begin{equation}
  \label{eq:117}
  \P[\widehat Y_{T^-}>e^{3A/2}/a\,|\,G_{\mathscr U}] \le C e^{-A/2}.
\end{equation}
By the same arguments, and since $\widehat Z^{(1)}_\emptyset \le \ep e^A\Ind_{(T=T^{(1)})}$ by the definition of a breakout, we have
\begin{align}
\label{eq:121}
 &\E[e^{-A}\widecheck Z_{\emptyset\vee T^-}\,|\,\F_{\mathscr U}]\Ind_{G_{\mathscr U}} \le C \ep\Ind_{(T=T^{(1)})}\\
  \label{eq:119}
  &\E[(e^{-A}\widecheck Z_{\emptyset\vee T^-})^2\,|\,\F_{\mathscr U}]\Ind_{G_{\mathscr U}} \le C \ep^2\Ind_{(T=T^{(1)})}\\
\label{eq:122}
&\E[\widecheck Y_{\emptyset\vee T^-}\,|\,\F_{\mathscr U}]\Ind_{G_{\mathscr U}} \le e^{A}\eta \le e^{-A},
\end{align}
be \eqref{eq:eta}. Equations \eqref{eq:117} and \eqref{eq:122} and Chebychev's inequality applied to \eqref{eq:716} and \eqref{eq:119} together with \eqref{eq:ep_upper} now give
\begin{equation}
  \label{eq:118}
  \P(\widehat G\cup\widecheck G\,|\,G_{\mathscr U}) \le 1-\ep^{5/4}.
\end{equation}
The lemma now follows from \eqref{eq:32}, \eqref{eq:Gfug} and \eqref{eq:118}.
\end{proof}

\paragraph{The probability of $G_{\mathrm{nbab}}$.}

\begin{lemma}
  \label{lem:Gnbab}
We have $\P((G_{\mathrm{nbab}})^c\,|\,\F_\Delta)\Ind_{G_\Delta} \le C\ep^2$ for large $A$ and $a$.
\end{lemma}
\begin{proof}
Define $R_1 = \widehat R_{[T^-,\Theta_1+e^Aa^2]} + R^{\mathrm{fug}}_{[T,\Theta_1+e^Aa^2]} +\widecheck R_{[T^-,\Theta_1+e^Aa^2]}$. By Markov's inequality, Lemmas~\ref{lem:Rt},  and \ref{lem:R_f} and Corollary~\ref{cor:hat_many_to_one} we have
\begin{equation*}
  \P(R_1 > 0\,|\,\F_\Delta)\Ind_{G_\Delta}\le C(Z_\Delta e^A/a + Y_\Delta) \le C\ep^2,
\end{equation*}
by the definition of $G_\Delta$, \eqref{eq:ep_lower} and \eqref{eq:eta}. Now, by  Corollary~\ref{cor:hat_quantities},
\begin{align*}
 \E[\widebar Z_T\,|\,\F_\Delta]\Ind_{G_\Delta} \le CA\mathscr E_{\mathscr U} \le CAe^{A/3}\quad\tand\quad\E[\widebar Y_T\,|\,\F_\Delta]\Ind_{G_\Delta} \le C \mathscr E_{\mathscr U} \le Ce^{A/3},
\end{align*}
where we used the fact that $\widecheck Z^{(1)}_\emptyset \le \ep e^A$ by the definition of a breakout. By Proposition~\ref{prop:T} and the inequality $1-e^{-x} \le x$, the probability conditioned on $\F_\Delta$ and on $G_\Delta$ that a check- or bar-particle breaks out between $T$ and $\Theta_1$ is now bounded by $p_B e^{A/3} \le C\ep^2$ by \eqref{eq:pB} and \eqref{eq:ep_lower}. This yields the lemma.
\end{proof}

\paragraph{The particles at the time $\Theta_1$. The probability of $G_1$.}

Recall that from the time $T^+ = T+\tau^{(\mathscr U,T)}_{\mathrm{max}}$ on, we move the barrier according to the function $f_\Delta$, which is equivalent to having the variable drift $-\mu_t = -\mu-f_\Delta(t/a^2)/a^2$. On $G_{\mathrm{fug}}$, we have $T^+ \le T+\zeta$ and $\Theta_1 = T+e^Aa^2$. Furthermore, $0 < \Delta \le C\log(1/\ep)$ on $G_\Delta$ and by the hypotheses on the functions $(f_x)$ and \eqref{eq:ep_lower} it follows,
\begin{equation}
 \label{eq:Delta_f}
\ton G_\Delta: \|f_\Delta\| \le \sqrt
a\quad\tand\quad \Delta - X^{[1]}_{\Theta_1}= O(e^{-A/2}),
\end{equation}
where $\|\cdot\|$ is defined in \eqref{eq:def_f_norm}.

\begin{lemma}
\label{lem:G1}
Suppose (HB$_0$). Then $\P(G_1) \ge 1-C\ep^{5/4}$ for large $A$ and $a$.
\end{lemma}

\begin{proof}
By Proposition~\ref{prop:quantities}, Lemmas~\ref{lem:mu_t_density} and \ref{lem:Gnbab}, Corollaries~\ref{cor:hat_many_to_one} and \ref{cor:hat_quantities} and \eqref{eq:Delta_f},
 \begin{align*}
   &\E_{\mathrm{nbab}}[e^{-A}Z_{\Theta_1}\,|\,\F_\Delta]\Ind_{G_\Delta} = e^{-\Delta - A}(Z_\Delta  + A\mathscr E_{\mathscr U})(1+O(e^{-A/2})) = 1+O(e^{-A/2})\\
&\Var_{\mathrm{nbab}}(e^{-A} Z_{\Theta_1}\,|\,\F_\Delta)\Ind_{G_\Delta} \le C \ep e^{-\Delta - A} \left(Z_\Delta e^A/a + Y_\Delta + \mathscr E_{\mathscr U}\right) \le C \ep e^{-2A/3},
 \end{align*}
such that by Chebychev's inequality and \eqref{eq:ep_lower},
\begin{equation}
  \label{eq:123}
  \P_{\mathrm{nbab}}(|e^{-A}Z_{\Theta_1} - 1| \ge \ep^{3/2}) \le C\ep^2.
\end{equation}
Furthermore, by Lemma~\ref{lem:hat_quantities}, Corollary~\ref{cor:hat_quantities} and Markov's inequality, we have
\begin{equation}
  \label{eq:124}
  \P_{\mathrm{nbab}}(\supp \nu_{\Theta_1} \subset (0,a),\ Y_{\Theta_1} \le \eta\,|\,\F_\Delta)\Ind_{G_\Delta} \le e^A/(a\eta) + C\eta e^A \le C\ep^2,
\end{equation}
by \eqref{eq:eta} and \eqref{eq:ep_lower}. The lemma now follows from \eqref{eq:123} and \eqref{eq:124} together with Lemmas~\ref{lem:GDelta} and \ref{lem:Gnbab}.
\end{proof}

\begin{remark}
  \label{rem:piece}
Lemma~\ref{lem:G1} obviously still holds if one replaces $G_1$ by $G_1\cap\{|e^{-A}Z_{\Theta_1}-1| < \tfrac 1 2 \ep^{3/2}\}$. This will be needed in Sections~\ref{sec:Bflat} and \ref{sec:Bsharp}.
\end{remark}

\paragraph{The Fourier transform of the barrier process.}
Define $\Delta_{\mathrm{drift}} = e^{-A} (\widehat Z_{T^-} + \widecheck Z_{\emptyset\vee T^-}) - 1$ and $\Delta_{\mathrm{jump}} = e^{-A} Z^{(\mathscr U,T)}$, which are independent random variables. Recall that $\Delta > 0$ on $G_\Delta$. By  \eqref{eq:121} and Lemmas~\ref{lem:T1_T2} and \ref{lem:GDelta}, we have
\begin{equation*}
  \E[e^{-A}\widecheck Z_{\emptyset\vee T^-}\,|\,G_{\mathscr U}]\le C\ep \P(T=T^{(1)}\,|\,G_{\mathscr U}) \le CA\ep^2,
\end{equation*}
such that with \eqref{eq:714}, we get,
\begin{equation}
 \label{eq:Delta_drift}
\E\left[\Delta_{\mathrm{drift}}\,|\,G_{\mathscr U}\right] = \pi^2\gamma_0(A+\log\ep+c+o(1)).
\end{equation}
Furthermore, \eqref{eq:716} and \eqref{eq:119} and the inequality $(x+y)^2 \le 2(x^2+y^2)$ yield
\begin{equation}
 \label{eq:Delta_drift_2}
\E\left[(\Delta_{\mathrm{drift}})^2\,|\,G_{\mathscr U}\right] = O(\ep^2A^2).
\end{equation}
Equation \eqref{eq:Delta_drift_2} and \eqref{eq:ep_upper} imply
\begin{equation}
 \label{eq:Delta_drift_bound}
\P(|\Delta_{\mathrm{drift}}| \ge \ep^{1/3},\ G_{\mathscr U}) = O(\ep^{4/3}A^2) = O(\ep^{7/6}).
\end{equation}

Fix $\lambda\in\R$. Since $\log(1+x+y) = \log(1+x)+\log((1+y)/(1+x))$, we have by \eqref{eq:Delta_f},  \eqref{eq:Delta_drift_bound} and Lemma~\ref{lem:GDelta},
\begin{equation}
\label{eq:782}
\begin{split}
 \E[e^{i\lambda X^{[\infty]}_{\Theta_1}}] &= \E[e^{i\lambda \Delta}\Ind_{(|\Delta_{\mathrm{drift}}| < \ep^{1/3})}\Ind_{G_\Delta}] + O(\ep^{9/8})\\
&= \E[e^{i\lambda \log(1+\Delta_{\mathrm{drift}})}\Ind_{(|\Delta_{\mathrm{drift}}| < \ep^{1/3})}\Ind_{G_\Delta}e^{i\lambda\log(1+\frac{\Delta_{\mathrm{jump}}}{1+\Delta_{\mathrm{drift}}})}] + O(\ep^{9/8}),
\end{split}
\end{equation}
for any $\lambda \in \R$. We will first study the term concerning $\Delta_{\mathrm{jump}}$. Write $Z = Z^{(\mathscr U,T)}$ and let $\rho$ be a real-valued constant with $|\rho| < \ep^{1/3}$. Then,
\begin{equation}
 \label{eq:785}
\begin{split}
\E[e^{i\lambda \log(1+\frac{\Delta_{\mathrm{jump}}}{1+ \rho})}] & = \E[e^{i\lambda \log(1+\frac{e^{-A}Z}{1+ \rho})}\,|\,Z > \ep e^A] + O(\P^a((G^{(\emptyset,0)})^c))\\
& = \int_\ep^\infty g(x)\,\P(e^{-A}Z \in \dd x\,|\,Z>\ep e^A) + O(\ep^2),
\end{split}
\end{equation}
where
\[
 g(x) = \exp\Big(i\lambda \log\Big(1+ \frac{x}{1+ \rho}\Big)\Big).
\]
By definition of $p_B$ and $\eta$, we have
\[
\begin{split}
 \int_\ep^1 x \P(e^{-A}Z \in \dd x\,|\,Z>\ep e^A) &= (p_B+O(\eta))^{-1}\E[e^{-A}Z\Ind_{(\ep e^A<Z\le  e^A)}]\\
 &= \pi^2\gamma_0(-\log\ep + o(1)),
\end{split}
\]
by \eqref{eq:W_expec}, \eqref{eq:ep_lower}, \eqref{eq:eta}, \eqref{eq:pB} and \eqref{eq:Z_W}. It follows that
\begin{multline}
\label{eq:786}
\int_\ep^\infty g(x)\,\P(e^{-A}Z \in \dd x\,|\,Z>\ep e^A) \\ = 1 + i\lambda \pi^2\gamma_0(1+\rho)^{-1}(-\log\ep + o(1))
 + \int_\ep^\infty h(x) \P(e^{-A}Z \in \dd x\,|\,Z>\ep e^A),
\end{multline}
where we define $h(x) =  g(x) - 1 - i\lambda (1+ \rho)^{-1}x \Ind_{(x\le 1)}$ for $x\ge 0$. Denote by $h^-(x)$ its left-hand derivative. Note that  $|h(x)|\le C(1\wedge x^2)$ and $|h^-(x)| \le C(x^{-1}\wedge x^2)$ for $x\ge 0$. By integration by parts,  \eqref{eq:W_tail} and \eqref{eq:Z_W}, we have
\begin{equation}
\label{eq:786a}
\begin{split}
  &\int_\ep^\infty h(x)\P(e^{-A}Z \in \dd x\,|\,Z>\ep e^A)\\
&= h(\ep) + p_B^{-1} (1+o(1))\Big(\int_\ep^\infty h^-(x)\P(Z > xe^A)\,\dd x + (h(1+)-h(1)) \P(Z > e^A)\Big)\\
&= \pi^2\gamma_0 \Big(\int_0^\infty h^-(x)\frac 1 x\,\dd x + (h(1+)-h(1)) + o(1)\Big)\\
&= \pi^2\gamma_0 (1+o(1)) \int_0^\infty h(x)\frac 1 {x^2}\,\dd x.
\end{split}
\end{equation}
Now, one readily sees that
\begin{equation}
\label{eq:786b}
 \int_0^\infty h(x)\frac 1 {x^2}\,\dd x =  i\lambda (c' + o(1)+O(\rho))\\
 + \int_0^\infty e^{i\lambda x}-1-i\lambda x\Ind_{(x\le 1)}\,\Lambda(\dd x),
\end{equation}
where $\Lambda(\dd x)$ is as in the statement of Proposition \ref{prop:piece} and $c'$ is a numerical constant. Equations \eqref{eq:785}, \eqref{eq:786}, \eqref{eq:786a} and \eqref{eq:786b} and the Taylor expansion of $e^{-x}$ at $x=0$ now yield
\begin{multline}
 \label{eq:789}
\E[e^{i\lambda\log(1+\frac{\Delta_{\mathrm{jump}}}{1+ \rho})}] \\
= \exp \Big[\gamma_0  \Big(i\lambda \pi^2 (-\log\ep  + c' + o(1) + \pi^2\int_0^\infty e^{i\lambda x}-1-i\lambda x\Ind_{(x\le 1)}\,\Lambda(\dd x)\Big) +O(\ep |\log\ep|\rho)\Big].
\end{multline}
Coming back to \eqref{eq:782}, we have by the Taylor expansion of $(1+x)^{i\lambda}$ at $x=0$,
\begin{align*}
&\E[e^{i\lambda (\log(1+\Delta_{\mathrm{drift}})+O(\ep |\log\ep| \Delta_{\mathrm{drift}}))}\Ind_{(|\Delta_{\mathrm{drift}}| < \ep^{1/3})}\Ind_{G_\Delta}]\\
& = \E[(1+i\lambda \Delta_{\mathrm{drift}} + O(\Delta_{\mathrm{drift}}^2))\Ind_{(T\le \sqrt\ep a^3)}] + O(\ep^{9/8}) && \text{by Lemma~\ref{lem:GDelta} and \eqref{eq:Delta_drift_bound}}\\
& = 1 + i\lambda \pi^2\gamma_0(A+\log\ep+c+o(1)) + \P(T> \sqrt\ep a^3) + O(\ep^{9/8}) && \text{by \eqref{eq:Delta_drift} and \eqref{eq:Delta_drift_2}}\\
& = \exp i\lambda \pi^2\gamma_0(A+\log\ep+c+o(1)+O(\ep^{1/8})),
\end{align*}
where the last equation follows from Lemma~\ref{lem:T0} and the Taylor expansion of $e^x$ at $x=0$. This equation, together with
\eqref{eq:782} and \eqref{eq:789} and the fact that $\Delta_{\mathrm{jump}}$ is independent from $\Delta_{\mathrm{drift}}$, $T$ and $G_{\mathrm{bulk}}$, yields \eqref{eq:piece_fourier} in the case $n=1$. 

\subsection{Proof of Theorems \ref{th:barrier} and \ref{th:barrier2}}
\label{sec:BBBM_proofs}
We define the process $(X_t'')_{t\ge 0}$ by
\[
X_t'' = X^{[\infty]}_{\Theta_{\lfloor t\gamma_0^{-1}a^3\rfloor}} - \pi^2A t.
\] 
\begin{proposition}
 \label{prop:Xt_second}
Suppose (HB). Then, as $A$ and $a$ go to infinity, the process $(X_t'')_{t\ge 0}$ converges in law (with respect to Skorokhod's topology) to the L\'evy process $(L_t)_{t\ge 0}$ defined in Theorem~\ref{th:1}.
\end{proposition}
\begin{proof}
Conditioning on $\F_0$, we can assume without loss of generality that the initial configuration $\nu_0$ is deterministic and that $G_0$ is verified. Denote by $(\F_t'')_{t\ge0}$ the natural filtration of the process $X_t''$ and note that $\F_t'' = \F_{\gamma_0 \lfloor t\gamma_0^{-1}\rfloor}'' \subset \F_{T_{\lfloor t\gamma_0^{-1}\rfloor}}$. In order to show convergence of the finite-dimensional distributions, it is enough to show (see Proposition~3.1 in \cite{Kurtz1975} or Lemma~8.1 in \cite{Ethier1986}, p.\ 225), that for every $\lambda\in\R$ and $t,s\ge 0$,
\begin{equation}
\begin{split}
\label{eq:910}
 \E\Big[\Big|\E[e^{i\lambda X_{t+s}''}\,|\,\F_t''] - e^{i\lambda X_t''}e^{s\kappa(\lambda)}\Big|\Big] \to 0,
\end{split}
\end{equation} 
as $A$ and $a$ go to infinity. Now, define $n := \lfloor t\gamma_0^{-1}\rfloor$ and $m := \lfloor (t+s)\gamma_0^{-1}\rfloor$, such that $(m-n)\gamma_0 = s+A^{-1}o(1)$, by \eqref{eq:ep_upper} and \eqref{eq:pB}. Then we have by Proposition \ref{prop:piece},
\begin{equation}
\label{eq:912}
 \begin{split}
  \E[e^{i\lambda (X_{t+s}''-X_t'')}\,|\,\F_{\Theta_m}]\Ind_{G_m} &= e^{-i\lambda As} \E[e^{i\lambda (X^{[\infty]}_{\Theta_m} - X^{[\infty]}_{\Theta_n})}\,|\,\F_{\Theta_m}]\Ind_{G_n}\\
&= \exp\Big(s\big(\kappa(\lambda) + o(1) + O(\ep^\numcst)\big)\Big)\Ind_{G_n}.
 \end{split}
\end{equation}
In total, we get for $A$ and $a$ large enough,
\[
 \E\Big[\Big|\E[e^{i\lambda (X_{t+s}''-X_t'')}\,|\,\F_t''] - e^{s\kappa(\lambda)}\Big|\Big] \le e^{s \kappa(\lambda)}\E[|e^{s(o(1) + O(\ep^\numcst))} - 1|] + \P(G_m^c).
\]
By Proposition \ref{prop:piece}, this goes to $0$ as $A$ and $a$ go to infinity, which proves \eqref{eq:910}.

In order to show tightness in Skorokhod's topology, we use Aldous' famous criterion \cite{Aldous1978} (see also \cite{Billingsley1999}, Theorem 16.10): If for every $M > 0$, every family of $(\F_t'')$-stopping times $\tau = \tau(A,a)$ taking only finitely many values, all of which in $[0,M]$ and every $h = h(A,a)\ge 0$ with $h(A,a) \to 0$ as $A$ and $a$ go to infinity, we have
\begin{equation}
\label{eq:920}
 X_{\tau+h}'' - X_{\tau}'' \to 0,\quad\text{ in probability as $A$ and $a$ go to infinity},
\end{equation}
then tightness follows for the processes $X_t''$ (note that the second point in the criterion, namely tightness of $X_t''$ for every fixed $t$, follows from the convergence in finite-dimensional distributions proved above). Now let $\tau$ be such a stopping time and let $V_\tau$ be the (finite) set of values it takes. We first note that since $G_n \supset G_{n+1}$ for every $n\in \N$, we have for every $t\in V_\tau$ and every $A$ and $a$ large enough,
\begin{equation}
\label{eq:922}
\P(G^c_{\lfloor t \gamma_0^{-1}\rfloor}) \le \P(G^c_{\lfloor M \gamma_0^{-1}\rfloor}) = O(M \ep^\numcst).
\end{equation}
by Proposition \ref{prop:piece}. Moreover, since $\F_t'' \subset \F_{T_{\lfloor t\gamma_0^{-1}\rfloor}}$ for every $t\ge 0$, we have for every $\lambda > 0$,
\begin{align*}
 \E[e^{i\lambda(X_{\tau+h}'' - X_{\tau}'')}] &= \sum_{t\in V_\tau}\E\Big[e^{i\lambda(X_{t+h}'' - X_t'')}\Ind_{(\tau = t)}\Big]\\
&= \sum_{t\in V_\tau}\E\Big[\E[e^{i\lambda(X_{t+h}'' - X_t'')}\,|\,\F_{T_{\lfloor t\gamma_0^{-1}\rfloor}}]\Ind_{(\tau = t)}\Ind_{G_{\lfloor t\gamma_0^{-1}\rfloor}}\Big] + O(M \ep^\numcst) && \text{by \eqref{eq:922}}\\
&= e^{h (\kappa(\lambda)+o(1) + O(\ep^\numcst))}(1-O(M\ep^\numcst)) + O(M \ep^\numcst), && \text{by \eqref{eq:912}},
\end{align*}
which converges to $1$ as $A$ and $a$ go to infinity. This implies \eqref{eq:920} and therefore proves tightness in Skorokhod's topology, since $M$ was arbitrary. Together with the convergence in finite-dimensional distributions proved above, the lemma follows.
\end{proof}

\paragraph{A coupling with a Poisson process.} 
Let $(V_n)_{n\ge 0}$ be a sequence of independent exponentially distributed random variable with mean $\gamma_0$. In order to prove convergence of the processes $X_t'$ and $X_t$, we are going to couple the BBM with the sequence $(V_n)$ in the following way: Suppose we have constructed the BBM until time $\Theta_{n-1}$. Now, on the event $G_{n-1}$, by Lemma~\ref{lem:coupling_T}, the strong Markov property of BBM and the transfer theorem (\cite{Kallenberg1997}, Theorem 5.10), we can construct the BBM up to time $\Theta_n$ such that $\P(|(T_n-\Theta_{n-1})/a^3 - V_n| > \ep^{3/2}) = O(\ep^2)$ (recall that $T_n$ denotes the time of the first breakout after $\Theta_{n-1}$). On the event $G_{n-1}^c$, we simply let the BBM evolve independently of $(V_j)_{j\ge n}$. Setting $G_n' = G_n \cap \{\forall j\le n: |(T_j-\Theta_{j-1})/a^3 - V_j| \le \ep^{3/2}\}$, there is by Lemma~\ref{lem:coupling_T} and Proposition~\ref{prop:piece} a $\numcst > 0$, such that for large $A$ and $a$,
\begin{equation}
 \label{eq:prob_Gnprime}
\P(G_n') \ge \P(G_0) - nO(\ep^{1+\numcst})
\end{equation}
Furthermore, we have $\Theta_n = T_n + a^{5/2}$ on $G_n'$, whence for large $A$ and $a$,
\begin{equation}
\label{eq:930}
 \ton G_n': |(\Theta_n-\Theta_{n-1})/a^3 - V_n| \le 2\ep^{3/2}.
\end{equation} 
This construction now permits us to do the

\begin{proof}[Proof of Theorem~\ref{th:barrier2}]
Let $d$ denote the Skorokhod metric on the space of cadlag functions $D([0,\infty))$ (see \cite{Ethier1986}, Section 3.5). Let $\Phi$ be the space of strictly increasing, continuous, maps of $[0,\infty)$ onto itself. Let $x,x_1,x_2,\ldots$ be elements of $D([0,\infty))$. Then (\cite{Ethier1986}, Proposition 3.5.3), $d(x_n,x) \to 0$ as $n\to\infty$ if and only if for every $M > 0$ there exists a sequence $(\varphi_n)$ in $\Phi$, such that
\begin{equation}
 \label{eq:940}
\sup_{t\in [0,M]}|\varphi_n(t) - t| \to 0,
\end{equation}
and
\begin{equation}
 \label{eq:941}
\sup_{t\in [0,M]} |x_n(\varphi_n(t)) - x(t)| \to 0.
\end{equation}
If $(x_n')_{n\in \N}$ is another sequence of functions in $D([0,\infty))$, with $d(x_n',x)\to 0$, then by the triangle inequality and the fact that $\Phi$ is stable under the operations of inverse and convolution, we have $d(x_n,x)\to 0$ if and only if there exists a sequence $(\varphi_n)$ in $\Phi$, such that \eqref{eq:940} holds and
\begin{equation}
 \label{eq:942}
\sup_{t\in [0,M]} |x_n(\varphi_n(t)) - x_n'(t)| \to 0.
\end{equation}

For every $A$ and $a$, we define the (random) map $\varphi_{A,a} \in \Phi$ by
\[
 \varphi_{A,a}(\gamma_0(n+r)) = ((1-r) \Theta_n + r \Theta_{n+1})a^{-3}, \text{ with } n\in\N,\ r\in [0,1].
\]
Let $M > 0$ and define $n_M = \lceil M\gamma_0 \rceil$.
Then we have by the monotonicity of $(X^{[\infty]}_t)_{t\ge0}$,
\begin{equation}
 \label{eq:944}
\sup_{t\in [0,M]} |\varphi_{A,a}(t) - t| \le  \max_{n\in\{0,\ldots,n_M\}} \left|a^{-3}\Theta_n - \gamma_0 n\right|,
\end{equation}
and
\begin{equation}
 \label{eq:945}
\sup_{t\in [0,M]}|X''_t - X'_{\varphi_{A,a}(t)}| \le  \max_{n\in\{0,\ldots,n_M\}} A \left|a^{-3}\Theta_n - \gamma_0 n\right|.
\end{equation}
By Doob's $L^2$ inequality and the fact that $\gamma_0 = E[V_1]$, we get
\[
 \P\Big(\max_{n\in\{0,\ldots,n_M\}} \left|\sum_{i=1}^n V_i - n\gamma_0\right| > \ep^{1/3}\Big) \le C \ep^{-2/3} n_M \Var(V_i) = O(\ep^{1/3}).
\]
Furthermore, on the set $G_{n_M}'$, we have for every $n\le n_M$,
\[
\left|\Theta_n - \sum_{i=0}^n V_i\right| \le O(n_M \ep^{3/2}) = O(\ep^{1/2}).
\]
In total, we get with \eqref{eq:944} and \eqref{eq:945}, as $A$ and $a$ go to infinity,
\begin{equation}
\label{eq:948}
 \forall M > 0:\sup_{t\in [0,M]} |\varphi_{A,a}(t) - t| \vee |X''_t - X'_{\varphi_{A,a}(t)}| \to 0,\quad\text{in probability},
\end{equation}
which is equivalent to 
\begin{equation}
\label{eq:949}
 \sum_{M=1}^\infty 2^{-M}\Big[1\wedge\Big(\sup_{t\in [0,M]} |\varphi_{A,a}(t) - t| \vee |X''_t - X'_{\varphi_{A,a}(t)}|\Big)\Big] \to 0,\quad\text{in probability}.
\end{equation}
Now, suppose that $A$ and $a$ go to infinity along a sequence $(A_n,a_n)_{n\in\N}$ and denote by $X'_{A_n,a_n}$, $X''_{A_n,a_n}$ and $\varphi_{A_n,a_n}$ the processes corresponding to these parameters. By Proposition \ref{prop:Xt_second} and Skorokhod's representation theorem (\cite{Billingsley1999}, Theorem 6.7), there exists a probability space, on which the sequence $(X''_{A_n,a_n})$ converges almost surely as $n\to\infty$ to the limiting L\'evy process $L = (L_t)_{t\ge 0}$ stated in the theorem. Applying again the representation theorem as well as the transfer theorem, we can transfer the processes $X'_{A_n,a_n}$ and $\varphi_{A_n,a_n}$ to this probability space in such a way that the convergence in \eqref{eq:949} holds almost surely, which implies that the convergence in \eqref{eq:948} holds almost surely as well. By the remarks at the beginning of the proof, it follows that on this new probability space,
\[
 d(X'_{A_n,a_n},L) \le d(X'_{A_n,a_n},X''_{A_n,a_n}) + d(X''_{A_n,a_n},L)\to 0,
\]
almost surely, as $n\to\infty$. This proves the theorem.
\end{proof}

\begin{proof}[Proof of Theorem~\ref{th:barrier}]
Let $(t_i^{A,a})_{i=1}^n$ be as in the statement of the theorem. By the virtue of Theorem~\ref{th:barrier2}, it suffices to show that
\begin{equation}
\label{eq:950}
 \P\Big(\forall i: X^{[\infty]}_{t_i^{A,a}a^3} = J_{t_i^{A,a}a^3}\Big) \to 1.
\end{equation}
Suppose for simplicity that $t_i^{A,a}\equiv t_i$ for all $i$; the proof works exactly the same in the general case.
Let $n := \lceil 2(t_k+2)/\gamma_0\rceil$, such that $n = O(\ep^{-1})$, by \eqref{eq:pB}. By Chebychev's inequality, we then have
\begin{equation}
 \label{eq:952}
\P(\sum_{i=1}^n V_i \le t_k+2) \le \P(\sum_{i=1}^n (V_i - \gamma_0) \le -\frac n 2 \gamma_0) = O(n\Var(V_i)) = O(\ep).
\end{equation}
Furthermore, define the intervals $I_i = t_i + [-2n\ep^{3/2}-e^A/a,2n\ep^{3/2}]$, $i=1,\ldots,k$ and denote by $\mathscr P$ the point process on the real line with points at the positions $V_1,V_1+V_2,V_1+V_2+V_3,\ldots$. Then $\mathscr P$ is a Poisson process with intensity $\gamma_0^{-1} = O(\ep^{-1})$ and thus,
\begin{equation}
 \label{eq:954}
\P\Big(\mathscr P \cap \bigcup_{i=1}^k I_i \ne \emptyset\Big) = O(\ep^{1/2}).
\end{equation}
We now have
\begin{align*}
 \P\Big(\forall i: X^{[\infty]}_{t_ia^3} = J_{t_ia^3}\Big) &\ge \P\Big(\nexists (i,j): t_ia^3\in[\Theta_j-T_{j-1},\Theta_j]\Big) && \text{by definition}\\
&\ge \P\Big(G_{n}',\ \sum_{i=1}^{n}V_i > t_k+2,\ \mathscr P \cap \bigcup_{i=1}^k I_i = \emptyset\Big) && \text{by definition of $G_n'$}\\
&\ge \P(G_0) - O(\ep^\numcst) && \text{by \eqref{eq:prob_Gnprime}, \eqref{eq:952}, \eqref{eq:954}}.
\end{align*}
Letting $A$ and $a$ go to infinity and using the hypothesis (HB) yields \eqref{eq:950} and thus proves the theorem.
\end{proof}

\section{The \texorpdfstring{\Bfl}{Bb}-BBM}
\label{sec:Bflat}
In this section, we define and study the B$^\flat$-BBM, which is obtained from the B-BBM by killing some of its particles. It will be used in Section~\ref{sec:NBBM} to bound the $N$-BBM from below. 

\subsection{Definition of the model}
We will use throughout the notation from Section~\ref{sec:BBBM}. Furthermore, we fix $\delta \in (0,1/2)$ and $K\ge 1$ such that $E_K \le \delta/10$, where $E_K$ was defined in \eqref{eq:def_E}. Define $N = \lceil 2\pi  e^{A+\delta}a^{-3}e^{\mu a}\rceil$, where $\mu = \sqrt{1 - \pi^2/a^2}$. We will then use interchangeably the phrases ``as $A$ and $a$ go to infinity'' and ``as $N$ go to infinity''. Furthermore, in the definition of these phrases (see Section~\ref{sec:before_breakout_definitions}, $A_0$ and the function $a_0(A)$ may now also depend on $\delta$. The symbols $C_\delta$ and $C_{\delta,\alpha}$ have the same meaning as $C$ (see Section~\ref{sec:interval_notation}), except that they may depend on $\delta$ or $\delta$ and $\alpha$ as well, $\alpha$ being defined later.

\begin{figure}[h]
 \centering
 \def\svgwidth{8cm}
\executeiffilenewer{sec8.svg}{sec8.pdf}%
{inkscape -z -D --file=sec8.svg %
--export-pdf=sec8.pdf --export-latex --export-area-drawing}%
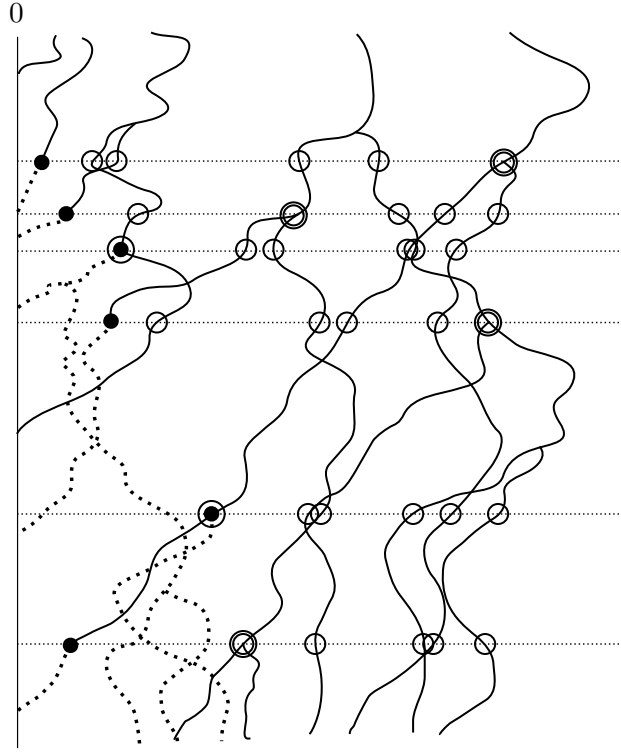%

 \caption{The \Bfl-BBM with parameter $N=6$ (no breakout is shown). White and red particles are drawn with solid and dotted lines, respectively. The blobs indicate when a white particle is colored red, the circles show the six white particles living at that time.}
 \label{fig:sec8}
\end{figure}

The B$^\flat$-BBM is then defined as follows: Given a possibly random initial configuration $\nu_0$ of particles in $(0,a)$, we let particles evolve according to B-BBM with barrier function given by \eqref{eq:wall_fn}, where, in addition, we colour the particles white and red as follows: Initially, all particles are coloured white. As soon as a white particle has $N$ or more white particles to its right, it is coloured red\footnote{This can be ambiguous if there are more than one of the particles at the same position, for example when the left-most particle branches. In order to eliminate this ambiguity, induce for every $t\ge0$ a total order on the particles in $\mathscr N(t)$ by $u<v$ iff $X_u(t) < X_v(t)$ or $X_u(t)=X_v(t)$ and $u$ precedes $v$ in the \emph{lexicographical order on $U$}. Whenever there are more than $N$ white particles, we then colour the particles in this order, which is well defined.}. Children inherit the colour of their parent. At each time $\Theta_n$, $n\ge1$, 
all the red particles are killed immediately and the process goes on with the remaining particles. See Figure~\ref{fig:sec8} for a graphical description.

For an interval $I\subset\R_+$, we define the stopping line $\mathscr L^{\mathrm{red}}_{\emptyset,I}$ by $(u,t)\in \mathscr L^{\mathrm{red}}_{\emptyset,I}$ if and only if the particle $u$ gets coloured red at the time $t$ and has been white up to the time $t$, with $t\in I$. We then set $Z^{\mathrm{red}}_{\emptyset,I}$ and  $Y^{\mathrm{red}}_{\emptyset,I}$ by summing respectively $w_Z$ and $w_Y$ over the particles of this stopping line. Furthermore, we define $\mathscr N^{\mathrm{red}}_t$ and $\mathscr N^{\mathrm{white}}_t$ to be the subsets of $\mathscr N_t$ formed by the red and white particles, respectively, and define $Z^{\mathrm{red}}_t$, $Y^{\mathrm{red}}_t$, $Z^{\mathrm{white}}_t$ and $Y^{\mathrm{white}}_t$ accordingly. 

Let $\nu^\flat_t$ be the configuration of white particles at the time $t$ and abuse notation by setting $\nu^\flat_n =\nu^\flat_{\Theta_n}$. We set $G^{\flat}_{-1}=\Omega$ and for each $n\in\N$, we define the event $G^\flat_n$ to be the intersection of $G^\flat_{n-1}$ with the following events (we omit the braces).
\begin{itemize}
 \item $\supp\nu^\flat_n\subset (0,a)$,
 \item $\mathscr N^{\mathrm{white}}_{\Theta_n} \subset U\times \{\Theta_n\}$ and $\Theta_n > T_n^+$ (for $n>0$),
 \item $|e^{-A}Z^{\mathrm{white}}_{\Theta_n} -1| \le \ep^{3/2}$ and $Y^{\mathrm{white}}_{\Theta_n} \le \eta$.
 \item $\P\Big(Z^{\mathrm{red}}_{\emptyset,[\Theta_n,\Theta_n+Ka^2]} + Y^{\mathrm{red}}_{\emptyset,[\Theta_n,\Theta_n+Ka^2]} \le a^{-1/2}\,\Big|\,\F_{\Theta_n}\Big) \ge 1-\ep^2$.
\end{itemize}
The last event is of course uniquely defined up to a set of probability zero. Note that $G^\flat_n\in\F_{\Theta_n}$ for each $n\in\N$.
Furthermore, we define the predicates
\begin{itemize}
 \item[(HB$^\flat$)] (The law of) $\nu_0$ is such that $\P(G^\flat_0)\to 1$ as $A$ and $a$ go to infinity.
 \item[(HB$^\flat_0$)] $\nu_0$ is deterministic and such that $G^\flat_0$ holds.
 \item[(HB$^\flat_{\perp}$)] $\nu_0$ is obtained from $\lfloor e^{-\delta}N\rfloor$ particles distributed independently according to the density proportional to $\sin(\pi x/a)e^{-\mu x}\Ind_{(0,a)}(x)$.
\end{itemize}
We now state the important results on the \Bfl-BBM.
\begin{lemma}
\label{lem:HBflat}
(HB$^\flat_{\perp}$) implies (HB$^\flat$) for large $A$ and $a$.
\end{lemma}

\begin{proposition}
  \label{prop:Bflat}
Proposition~\ref{prop:piece} still holds for the B$^\flat$-BBM, with $G_n$ replaced by $G_n^\flat$. The same is true for Theorems~\ref{th:barrier} and \ref{th:barrier2}, with (HB) replaced by (HB$^\flat$).
\end{proposition}

Recall the definition of $\med^N_\alpha$ and $x_\alpha$ from the introduction and of (Ht) from Section~\ref{sec:BBBM}.
\begin{proposition}
\label{prop:Bflat_med}
Suppose (HB$^\flat_\perp$). Let $(t^{A,a}_j)$ satisfy (Ht) and let $\alpha \in (0,e^{-2\delta})$. As $A$ and $a$ go to infinity,
\[
\P(\forall j: \med^N_\alpha(\nu^{\flat}_{a^3t^{A,a}_j}) \ge x_{\alpha e^{2\delta}})\to 1.
\]
\end{proposition}

\subsection{Preparatory lemmas}

In this section, we will establish some fine estimates for the number of particles  of the process, which will be used later to bound the number of creations of red particles. If a particle $u\in\mathscr N_0(t)$ satisfies $\tau_l(u) \le t < \sigma_{l+1}(u)$ for some $l\ge 1$, we call it an ``in between'' particle (because it is in some kind of intermediary stage) and otherwise a ``regular'' particle. We then define
\[
 N_t(r) = \sum_{(u,s)\in\mathscr N_t}\Ind_{(X_u(s) \ge r)}.
\]
The main lemma in this section is the following:
\begin{lemma}
\label{lem:Bflat_N} Suppose (HB$_0$). Let $Ka^2\le t \le \tconst{eq:tcrit}$.  Then, for $0\le r\le 9a/10$ and every $\alpha > 0$ there exists $C_{\delta,\alpha}$, such that for large $A$ and $a$,
\[
\P(N_t(r) \ge e^{-(\delta/2+ r/3)} N\,|\,T > t) \le C_{\delta,\alpha} A^2\ep\Big(\frac t {a^3}+\eta\Big)e^{- (\frac 4 3 - \alpha)  r}
\]
Furthermore, conditioned on $\F_\Delta$, for $t\le \Theta_1+e^Aa^2$, for large $A$ and $a$,
\[
\P_{\mathrm{nbab}}(N_t(r) \ge e^{-(\delta/2+ r/3)} N\,|\,\F_\Delta)\Ind_{G_\Delta}\le C_{\delta,\alpha} A^2\ep^2e^{- (\frac 4 3 - \alpha)  r}.
\]
\end{lemma}
The following lemma about BBM conditioned not to break out will be used many times in the proof:

\begin{lemma}
  \label{lem:N_with_tiers}
Let $f$ be a barrier function, $0\le t\le t_0\le \tconst{eq:tcrit}$ and suppose that $\Err(f,t)\le C$. For $x\in(0,a)$, define $\Phat_f^x = \P_f^x(\cdot\,|\,T>t_0)$.   Then, for all $r\le (9/10)a$ and for large $A$ and $a$,
\begin{align}
  \label{eq:66}
&\Ehat_f^x[N^{(1+)}_{t}(r)] \le C Ae^{-A}N(1+\mu r)e^{-\mu r}\big(1+(\eta r^2/a)\big)\Big(\frac {t} {a^3}w_Z(x)+w_Y(x)\Big)\\
  \label{eq:67}
&\Ehat_f^x[(N_{t}(r))^2] \le C (1+r^4)e^{- 2  r} A^2 \ep e^{-A} N^2 \Big(\frac {t} {a^3}w_Z(x)+w_Y(x)\Big).
\end{align}
\end{lemma}
\begin{proof}
By Lemma~\ref{lem:R_f}, the upper bound of Lemma~\ref{lem:Zl_exp_upperbound} is still valid with varying drift (see also Remark~\ref{rem:hat_f}). This gives,
\begin{align}
  \label{eq:28}
\Ehat^x_f[Z^{(1+)}_{\emptyset,[0,t]}] \le CA\Big(\frac {t} {a^3}w_Z(x)+w_Y(x)\Big),\quad\tand\quad  \Ehat^x_f[Y^{(1+)}_{\emptyset,[0,t]}] \le \eta \Ehat^x_f[Z^{(1+)}_{\emptyset,[0,t]}].
\end{align}
where the second inequality follows from Lemma~\ref{lem:critical_line}.
Equation \eqref{eq:66} then follows from Lemma~\ref{lem:N_expec_upper_bound} and \eqref{eq:28}, noting that $a-y-f(\zeta/a^2) \ge r+a/20$ for large $a$, by the hypotheses on $r$ and $f$.

For the second equation, suppose for simplicity that $f\equiv 0$. By Lemma~\ref{lem:many_to_two_with_tiers}, we have for all $x\in(0,a)$,
  \begin{multline}
    \label{eq:27}
    \Ehat^x[(N_{t}(r))^2] \le \Ehat^x[N_{t}(r)] + \beta_0 m_2\int_0^t\int_0^a
    \widehat{\p}^{f,(0+)}_s(x,z) \Ehat^{(z,s)}[N_{t}(r)]^2\,\dd z\,\dd s\\
+ \Ehat^x\Big[\sum_{(u,s)\in\mathscr R_{t}} \sum_{(v_1,s_1),(v_2,s_2)\in\mathscr
  S^{(u,s)},\,v_1\ne v_2} \Ehat^{(X_{v_1}(s_1),s_1)}[N_{t}(r)]\Ehat^{(X_{v_2}(s_2),s_2)}[N_{t}(r)]\Big].
  \end{multline}
Now, for large $A$ we have by \eqref{eq:ptilde_estimate_f} and Lemma~\ref{lem:mu_t_density},
\begin{multline}
  \label{eq:30}
\int_0^t\int_0^a \widehat{\p}^{f,(0+)}_s(x,z)\Ehat^{(z,s)}[N_{t}(r)]^2\,\dd z\,\dd s \\
 \le C \Ehat\Big[\sum_{(v,s_0)\in\mathscr S_{t}} \int_{s_0}^t\int_0^a \p_{s-s_0}(X_v(s_0),z)\Ehat^{(z,s)}[N_{t}(r)]^2\,\dd z\,\dd s\Big]
\end{multline}
With Lemmas~\ref{lem:N_expec_upper_bound} and \ref{lem:Zl_exp_upperbound} and \eqref{eq:66}, we have as in the proof of Lemma~\ref{lem:N_2ndmoment}, for any $s_0\le t$,
\begin{equation}
  \label{eq:29}
  \int_{s_0}^t\int_0^a \p_{s-s_0}(x,z)  \Ehat^{(z,s)}[N_{t}(r)]^2\,\dd z\,\dd s \le CA^2e^{-2A}N^2(1+r^4)e^{-2\mu r}\Big(\frac t {a^3} w_Z(x)+w_Y(x)\Big).
\end{equation}
Equations \eqref{eq:30} and \eqref{eq:29} together with Lemma~\ref{lem:Zl_exp_upperbound} now give
\begin{equation}
\label{eq:29b}
\int_0^t\int_0^a \widehat{\p}^{f,(0+)}_s(x,z)\Ehat^{(z,s)}[N_{t}(r)]^2\,\dd z\,\dd s \le 
CA^2e^{-2A}N^2(1+r^4)e^{-2\mu r}\Big(\frac t {a^3} w_Z(x)+A w_Y(x)\Big).
\end{equation}
Moreover, we have
\begin{equation}
  \label{eq:31}
\begin{split}
  \Ehat^x\Big[\sum_{(u,s)\in\mathscr R_{t}} &\sum_{(v_1,s_1),(v_2,s_2)\in\mathscr
  S^{(u,s)},\,s_1<s_2} \Ehat^{X_{v_1}(s_1)}[N_{t}(r)]\Ehat^{X_{v_2}(s_2)}[N_{t}(r)]\Big]\\
&\le \Ehat^x\Big[\sum_{(u,s)\in\mathscr R_{t}} \Big(\sum_{(v,s_0)\in\mathscr
  S^{(u,s)}} \Ehat^{X_v(s_0)}[N_{t}(r)]\Big)^2\Big]\\
&\le C (1+r^4)e^{-2\mu r} A^2e^{-2A}N^2 \Ehat^x\Big[\sum_{(u,s)\in\mathscr R_{t}} (Z^{(u,s)})^2\Big]\\
&\le C (1+r^4)e^{-2\mu r} \ep e^{-A} A^2N^2 \Big(\frac {t} {a^3} w_Z(x) + w_Y(x)\Big),
\end{split}
\end{equation}
by \eqref{eq:QaZsquared} and Lemma~\ref{lem:Rt}. Equation \eqref{eq:67} now follows from \eqref{eq:ep_lower}, \eqref{eq:66}, \eqref{eq:27}, \eqref{eq:29b} and \eqref{eq:31}.
\end{proof}

\begin{corollary}
\label{cor:N_with_tiers}
With the hypotheses from Lemma~\ref{lem:N_with_tiers}, we have for large $A$ and $a$,
\[
\Ehat_f^x[N_{t}(r)] \le C Ae^{-A}N(1+ r^2)e^{-\mu r}e^{-(a-x)/2}.
\]
\end{corollary}
\begin{proof}
  This follows immediately from Lemmas~\ref{lem:mu_t_density}, \ref{lem:N_expec_upper_bound} and \ref{lem:N_with_tiers} and Corollary~\ref{cor:hat_many_to_one}.
\end{proof}

\begin{lemma}
  \label{lem:N_bar_with_tiers}
Suppose (HB$_0$). Let $\alpha > 0$. For large $A$ and $a$, we have for every $0\le r\le 9a/10$ and $t\le \Theta_1+e^Aa^2$, 
  \begin{equation}
    \label{eq:65}
  \P_{\mathrm{nbab}}(\overline N_t(r)>e^{-A/6- r/3}N\,|\,\F_\Delta)\Ind_{G_\Delta}\le C_\alpha A^2\ep e^{-2A/3}e^{-(\frac 4 3 -\alpha) r}.
  \end{equation}
\end{lemma}
\begin{proof}
Conditioned on $\F_\Delta$, define $\Phat = \P(\cdot\,|\,T>\Theta_1+e^Aa^2)$. Then,
\begin{equation}
\label{eq:63}
\Var_{\mathrm{nbab}}(\overline N_t(r)\,|\,\F_\Delta)\Ind_{G_\Delta}\le C  \sum_{i=0}^l \sum_{u\preceq\mathscr U} \Ind_{(\sigma_i\le d_u < \tau_i)} (k_u-1)\Ehat^{(X_u(d_u),d_u)}_{X^{[1]}}[N_t(r)^2],
\end{equation}
and Lemma \ref{lem:N_with_tiers} now implies
\begin{equation}
  \label{eq:68}
 \Var_{\mathrm{nbab}}(\overline N_t(r)\,|\,\F_\Delta)\Ind_{G_\Delta}\le C (1+r^4)e^{- 2  r} A^2 \ep e^{-A} N^2 \mathscr E_{\mathscr U}. 
\end{equation}
Furthermore, by Corollary~\ref{cor:N_with_tiers},
\begin{equation}
  \label{eq:70}
  \E_{\mathrm{nbab}}[\overline N_t(r)\,|\,\F_\Delta]\Ind_{G_\Delta} \le  C Ae^{-A}N(1+r^2)e^{-\mu r} \mathscr E_{\mathscr U}.
\end{equation}
Equation \eqref{eq:65} now follows from \eqref{eq:68} and \eqref{eq:70} together with the conditional Chebychev inequality and the fact that $\mathscr E_{\mathscr U}\le e^{A/3}$ on $G_{\mathscr U}\subset G_\Delta$.
\end{proof}


\begin{proof}[Proof of Lemma \ref{lem:Bflat_N}]
  Assume  $Ka^2 \le t \le  \tconst{eq:tcrit}$. By the hypothesis on $Z_0$, the definition of $K$, Proposition~\ref{prop:N_expec} and Corollary~\ref{cor:hat_many_to_one}, we have for large $A$ and $a$,
\begin{equation}
  \label{eq:22}
\E[N^{(0)}_{t}(r)\,|\,T>t] \le (1-7\delta/8)(1+\mu r)e^{-\mu r}N, 
\end{equation}
Moreover, we have by Lemma~\ref{lem:N_with_tiers},
 \begin{equation}
   \label{eq:23}
   \begin{split}
\E[N^{(1+)}_{t}(r)\,|\,T>t] &\le C(1+\mu r) e^{-\mu r}N (1+\eta r) A e^{-A}(\frac t {a^3}Z_0 + Y_0)\\
& \le C\sqrt{\ep}A(1+r^2)e^{-\mu r}N,
   \end{split}
 \end{equation}
by the hypotheses on $Z_0$ and $Y_0$. Equations \eqref{eq:22} and \eqref{eq:23} now give for large $A$ and $a$, for all $r\in[0,a]$,
\begin{equation}
  \label{eq:33}
   \E[ N_t(r)\,|\,T>t] \le (e^{-\delta/2} - \delta/8) e^{- r/3} N,
\end{equation}
Furthermore, by Lemma~\ref{lem:N_with_tiers} and the hypotheses on $Z_0$ and $Y_0$, we have
\begin{equation}
  \label{eq:75}
\Var( N_t(r)\,|\,T>t) \le CA^2 \ep (1+r^4)e^{-2r} (t/a^3 +\eta),
\end{equation}
Chebychev's inequality, \eqref{eq:33} and \eqref{eq:75} yield the first equation of the lemma.

Conditioned on $\F_\Delta$, let $T \le t \le \Theta_1+e^Aa^2$. Define $\widetilde N^{\mathrm{bulk}}_t(r)$ to be the number of hat- and check-particles to the right of $r$ at time $t$ that have not hit $a$ between $T^-$ and $t$ and likewise $\widetilde N^{\mathrm{fug}}_t(r)$ the number of fug-particles with the same properties. Set $M_t = e^{- X^{[1]}_{t}} = (1+\Delta_t)^{-1}$, where $\Delta_t =  \thbar((t-T^+)/a^2)\Delta$. Note that we have on $G_\Delta$:  $|\Delta - e^{-A}Z^{(\mathscr
     U,T)}| \le \ep^{1/8}$.


Define $\alpha_r = (1+\mu r)e^{-\mu r}$. By Lemma~\ref{lem:mu_t_density} and Proposition~\ref{prop:N_expec}, we have for large $A$ and $a$,
\begin{equation}
  \label{eq:48}
  \begin{split}
  \E[\widetilde N^{\mathrm{bulk}}_t(r)\,|\,\F_\Delta]\Ind_{G_\Delta} & \le \alpha_re^{o(1)}NM_t(1+E_K)e^{-A}(\widehat Z_{T^-}+\widecheck Z_{\emptyset,T^-})\Ind_{G_\Delta}\\ 
&\le \alpha_rNM_t(1-3\delta/4),
  \end{split}
\end{equation}
by the definition of $\widehat G$ and $\widecheck G$.
Furthermore, we have by Proposition~\ref{prop:N_expec} and Lemma~\ref{lem:mu_t_density}, for large $a$,
\begin{equation}
  \label{eq:51}
  \begin{split}
  \E[\widetilde N^{\mathrm{fug}}_t(r)\,|\,\F_\Delta] \Ind_{G_\Delta} &\le \alpha_re^{o(1)} N M_t \left(\thbar\left((t-T)/a^2\right) + O\left((y+\Delta+r)^2\eta/a\right)\right) e^{-A}Z^{(\mathscr U,T)}\Ind_{G_\Delta}\\
&\le \alpha_r N M_t \left(\Delta_t+ O(\ep^{1/8}+a^{-1}r^2)\right)\Ind_{G_\Delta},
  \end{split}
\end{equation}
Equations \eqref{eq:48} and \eqref{eq:51} now give,
\begin{equation}
  \label{eq:55}
\E[\widetilde N^{\mathrm{bulk}}_t(r)+\widetilde N^{\mathrm{fug}}_t(r)\,|\,\F_\Delta]\Ind_{G_\Delta} \le \alpha_r  N \left(1-\delta/2 + O(a^{-1}r^2)\right).
\end{equation}
Similarly, one has by Lemma~\ref{lem:N_2ndmoment},
\begin{equation}
  \label{eq:40}
\Var(\widetilde N^{\mathrm{bulk}}_t(r)+\widetilde N^{\mathrm{fug}}_t(r)\,|\,\F_\Delta) \Ind_{G_\Delta}\le Ce^{-A} (1+r^4) e^{-2\mu r}.
\end{equation}
Equations \eqref{eq:55} and \eqref{eq:40} and the conditional Chebychev inequality now yield for large $A$ and $a$,
\[
\P(\widetilde N^{\mathrm{bulk}}_t(r)+\widetilde N^{\mathrm{fug}}_t(r) \ge e^{-2\delta/5- r/3}N\,|\,\F_\Delta)\Ind_{G_\Delta} \le C_\delta (1+r^4)e^{-A}e^{-\frac 4 3 \mu r}.
\]
This, together with Lemma \ref{lem:N_bar_with_tiers} and the fact that the hat-, fug- and check-particles do not hit $a$ on $G_{\mathrm{nbab}}$ finishes the proof of the lemma.
\end{proof}

We finish this section with a result which fills the gap in Lemma~\ref{lem:Bflat_N} under the hypothesis (HB$^\flat_\perp$).

\begin{lemma}
  \label{lem:Bflat_start}
Suppose (HB$^\flat_\perp$). Then $\P(|e^{-A}Z_0 - 1| \le \ep^{3/2},\ Y_0 \le \eta) \le o(1)$. Moreover, let $\alpha > 0$. Then for large $A$ and $a$, we have for every $t\le Ka^2$ and $0\le r\le 9a/10$,
\[
\P(N_t(r) \ge e^{-(\delta/2+ r/3)} N,\,R_{Ka^2} = 0) \le C_{\alpha,\delta} e^{-(\frac 4 3 -\alpha) r}/a,
\]
and $\P(R_{Ka^2} = 0) \ge 1-C_\delta e^A/a$.
\end{lemma}
\begin{proof}
Recall that initially,  there are $\lfloor e^{-\delta}N\rfloor$ particles distributed independently according to the density $\phi(x)$ proportional to $\sin(\pi x/a)e^{-\mu x}\Ind_{(0,a)}(x)$. An elementary calculation yields that
\begin{equation}
  \label{eq:49}
  \E[Z_0] = e^A(1+o(1)),\quad\Var(Z_0) \le Ce^A/a^3,\quad \E[Y_0] \le Ce^A/a.
\end{equation}
This immediately yields the first statement, by Chebychev's and Markov's inequalities. Moreover, \eqref{eq:49} with Lemmas~\ref{lem:Rt} and \ref{lem:R_f} yields
  \begin{equation}
\label{eq:RKa2}
\P(R_{Ka^2} \ge 1) \le \E[R_{Ka^2}] \le C \Big(\frac K a \E[Z_0] + \E[Y_0]\Big) \le C \frac {Ke^A} a,
  \end{equation}
which gives the last statement. Note that on the event $\{R_{Ka^2} = 0\}$, the B$^\flat$-BBM equals BBM with absorption at $0$ and $a$. Since the density $\phi(x)$ is stationary w.r.t. this process, a quick calculation shows that
\begin{equation}
  \label{eq:69}
  \E[N_t^{(0)}(r)] =  \E[N_0^{(0)}(r)] = \lfloor e^{-\delta}N\rfloor \int_r^a \phi(x)\,\dd x \le e^{-\delta}N (1+\mu r) e^{-\mu r}.
\end{equation}
 Furthermore, by the independence of the initial particles and the law of total variance,
\begin{equation}
  \label{eq:83}
  \Var(N_t^{(0)}(r)) = \lfloor e^{-\delta}N\rfloor \Big(\Var(\E^X[N_t^{(0)}(r)]) + \E[\Var^X(N_t^{(0)}(r))]\Big),
\end{equation}
where $X$ is a random variable distributed according to the density $\phi$ and the outer variance and expectation are with respect to $X$.  
By Lemma~\ref{lem:N_2ndmoment}, we have for every $x\in [0,a]$,
\[
  \Var^x(N_t^{(0)}(r)) \le \E^x[N_t^{(0)}(r)^2] \le Ce^{-2A}N^2(1+r^4)e^{-2\mu r}((t/a^3)w_Z(x)+w_Y(x)),
\]
and a simple calculation then yields for $t\le Ka^2$,
\begin{equation}
  \label{eq:84}
 \E[\Var^X(N_t^{(0)}(r))] = \int_0^a  \Var^x(N_t^{(0)}(r)) \phi(x)\,\dd x \le C(K/a)e^{-A} N (1+r^4)e^{-2\mu r}.
\end{equation}
Moreover, by Lemma~\ref{lem:N_expec_upper_bound} and the hypothesis on $r$, we have for every $x\in[0,a]$,
\[
\E^x[N_t^{(0)}(r)] \le C(1+r^2)e^{-\mu r}\Big(e^{\mu x}\Ind_{(x\le \frac {19a} {20})} + e^{-A}N(w_Z(x)+w_Y(x))\Ind_{(x > \frac {19a} {20})}\Big),
\]
which yields
\begin{equation}
  \label{eq:50}
\begin{split}
 \E\left[ (\E^X[N_t^{(0)}(r)])^2\right] &\le C(1+r^4)e^{-2\mu r}\Big(\int_0^{\frac {19a} {20}} a e^{\mu x}\,\dd x + e^{-2A}N^2\int_{\frac {19a} {20}}^a (a-x+1)^3 e^{\mu(x-2a)}\,\dd x\Big) \\
& \le Ca^{-3}e^{-A}N(1+r^4)e^{-2\mu r}.
\end{split}
\end{equation}
Equations \eqref{eq:83}, \eqref{eq:84} and \eqref{eq:50} now yield for large $a$,
\begin{equation}
  \label{eq:85}
  \Var(N_t^{(0)}(r)) \le C(K/a)e^{-A} N (1+r^4)e^{-2\mu r}.
\end{equation}
The lemma now follows from \eqref{eq:69} and \eqref{eq:85}, together with Chebychev's inequality.
\end{proof}

\subsection{The probability of \texorpdfstring{$G^{\flat}_1$}{Gb1}}
Most of the work in this section will be devoted to bounding the number of creations of red particles. For this, we will discretize time: We set $t_n = Ka^2+n\delta/4$ and $I_n = [t_n,t_{n+1})$ for all $n\in\N$. Furthermore, define the $\F_{\mathscr U}$-measurable random variables $n_1 := \max\{n:t_n<T\}$ and $n_2:=\max\{n:t_n<\Theta_1+e^A/a^2\}$.

For an interval $I\subset \R_+$ and $r\in[0,a]$, we now define 
\[
U_I(r) = \# \{(u,t)\in\mathscr L^{\mathrm{red}}_{\emptyset,I}:X_u(t)\ge r\},
\]
the number of red particles created to the right of $r$ during the time interval $I$ and set $U_n(r) = U_{I_n}(r)$.
Furthermore, we denote by $N^{\mathrm{white}}_{t}(r)$ the number of white particles with positions $\ge r$ at time $t$, \emph{including} the ``in between particles''. 

\begin{lemma}
  \label{lem:Unr} Suppose (HB$_0$). Let $I=[t_l,t_r]$, with $d:=t_r-t_l \le \delta/4$. Then for every $\alpha > 0$, for sufficiently small $\delta$, there exists $C_{\delta,\alpha}$, such that for large $A$ and $a$, for every $r\in [0,9a/10]$,  if $Ka^2 \le t_l < \tconst{eq:tcrit}$,
  \begin{eqnarray}
    \label{eq:Unr}
\E[U_I(r)\,|\,T>t_l] \le C_{\delta,\alpha}A^2\ep\Big(\frac {t_l} {a^3} +\eta\Big)e^{-(\frac 4 3-\alpha) r}N,\\
 \label{eq:tau}
\P(\exists t\in I: N^{\mathrm{white}}_{t}(r)\ge N\,|\,T>t_l) \le C_{\delta,\alpha}A^2\ep\Big(\frac {t_l} {a^3} +\eta\Big)e^{-(\frac 4 3-\alpha) r},\\
\label{eq:Unr0}
\E[U_I(r)\Ind_{(R_{[t_l-\zeta,t_r]}=0)}\,|\,T>t_l] \le C_{\delta,\alpha}dA^2\ep\Big(\frac {t_l} {a^3} +\eta\Big)e^{-(\frac 4 3-\alpha) r}N.
 \end{eqnarray}
Likewise, conditioned on $\F_\Delta$, we have for $T\le t_l<\Theta_1+e^Aa^2$,
\begin{eqnarray}
    \label{eq:Unr1}
\E_{\mathrm{nbab}}[U_I(r)\,|\,\F_\Delta]\Ind_{G_\Delta} \le C_{\delta,\alpha}A^2\ep^2e^{-(\frac 4 3-\alpha) r}N,\\
 \label{eq:tau1}
\P_{\mathrm{nbab}}(\exists t\in I: N^{\mathrm{white}}_{t}(r)\ge N\,|\,\F_\Delta) \Ind_{G_\Delta} \le C_{\delta,\alpha}A^2\ep^2e^{-(\frac 4 3-\alpha) r}.
\end{eqnarray}
Finally, suppose (HB$^\flat_\perp$). Then
\begin{equation}
  \label{eq:44}
\E[U_{[0,Ka^2]}(r)\Ind_{(R_{Ka^2} = 0)}] \le C_{\delta,\alpha}ae^{-(\frac 4 3-\alpha) r}N.
\end{equation}
\end{lemma}
\begin{proof}
Define the stopping time $\tau = \inf\{t\in I: N^{\mathrm{white}}_{t}(r)\ge N\}$, with $\inf
\emptyset = +\infty$. Define further the law $\Phat = \P(\cdot\,|\,T>t_l)$ and  the event $E = \{N_{t_l}({r'}) < e^{-\delta/2-r'/3}N\}$, where  $r' = r - \sqrt r \vee 0$. Then $e^{-4r'/3} \le C_\alpha e^{-r(4/3-\alpha)}$ for every $\alpha > 0$, such that by Lemma~\ref{lem:Bflat_N}, we have 
\begin{equation}
  \label{eq:18}
\begin{split}
\Phat(\tau < \infty)   &\le \Phat(\tau < \infty,\, E)  + \Phat(E^c) \\
 &\le  \Phat(\tau < \infty\,|\,E)+  C_{\delta,\alpha} A^2\ep\Big(\frac {t_l} {a^3}+\eta\Big)e^{-(\frac 4 3 -\alpha) r}.
\end{split}
\end{equation}
For $t\in I$, we can bound $N^{\mathrm{white}}_{t}(r)$ by the sum of
\begin{enumerate}[nolistsep]
\item $N^{[1]}_{t}$: the number of descendants at time $t$ of \{the particles in $\mathscr N_{t_l}$ which are to the \emph{right} of $r'$\},
\item $N^{[2]}$: the number of descendants of \{the particles in $\mathscr N_{t_l}$ which are to the \emph{left} of $r'$\} and which reach the point $r$ before the time $t_r$ and
\item $N^{[3]}_{t}$: the number of ``in between'' particles at time $t$.
\end{enumerate}
Conditioned on $\F_{t_l}$ and the event $E$, $N^{[1]}_{t}$ is stochastically 
bounded by the number of individuals at time $t-t_l$ in a Galton--Watson process with
reproduction law $q$ and branching rate $\beta_0$, starting from $ \lfloor e^{-\delta/2-r'/3}N\rfloor$
individuals. Let $(\mathrm{GW}_t)_{t\ge 0}$ be such a process, then $e^{-\beta_0 mt}\mathrm{GW}_t = e^{-t/2}\mathrm{GW}_t$ is a
martingale and by Doob's $L^2$-inequality, we get
\begin{equation}
  \label{eq:87}
E[\sup_{t\in[0,d]}(e^{-t/2}\mathrm{GW}_t-\mathrm{GW}_0)^2] \le 4 \Var(e^{-d/2}\mathrm{GW}_d) \le C d Ne^{-r'/3}.
\end{equation}
By \eqref{eq:87}, Chebychev's inequality and the hypothesis on $d$, we get,
\begin{equation}
\label{eq:26}
\begin{split}
\Phat(\sup_{t\in I} N^{[1]}_{t} \ge Ne^{-\delta/4}) &\le P(\sup_{t\in[0,d]} \mathrm{GW}_t \ge Ne^{-\delta/4})\\
&\le P(\sup_{t\in[0,d]}e^{-t/2}\mathrm{GW}_t-\mathrm{GW}_0\ge Ne^{-3\delta/8}(1-e^{-\delta/8}))\\
&\le C_\delta d e^{-r'/3}/N.
\end{split}
\end{equation}
Now, if $r < 1$, we have $N^{[2]} = 0$, so suppose that $r\ge 1$. Let the random variable $\Gamma_r$ denote the number of particles in BBM which reach $r$ before the time $d$. Since there are at most $N$ white particles to the left of $r'$ at time $t_l$, we have by Lemma~\ref{lem:many_to_one},
\begin{equation}
  \label{eq:17}
  \Ehat[N^{[2]}\,|\,\F_{t_l}] \le N \sup_{x\in [0,r']} \E^x[\Gamma_r] \le C N W(\sup_{t\in[0,d]} B_t \ge \sqrt r) \le CNd^{-1/2} e^{-r/(2d)},
\end{equation}
where $(B_t)$ follows standard Brownian motion under $W$. For each $(u,t)\in\mathscr N_{t_l}$, let $N^{[2]}_u$ be the contribution of the descendants of $u$ to $N^{[2]}$. By Lemma~\ref{lem:many_to_two_fixedtime},
\begin{equation}
  \label{eq:34}
  \begin{split}
  \Var[N^{[2]}\,|\,\F_{t_l}] &\le \sum_{(u,t)\in\mathscr N_{t_l}} \E[(N_u^{[2]})^2] \le N \sup_{x\in[0,r']}\E^x[\Gamma_r^2] \\
&\le N \Big(\sup_{x\in[0,r']}\E^x[\Gamma_r] + C\int_0^{d} \int_0^a \p_t(x,z)\E^{(z,t)}[\Gamma_r]^2\,\dd z\,\dd t\Big).
  \end{split}
\end{equation}
Trivially, $\E^x[\Gamma_r] \le e^{d/2}$ and $\int_0^a \p_t(x,z)\,\dd z \le e^{d/2}$ for any $x\in\R$, $t\le d$. Together with  \eqref{eq:17} and \eqref{eq:34}, this yields,
\begin{equation}
  \label{eq:76}
  \Var[N^{[2]}\,|\,\F_{t_l}] \le CN \sup_{x\in[0,r']}\E^x[\Gamma_r] \le CNd^{-1/2} e^{-r/(2d)}.
\end{equation}
Chebychev's inequality and the inequality $t^{-1/2}e^{-1/(4t)} \le Ct$ now yields
\begin{equation}
  \label{eq:78}
  \P(N^{[2]} \ge N(1-e^{-\delta/8})\,|\,\F_{t_l}) \le C_\delta d e^{-r/(4d)}/N.
\end{equation}
As for $N^{[3]}$, we have  $\Ehat[R_{t_l}^2] = \Varhat(R_{t_2}) + \Ehat[R_{t_l}]^2 \le Ce^{2A}$ by Lemma~\ref{lem:Zt_variance} as well as Lemmas~\ref{lem:Rt}, \ref{lem:Zl_exp_upperbound} and Corollary~\ref{cor:hat_many_to_one}, such that by the definition of a breakout event,
\begin{equation*}
\Ehat[\sup_{t\in I}(N^{[3]}_{t})^2] \le Ce^{2A}\zeta^2.
\end{equation*}
Markov's inequality then gives
\begin{equation}
  \label{eq:88}
\Phat\big(\sup_{t\in I}N^{[3]}_{t} \ge (e^{-\delta/4}-e^{-\delta/8})N\big) \le C_\delta e^{2A}\zeta^2/N^2.
\end{equation}

Equations \eqref{eq:18}, \eqref{eq:26}, \eqref{eq:78} and \eqref{eq:88} now yield \eqref{eq:tau}. As for \eqref{eq:Unr} and \eqref{eq:Unr0}, we first note that the expected number of individuals created until time $t$ in a Galton--Watson process with reproduction law $(q(k))_{k\ge0}$ is bounded by the number of individuals at time $t$ of a Galton--Watson process with reproduction law $(1-q(0))^{-1}q(k)_{k\ge 1}$, whose expectation is bounded by $Ct$. This gives, 
\begin{equation}
  \label{eq:133}
\Ehat[U_{(\tau,t_r]}\,|\,\tau < \infty] \le C(dN + \zeta e^A).
\end{equation}
because there are at most $N$ white particles at the time $t_l$ and the total number of in-between particles is bounded by $\zeta \E[R_{t_l}] \le C\zeta e^A$ by Lemmas~\ref{lem:Rt} and \ref{lem:Zl_exp_upperbound} and Corollary~\ref{cor:hat_many_to_one}. 
As for the overshoot, i.e.\ the particles that have been coloured red at the time $\tau$, we first note that 
\begin{equation}
  \label{eq:80}
\Ehat[U_{\{\tau\}} \Ind_{(\tau<\infty)}] \le \Ehat[U_{\{\tau\}}\Ind_{(\tau < \infty,\,U_{\{\tau\}} \ge N)}] + N \Phat(\tau < \infty).
\end{equation}
Now, 
\begin{equation}
  \label{eq:81}
U_{\{\tau\}}\Ind_{(\tau < \infty,\,U_{\{\tau\}} \ge N)} \le \max\{k_u: u\in U, d_u \in I,\,k_u \ge N,\, X_u(d_u-) \ge r\}
\end{equation}
By the definition of a breakout event, the set on the right-hand side of \eqref{eq:81} contains no ``in between'' particles for $N>\zeta$. Furthermore, the expected number of branching events of ``regular'' particles to the right of $r$ is bounded by 
\begin{equation}
  \label{eq:24}
C\int_{t_l}^{t_r}\Ehat[N_{t}(r)]\,\dd t \le Cd N(1+r^2)e^{- r},
\end{equation}
where the inequality follows from \eqref{eq:22}. This yields,
\begin{equation}
  \label{eq:82}
 \Ehat[U_{\{\tau\}}\Ind_{(\tau < \infty,\,U_{\{\tau\}} \ge N)}]  \le C d N (1+r^2) e^{- r} E[L\Ind_{(L \ge N)}] \le C d (1+r^2) e^{- r},
\end{equation}
where the last inequality follows from the Cauchy--Schwarz and Chebychev inqualities. 
By \eqref{eq:tau}, \eqref{eq:133}, \eqref{eq:80} and \eqref{eq:82} we finally get for large $A$ and $a$,
\begin{equation}
  \label{eq:20}
\Ehat[U_I(r)] = \Ehat[(U_{\{\tau\}}(r)+U_{(\tau,t_r]}(r))\Ind_{(\tau <\infty)}] \le C_{\delta,\alpha} N A^2\ep\Big(\frac {t_l} {a^3} +\eta\Big)e^{-(\frac 4 3-\alpha) r}.
\end{equation}
For the  proof of \eqref{eq:Unr0}, we note that on the set $R_{[t_l-\zeta,t_r]}=0$ there are no in-between particles during $I$. Taking away the corresponding terms in the above proof yields \eqref{eq:Unr0}.

The proof of the other three equations is very similar. The proof of \eqref{eq:Unr1} and \eqref{eq:tau1} uses the second half of Lemma~\ref{lem:Bflat_N} instead of the first and \eqref{eq:70} in addition to \eqref{eq:22} for the proof of \eqref{eq:24}. The proof of the last equation draws on Lemma~\ref{lem:Bflat_start} instead of Lemma~\ref{lem:Bflat_N} and requires covering the interval $[0,Ka^2]$ by pieces of length $d$.
\end{proof}



\begin{corollary}
  \label{cor:Ga2}
Suppose (HB$_0$). For an interval $I\subset \R_+$, let $G_{a/2,I}$ be the event that no particle is coloured red during the time interval $I$ and to the right of $a/2$. Then,
\[
\P(G_{\mathscr U}\backslash G_{a/2,[Ka^2,T]}) \le C_\delta o(1)\quad\tand\quad \P_{\mathrm{nbab}}\left((G_{a/2,[T,\Theta_1+e^Aa^2]})^c\,|\,\F_\Delta\right)\Ind_{G_\Delta}\le C_\delta o(1).
\]
\end{corollary}
\begin{proof}
Follows immediately from \eqref{eq:tau} and \eqref{eq:tau1} by summing over the intervals $I_n$ (note that $T\le \sqrt\ep a^3$ on $G_{\mathscr U}$).
\end{proof}


\begin{lemma} Suppose (HB$_0$). Set $T'=T\wedge \sqrt\ep a^3$ and define the intervals $I = [Ka^2,T']$ and $J = [T,\Theta_1+e^Aa^2]$ For large $A$ and $a$,
\label{lem:ZY_red}
\begin{align*}
&\E[(Z^{\mathrm{red}}_{\emptyset,I} + Y^{\mathrm{red}}_{\emptyset,I})\Ind_{G_{a/2,I}}] \le C_{\delta} A^2\ep^3 e^A\\ 
&\E_{\mathrm{nbab}}[(Z^{\mathrm{red}}_{\emptyset,J} + Y^{\mathrm{red}}_{\emptyset,J})\Ind_{G_{a/2,J}}\,|\,\F_\Delta]\Ind_{G_\Delta} \le C_{\delta} e^{2A}/a. 
\end{align*}
Now suppose (HB$^\flat_\perp$). Then
\begin{equation*}
\E[(Z^{\mathrm{red}}_{\emptyset,[0,Ka^2]} + Y^{\mathrm{red}}_{\emptyset,[0,Ka^2]})\Ind_{(R_{Ka^2}=0)}] \le C_{\delta}e^A/a^2.
\end{equation*}
\end{lemma}
\begin{proof}
By Lemma~\ref{lem:Unr}, we have for every $r\in [0,a/2]$ and $n\le n_0:=\lfloor(\sqrt\ep a^3-Ka^2)/(t_1-t_0)\rfloor$,
\begin{equation}
  \label{eq:37}
\begin{split}
  \E[U_I(r)] &\le \sum_{n=0}^{n_0}\E[U_n(r)\,|\,T>t_n]\P(T>t_n)\\
 &\le C_{\delta}A^2\ep e^{-\frac 8 7 r}N \sum_{n=0}^{n_0}\Big(\frac {t_n} {a^3} +\eta\Big)\P(T>t_n)\\
&\le C_\delta A^2 \ep e^{-\frac 8 7 r}e^A \E\Big[\Big(\frac {T\wedge \sqrt\ep a^3} {a^3}+\eta\Big)^2 \Big] \le C_\delta A^2 \ep^3 e^{-\frac 8 7 r},
\end{split}
\end{equation}
by Lemma~\ref{lem:moments_T} and \eqref{eq:eta}. 
Now, by integration by parts, we have for every (possibly random) interval $I$,
\begin{equation}
  \label{eq:71}
 Z^{\mathrm{red}}_{\emptyset,I} = \int_0^a w_Z'(r) U_I(r)\,\dd r\quad\tand\quad Y^{\mathrm{red}}_{\emptyset,I} = U_I(0) + \int_0^a w_Y'(r) U_I(r)\,\dd r.
\end{equation}

The first equation in the statement of the lemma now follows from \eqref{eq:37} and \eqref{eq:71}, since $w_Z'(x)\le C(1+x)e^{x-a}$ and $w_Y'(x) \le Ce^{x-a}$ for $x\in[0,a]$. The remaining equations follow similarly.
\end{proof}

The following lemma will only be needed in the proof of Proposition~\ref{prop:Bflat_med}.
\begin{lemma}
\label{lem:Nred}
Suppose (HB$^\flat_0$). Let $r\in[0,9a/10]$ and $(K+1)a^2\le t\le \tconst{eq:tcrit}$. Then, for every $\alpha > 0$, for large $A$ and $a$,
  \begin{equation*}
  \P[N^{\mathrm{red},{(0)}}_{t}(r)>\alpha N,\,T>t] \le C_\delta \Big(\alpha^{-1}(1+ \mu r)e^{-\mu r}A^2\ep \Big(\frac t {a^3} +\eta\Big)+ \ep^2\Big).
  \end{equation*}
\end{lemma}
\begin{proof}
For an interval $I$, denote by $N^{\mathrm{red},{(0)}}_t(r,I)$ the number of red particles to the right of $r$ at time $t$ which have turned red during $I$. Define the event 
\[G = G_{a/2,[0,t]}\cap \{R_{[t-a^2-\zeta,t]}=0\}\cap \{Z^{\mathrm{red}}_{\emptyset,[0,Ka^2]} + Y^{\mathrm{red}}_{\emptyset,[0,Ka^2]} \le a^{-1/2}\},\]
and note that $\P(G^c,\,T>t) \le C\ep^2$ for large $A$ and $a$, by Corollary~\ref{cor:Ga2}, Lemma~\ref{lem:hat_quantities} and the hypothesis.

Write $\Ehat = \E[\cdot\,|\,T>t]$. If $\Ehat^{(x,s)}[N^{(0)}_{t}(r)]\le f(x,s)$ for some function $f$ with $f(0,s)\equiv f(a,s)\equiv 0$, then by integration by parts,
  \begin{align*}
  \Ehat[N^{\mathrm{red},{(0)}}_t(r,I)\Ind_{G}] 
\le\int_I \int_0^{a/2} \Ehat[U_{\dd s}(x)\Ind_{(R_{[t-\zeta,t]}=0)}]  \frac \dd {\dd x}f(x,s)\,\dd x\,\dd s.
  \end{align*}
Define $I_1 = [0,t-a^2]$ and $I_2 = [t-a^2,t]$. By \eqref{eq:N_expec_large_t}, \eqref{eq:250} and Corollary~\ref{cor:hat_many_to_one}, we have for $x\in[0,a/2]$,
\begin{equation*}
  \Ehat^{(x,s)}[N^{(0)}_{t}(r)] \le \begin{cases} C(1+\mu r)e^{-\mu r}N e^{-A}w_Z(x), &\tif s\in I_1,\\ Ce^{\mu x}\int_r^\infty e^{-\mu z} p^a_{t-s}(x,z),& \tif s\in I_2, \end{cases}
\end{equation*}
Similarly to the proof of Lemma~\ref{lem:ZY_red}, we now have by the inequality $w_Z'(x)\le C(1+x)e^{x-a}$, for large $A$ and $a$,
\begin{align*}
\Ehat[N^{\mathrm{red},{(0)}}_t(r,I_1)\Ind_{G}]  &\le Ce^{-A}N(1+\mu r)e^{-\mu r} \int_0^{a/2}\Ehat[U_{I_1}(x)\Ind_{G}](1+x)e^{x-a}\,\dd x\\
&\le C_\delta A^2\ep ((t/a^{3})+\eta)(t/a^{3})(1+\mu r)e^{-\mu r}N,
\end{align*}
by Lemma~\ref{lem:Unr} and the definition of $G$. As for the particles created during $I_2$, we have by Lemma~\ref{lem:Unr}, the definition of $G$ and the hypothesis on $t$, for all $s\in I_2$,
\begin{equation*}
  \Ehat[U_{\dd s}(x)\Ind_{G}] \le C_{\delta}A^2\ep\Big(\frac t {a^3} +\eta\Big)Ne^{-\frac 5 4 x}\dd s.
\end{equation*}
Now, with \eqref{eq:p_sin}, one easily sees that for large $a$, $(\dd/\dd x) (e^{\mu x}p_s^a(x,z)) \ge 0$ for all $x\le a/2$ and $s\ge a^2$, whence, by \eqref{eq:p_potential} and Fubini's theorem, we get
\begin{align*}
 \Ehat[N^{\mathrm{red},{(0)}}_t(r,I_2)\Ind_{G}] &\le C \int_0^{a^2} \int_0^{a/2} \Ehat[U_{\dd s}(x)\Ind_{G}]\int_r^\infty e^{-\mu z} \frac \dd {\dd x} \left(e^{\mu x}p^a_s(x,z)\right)\,\dd z\,\dd x\,\dd s\\
&\le  C \int_0^{a/2} \Ehat[U_{\dd s}(x)\Ind_{G}]\int_r^\infty e^{-\mu z} \frac \dd {\dd x}\big(e^{\mu x}a^{-1}(x\wedge z)(a-x\vee z)\big)\,\dd z\,\dd x\\
& \le C_{\delta}A^2\ep\Big(\frac t {a^3} +\eta\Big)N\int_r^\infty \int_0^{a/2}  e^{-\mu (x/4+ z)} (1+x)\,\dd x\,\dd z\\
&\le C_\delta A^2\ep ((t/a^{3})+\eta) e^{-\mu r}N.
\end{align*}
In total, we have
\begin{equation*}
 \Ehat[N^{\mathrm{red},{(0)}}_t(r)\Ind_{G}] \le  C_{\delta}A^2\ep\Big(\frac t {a^3} +\eta\Big)(1+\mu r)e^{-\mu r}N.
\end{equation*}
An application of Markov's inequality finishes the proof.
\end{proof}

\begin{lemma} Suppose (HB$^\flat_0$). Let $G_{\mathrm{fred}}$ be the event that the fugitive does not get coloured red. Then, for large $A$ and $a$,
  \label{lem:fred}
\[
\P(G_{\mathscr U}\cap G_{a/2,[Ka^2,T]}\backslash G_{\mathrm{fred}}) \le C_\delta \ep^2.
\]
\end{lemma}
\begin{proof}
Set $G_{a/2} = G_{a/2,[Ka^2,T]}$. Let $\P_0$ be the law of BBM with absorption at zero from Section~\ref{sec:before_breakout} (i.e.\ we do not move the barrier when a breakout occurs). Let $T_{\mathrm{red}}$ be the first breakout of a red particle. Then
\begin{equation}
  \label{eq:92}
  \P(G_{\mathscr U}\cap G_{a/2}\backslash G_{\mathrm{fred}}) = \P_0(T = T_{\mathrm{red}},\,\G_{\mathscr U}\cap G_{a/2}) \le \P_0(T_{\mathrm{red}} \le \sqrt\ep a^3,\,G_{a/2}).
\end{equation} 
Define $T' = T\wedge \sqrt\ep a^3$. Then, since $\mathscr L^{\mathrm{red}}_{\emptyset,T'}$ is a stopping line, we have by Proposition~\ref{prop:T} and Corollary~\ref{cor:Ga2},
\begin{equation}
\label{eq:93}
\P_0(T_{\mathrm{red}} \le \sqrt\ep a^3\,|\,\F_{\mathscr L^{\mathrm{red}}_{\emptyset,T'}})\Ind_{G_{a/2}} \le C \ep^{-1/2}e^{-A} (Z^{\mathrm{red}}_{\emptyset,T'}+Y^{\mathrm{red}}_{\emptyset,T'}).
\end{equation}
Recall that $\P(Z^{\mathrm{red}}_{\emptyset,[0,Ka^2]} + Y^{\mathrm{red}}_{\emptyset,[0,Ka^2]} \le a^{-1/2}) \ge 1-\ep^2$ by hypothesis.
By the tower property of conditional expectation and \eqref{eq:93}, this gives
\begin{align*}
  \P_0(T_{\mathrm{red}} \le \sqrt\ep a^3,\,G_{a/2}) &\le \E[\P_0(T_{\mathrm{red}} \le \sqrt\ep a^3\,|\,\F_{\mathscr L^{\mathrm{red}}_{\emptyset,T'}})\Ind_{G_{a/2}}]+ \ep^2 + o(1)\\
&\le C\E[C \ep^{-1/2}e^{-A} (Z^{\mathrm{red}}_{\emptyset,[Ka^2,T']}+Y^{\mathrm{red}}_{\emptyset,[Ka^2,T']})\Ind_{G_{a/2}}] + C\ep^2.
\end{align*}
The lemma now follows from \eqref{eq:92} together with Lemma~\ref{lem:ZY_red},  \eqref{eq:ep_upper} and the hypothesis.
\end{proof}

\begin{lemma}
\label{lem:G_Bflat}
Suppose (HB$^\flat_0$). There exists a numerical constant $\numcst > 0$, such that for large $A$ and $a$,
\[
\P(G^{\flat}_1)\ge 1- \ep^{1+\numcst}.
\]
\end{lemma}
\begin{proof}
Recall the event $G_1$ from Section~\ref{sec:BBBM} and change its definition slightly by requiring that $|e^{-A}Z_{\Theta_1}-1| \le \ep^{3/2}/2$, call the new event $G_1'$. Define the random variable 
\[
X = \P(Z^{\mathrm{red}}_{\emptyset,[\Theta_1,\Theta_1+Ka^2]} + Y^{\mathrm{red}}_{\emptyset,[\Theta_1,\Theta_1+Ka^2]} \le a^{-1/2}\,|\,\F_{\Theta_1}).
\]
By Proposition~\ref{prop:piece} and Remark~\ref{rem:piece}, it suffices to show that
\begin{equation*}
  \P\left(G_1'\backslash \left\{Z^{\mathrm{red}}_{\Theta_1} \le \ep^{3/2}/2,\,X>1-\ep^2\right\}\right) \le \ep^{1+\numcst},
\end{equation*}
for some numerical constant $\numcst$.

Let $T'$, $I$ and $J$ as in Lemma~\ref{lem:ZY_red}. We then have by Markov's inequality, Lemma~\ref{lem:ZY_red} and Corollary~\ref{cor:Ga2},
\begin{equation}
  \label{eq:35}
  \E_{\mathrm{nbab}}[1-X\,|\,\F_\Delta]\Ind_{G_\Delta} \le C_\delta o(1).
\end{equation}
Markov's inequality applied to \eqref{eq:35} then yields
\begin{equation}
  \label{eq:52}
  \P(X>1-\ep^2,\ G_\Delta\cap G_{\mathrm{nbab}}) \le C_\delta o(1).
\end{equation}
Now define $G_{\mathrm{red}} = \{T\le\sqrt\ep a^3,\,Z^{\mathrm{red}}_{\emptyset,[0,Ka^2]} + Y^{\mathrm{red}}_{\emptyset,[0,Ka^2]} \le a^{-1/2}\}\cap G_{a/2,I}\cap G_{\mathrm{fred}}$. Conditioned on $T$ and $\F_{\mathscr L^{\mathrm{red}}_{\emptyset,T}}$, on the set $G_{\mathrm{red}}$, the particles from the stopping line $\mathscr L^{\mathrm{red}}_{\emptyset,T}$ then all spawn BBM conditioned not to break out before $T$ (because the neither the fugitive nor any in-between particles are on the stopping line). By Lemma~\ref{lem:ZY_red}, Lemma~\ref{lem:hat_quantities} and the tower property of conditional expectation, we then have
\begin{align*}
  \E[Z^{\mathrm{red}}_{T}\Ind_{G_{\mathrm{red}}}] &\le C\Big(Aa^{-1/2}+\E[(Z^{\mathrm{red}}_{\emptyset,I} + AY^{\mathrm{red}}_{\emptyset,I})\Ind_{G_{a/2,I}}]\Big) \le C_\delta (A\ep)^3e^A\\
  \E[Y^{\mathrm{red}}_{T}\Ind_{G_{\mathrm{red}}}] &\le C\Big(Aa^{-1/2}+\E[(Z^{\mathrm{red}}_{\emptyset,I} + Y^{\mathrm{red}}_{\emptyset,I})\Ind_{G_{a/2,I}}]\Big) \le C_\delta A^2\ep^3e^A.
\end{align*}
By Lemma~\ref{lem:hat_quantities}, we now have
\begin{align*}
\E_{\mathrm{nbab}}[Z^{\mathrm{red}}_{\Theta_1}\,|\,\F_\Delta\vee \F_T]\Ind_{G_\Delta} \le C_\delta \Big(Z^{\mathrm{red}}_T + Z^{\mathrm{red}}_{\emptyset,J}+ A(Y^{\mathrm{red}}_T + Y^{\mathrm{red}}_{\emptyset,J})\Big),
\end{align*}
and Lemma~\ref{lem:ZY_red} and the tower property of conditional expectation give
\begin{equation}
 \label{eq:96}
  \E[Z^{\mathrm{red}}_{\Theta_1}\Ind_{G_{\Delta}\cap G_{\mathrm{nbab}} \cap G_{\mathrm{red}}}] \le C_\delta (A\ep)^3 e^A.
\end{equation}
Furthermore, by the hypothesis, Corollary~\ref{cor:Ga2} and Lemmas~\ref{lem:GDelta}, \ref{lem:Gnbab} and \ref{lem:fred}, we have
\begin{equation}
\label{eq:979}
\P(G_\Delta\cap G_{\mathrm{nbab}}\cap G_{\mathrm{red}}) \ge 1-\ep^{1+\numcst}. 
\end{equation}
The lemma now follows from \eqref{eq:52}, \eqref{eq:979} and Markov's inequality applied to \eqref{eq:96}, together with \eqref{eq:ep_upper}.
\end{proof}

\subsection{Proofs of the main results}

\begin{proof}[Proof of Lemma~\ref{lem:HBflat}]
By Lemma~\ref{lem:Bflat_start}, it remains to estimate the probability of the last event in the definition of $G^\flat_0$. Define the random variable $X = \P(Z^{\mathrm{red}}_{\emptyset,[0,Ka^2]} + Y^{\mathrm{red}}_{\emptyset,[0,Ka^2]} \le a^{-1/2}\,|\,\F_0)$. By Lemma~\ref{lem:Bflat_start} and Markov's inequality applied to the last equation of Lemma~\ref{lem:ZY_red}, we have $\E[(1-X)\Ind_{(R_{Ka^2}=0)}] \le C_\delta e^A/a$. By Markov's inequality, we then have 
\[
\P(X\le 1-\ep^2,\ R_{Ka^2}=0) \le \ep^{-2}\E[(1-X)\Ind_{(R_{Ka^2}=0)} ] \le C_\delta o(1).
\]
Together with \eqref{eq:RKa2}, this finishes the proof.
\end{proof}

\begin{proof}[Proof of Proposition~\ref{prop:Bflat}]
Let $n\ge 0$. Conditioned on $\F_{\Theta_n}$,  the barrier process until $\Theta_{n+1}$ is by definition the same in B$^\flat$-BBM and B-BBM. Furthermore, $G_n\subset G^\flat_n$ and $G^\flat_n\in\F_{\Theta_n}$. The first statement then follows by induction from Lemma~\ref{lem:G_Bflat}, as in the beginning of the proof of Proposition~\ref{prop:piece}.

As for the second statement, inspection of the proofs of Theorems~\ref{th:barrier} and \ref{th:barrier2} shows that they only rely on Proposition~\ref{prop:piece} and on the existence of the coupling with a Poisson process constructed in Section~\ref{sec:BBBM_proofs}. But this construction only relied on the law of $T_{n+1}$ conditioned on $\F_{\Theta_n}$ and $G_n$ and thus readily transfers to the \Bfl-BBM.
\end{proof}

\begin{proof}[Proof of Proposition~\ref{prop:Bflat_med}]
Define $\rho = x_{\alpha e^{2\delta}}$, such that $(1+\mu\rho)e^{-\mu\rho} \ge \alpha e^{2\delta}(1+o(1))$. Recall from Section~\ref{sec:BBBM_proofs} the coupling of the process $(\Theta_n/a^3)_{n\ge 0}$ with a Poisson process $(V_n)_{n\ge 0}$ of intensity $\gamma_0^{-1} =  p_B e^A \pi$. As indicated above, this coupling can be constructed for the \Bfl-BBM as well. In particular, setting  $\widetilde G^\flat_n = G^\flat_n \cap \{\forall j\le n: |(T_j-\Theta_{j-1})/a^3 - V_j| \le \ep^{3/2}\}$, we have as in \eqref{eq:prob_Gnprime}, for every $n\ge 0$,
\begin{equation}
  \label{eq:64}
\P(\widetilde G^\flat_n) \ge \P(G^\flat_n)-nO(\ep^{2}) \ge \P(G^\flat_0) - nO(\ep^{1+\numcst}),
\end{equation}
by Proposition~\ref{prop:Bflat} (where $\numcst>0$ is a numerical constant).
Now let $(t_i^{A,a})$ satisfy (Ht). For simplicity, suppose that $t_i^{A,a}\equiv t_i$ for all $i$, the proof of the general case is exactly the same. We want to show that $\P(\forall 1\le i\le k: N^{\mathrm{white}}_{t_ia^3}(\rho) \ge \alpha N) \to 1$ as $A$ and $a$ go to infinity. By \eqref{eq:954}, the probability that there exists $(i,j)$, such that $t_i$ and $t_{i+1}$ are in the same interval $[\Theta_j,\Theta_{j+1}]$ is bounded by $O(\sqrt\ep)$ for large $A$ and $a$. We can therefore suppose that $k=1$, the general case is a straightforward extension.

It thus remains to show that for every $t>0$, $\P(N^{\mathrm{white}}_{ta^3}(\rho) \ge \alpha N) \to 1$ as $A$ and $a$ go to infinity (the case $t= 0$ is a straightforward calculation). Let $t>0$ and $n := \lceil \gamma_0^{-1}(t+1)\rceil$. Then $\E[V_n] \ge t+1$ and $\Var(V_n) \le 2\gamma_0(t+1)$, such that $\P(V_n < t+1/2) \le C\gamma_0 = O(\ep)$ (for this proof, we allow the constants $C$ and the expression $O(\cdot)$ depend on $t$). Defining the event
\[
E_n(s) = \{\Theta_n> s\} \cap \{\nexists j\in\{0,\ldots,n\}: s\in [T_{j-1},\Theta_{j}+e^Aa^2]\} \cap \widetilde G^\flat_n,
\]
(define $T_{-1} = 0$) we then have by \eqref{eq:954} and \eqref{eq:64}, $\P(E_n) \ge 1-O(\ep^{1+\numcst})$. 

We now prove by induction that there exists a numerical constant $\numcst > 0$, such that for every $1\le j\le n$, 
\begin{equation*}
S_j = \sup_{\nu\in G^\flat_0} \sup_{s\in[0,ta^3]} \P^{\nu}(N^{\mathrm{white}}_{s}(\rho) < \alpha N,\,E_j(s)) \le C_{\delta,\alpha} j\ep^{1+\numcst},
\end{equation*}
 as $A$ and $a$ go to infinity (we accept the abuse of notation $\nu\in G^\flat_0$). By definition, $S_0 = 0$. 
Now, suppose the statement is true for some $j\ge 0$. We then have for every $\nu\in G^\flat_0$ and $s\in[0,ta^3]$,
\begin{multline*}
  \P^\nu(N^{\mathrm{white}}_{s}(\rho) < \alpha N,\,E_{j+1}(s))  \le \P^\nu(N^{\mathrm{white}}_{s}(\rho) < \alpha N,\,E_{j+1}(s),\,s > \Theta_1)\\
 + \P^\nu(N^{\mathrm{white}}_{s}(\rho) < \alpha N,\,E_{j+1}(s),\,s \le \Theta_1)=: P_1+P_2.
\end{multline*}
We first have, by the definition of $E_{j+1}(s)$ and the fact that the process starts afresh at the stopping time $\Theta_1$,
\begin{align*}
  P_1 &\le \E^\nu[\P(N^{\mathrm{white}}_{s}(\rho) < \alpha N,\,E_{j+1}(s)\,|\,\F_{\Theta_1})\Ind_{G^\flat_1\cap\{\Theta_1\le s\}}] \\
&\le \E^\nu[\sup_{\nu_1\in G^\flat_0}\P^{\nu_1}(N^{\mathrm{white}}_{s-\Theta_1}(\rho) < \alpha N,\,E_j(s-\Theta_1))\Ind_{(\Theta_1\le s)}] 
\le S_j.
\end{align*}
Second, we have again by the definition of $E_{j+1}(s)$,
\begin{equation}
\label{eq:Nalpha}
  P_2 \le \P^\nu(N^{\mathrm{white}}_{s}(\rho) < \alpha N,\,T> s)\Ind_{(s\ge e^Aa^2)}.
\end{equation}
Write $\Phat^\nu = \P^\nu(\cdot\,|\,T>s)$.
For $e^Aa^2\le s\le \tconst{eq:tcrit}$, we now have by Proposition~\ref{prop:N_expec} and Corollary~\ref{cor:hat_many_to_one} for the first inequality and Lemma~\ref{lem:N_2ndmoment} and Corollary \ref{cor:hat_many_to_one} for the second,
\begin{align*}
  \label{eq:130}
  &\Ehat^\nu[N^{(0)}_s(\rho)] \ge (1-O(\ep^{3/2}))e^{2\delta} \alpha e^{-\delta} N \ge (1-O(\ep^{3/2}))(1+\delta) \alpha,\\
  &\Varhat^\nu(N^{(0)}_s(\rho)) \le C_\alpha \alpha^2 e^{-2A} N^2,
\end{align*}
such that 
\begin{align*}
  \P^\nu(&N^{\mathrm{white},{(0)}}_{s}(\rho) < \alpha N,\,T>s) \\
&\le \Phat^\nu(N^{(0)}_{s}(\rho) < (1+\delta/2)\alpha N) +  \P^\nu(N^{\mathrm{red},{(0)}}_{s}(\rho) \ge (\delta/2)\alpha N,\,T>s)\\
&\le C_{\delta,\alpha} \left(e^{-2A} + A^2\ep(s/a^3+\eta)\right).
\end{align*}
by Chebychev's inequality applied to the first term and Lemma~\ref{lem:Nred} to the second.
Altogether, this gives,
\begin{align*}
  P_2 \le C_{\delta,\alpha} \Big(\P(T>\sqrt\ep a^3) + e^{-2A} +  A^2\ep^{3/2}\Big) \le C_{\delta,\alpha} A^2\ep^{3/2},
\end{align*}
for large $A$ and $a$.
Finally, we get $S_{j+1} \le S_j + C_{\delta,\alpha} A^2\ep^{3/2}$ for all $0\le j\le n$, which yields $S_j \le C_{\delta,\alpha} j\ep^{5/4}$ for large $A$ and $a$, by \eqref{eq:ep_upper}. This gives
\[
\P(N^{\mathrm{white}}_{ta^3}(\rho) \le \alpha N) \le S_n + \P(V_n < t+1/2) + (1-\P(E_n(s)))\to 0,
\]
as $A$ and $a$ go to infinity. The statement follows.
%
\end{proof}

\section{The \texorpdfstring{\Bsh}{B\#}-BBM}
\label{sec:Bsharp}

In this section, we define and study the \Bsh-BBM, a model which will be used in Section~\ref{sec:NBBM} to bound the $N$-BBM from \emph{above}.

\subsection{Definition of the model}
As in the previous section, we will use throughout the notation from Section~\ref{sec:BBBM} and furthermore fix $\delta \in (0,1/2)$ and define $K$ to be the smallest number, such that $K \le 1$ and $E_K \le \delta/10$. We now define however  $N = \lfloor 2\pi  e^{A-\delta}a^{-3}e^{\mu a} \rfloor $. Again, we will use interchangeably the phrases ``as $A$ and $a$ go to infinity'' and ``as $N$ go to infinity'' and $A_0$ and the function $a_0(A)$ in the definition of these phrases may now also depend on $\delta$. The symbols $C_\delta$ and $C_{\delta,\alpha}$ have the same meaning as in the last section.

\begin{figure}[h]
 \centering
\def\svgwidth{10cm}
\executeiffilenewer{sec9.svg}{sec9.pdf}%
{inkscape -z -D --file=sec9.svg %
--export-pdf=sec9.pdf --export-latex --export-area-drawing}%
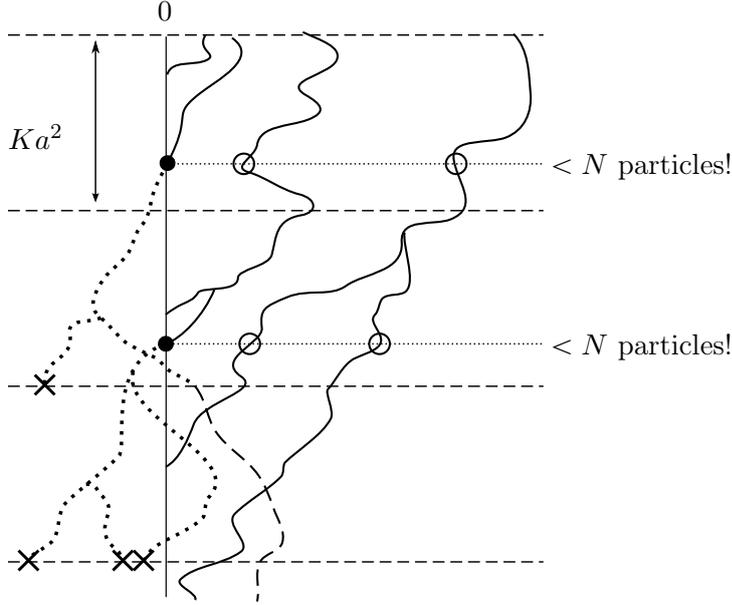%

 \caption{The \Bsh-BBM with parameter $N=3$ (no breakout is shown). White, blue and grey particles are drawn with solid, dotted and dashed lines, respectively. A cross indicates that a blue particle is killed.}
 \label{fig:sec9}
\end{figure}

The B$^\sharp$-BBM is then defined as follows: Given a possibly random initial configuration $\nu_0$ of particles in $(0,a)$, we let particles evolve according to B-BBM with barrier function given by \eqref{eq:wall_fn} and with the following changes: Define $t_n = n(K+3)a^2$ and $I_n = [t_n,t_{n+1})$ (note that these definitions differ from those of Section~\ref{sec:Bflat}). Colour all initial particles white. When a white particle hits $0$ during the time interval $I_n$ and has at least $N$ particles to its right, it is killed immediately. If less than $N$ particles are to its right, it is coloured blue and survives until the time $t_{n+2}\wedge \Theta_1$, where all of its descendants \emph{to the left of $0$} are killed and the remaining survive and are coloured white again. At the time $\Theta_1$, the process starts afresh. 
See Figure~\ref{fig:sec9} for a graphical description.

For bookkeeping, we add a shade of grey to the white particles which have hit 0 at least once (and call them hence the grey particles). We then add the superscripts ``nw'', ``gr'', ``blue'' or ``tot'' to the quantities referring respectively to the non-white, grey, blue or all the particles. Quantities without this superscript refer to the white particles.

In particular, we define $B_I$ and $B^{\mathrm{tot}}_I$ to be the number of white, respectively, white and grey particles touching the left barrier during the time interval $I$ with less than $N$ particles to their right (i.e.\ those which are coloured blue). We set $B_n= B_{I_n}$ and $B^{\mathrm{tot}}_n = B^{\mathrm{tot}}_{I_n}$.

Let $\nu^\sharp_t$ be the configuration of (all) the particles at the time $t$ and abuse notation by setting $\nu^\sharp_n =\nu^\sharp_{\Theta_n}$. We set $G^{\sharp}_{-1}=\Omega$ and for each $n\in\N$, we define the event $G^\sharp_n$ to be the intersection of $G^\sharp_{n-1}$ with the following events (we omit the braces).
\begin{itemize}
 \item $\supp\nu^\sharp_n\subset (0,a)$,
 \item $\mathscr N^{\mathrm{tot}}_{\Theta_n} \subset U\times \{\Theta_n\}$ and $\Theta_n > T_n^+$ (for $n>0$),
 \item $|e^{-A}Z^{\mathrm{tot}}_{\Theta_n} -1| \le \ep^{3/2}$ and $Y^{\mathrm{tot}}_{\Theta_n} \le \eta$,
 \item $N(\Theta_n) \ge N$ and for all $j\ge 0$ with $\Theta_{n-1}+t_j<\Theta_n$: $N( \Theta_{n-1}+t_j) \ge N$,
 \item $\P\Big(B^{\mathrm{tot}}_{[\Theta_n,\Theta_n+t_1]} \le e^{-A} e^{\mu a}/a\,\Big|\,\F_{\Theta_n}\Big) \ge 1-\ep^2$.
\end{itemize}
The last event is of course uniquely defined up to a set of probability zero. Note that $G^\sharp_n\in\F_{\Theta_n}$ for each $n\in\N$.
Furthermore, we define the predicates
\begin{itemize}
 \item[(HB$^\sharp$)] (The law of) $\nu_0$ is such that $\P(G^\sharp_0)\to 1$ as $A$ and $a$ go to infinity.
 \item[(HB$^\sharp_0$)] $\nu_0$ is deterministic and such that $G^\sharp_0$ holds.
 \item[(HB$^\sharp_{\perp}$)] $\nu_0$ is obtained from $\lceil e^{\delta}N\rceil$ particles distributed independently according to the density proportional to $\sin(\pi x/a)e^{-\mu x}\Ind_{(0,a)}(x)$.
\end{itemize}

We now state the important results on the \Bsh-BBM.
\begin{lemma}
\label{lem:HBsharp}
(HB$^\sharp_{\perp}$) implies (HB$^\sharp$) for large $A$ and $a$.
\end{lemma}

\begin{proposition}
  \label{prop:Bsharp}
Proposition~\ref{prop:piece} still holds for the B$^\sharp$-BBM, with $G_n$ replaced by $G_n^\sharp$. The same is true for Theorems~\ref{th:barrier} and \ref{th:barrier2}, with (HB) replaced by (HB$^\sharp$).
\end{proposition}

Recall the definition of $\med^N_\alpha$ and $x_\alpha$ from the introduction and of (Ht) from Section~\ref{sec:BBBM}.
\begin{proposition}
\label{prop:Bsharp_med}
Suppose (HB$^\sharp_\perp$). Let $(t^{A,a}_j)$ satisfy (Ht). Let $\alpha\in (0,1)$. As $A$ and $a$ go to infinity,
\[
\P(\forall j: \med^N_\alpha(\nu^{\sharp}_{a^3t^{A,a}_j}) \le x_{\alpha e^{-2\delta}})\to 1.
\]
\end{proposition}

\begin{lemma}
  \label{lem:Csharp}
Define a variant called C$^\sharp$-BBM of the B$^\sharp$-BBM by killing blue particles only if there are at least $N$ particles to their right. Then Propositions~\ref{prop:Bsharp} and \ref{prop:Bsharp_med} hold for the C$^\sharp$-BBM as well.
\end{lemma}

\subsection{Preparatory lemmas}
\label{sec:Bsharp-preparatory-lemmas}

We first derive upper bounds on the probability that the number of white particles is less than $N$ at a given time $t$. We do not impose any particular condition on the initial configuration..

\begin{lemma}
\label{lem:Nt_prob_less_N}
Let $\Phat = \P(\cdot\,|\,T>s)$ with $s\le \tconst{eq:tcrit}$ and let $Ka^2\le t\le s$. We have for small $\delta$ and for $A$ and $a$ large enough,
\[
 \Phat(N^{(0)}_{t} < N\,|\,\F_0)\Ind_{(Z_0 \ge (1-\delta/2)e^A)} \le C_\delta e^{-A} \Big(\frac K a + e^{-A}Y_0\Big).
\]
\end{lemma}
\begin{proof}
By \eqref{eq:N_expec_large_t}, Corollary~\ref{cor:hat_many_to_one} and the definition of $K$, we have for  $A$ and $a$ large enough,
\begin{equation*}
\Ehat[N^{(0)}_{t}\,|\,\F_0] \ge (1-Cp_B)\E[N^{(0)}_{t}\,|\,\F_0] \ge (1+\delta/4) N e^{-A}Z_0.
\end{equation*}
By the conditional Chebychev inequality, we then have 
\begin{equation*}
\Phat(N^{(0)}_{t} < N\,|\,\F_0)\Ind_{(Z_0\ge (1-\delta/2)e^A)} \le C \frac{\Varhat(N^{(0)}_{t}\,|\,\F_0)}{(e^{-A} N\delta Z_0)^2}.
\end{equation*}
The lemma now follows from Lemma~\ref{lem:N_2ndmoment} and Corollary~\ref{cor:hat_many_to_one}
\end{proof}
\begin{corollary}
\label{cor:Nt_prob_less_N}
  Under the conditions of Lemma~\ref{lem:Nt_prob_less_N}, suppose furthermore that $t\ge (K+1)a^2$. Then, for small $\delta$ and for $A$ and $a$ large enough,
\[
 \Phat(N^{(0)}_t < N,\,Z^{(0)}_{t-Ka^2} \ge (1-\delta/2)e^A\,|\,\F_0) \le C_\delta e^{-A}(1+e^{-A}Z_0)/a.
\]
\end{corollary}
\begin{proof}
 First condition on $\F_{t-Ka^2}$ and apply Lemma~\ref{lem:Nt_prob_less_N}. Then condition on $\F_0$ and apply \eqref{eq:Yt} and Corollary~\ref{cor:hat_many_to_one}.
\end{proof}


For $s\ge 0$, define the event
\[
G_{Z,s} = \left\{\sup_{0\le t\le s}\left|e^{-A}Z^{(0)}_t - 1\right| < \delta/4\right\},
\]
and set $G_{Z,n} = G_{Z,t_n}$.
In using the last two results, the following lemma will be crucial:
\begin{lemma}
\label{lem:G_Z}
Suppose $|e^{-A}Z_0-1|\le \delta/8$. Let $s\le \tconst{eq:tcrit}$ and set $\Phat = \P(\cdot\,|\,T>s)$. Then $\Phat(G_{Z,s}^c) \le C_\delta e^{-A}$  for large $A$ and $a$.
\end{lemma}
\begin{proof}
Let $h(x,t)=\P^{(x,t)}(T>s)$ and note that $h(x,t) \ge 1-Cp_B$ by \eqref{eq:hbar}. Define $Z^h_t = \sum_{u\in\widehat{\mathscr N}^{(0)}_t} (h(X_u(t),t))^{-1}w_Z(X_u(t))$, which is a supermartingale under $\Phat$ by Lemma~\ref{lem:many_to_one_conditioned}, with $\Ehat[Z^h_s] = (1+O(p_B))Z^h_0$ by Corollary~\ref{cor:hat_many_to_one}. For large $A$ and $a$, we then have by Doob's $L^2$ inequality,
\begin{equation*}
\Phat(G_{Z,s}^c) \le \Phat(\sup_{0\le t\le s}\big|Z^h_t - \E[Z^h_s]\big| \ge e^A\delta/16)\le C_\delta e^{-2A}\Varhat(Z^h_s).
\end{equation*}
The lemma now follows from the previous equation, Proposition~\ref{prop:quantities} and Corollary~\ref{cor:hat_many_to_one}.
\end{proof}

The following lemma is the analogue of Lemma~\ref{lem:Nt_prob_less_N} for the system after the breakout. 
\begin{lemma}
  \label{lem:Nt_prob_less_N_moving_wall}
Let $t\in[T,\Theta_1]$ and let $(u,t_u)\in \widehat{\mathscr N}_{T^-}\cup \mathscr S^{(\mathscr U,T)}$. Then, for large $A$ and $a$,
\[
 \P_{\mathrm{nbab}}(\widehat N^{\neg u}_t + N^{\mathrm{fug},\neg u}_t<N\,|\,\F_\Delta)\Ind_{G_\Delta} \le C_{\delta} \eta e^{-A}/\ep.
\]
\end{lemma}
\begin{proof}
As in the proof of Lemma~\ref{lem:Bflat_N}, set $M_t = e^{- X^{[1]}_{t}} = (1+\Delta_t)^{-1}$, where $\Delta_t =  \thbar((t-T^+)/a^2)\Delta$. Recall that on $G_\Delta$:  $|\Delta - e^{-A}Z^{(\mathscr
     U,T)}| \le \ep^{1/8}$. For a particle $(u,t_u)\in \widehat{\mathscr N}_{T^-}\cup \mathscr S^{(\mathscr U,T)}$, we define $\widehat N'^{\neg u}_t$ and $N'^{\mathrm{fug},\neg u}_t$ to be the number of hat-, respectively, fug-particles which are not descendants of $u$ and which have not hit $a$ after the time $T^-$.

Now, by Lemma~\ref{lem:mu_t_density} and Proposition~\ref{prop:N_expec}, we have for large $A$ and $a$,
\[
 \E[\widehat N'^{\neg u}_t\,|\,\F_\Delta]\Ind_{G_\Delta} \ge N e^{-A}M_t(\widehat Z_{T^-}-\|w_Z\|_\infty)(1-\delta/8)\Ind_{G_\Delta}\ge N M_t(1+\delta/2)\Ind_{G_\Delta}.
\]
Moreover, by Lemma~\ref{lem:critical_line}, we have for large $A$ and $a$,
\[
\begin{split}
  \E[N'^{\mathrm{fug},\neg u}_t\,|\,\F_\Delta]\Ind_{{G_\Delta}} &=
  Ne^{-A}M_t Z^{(\mathscr U,T)}(\thbar((t-T)/a^2) + O(y/a))\Ind_{{G_\Delta}}\\
 &\ge
  NM_t(\Delta_t(1+\delta) -\delta/8)\Ind_{{G_\Delta}}.
\end{split}
\]
In total, this gives for large $A$ and $a$,
\begin{equation}
  \label{eq:7}
  \E[\widehat N'^{\neg u}_t + N'^{\mathrm{fug},\neg u}_t\,|\,\F_\Delta]\Ind_{{G_\Delta}} \ge N(1+\delta/4)\Ind_{{G_\Delta}}.
\end{equation}
Moreover, by Lemma~\ref{lem:N_2ndmoment}, we have for large $A$ and $a$,
\begin{multline}
  \label{eq:3}
  \Var(\widehat N'^{\neg u}_t+N'^{\mathrm{fug},\neg u}_t\,|\,\F_\Delta) \\
\le C(Ne^{-A})^2 ((\widehat Z_{T^-}+Z^{(\mathscr U,T)})(t-T^-)/a^3 +(\widehat Y^{(0)}_{T^-}+Y^{(\mathscr U,T)}))\le CN^2 \eta e^{-A}/\ep,
\end{multline}
by the definition of $G_\Delta$. Lemma~\ref{lem:Gnbab}, \eqref{eq:7} and \eqref{eq:3} and the conditional Chebychev
inequality now yield the lemma.
\end{proof}

For a barrier function $f$, let $\P_f^0$ denote the law of BBM with drift $-\mu_s$ defined in \eqref{eq:mu_t} starting from a single particle at $0$. Let $R_t$ be the number of particles hitting $a$ before the time $t$ and define $Z_t$, $Y_t$ by summing $w_Z$, respectively $w_Y$ over the particles at time $t$.
\begin{lemma}
\label{lem:BBM_0}
Let $t \ge 0$ and $f$ be a barrier function. Then,
\begin{equation*}
\E^0_f[Z_t] \le a e^{\frac{\pi^2}{2a^2}t - \mu a},\quad
\E^0_f[Y_t] \le e^{\frac{\pi^2}{2a^2}t - \mu a}\quad\tand\quad
\E^0_f[R_t] \le e^{\frac{\pi^2}{2a^2}t - \mu a}.
\end{equation*}
\end{lemma}
\begin{proof}
By Lemma~\ref{lem:many_to_one_simple} and Girsanov's theorem, we have
\[
 \E^0_f[Y_t] = e^{\beta_0 mt}e^{-\mu a}W^0_{-\mu_s}\Big[e^{\mu X_t}\Ind_{(X_t\ge 0)}\Big] \le  e^{\beta_0 mt}e^{-\mu a}W^0_{-\mu}\Big[e^{\mu (X_t-f(t/a^2))}\Big] \le e^{\frac{\pi^2}{2a^2}t - \mu a},
\]
which implies the first two inequalities, since $Z_t \le a Y_t$. As for the third one, by Lemma~\ref{lem:many_to_one} and again Girsanov's theorem,
\[
 \E^0[R_t] = W^0_{-\mu}\Big[e^{\beta_0 mH_a}\Ind_{(H_a \le t)}\Big] \le e^{-\mu a} W^0\Big[e^{\frac{\pi^2}{2a^2}H_a}\Ind_{(H_a \le t)}\Big] \le e^{\frac{\pi^2}{2a^2}t - \mu a}.
\]
Furthermore, we have $\E^0_f[R_t]\le  \E^0[R_t]$ as in Lemma~\ref{lem:R_f}. This finishes the proof of the lemma. 
\end{proof}

\subsection{The probability of \texorpdfstring{$G^\sharp_1$}{G\#1}}
\label{sec:blue-particles}
For an interval $I$, we define the random variable $L_I$ counting the number of particles hitting the origin during the time interval $I$:
\[
L_I = \sum_{(u,s)\in\mathscr L_{H_0}}\Ind_{(s\in I)},
\]
and set $L_n = L_{I_n}$.

\begin{lemma}
\label{lem:Bsharp_left_border}
  Let $n\ge0$, $s\le \tconst{eq:tcrit}$, $I = [t_l,t_r] \subset [0,s]$ and $x\in[0,a]$. Furthermore, let $f$ be a barrier function with $\Err(f,t_r)\le C$. Then, for large $A$ and $a$, we have
\[
\E^{x}_f[L_I\,|\,T>s] \le \begin{cases}
C e^{\mu a}a^{-3}(t_r-t_l)(Aw_Y(x) + w_Z(x)), &\tif t_l \ge a^2 \tor x\ge a/2\\
C e^{\mu a} (w_Y(x) +a^{-3}(t_r-t_l)w_Z(x)), &\text{ for any $t,x$.}
\end{cases}
\]
\end{lemma}
\begin{proof}
Write $\Ehat[\cdot] = \E[\cdot\,|\,T>s]$. Define $L^{(0)}_I = \sum_{(u,s)\in\mathscr L_{H_0}}\Ind_{(H_a(X_u) > s\in I)}$. As in the proof of Lemma~\ref{lem:mu_t_density}, we have for every interval $I$, $x\in[0,a]$ and $t\le s$,
\begin{equation}
\label{eq:left_border}
\Ehat^{x}_f[L^{(0)}_I \Big] \le e^{O(\Err(f,t))}\Ehat^{x}\Big[L^{(0)}_I ] \le C e^{\mu x} I^a(a-x,I),
\end{equation}
by \eqref{eq:Ia} and Corollary~\ref{cor:hat_many_to_one}. Lemma~\ref{lem:I_estimates} and \eqref{eq:left_border} now give,
\begin{equation}
  \label{eq:99}
\begin{split}
\Ehat^{x}_f[L_I] &= \Ehat^{x}_f[L^{(0)}_I ] + \Ehat^{x}_f\Big[\sum_{(u,s)\in  \mathscr S^{(1+)}_{t_r}} \Ehat^{(X_u(s),s)}_f[L^{(0)}_I] \Big]\\
&\le Ce^{\mu a}\left(w_Y(x) I^a(a-x,I) + \frac K a \Ehat^{x}_f[Z^{(1+)}_{\emptyset,[0,t_r]}]\right).
\end{split}
\end{equation}
Lemma~\ref{lem:I_estimates} now gives for large $a$,
\begin{equation}
  \label{eq:100}
I^a(a-x,I) \le C\left(a^{-3}(t_r-t_l)\sin(\pi x/a)+\Ind_{(t_l< a^2\tand x\le a/2)}\right),
\end{equation}
and by Lemmas~\ref{lem:R_moments}, \ref{lem:R_f} and \ref{lem:Zl_exp_upperbound}, Corollary~\ref{cor:hat_many_to_one} and \eqref{eq:QaZ}, we have 
\begin{equation}
  \label{eq:98}
  \Ehat^{x}_f[Z^{(1+)}_{\emptyset,[0,t_r]}] \le CAw_Y(x)\left(1+\frac{t_r} {a^2} \sin(\pi x/a)\right).
\end{equation}
The lemma now follows from \eqref{eq:99}, \eqref{eq:100} and \eqref{eq:98}, together with the hypothesis $t_r \le s \le \tconst{eq:tcrit} \le Ca^3/A$. 
\end{proof}

The following lemma is crucial. It will permit to estimate the number of particles turning blue upon hitting the origin.

\begin{lemma}
  \label{lem:blue_particles}
Suppose that $\nu_0$ is deterministic with $Z_0 \ge (1-\delta/4)e^A$.
Let $s\le\tconst{eq:tcrit}$ and $I=[t_l,t_r] \subset [Ka^2,s]$. Write $\Phat=\P(\cdot\,|\,T>s).$ Then, for large $A$ and $a$,
\begin{equation*}
  \Ehat\Big[\sum_{(u,t)\in\mathscr L_{H_0}} \Ind_{(t\in I,\,N^{(0)}_t < N)}\Big] \le C_\delta e^{-A}(a^{-1} + e^{-A}Y_0)\Ehat[L_I].
\end{equation*}
\end{lemma}
\begin{proof}
For an individual $u\in U$ and $t\ge 0$, define $\widehat N^{(0),\neg u}_t = \sum_{(v,s)\in \mathscr N^{(0)}_t}\Ind_{(v\nsucceq u_0)}$, where $u_0$ is the ancestor of $u$ at the time 0. By the trivial inequality $N^{(0)}_t \ge \widehat N^{(0),\neg u}_t$ for every $u$ and by the independence of the initial particles,
  \begin{equation*}
  \Ehat\Big[\sum_{(u,t)\in\mathscr L_{H_0}} \Ind_{(t\in I,\,N^{(0)}_t < N)}\Big] \le \Ehat\Big[\sum_{(u,t)\in\mathscr L_{H_0}} \Ind_{(t\in I)} \P(\widehat N^{(0),\neg u}_t < N)\Big].
 \end{equation*}
The lemma now follows from Lemma~\ref{lem:Nt_prob_less_N}, since for every $u\in\mathscr N(0)$, we have $Z_0 - w_Z(X_u(0)) \ge Z_0 - C \ge (1-\delta/2)e^A$ for large $A$, by hypothesis.
\end{proof}

\paragraph{Before the breakout.}

For every $n\ge0$, define the (sub-probability) measure $\P_n = \P(\cdot,\ T>t_{n+1})$. The following lemma gives an estimate on the number of \emph{white} particles (not counting the grey ones) turning blue during the interval $I_n$. 

\begin{figure}[ht]
 \centering
\executeiffilenewer{sec93.svg}{sec93.pdf}%
{inkscape -z -D --file=sec93.svg %
--export-pdf=sec93.pdf --export-latex --export-area-drawing}%
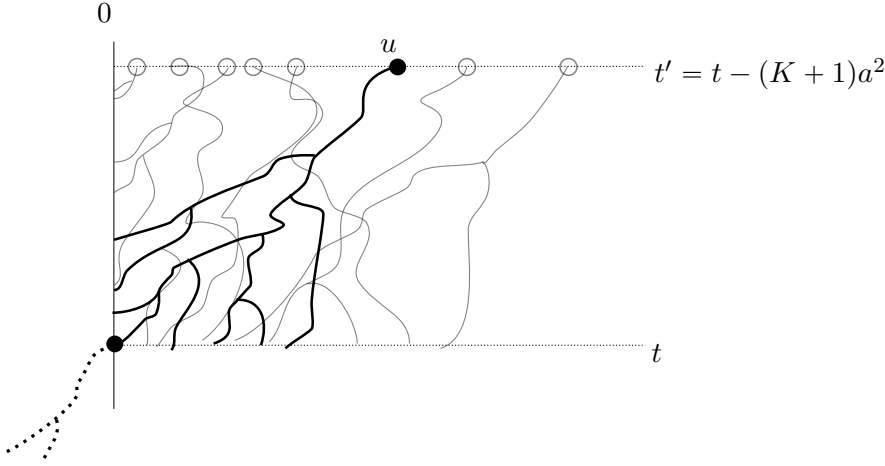%

 \caption{The key step in the proof of Lemma~\ref{lem:B_hat}: Conditioned on $\F_{t'}$, the descendants of $u$ are independent of the descendants of the other particles by the strong branching property. Furthermore, conditioned on $\F_{t'}$, the probability that at time $t$ there are less than $N$ particles not descending from $u$ is bounded by $C(a^{-1} + Y^{(0),\neg u}_{t'})$. These facts imply that the expected proportion of white particles turning blue upon hitting $a$ is bounded by $C_\delta e^{-A}/a$.}
\end{figure}

\begin{lemma}
\label{lem:B_hat}
For every $n\ge 1$ with $t_{n+1}\le \tconst{eq:tcrit}$, we have for large $A$ and $a$,
 \[
\E_n[B_n\Ind_{G_{Z,n}}] \le C_\delta \frac{e^{\mu a}}{a^2}.
 \]
\end{lemma}
\begin{proof}
Let $B_n^{(0)}$ be the number of particles turning blue during $I_n$ and which have not hit $a$ before $t'=t_n-(K+1)a^2 = t_{n-1}+2a^2$ and let $B_n^{(1+)} = B_n-B_n^{(0)}$. 
Lemmas~\ref{lem:Bsharp_left_border} and \ref{lem:blue_particles} now yield,
\begin{equation*}
  \begin{split}
 \E_n[B^{(0)}_n\Ind_{G_{Z,n}}\,|\,\F_{t'}] &\le  C_\delta e^{ - A} \Big(\frac 1 a + e^{-A}Y_{t'}^{(0)}\Big)\Ind_{G_{Z,t'}} \sum_{(u,s)\in\mathscr N^{(0)}_{t'}} \E_n^{(X_u(s),s)}[L_n(\emptyset,s)]\\
&\le C_\delta e^{\mu a - A} \Big(\frac 1 a Y_{t'}^{(0)} + e^{-A}(Y_{t'}^{(0)})^2\Big)\Ind_{G_{Z,t'}},
  \end{split}
\end{equation*}
Proposition~\ref{prop:quantities} now gives
\begin{equation}
  \label{eq:132}
  \E_n[B^{(0)}_n\Ind_{G_{Z,n}}] \le C_\delta e^{\mu a} \Big(a^{-2} + e^{-2A}a^{-1}\E[Y_{t_{n-1}+a^2}^{(0)}\Ind_{G_{Z,t_{n-1}}}]\Big) \le C_\delta e^{\mu a} /a^2.
\end{equation}
As for the remaining particles, by the strong branching property,
\begin{equation*}
\begin{split}
   \E_n[B^{(1+)}_n\Ind_{G_{Z,n}}\,|\,\mathscr L_{H_a\wedge t'}] &= \sum_{(u,s)\in\mathscr R^{(0)}_{t'}}\int_{I_n} \E_n^{(a,s)}[L_{[r,r+\dd r)}]\P_n(N^{(0)}_r< N,\,G_{Z,n}\,|\,\mathscr N_{t'})\\
&\le C_\delta (e^{-A}/a) R^{(0)}_{t'} (e^{\mu a}/a) \Qhat^a[Z],
\end{split}
\end{equation*}
by Corollary~\ref{cor:Nt_prob_less_N} and  Lemma~\ref{lem:Bsharp_left_border}. Lemma~\ref{lem:Rt} and \eqref{eq:QhatZ} then give
\begin{equation}
  \label{eq:134}
  \E_n[B^{(1+)}_n\Ind_{G_{Z,n}}] \le C_\delta (e^{\mu a}/a^2) A t_{n+1}/a^3 \le C_\delta e^{\mu a}/a^2,
\end{equation}
by the hypothesis on $t_{n+1}$. The lemma now follows from \eqref{eq:132} and \eqref{eq:134}.
\end{proof}

Up to now, we have only considered the white particles turning blue and will now turn to $B_n^{\mathrm{tot}}$, $n\ge0$.  For $n\ge 0$ and $0\le k \le n-2$, we define $B_{k,n}$ to be the number of particles that turn blue at a time $t\in I_n$, have an ancestor that turned blue at a time $t'\in I_k$, have none that has hit 0 between $t_{k+2}$ and $t_n$ and have never hit $a$ between $t'$ and $t$.

Let $G_{\mathrm{b,n}}$ the event that no blue particle hits $a$ before $t_n$ and $G_{\mathrm{g},n}$ the event that among the descendants of the particles counted by $B_k$, $k=0,\ldots,n-2$, no grey particle breaks out before $\tconst{eq:tcrit}$ (i.e.\ that no particle that turned blue before $t_{n-1}$ has a grey descendent that breaks out before $\tconst{eq:tcrit}$). Then set 
\[
G_n^{\mathrm{tot}} = G_{Z,n}\cap G_{\mathrm{b,n}} \cap G_{\mathrm{g},n}\cap 
\{B_0 \le e^Ae^{\mu a}/a^2,\,Y_0\le \eta\}.
\]

\begin{lemma}
 \label{lem:B_kn}
For every $n\ge 0$ with $t_{n+1}\le \tconst{eq:tcrit}$ and $k\le n-2$, we have for large $A$ and $a$,
\[
 \E_n[B_{k,n}\Ind_{G_n^{\mathrm{tot}}}] \le C_{\delta} \E_k[B^{\mathrm{tot}}_k\Ind_{G_k^{\mathrm{tot}}}] \frac{e^{-A}}{a}.
\]
\end{lemma}
\begin{proof}
Let $\mathscr B_k$ be the stopping line consisting of the particles that turn blue during $I_k$, at the moment at which they turn blue (hence, $B^{\mathrm{tot}}_k = \#\mathscr B_k$). Note that $\mathscr B_k\wedge \mathscr N^{(0)}_{t_{n-1}} = \mathscr B_k \cup \mathscr N^{(0)}_{t_{n-1}}$, such that the descendants of $\mathscr B_k$ and of $\mathscr N^{(0)}_{t_{n-1}}$ are independent, given their past, by the strong branching property. Note also that $G_k^{\mathrm{tot}}\in\F_{\mathscr B_k\wedge \mathscr N^{(0)}_{t_{n-1}}}$, since by definition this $\sigma$-field contains all the information about the descendents of the particles in $\mathscr B_j$, $j<k$. By Corollary~\ref{cor:Nt_prob_less_N}, we then have
\begin{multline}
  \label{eq:103}
 \E_n[B_{k,n}\Ind_{G_n^{\mathrm{tot}}}\,|\,\F_{\mathscr B_k\wedge \mathscr N^{(0)}_{t_{n-1}}}] \\
\le C_{\delta} \frac{e^{-A}}{a}\Ind_{G_k^{\mathrm{tot}}}\Ind_{(T>t_{k+1})}\sum_{(u,s)\in\mathscr B_k}\E_n^{(0,s)}\Big[\sum_{(v,t)\in\mathscr L_{H_0^{k+2}}}\Ind_{(H_a(X_v) > t_{k+2},\,X_v(t_{k+2})\in (0,a),\,t\in I_n)}\Big],
\end{multline}
where here $H_0^{k+2}(X) = \inf\{t\ge t_{k+2}: X_t = 0\}$. By Lemmas~\ref{lem:BBM_0} and \ref{lem:Bsharp_left_border}, each summand in the right-hand side of the above inequality is bounded by $C_\delta e^{\mu a} \E_n^{0}[Y_{t_{k+2}-s}] \le C_\delta$.
The lemma follows.
\end{proof}
\begin{lemma}
\label{lem:B_n}
For all $n\ge 1$ with $t_{n+1}\le \tconst{eq:tcrit}$, we have for large $A$ and $a$,
\begin{equation}
\label{eq:B_n}
\E_n[B^{\mathrm{tot}}_n\Ind_{G_n^{\mathrm{tot}}}] \le C_\delta \frac{e^{\mu a}}{a^2}.
\end{equation}
\end{lemma}
\begin{proof}
By Lemma~\ref{lem:B_hat}, the statement is true for $n=1$, because $B^{\mathrm{tot}}_1 = B_1$ by definition. Now we have for every $n$, by Lemmas~\ref{lem:B_hat} and \ref{lem:B_kn} and the fact that $B^{\mathrm{tot}}_0 \le e^{-A}e^{\mu a}/a$ on $G_n^{\mathrm{tot}}$, 
\begin{equation}
 \label{eq:840}
 \begin{split}
\E_n[B^{\mathrm{tot}}_n\Ind_{G_n^{\mathrm{tot}}}] &= \sum_{k=0}^{n-2} \E_n[B_{k,n}\Ind_{G_n^{\mathrm{tot}}}] + \E_n[B_n\Ind_{G_n^{\mathrm{tot}}}]\\
 &\le C_\delta\Big(\frac{e^{-A}}{a}\sum_{k=1}^{n-2} \E_k[B^{\mathrm{tot}}_k\Ind_{G_k^{\mathrm{tot}}}] + \frac{e^{\mu a}}{a^2}\Big).
 \end{split}
\end{equation}
The lemma now follows easily by induction over $n$, since $n\le a$ by hypothesis.
\end{proof}

\begin{lemma}
\label{lem:G_n_tot}
 Suppose (HB$^\sharp_0$). For large $A$ and $a$, we have $\P_n((G_n^{\mathrm{tot}})^c) \le C_\delta e^{-2A/3}(n/a)+ \ep^2$ for every $n\ge 0$.
\end{lemma}
\begin{proof}
Let $n\ge 1$. By definition, the event $G_{n-1}^{\mathrm{tot}}\backslash G_{\mathrm{b},n}$ implies that a descendant of a particle in $\mathscr B_{n-2}$ hits $a$ before $t_n$. Markov's inequality and Lemma~\ref{lem:BBM_0} then imply,
\begin{equation}
  \label{eq:90}
  \P_n(G_{n-1}^{\mathrm{tot}}\backslash G_{\mathrm{b},n}) \le \E_n[B^{\mathrm{tot}}_{n-2}\Ind_{G_{n-1}^{\mathrm{tot}}}] \Big(\sup_{t\le t_2} \E^0[R_t]\Big)\le C_\delta a^{-2},
\end{equation}
by Lemma~\ref{lem:B_n}.
Furthermore, by Proposition~\ref{prop:T} and Lemmas~\ref{lem:BBM_0} and \ref{lem:B_n}, we have
\begin{equation}
  \label{eq:106}
  \P_n((G_{n-1}^{\mathrm{tot}}\cap G_{\mathrm{b},n}) \backslash G_{\mathrm{g},n}) \le \E_n[B^{\mathrm{tot}}_{n-2}\Ind_{G_{n-1}^{\mathrm{tot}}}]\Big(\sup_{t\le t_2}p_B\E^0[ \tconst{eq:tcrit}Z_t+Y_0]\Big) \le C_\delta e^{-2A/3} a^{-1},
\end{equation}
by \eqref{eq:pB} and \eqref{eq:ep_lower}. By \eqref{eq:90} and \eqref{eq:106}, we have for large $A$ and $a$,
\[
\P_n(G_{Z,n}\backslash G_n^{\mathrm{tot}}) \le \P_n(G_{Z,n}\backslash G_{n-1}^{\mathrm{tot}}) + C_\delta e^{-2A/3}a^{-1} \le \P_{n-1}(G_{Z,n-1}\backslash G_{n-1}^{\mathrm{tot}}) + C_\delta e^{-2A/3}a^{-1}.
\]
By induction over $n$ and the hypothesis, this gives for every $n$,
\[
\P_n(G_{Z,n}\backslash G_n^{\mathrm{tot}}) \le C_\delta e^{-2A/3}(n/a)+ \ep^2.
\]
Together with Lemma~\ref{lem:G_Z}, the statement follows.
\end{proof}

Define the random variable $n_0:= \lfloor T/t_1\rfloor $ and note that $t_{n_0-1}>T^-$ for large $A$. Define the events $G_{\mathrm{nw}} = \{Z^{\mathrm{gr}}_{t_{n_0}}\le e^{A/2},\,Y^{\mathrm{gr}}_{t_{n_0}} \le e^{A/2}/a,\,B^{\mathrm{gr}}_{n_0-1}\le \eta e^{\mu a}/a\}$ and $G_{n_0}^{\mathrm{tot}} = \bigcup_n (G_n^{\mathrm{tot}}\cap \{n_0=n\})$.

\begin{lemma}
  \label{lem:G_ex} We have for large $A$ and $a$,
\[
\P((G_{n_0}^{\mathrm{tot}}\cap G_{\mathscr U})\backslash G_{\mathrm{nw}}) \le C_\delta e^{-A/2}.
\]
\end{lemma}
\begin{proof}
By definition, $\F_{\mathscr U}\subset \F_{\mathscr L_{H_0}}$. Lemmas~\ref{lem:hat_quantities}, \ref{lem:BBM_0} and \ref{lem:B_n} then give
\begin{equation}
  \label{eq:2}
  \begin{split}
  \E[Z^{\mathrm{gr}}_{t_{n_0}} \Ind_{G_{n_0}^{\mathrm{tot}}\cap G_{\mathscr U}}]&\le C\sup_{t\le t_2}\E^0[Z_t + AY_t] \Big(\sum_{n=0}^\infty \E\Big[\E_n[B^{\mathrm{tot}}_n\Ind_{G_n^{\mathrm{tot}}}\,|\,\F_{\mathscr L_{H_0}}]\Ind_{G_{\mathscr U}}\Big]\Big)\\
&\le C_\delta a e^{\mu a} \sum_{n=0}^{\lfloor \sqrt\ep a^3/t_1\rfloor} \E_n[B^{\mathrm{tot}}_n\Ind_{G_n^{\mathrm{tot}}}]\le C_\delta.
  \end{split}
\end{equation}
Similarly, we have

\begin{align}
  \label{eq:36}
 \P(R^{\mathrm{gr}}_{I_{n_0-1}} > 0,\ G_{n_0}^{\mathrm{tot}},\ G_{\mathscr U}) \le C_\delta \eta,\quad\tand \quad\E[Y^{\mathrm{gr}}_{t_{n_0}}\Ind_{(R^{\mathrm{gr}}_{I_{n_0-1}} = 0)}\Ind_{G_{n_0}^{\mathrm{tot}}\cap G_{\mathscr U}}] \le C_\delta/a.
\end{align}
Finally, we have as in \eqref{eq:2}, by Lemmas~\ref{lem:B_kn} and \ref{lem:B_n},
\begin{multline}
  \label{eq:38}
   \E[B^{\mathrm{gr}}_{n_0-1} \Ind_{G_{n_0-1}^{\mathrm{tot}}\cap G_{\mathscr U}}] = \E\Big[\sum_{k=0}^{n_0-2}\E[B_{k,{n_0-1}}\Ind_{G_{n_0-1}^{\mathrm{tot}}}\,|\,\F_{\mathscr L_{H_0}}]\Ind_{G_{\mathscr U}}\Big] \\
\le C_\delta (e^{-A}/a) \sum_{k=0}^{\lfloor \sqrt\ep a^3/t_1\rfloor} \E_k[B^{\mathrm{tot}}_k\Ind_{G_k^{\mathrm{tot}}}] \le C_\delta e^{-A+\mu a}/a^2.
\end{multline}
The lemma now follows from \eqref{eq:2}, \eqref{eq:36} and \eqref{eq:38}, together with Markov's inequality.
\end{proof}

\paragraph{After the breakout.}

We study now the system after the breakout. Recall that $n_0 = \lfloor T/t_1\rfloor$ by definition and define $n_1 := \lfloor \Theta_1/t_1 +
 1\rfloor$. 
We define $\F^\sharp_\Delta = \F_\Delta\cap \F_{\mathscr N^{\mathrm{gr}}_{t_{n_0}}\wedge \mathscr B_{n_0-1}}$ and the event $G^\sharp_\Delta = G_\Delta\cap G_{\mathrm{nw}}\in\F^\sharp_\Delta$. We denote by the superscript ``gr$\preceq$" the quantities relative to particles the descending from those that were grey before or at time $t_{n_0}$. Then let $G^\sharp_{\mathrm{nbab}}$ be the intersection of $G_{\mathrm{nbab}}$ with the event that none of these particles hits $a$ between $t_{n_0}$ and $\Theta_1+e^Aa^2$ before hitting 0. Define then the (sub-probability) measure $\P^\sharp_{\mathrm{nbab}} = \P(\cdot,\,G^\sharp_{\mathrm{nbab}})$.

  \begin{lemma}
\label{lem:Gsharp_Delta_nbab}
For large $A$ and $a$, $\P(G^\sharp_\Delta)\ge 1-C\ep^{5/4}$ and $\P((G^\sharp_{\mathrm{nbab}})^c\,|\,\F^\sharp_\Delta)\Ind_{G^\sharp_\Delta} \le C\ep^2$.
  \end{lemma}
  \begin{proof}
On $G_\Delta$, we have $t_{n_0}\le \tconst{eq:tcrit}$ for large $A$ and $a$, in particular, $n_0 \le a$. By Lemmas~\ref{lem:G_n_tot} and \ref{lem:G_ex} and \eqref{eq:ep_lower}, we then have for large $A$ and $a$,
\begin{equation*}
  \P(G_\Delta\backslash G^\sharp_\Delta) \le \P(G_\Delta \backslash G^{\mathrm{tot}}_{n_0}) + \P((G_\Delta \cap G^{\mathrm{tot}}_{n_0})\backslash G_{\mathrm{nw}}) \le C\ep^2.
\end{equation*}
With Lemma~\ref{lem:GDelta}, this yields the first statement. For the second statement, we have by Lemmas~\ref{lem:Rt} and \ref{lem:R_f},
\begin{equation*}
  \P(R^{\mathrm{gr}\preceq}_{[t_{n_0},\Theta_1+e^Aa^2]} > 0\,|\,\F^\sharp_\Delta)\Ind_{G^\sharp_\Delta} \le C((e^A/a)Z^{\mathrm{gr}}_{t_{n_0}}+Y^{\mathrm{gr}}_{t_{n_0}})\Ind_{G^\sharp_\Delta} \le Ce^{A/2}\eta.
\end{equation*}
Together with Lemma~\ref{lem:Gnbab}, this implies the lemma.
  \end{proof}

We change the notation of $B_n$ a bit: it is defined to be the number of particles hitting $0$ during $I_n$ for the first time \emph{after $t_{n_0-1}$}. In particular, we also count the descendants of the grey particles at that time.

\begin{lemma}
 For large $A$ and $a$, we have for every $n\ge n_0-1$,
\[
\E^\sharp_{\mathrm{nbab}}[B_n\,|\,\F^\sharp_\Delta]\Ind_{G^\sharp_\Delta} \le C_{\delta}e^{-5A/3}\frac{e^{\mu a}}{a}.
\]
\end{lemma}
\begin{proof}
Define $\mathscr L := \widehat{\mathscr N}_{T^-}\cup \mathscr S^{(\mathscr U,T)}$. We have for $t\in I_n$, as in Lemma~\ref{lem:blue_particles},
\begin{equation*}
  \begin{split}
  \E^\sharp_{\mathrm{nbab}}[\widehat B_n + B^{\mathrm{fug}}_n\,|\,\F^\sharp_\Delta]\Ind_{G^\sharp_\Delta} &\le C\sum_{(u,t)\in \mathscr L} \Ehat_{X^{[1]}}^{(X_u(t),t)}[L_n]\P(\widehat N'^{\neg u}_t + N'^{\mathrm{fug},\neg u}_t < N\,|\,\F_\Delta)\Ind_{G^\sharp_\Delta}\\
&\le C_\delta \eta (\ep e^{A})^{-1} \sum_{(u,t)\in \mathscr L}\Ehat_{X^{[1]}}^{(X_u(t),t)}[L_n]\Ind_{G^\sharp_\Delta}
  \end{split}
\end{equation*}
by Lemma~\ref{lem:Nt_prob_less_N_moving_wall}. By Lemma~\ref{lem:Bsharp_left_border}, \eqref{eq:ep_upper} and the definition of $\widehat G$ and $G_{\mathrm{fug}}$, this gives for large $A$ and $a$,
\begin{equation}
  \label{eq:109}
  \begin{split}
  \E^\sharp_{\mathrm{nbab}}[\widehat B_n  + B^{\mathrm{fug}}_n\,|\,\F^\sharp_\Delta]\Ind_{G^\sharp_\Delta} &\le C_\delta e^{\mu a}\eta (\ep e^Aa)^{-1}(\widehat Z_{T^-}  + A\widehat Y_{T^-} +Z^{(\mathscr U,T)}) \Ind_{G^\sharp_\Delta} \\
&\le C_\delta e^{\mu a}\eta\ep^{-2}/a.
  \end{split}
\end{equation}
Furthermore, by the independence of the check- and bar particles from the others, we have by Lemmas~\ref{lem:Nt_prob_less_N_moving_wall} and \ref{lem:Bsharp_left_border} and by the inequality $w_Y(x) \le a^{-1}e^{-(a-x)/2}$, valid for every $x\le a/2$ and $a$ large enough,
\begin{multline}
  \label{eq:110}
 \E^\sharp_{\mathrm{nbab}}[\widecheck B_n+\widebar B_n\,|\,\F^\sharp_\Delta]\Ind_{G^\sharp_\Delta}\\
 \le C_\delta e^{\mu a}\eta (\ep e^{A})^{-1}a^{-1}(\widecheck Z_{\emptyset\vee T^-}+A\widecheck Y_{\emptyset\vee T^-} + A \mathscr E_{\mathscr U}+Z^{\mathrm{gr}}_{t_{n_0}}+AY^{\mathrm{gr}}_{t_{n_0}})\Ind_{G^\sharp_\Delta}\le C_\delta e^{\mu a} \eta /a,
\end{multline}
by \eqref{eq:ep_upper} and the definition of $G^\sharp_\Delta$. 
Similarly, we have for $n\ge n_0$,
\begin{equation}
  \label{eq:39}
  \E^\sharp_{\mathrm{nbab}}[B^{\mathrm{gr}\preceq}_n\,|\,\F^\sharp_\Delta]\Ind_{G^\sharp_\Delta} \le C_\delta e^{\mu a}\eta (\ep e^{A})^{-1}Y^{\mathrm{gr}}_{t_{n_0}}\Ind_{G^\sharp_\Delta} \le C_\delta e^{\mu a} \eta /a,
\end{equation}
and $B^{\mathrm{gr}\preceq}_{n_0-1}\le \eta e^{\mu a}/a$ on the event $G^\sharp_\Delta$. 
The lemma now follows from \eqref{eq:109}, \eqref{eq:110} and \eqref{eq:39}, together with \eqref{eq:eta} and \eqref{eq:ep_lower}.
\end{proof}

For $n\ge n_0-1$, define $G'^{\mathrm{tot}}_n$ to be the event that no descendant of a particle which has been coloured blue between $t_{n_0-1}$ and  $t_n$ hits $a$ before $\Theta_1+e^Aa^2$.
\begin{lemma}
\label{lem:B_n_moving_wall}
 We have for every $n\ge n_0$ and for large $A$ and $a$,
\[
\E^\sharp_{\mathrm{nbab}}[B^{\mathrm{tot}}_n\Ind_{G'^{\mathrm{tot}}_n}\,|\,\F^\sharp_\Delta]\Ind_{G^\sharp_\Delta} \le C_{\delta}e^{-5A/3}\frac{ e^{\mu a}}{a}. 
\]
\end{lemma}
\begin{proof}
  The proof is similar to the proof of Lemma~\ref{lem:B_n}: Substituting Lemma~\ref{lem:Nt_prob_less_N} by Lemma~\ref{lem:Nt_prob_less_N_moving_wall}, one first shows that for every $n_0-1 \le k\le n-2$, one has
  \begin{equation}
    \label{eq:107}
    \E^\sharp_{\mathrm{nbab}}[B_{k,n}\Ind_{G'^{\mathrm{tot}}_n}\,|\,\F_{\mathscr B_k}]\Ind_{G^\sharp_\Delta} \le C_\delta \eta e^{-A} \E^\sharp_{\mathrm{nbab}}[B^{\mathrm{tot}}_k\Ind_{G'^{\mathrm{tot}}_k}\,|\,\F^\sharp_\Delta]\Ind_{G^\sharp_\Delta},
  \end{equation}
which, by a recurrence similar to \eqref{eq:840}, yields the lemma.
\end{proof}
\begin{lemma}
  \label{lem:G_n_prime}
 For large $A$ and $a$, $\P^\sharp_{\mathrm{nbab}}((G'^{\mathrm{tot}}_{n_1})^c\,|\,\F^\sharp_\Delta)\Ind_{G^\sharp_\Delta} \le 1/a.$
\end{lemma}
\begin{proof}
  Similarly to the proof of Lemma~\ref{lem:G_n_tot}, but using Lemmas~\ref{lem:Rt} and \ref{lem:R_f} instead of Proposition~\ref{prop:T}, we have for large $A$ and $a$,
\[
\P^\sharp_{\mathrm{nbab}} (G'^{\mathrm{tot}}_{n-1}\backslash G'^{\mathrm{tot}}_n\,|\,\F^\sharp_\Delta)\Ind_{G^\sharp_\Delta} \le C_\delta \E^\sharp_{\mathrm{nbab}}[B^{\mathrm{tot}}_{n-2}\Ind_{G'^{\mathrm{tot}}_{n-1}}\,|\,\F^\sharp_\Delta]\Ind_{G^\sharp_\Delta} e^{A-\mu a}/a \le C_\delta a^{-2},
\]
where the last inequality follows from Lemma~\ref{lem:B_n_moving_wall} and \eqref{eq:eta}. The lemma now follows by induction over $n$.
\end{proof}

\begin{lemma}
\label{lem:G_sharp}
Suppose (HB$^\sharp_0$). Then $\P(G^\sharp_1) \ge 1-\ep^{4/3}$  for large $A$ and $a$.
\end{lemma}
\begin{proof}
By Lemma~\ref{lem:Nt_prob_less_N}, Corollary~\ref{cor:Nt_prob_less_N} and the union bound we have
\begin{equation}
  \label{eq:113}
\P(\exists n \le n_0:N(t_n) < N,\ G_{Z,n_0},\ G_{\mathscr U}) \le C_{\delta} e^{-A}.
\end{equation}
Furthermore, by Lemma~\ref{lem:Nt_prob_less_N_moving_wall}, we have
\begin{equation}
  \label{eq:9}
\P(\exists n_0 < n < n_1:N(t_n) < N\tor N(\Theta_1) < N,\ G_\Delta) \le C_{\delta} \eta/\ep \le C_\delta e^{-A},
\end{equation}
by \eqref{eq:eta} and \eqref{eq:ep_lower}.
Now, by Proposition~\ref{prop:quantities} and Lemmas~\ref{lem:BBM_0} and \ref{lem:Gsharp_Delta_nbab},
\begin{equation}
  \label{eq:4}
  \begin{split}
 \E^\sharp_{\mathrm{nbab}}[Z^{\mathrm{tot}}_{\Theta_1}\Ind_{G'^{\mathrm{tot}}_{n_1}}\,|\,\F^\sharp_\Delta]\Ind_{G^\sharp_\Delta} &\le C\Big(Z^{\mathrm{nw}}_{t_{n_0}}+ C_\delta a e^{\mu a}\sum_{n=n_0-1}^{n_1-1}  \E^\sharp_{\mathrm{nbab}}[B^{\mathrm{tot}}_n\Ind_{G'^{\mathrm{tot}}_n}\,|\,\F^\sharp_\Delta]\Big)\Ind_{G^\sharp_\Delta}\\
&\le Ce^{A/2} + C_\delta  \le Ce^{A/2},
  \end{split}
\end{equation}
by Lemma~\ref{lem:B_n_moving_wall}. Similarly, we get
\begin{equation}
  \label{eq:5}
  \E^\sharp_{\mathrm{nbab}}[Y^{\mathrm{tot}}_{\Theta_1}\Ind_{G'^{\mathrm{tot}}_{n_1}}\,|\,\F^\sharp_\Delta]\Ind_{G^\sharp_\Delta}  \le Ce^{A/2}/a.
\end{equation}
Moreover, we have by Lemma~\ref{lem:B_n_moving_wall},
\begin{equation*}
\E^\sharp_{\mathrm{nbab}}[(B^{\mathrm{tot}}_{n_1-1}+B^{\mathrm{tot}}_{n_1}) \Ind_{G'^{\mathrm{tot}}_{n_1}}\,|\,\F^\sharp_\Delta]\Ind_{G^\sharp_\Delta}  \le C_\delta e^{\mu a -5A/3}/a.
\end{equation*}
Setting $X = \P\Big(B^{\mathrm{tot}}_{[\Theta_1,\Theta_1+t_1]} \le e^{-A} e^{\mu a}/a\,\Big|\,\F_{\Theta_1}\Big)$, we then have $\E[(1-X) \Ind_{G^\sharp_{\mathrm{nbab}}\cap G'^{\mathrm{tot}}_{n_1} \cap G^\sharp_\Delta}] \le C_\delta e^{-2A/3}$ by Markov's inequality.
Applying the Markov inequality once more to $X$ as in the proof of Lemma~\ref{lem:HBflat} yields 
\begin{equation}
  \label{eq:8}
\P(X \le 1-\ep^2,\,G^\sharp_{\mathrm{nbab}}\cap G'^{\mathrm{tot}}_{n_1} \cap G^\sharp_\Delta) \le C_\delta\ep^{-2}e^{-2A/3} \le C_\delta \ep^2,
\end{equation}
by \eqref{eq:ep_lower}.
The lemma now follows from \eqref{eq:113}, \eqref{eq:9} and \eqref{eq:8}, Lemmas~\ref{lem:Gsharp_Delta_nbab} and \ref{lem:G_n_prime} and Markov's inequality applied to \eqref{eq:4} and \eqref{eq:5}.
\end{proof}

\subsection{Proofs of the main results}

\begin{proof}[Proof of Lemma~\ref{lem:HBsharp}]
The property about $Z_0$ and $Y_0$ follows as in the proof of Lemma~\ref{lem:Bflat_start} from \eqref{eq:49}. Moreover, again as in the proof of Lemma~\ref{lem:Bflat_start}, we have by \eqref{eq:49}, \eqref{eq:69} and \eqref{eq:84}, for every $i\le \lceil e^\delta N\rceil$,
\begin{equation*}
\P(N^{\neg i}_{t_1} < N,\,R_{t_1}=0) \le C_\delta e^{-A}/a\quad\tand\quad\P(R_{t_1}> 0) \le C_\delta e^A/a,
\end{equation*}
where we denote the initial particles by $1,\dots,\lceil e^\delta N\rceil$.
By the independence of the particles, we now have,
\begin{equation*}
  \begin{split}
  \E[B_0\Ind_{(R_{t_1}=0)}] &\le \sum_{i=1}^{\lceil e^\delta N\rceil} \P(N^{\neg i}_{t_1} < N,\,R_{t_1}=0) \E^{X_i(0)}[L_0\Ind_{(R_{t_1}=0)}]\\
& \le C_\delta (e^{-A}/a) e^{\mu a} \E[Y_0] \le C_\delta e^{\mu a}/a^2,
  \end{split}
\end{equation*}
by \eqref{eq:49}. Setting $X = \P\Big(B_{0} \le e^{-A} e^{\mu a}/a\,\Big|\,\F_0\Big)$, we then have $\E[1-X] \le C_\delta e^{A}/a$ by Markov's inequality. Applying Markov's inequality once more to $X$ as in the proof of Lemma~\ref{lem:HBflat} yields the lemma.
\end{proof}

\begin{proof}[Proof of Proposition~\ref{prop:Bsharp}]
Exactly the same reasoning as in the proof of Proposition~\ref{prop:Bflat}, but using Lemma~\ref{lem:G_sharp} instead of Lemma~\ref{lem:G_Bflat}.
\end{proof}

\begin{proof}[Proof of Proposition~\ref{prop:Bsharp_med}]
Define $\rho = x_{\alpha e^{-2\delta}}$, such that $(1+\mu\rho)e^{-\mu\rho} \le \alpha e^{-2\delta}(1+o(1))$. As in the proof of Proposition~\ref{prop:Bflat_med} (see \eqref{eq:Nalpha}), it is enough to show that
\begin{equation}
\label{eq:97}
\sup_{\nu\in G^\sharp_0}\sup_{t\ge e^Aa^2} \P^\nu(N^{\mathrm{tot}}_t(\rho) > \alpha N,\,T>t) \le C_\delta \ep^{5/4},
\end{equation}
for large $A$ and $a$. Let $\nu\in G^\sharp_0$ and $t\ge e^Aa^2$. If $t> \sqrt\ep a^3$, then the above probability is bounded by $C\ep^2$ by Proposition~\ref{prop:T}. Suppose therefore that $e^Aa^2 \le t\le \sqrt\ep a^3$ and write $\Phat = \P(\cdot\,|\,T>t)$. By definition, we have $(1+ x_\alpha)e^{- x_\alpha} = \alpha$. Let $n$ be the largest integer such that $t_{n+1} < t$; note that $n\ge1$ for large $A$. 
By Lemmas~\ref{lem:hat_quantities} and \ref{lem:Zt_variance} together with Chebychev's and Markov's inequalities, \eqref{eq:ep_lower} and \eqref{eq:ep_upper}, we have for large $A$ and $a$,
\begin{equation}
  \label{eq:16}
\Phat(|e^{-A}Z_{t_n}-1|\le \ep^{1/8},\,Y_{t_n}\le e^{4A/3}/a) \ge 1-C_\delta  \ep^{5/4}.
\end{equation}
And as in the proof of Lemma~\ref{lem:G_ex}, the same holds for $Z^{\mathrm{gr}}_{t_n}$ and $Y^{\mathrm{gr}}_{t_n}$ as well, which yields
\begin{equation}
  \label{eq:89}
\Phat(|e^{-A}Z^{\mathrm{wg}}_{t_n}-1|\le \ep^{1/8},\,Y^{\mathrm{wg}}_{t_n}\le e^{4A/3}/a) \ge 1-C_\delta \ep^{5/4},
\end{equation}
Lemmas~\ref{lem:Rt} and \ref{lem:BBM_0} and Corollary~\ref{cor:hat_many_to_one} then show that $\P(R^{\mathrm{tot}}_{[t_n,t]} > 0) \le o(1)$. 

Let $N'^{\mathrm{wg}}_t(r)$ be the white and grey particles to the right of $r$ at time $t$ which descend from the white and grey particles at time $t_n$ and which have not hit $a$ between $t_n$ and $t$. By Proposition~\ref{prop:N_expec} and Corollary~\ref{cor:hat_many_to_one}, we have for large $A$ and $a$,
\begin{align}
  \Ehat[N'^{\mathrm{wg}}_t(\rho)\,|\,\F_{t_n}] &\le e^{-\delta/2} \alpha N e^{-A} Z^{\mathrm{wg}}_{t_n},\\
  \Varhat(N'^{\mathrm{wg}}_t(\rho)\,|\,\F_{t_n}) &\le C_{\delta,\alpha} \alpha^2N^2 e^{-2A} (a^{-1}Z^{\mathrm{wg}}_{t_n}+Y^{\mathrm{wg}}_{t_n}).
\end{align}
Chebychev's inequality and \eqref{eq:89} then give for large $A$ and $a$,
\begin{equation}
  \label{eq:94}
  \P(N'^{\mathrm{wg}}_t(\rho) > e^{-\delta/4} \alpha N) \le C\ep^{5/4}.
\end{equation}
Furthermore, denote by $N'^{\mathrm{blue}}_t(r)$ the number of blue particles to the right of $r$ at time $t$ which have turned blue after $t_{n}+a^2$ and which have not hit $a$ between $t_n$ and $t$. Then by Lemma~\ref{lem:blue_particles} and \eqref{eq:left_border}, we have
\begin{equation}
\label{eq:73}
  \begin{split}
 & \Ehat[N'^{\mathrm{blue}}_t(\rho)\,|\,\F_{t_n}]\\
 \le\ & C_\delta e^{-A}(K/a + e^{-A}Y^{\mathrm{wg}}_{t_n})\sum_{(u,s)\in\mathscr N^{\mathrm{wg}}_{t_n}}\int_{t_{n}+a^2}^t \Ehat[L_{[\tau,\tau+\dd \tau]}(u,s)]\E^0[N_{{t-\tau}}({\rho})]\\
\le\ & C_\delta e^{-A}(K/a + e^{-A}Y^{\mathrm{wg}}_{t_n})\sum_{(u,s)\in\mathscr N^{\mathrm{wg}}_{t_n}}e^{\mu X_u(s)}\int_{a^2}^{d} I^a(a-X_u(s),\dd \tau)\E^0[N_{t-t_n-\tau}(\rho)],
  \end{split}
\end{equation}
where we set $d = t-t_n$. Note that $(K+3)a^2\le d \le 2(K+3)a^2$. By Lemma~\ref{lem:many_to_one_simple} and Girsanov's theorem, we have now for every $\tau\ge 0$,
\begin{equation}
\label{eq:53}
  \E^0[N_{\tau}({\rho})] = e^{\beta_0 m\tau}W^0_{-\mu}(X_\tau > \rho) = e^{\frac{\pi^2}{2a^2}\tau}W^0[e^{-\mu \rho}\Ind_{(X_\tau > \rho)}] = e^{\frac{\pi^2}{2a^2}\tau}\int_\rho^\infty e^{-\mu z}g_{d-\tau}(z)\,\dd z,
\end{equation}
where $g_\tau(x)=(2\pi\tau)^{-1/2}e^{-x^2/(2\tau)}$ is the Gaussian density with variance $\tau$. If $\tau\ge a^2$, then $\sup_z g_\tau(z) \le C/a$, such that for every $x\in[0,a]$ and $z\ge\rho$,
\begin{equation}
\label{eq:45}
  \int_{a^2}^{d-a^2}I^a(a-x,\dd\tau) g_{d-\tau}(z) \le Ca^{-1} I^a(a-x,[a^2,d-a^2]) \le CKa^{-1}\sin(\pi x/a),
\end{equation}
by Lemma~\ref{lem:I_estimates}. Moreover, by Lemma~\ref{lem:I_estimates}, we have $I^a(a-x,\dd \tau) \le Ca^{-2} \sin(\pi x/a)\,\dd\tau$ for every $\tau\ge a^2$, such that
\begin{equation}
  \label{eq:72}
  \int_{d-a^2}^d I^a(a-x,\dd\tau)g_{d-\tau}(z) \le Ca^{-2}\sin(\pi x/a)\int_0^{a^2} \tau^{-1/2}\,\dd\tau \le Ca^{-1}\sin(\pi x/a).
\end{equation}
Equations \eqref{eq:73}, \eqref{eq:53}, \eqref{eq:45} and \eqref{eq:72} now yield
\begin{equation}
  \label{eq:74}
  \begin{split}
 \Ehat[N'^{\mathrm{blue}}_t(\rho)\,|\,\F_{t_n}] &\le C e^{-\mu \rho}e^{\mu a-A}a^{-2}(Ka^{-1} + e^{-A}Y^{\mathrm{wg}}_{t_n}) Z^{\mathrm{wg}}_{t_n}\\
& \le C_\delta \alpha e^{-2A}N(1+ae^{-A}Y^{\mathrm{wg}}_{t_n}) Z^{\mathrm{wg}}_{t_n}.
  \end{split}
\end{equation}
Furthermore, if $N^{\mathrm{rest}}_t(\rho)$ denotes the  particles to the right of $r$ at time $t$ descending from those turning blue between $t_{n-1}$ and $t_n+a^2$, then by \eqref{eq:53} and the supremum bound on $g_\tau(z)$,
\begin{equation}
  \label{eq:91}
  \Ehat[N^{\mathrm{rest}}_t(\rho)G_n^{\mathrm{tot}}] \le C_\delta a^{-1}e^{-\mu \rho} \E[(B_n^{\mathrm{tot}} +B_{n-1}^{\mathrm{tot}})\Ind_{G_n^{\mathrm{tot}}}] \le C_\delta \alpha e^{-A}N,
\end{equation}
by Lemma~\ref{lem:B_n}. Markov's inequality applied to \eqref{eq:91} and \eqref{eq:74} together with \eqref{eq:89}, Lemmas~\ref{lem:G_n_tot} and \ref{lem:G_Z} give
\begin{equation}
  \label{eq:125}
  \P(N_t^{\mathrm{rest}}+N_t'^{\mathrm{blue}} > \ep \alpha N) \le C_\delta \ep^2,
\end{equation}
by \eqref{eq:ep_lower}. The statement now follows from \eqref{eq:94}, \eqref{eq:125} and the above-mentioned bound on $\P(R^{\mathrm{tot}}_{[t_n,t]} > 0)$.
\end{proof}

\begin{proof}[Proof of Lemma~\ref{lem:Csharp}]
By definition, the B$^\sharp$-BBM and C$^\sharp$-BBM coincide until $\Theta_n$ on the set $G^\sharp_n$. By Proposition~\ref{prop:Bsharp}, we have $\P(G^\sharp_n) \ge 1-n\ep^{1+\numcst}$ for some numerical constant $\numcst > 0$ and by the coupling of  $(\Theta_n)_{n\ge0}$ with a Poisson process of intensity $p_B e^A\pi \sim \ep^{-1}$, we have $\P(\Theta_{\ep^{-1-\numcst/2}}> \ep^{-\numcst/4})\to 1$, as $A$ and $a$ go to infinity. Propositions~\ref{prop:Bsharp} and \ref{prop:Bsharp_med} then readily transfer to the C$^\sharp$-BBM.
\end{proof}

\section{The \texorpdfstring{$N$}{N}-BBM: proof of Theorem \ref{th:1}}
\label{sec:NBBM}

We will first establish a monotone coupling between the $N$-BBM and a class of slightly more general BBM with selection which includes the B$^\flat$-BBM from Section~\ref{sec:Bflat} and the C$^\sharp$-BBM from Lemma~\ref{lem:Csharp}. In a second part, Theorem~\ref{th:1} is proven.
%
%

\subsection{A monotone coupling between \texorpdfstring{$N$}{N}-BBM and more general particle systems}
\label{sec:coupling}

A \emph{selection mechanism} for branching Brownian motion is by definition a stopping line $\mathscr L$, which has the interpretation that if $(u,s)\in\mathscr L$, we think of $u$ being \emph{killed} at the time $t$. The set of particles in the system at time $t$ then consists of all the particles $u\in\mathscr N(t)$, which do not have an ancestor which has been killed at a time $s\le t$, i.e.\ all the particles $u\in\mathscr N(t)$ with $\mathscr L \precneq (u,t)$.

Now suppose we have two systems of BBM with selection, the $N^+$-BBM and the
$N^-$-BBM, whose selection mechanisms satisfy the following rules.
\begin{enumerate}[nolistsep]
\item Only left-most particles are killed.
\item $N^+$-BBM: Whenever a particle gets killed, there are at least $N$
  particles to its right (but not necessarily all the particles which have $N$ particles
  to their right get killed). $N^-$-BBM: Whenever at least $N$ particles are to
  the right of a particle, it gets killed (but possibly more particles get
  killed).
\end{enumerate}

Let $\nu^+_t$, $\nu^-_t$ and $\nu^N_t$ be the counting measures of the
particles at time t in $N^+$-BBM, $N^-$-BBM and $N$-BBM, respectively. On the space of (finite) counting measures on $\R$ we denote by $\preceq$ the usual stochastic ordering: For two counting measures $\nu_1$ and $\nu_2$, we write $\nu_1\preceq\nu_2$ if and only if $\nu_1([x,\infty)) \le \nu_2([x,\infty))$ for every $x\in\R$. If $x_1,\ldots,x_n$ and $y_1,\ldots,y_m$ denote the atoms of $\nu_1$ and $\nu_2$ respectively, then this is equivalent to the existence of an injective map\footnote{We use the notation $[n] = \{1,\ldots,n\}$.}  $\phi:[n]\to[m]$ with $x_i \le y_{\phi(i)}$ for all $i\in[n]$. Furthermore, for two families of counting measures $(\nu_1(t))_{t\ge0}$ and $(\nu_2(t))_{t\ge0}$, we write $(\nu_1(t))_{t\ge0}\preceq (\nu_2(t))_{t\ge0}$ if $\nu_1(t)\preceq\nu_2(t)$ for every $t\ge 0$. If $(\nu_1(t))_{t\ge0}$ and $(\nu_2(t))_{t\ge0}$ are random, then we write $(\nu_1(t))_{t\ge0}\stackrel{\mathrm{st}}{\preceq} (\nu_2(t))_{t\ge0}$ if there exists a coupling between the two (i.e.\ a realisation of both 
on the same probability space), such that $(\nu_1(t))_{t\ge0}\preceq (\nu_2(t))_{t\ge0}$.

\begin{lemma}
\label{lem:coupling}
Suppose that $\nu^-_0 \stackrel{\mathrm{st}}{\preceq} \nu^N_0 \stackrel{\mathrm{st}}{\preceq} \nu^+_0$. Then $(\nu^-_t)_{t\ge0} \stackrel{\mathrm{st}}{\preceq} (\nu^N_t)_{t\ge0} \stackrel{\mathrm{st}}{\preceq} (\nu^+_t)_{t\ge0}$.
\end{lemma}
\begin{proof}
 We only prove the second inequality $\nu^N_t \stackrel{\mathrm{st}}{\preceq} \nu^+_t$, the proof of the first one is similar. By a coupling argument and conditioning on $\F_0$ it is enough to show it for deterministic $\nu^N_0$ and $\nu^+_0$. Let $n^+ = \nu^+_0(\R)$ and let $\Pi = (\Pi^{(1)},\Pi^{(2)},\ldots,\Pi^{(n^+)})$ be a forest of independent BBM trees with the atoms of $\nu^+_0$ as initial positions. We denote by $\mathscr N^\Pi(t)$ the set of individuals alive\footnote{The term ``alive'' has the same meaning here as in Section~\ref{sec:preliminaries_definition}.} at time $t$ and by $X^\Pi_u(t)$ the position of an individual $u\in\mathscr N^\Pi(t)$. Denote by $\mathscr N^+(t)\subset \mathscr N^\Pi(t)$ the subset of individuals which form the $N^+$-BBM (i.e.\ those which have not been killed by the selection mechanism of the $N^+$-BBM). We set $\nu^+_t = \sum_{u\in\mathscr N^+(t)} \delta_{X^\Pi_u(t)}$.

From the forest $\Pi$ we will construct a family of forests $\left(\Xi_T = (\Xi^{(1)}_T,\ldots,\Xi^{(N)}_T)\right)_{T\ge0}$ (not necessarily comprised of independent BBM trees), such that
\begin{itemize}[nolistsep]
 \item if $T_1\le T_2$, then the forests $\Xi_{T_1}$ and $\Xi_{T_2}$ agree on
   the time interval $[0,T_1]$,
\item the initial positions in the forest $\Xi_0$ are the atoms of $\nu^N_0$,
\item for every $T\ge 0$, the $N$-BBM is embedded in $\Xi_T$ up to the time $T$, i.e.\ for $0\le t\le T$, if $\mathscr N^\Xi(t)$ denotes the set of individuals\footnote{Note that this does not depend on $T$.} from $\Xi_T$ alive at time $t$ and $X^\Xi_u(t)$ the position of the individual $u\in\mathscr N^\Xi(t)$, then there is a subset $\mathscr N^N(t)\subset \mathscr N^\Xi(t)$ such that $(\nu^N_t)_{0\le t\le T} = \big(\sum_{u\in\mathscr N^N(t)} \delta_{X^\Xi_u(t)}\big)_{0\le t\le T}$ is equal in law to the empirical measure of $N$-BBM,
\item for every $t\ge0$, there exists a (random) injective map $\phi_t:\mathscr N^N(t) \to \mathscr N^+(t)$, such that $X^\Xi_u(t) \le X^\Pi_{\phi_t(u)}(t)$ for every $u\in\mathscr N^N(t)$.
\end{itemize}
We will say that the individuals $u$ and $\phi_t(u)$ are \emph{connected}. If at a time $t$ an individual $v\in \mathscr N^{+}(t)$ is not connected to another individual (i.e.\ $v\notin \phi_t(\mathscr N^N(t))$), we say that $v$ is \emph{free}.

The construction of the coupling goes as follows: Since $\nu^N_0 \preceq \nu^+_0$, we can construct $\Xi_0$, $\mathscr N^N(0)$ and $\phi_0$, such that for every $u\in \mathscr N^N(0)$ we have $X^\Xi_u(0) \le X^\Pi_{\phi_0(u)}(0)$ and the subtrees $\Xi^{(u,0)}_0$ and $\Pi^{(\phi_0(u),0)}$ are the same up to translation. We now define a sequence of random times $(t_n)_{n\ge0}$ recursively by $t_0 = 0$ and for each $n$, we define $t_{n+1}$ to be the first time after $t_n$ at which either
\begin{enumerate}[nolistsep]
 \item a particle of the $N$-BBM branches, or
 \item the left-most particle of the $N^+$-BBM dies without a particle of the $N$-BBM branching.
\end{enumerate}
We then set $\Xi_t = \Xi_0$ and $\phi_t = \phi_0$ for all $t\in[0,t_1)$. Now, let $n\in\N$ and suppose that
\begin{itemize}[nolistsep]
 \item[a)] $\Xi_t$ and $\phi_t$ have been defined for all $t<t_{n+1}$ and are equal to $\Xi_{t_n}$ and $\phi_{t_n}$, respectively,
 \item[b)] for each $u\in \mathscr N^N(t_n)$, the subtrees $\Xi_{t_n}^{(u,t_n)}$ and $\Pi_{t_n}^{(\phi_{t_n}(u),t_n)}$ are the same, up to translation.
\end{itemize}
Note that this is the case for $n=0$. We now distinguish between the two cases above, starting with the second:

\emph{Case 2:} The left-most particle $w'$ of the $N^+$-BBM gets killed without a particle of the $N$-BBM branching. If $w'$ is free, nothing has to be done. Suppose therefore that $w'$ is connected to a particle $w$ of the $N$-BBM. Then, since there are at most $N-1$ remaining particles in the $N$-BBM and there are at least $N$ particles to the right of $w'$ in the $N^+$-BBM (otherwise it would not have been killed), at least one of those particles is free. Denote this particle by $v'$. We then ``rewire'' the particle $w$ to $v'$ by setting $\phi_{t_{n+1}}(w) := v'$ and define $\Xi_{t_{n+1}}$ by replacing the subtree $\Xi_{t_n}^{(w,t_{n+1})}$ in $\Xi_{t_n}$ by $\Pi_{t_n}^{(v',t_{n+1})}$, properly translated. Note that we then have 
\[
X^\Pi_{\phi_{t_{n+1}}(w)}(t_{n+1}) = X^\Pi_{v'}(t_{n+1}-)\ge X^\Pi_{w'}(t_{n+1}-) \ge X^\Xi_w(t_{n+1}-),
\]
where the first inequality follows from the fact that $w'$ is the left-most individual in $N^+$-BBM at time $t_{n+1}$ and the second inequality holds by hypothesis. See Figure~\ref{fig:sec10} for a picture.

\begin{figure}[ht]
 \begin{center}
\executeiffilenewer{sec10.svg}{sec10.pdf}%
{inkscape -z -D --file=sec10.svg %
--export-pdf=sec10.pdf --export-latex --export-area-drawing}%
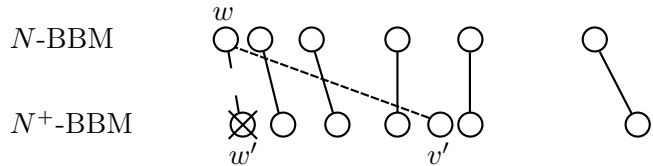%

 \caption{The connection between the particles $w$ of the $N$-BBM and $w'$ of the $N^+$-BBM breaks. By definition of the $N^+$-BBM, there exists a free particle $v'$ to the right of $w'$ and $w$ is rewired to that particle.}
 \label{fig:sec10}
 \end{center}
\end{figure}

If more than one particle of the $N^+$-BBM gets killed at the time $t_{n+1}$, we repeat the above for every particle, starting from the left-most.

\emph{Case 1:} A particle $u$ of the $N$-BBM branches at time $t_{n+1}$. By the hypothesis b), the particle $\phi_{t_n}(u)$ then branches as well into the same number of children. We then define $\phi_{t_n}(uk) = \phi_{t_n}(u)k$ for each $k\in[C_u]$ (recall that $C_u$ denotes the number of children of $u$), i.e.\ we connect each child of $u$ to the corresponding child of $\phi_{t_n}(u)$. Now first define $\phi_{t_{n+1}}'$ to be the restriction of $\phi_{t_n}$ to the surviving particles. Then continue as in Case~2, i.e.\ for each particle $w'$ of the $N^+$-BBM which gets killed and which is connected through  $\phi_{t_{n+1}}'$ to a particle $w$ of the $N$-BBM, rewire $w$ to a free particle $v'$. In the end, we get $\phi_{t_{n+1}}$.
 
In both cases, we then set $\Xi_t = \Xi_{t_{n+1}}$ and $\phi_t = \phi_{t_{n+1}}$ for each $t\in[t_{n+1},t_{n+2})$. Note that each time we are rewiring a particle, we rewire it to a particle whose subtree is independent of the others by the strong branching property, whence the particles from $\mathscr N^N(t)$ and $\mathscr N^+(t)$ still follow the law of $N$-BBM and $N^+$-BBM, respectively. Furthermore, we have for every $\omega\in\Omega$: $\nu^N_t(\omega) \preceq \nu^+_t(\omega)$ for every $t\ge 0$. This finishes the proof.
\end{proof}

\subsection{Proof of Theorem \ref{th:1}}
Let $(\nu_t^N)_{t\ge 0}$ be measure-valued $N$-BBM starting from the initial condition
\begin{itemize}
 \item[(H0)]at time 0, there are $N$ particles independently distributed according to the
density proportional to $\sin(\pi x/a_N) e^{- x}\Ind_{(x\in(0,a_N))}$, where $a_N = (\log N + 3 \log\log N)$.
\end{itemize}

Recall the definitions from the introduction. Let $\alpha \in (0,1)$. We wish to show that the finite-dimensional distributions of the process 
\begin{equation}
  \label{eq:126}
\left(M^N_\alpha\big(t \log^3 N\big)\right)_{t\ge0}
\end{equation}
converge weakly as $N\to\infty$ to those of the L\'evy process $(L_t)_{t\ge0}$ stated in Theorem~\ref{th:1} with $L_0 = x_\alpha$. We will do this by proving seperately a lower and an upper bound and show that in the limit these bounds coincide and equal the L\'evy process $(L_t)_{t\ge0}$.

\paragraph{Lower bound.} Fix $\alpha\in(0,1)$ and $\delta > 0$. We will let $N$ and in parallel $A$ and $a$ go to infinity (in the meaning of Section~\ref{sec:BBBM}) in such a way that $N = 2\pi  e^{A+\delta}a^{-3}e^{\mu a}$ and such that $A$ goes to infinity sufficiently slowly such that the results from Section~\ref{sec:Bflat} hold. We then have with $c=\log(2\pi )$,
\begin{equation}
  \label{eq:127}
  a = a_N - (A+\delta +c + o(1)),\quad \tand\quad \mu = \mu_N - \frac{\pi^2}{a^3}(A+\delta+c + o(1)).
\end{equation}
Let $(\nu^{\flat}_t)_{t\ge 0}$ be the measure-valued B$^\flat$-BBM starting from the initial configuration (HB$^\flat_{\perp}$), i.e.\ $\nu^\flat_0$ is obtained from $\lfloor e^{-\delta}N\rfloor$ particles distributed independently according to the density proportional to $\sin(\pi x/a)e^{-\mu x}\Ind_{(0,a)}(x)$. An easy calculation now shows that $\nu^\flat_0 \stackrel{\mathrm{st}}{\preceq} \nu^N_0$ for large $N$, by \eqref{eq:127}. Now, if $X_t^\flat$ denotes the barrier process of the B$^\flat$-BBM, then $(\nu^-_t)_{t\ge 0} = (\nu^\flat_t+\mu t + X^\flat_t)_{t\ge0}$  is by definition an instance of the $N^-$-BBM. Lemma~\ref{lem:coupling} now gives for large $N$, $(\nu^-)_{t\ge0} \stackrel{\mathrm{st}}{\preceq} (\nu^N_t)_{t\ge0}$, which by definition implies
\begin{equation}
  \label{eq:128}
 (\med^N_\alpha(\nu^N_t))_{t\ge0} \stackrel{\mathrm{st}}{\ge} (\med^N_\alpha(\nu^-_t))_{t\ge0}.
\end{equation}
Given $0\le t_1<\ldots<t_n$, we now define $t_i^N = a^{-3}t_i \log^3 N$, such that $t_i^N\to t_i$ for every $i$, as $N\to\infty$. We then have 
\begin{align*}
\Big(M^N_\alpha&\big(t_i \log^3 N\big)\Big)_{i=1,\ldots,n} \\
&= \big(\med^N_\alpha(\nu^N_{t_i \log^3 N})-\mu_Nt_i \log^3 N\big)_{i=1,\ldots,n} && \text{by definition}\\
&\stackrel{\mathrm{st}}{\ge} \big(\med^N_{\alpha}(\nu^-_{t_i \log^3 N})-\mu_Nt_i \log^3 N\big)_{i=1,\ldots,n} && \text{by \eqref{eq:128}}\\
&= \big(\med^N_{\alpha}(\nu^\flat_{a^3t_i^N}) + X_{a^3t_i^N}-\pi^2(A+\delta+c+o(1))t_i^N\big)_{i=1,\ldots,n} && \text{by \eqref{eq:127}}\\
&\stackrel{\mathrm{st}}{\ge} \big(x_{\alpha e^{2\delta}} + L_{t_i} - O(\delta) + c\big)_{i=1,\dots,n} &&\text{for large $N$},
\end{align*}
where the last inequality follows from Propositions~\ref{prop:Bflat} and \ref{prop:Bflat_med}, with $(L_t)_{t\ge 0}$ being the L\'evy process from the statement of Theorem~\ref{th:1} starting from $0$. Letting first $N\to\infty$, then $\delta\to 0$ yields the proof of the lower bound.

\paragraph{Upper bound.} The proof is analogous to the previous case, relying on Propositions~\ref{prop:Bsharp} and \ref{prop:Bsharp_med} instead of Propositions~\ref{prop:Bflat} and \ref{prop:Bflat_med}. There are only two differences to notice: First, the B$^\sharp$-BBM is not a realisation of the $N^+$-BBM. However, the C$^\sharp$-BBM (defined in Lemma~\ref{lem:Csharp}) is such a realisation and by that lemma, Propositions~\ref{prop:Bsharp} and \ref{prop:Bsharp_med} hold for the C$^\sharp$-BBM as well. Second, if $\nu^\sharp_0$ is distributed according to (HB$^\sharp_{\perp}$), we do not have $\nu^N_0\stackrel{\mathrm{st}}{\preceq}\nu^\sharp_0$. However, if $\widetilde \nu^N_0$ is obtained from $\nu^N_0$ by killing the particles in the interval $[a_N-A^2,a_N]$, then by \eqref{eq:127}, a quick calculation shows that $\widetilde\nu^N_0 = \nu^N_0$ with high probability and $\widetilde\nu^N_0 \stackrel{\mathrm{st}}{\preceq} \nu^\sharp_0$ as $N$ goes to infinity. This finishes the proof of the upper 
bound and of Theorem~\ref{th:1}.

\section*{Acknowledgments}
\addcontentsline{toc}{section}{Acknowledgements}
I wish to thank my PhD advisor Zhan Shi for his continuous support and encouragements and Julien Berestycki for numerous discussions and his interest in this work.



\chapter{A note on stable point processes occurring in branching Brownian motion}
\label{ch:stable}

\section{Introduction}

Brunet and Derrida \cite[p.\ 18]{Brunet2011} asked the following question, which arose during the study of the extremal particles in branching Brownian motion: Let $Z$ be a point process on $\R$, with the following property they called ``superposability'': $Z$ is equal in law to $T_\alpha Z+T_\beta Z'$, where $Z'$ is an independent copy of $Z$, $e^\alpha+e^\beta=1$ and $T_x$ is the translation by $x$. Is it true that $Z$ can be obtained from a Poisson process of intensity $e^{-x}\,\dd x$ on $\R$ by replacing each point by independent copies of an auxiliary point process $D$ (they called $D$ the ``decoration'')? More precisely, can $Z$ be written as
\begin{equation}
\label{eq_definition_dppp}
Z = \sum_{i=1}^\infty T_{\xi_i} D_i,
\end{equation}
where $(\xi_i)_{i\ge1}$ are the atoms of the above-mentioned Poisson process and $D_1,D_2,\ldots$ are independent copies of $D$ and independent of $\xi$?
This question was answered in the affirmative by the author \cite{Maillard2011}, and independently in the special case arising in branching Brownian motion by Arguin, Bovier, Kistler \cite{Arguin2011,Arguin2011a} and A\"{\i}d\'ekon, Berestycki, Brunet, Shi \cite{Aidekon2011c}. The representation \eqref{eq_definition_dppp} was also shown for the branching random walk by Madaule \cite{Madaule2011}, relying on the author's result. See also \cite{Kabluchko2011} for a related result concerning branching random walks.

Immediately after the article \cite{Maillard2011} was published on the arXiv, the author was informed by Ilya Molchanov that the superposability property had a classical interpretation in terms of \emph{stable point processes}, and the representation \eqref{eq_definition_dppp} was known in this field as the \emph{LePage series representation} of a stable point process.

The purpose of this note is two-fold: First, we want to outline how \eqref{eq_definition_dppp} can be obtained via the theory of stability in convex cones. Second, we give a succinct and complete proof of \eqref{eq_definition_dppp} and an extension to random measures for easy reference. 

\section{Stability in convex cones}
\label{sec:cones}

Let $Y$ be the image (in the sense of measures) of $Z$ by the map $x\mapsto e^x$ (this was suggested by Ilya Molchanov), such that $Y$ is a point process on $\R_+ = (0,\infty)$. By the superposability of $Z$, $Y$ has the following \emph{stability property}: $Y$ is equal in law to $aY + bY'$, where $Y'$ is an independent copy of $Y$, $a,b\ge 0$ with $a+b=1$ and $aY$ is the image of $Y$ by the map $x\mapsto ax$. If $Y$ is a simple point process, one can see the collection of points of the point process $Y$ as a random closed subset of $\R_+$, and the stability property is then also known as the \emph{union-stability} for random closed sets (see e.g.\ \cite[Ch. 4.1]{Molchanov2005}). 

Davydov, Molchanov and Zuyev \cite{Davydov2008} have introduced a very general framework for studying stable distributions in \emph{convex cones}, where a convex cone $\mathbb K$ is a topological space equipped with two continuous operations: addition (i.e.\ a commutative and associative binary operation $+$ with neutral element $\mathbf e$) and multiplication by positive real numbers. Furthermore, the two operations must distribute and  $\mathbb K\backslash\{\mathbf e\}$ be a complete separable metric space. For example, the space of compact subsets of $\R^d$ containing the origin is a convex cone, where the addition is the union of sets and the multiplication by $a>0$ is the image of the set by the map $x\mapsto ax$ (see Example~8.11 in \cite{Davydov2008}). Furthermore, it is a \emph{pointed cone}, in the sense that there exists a unique \emph{origin} $\mathbf 0$, such that for each compact set $K\subset \R^d$, $aK\to \mathbf 0$ as $a\to 0$ (the origin is of course $\mathbf 0 = \{0\}$). The existence of 
the 
origin permits to define a \emph{norm} by $\|K\| = d(\mathbf 0,K)$, with $d$ the Hausdorff distance in $\R^d$. An example of a convex cone without origin (Example~8.23 in \cite{Davydov2008}) is the space of measures on $\R^d\backslash\{0\}$ equipped with the usual addition of measures and multiplication by $a>0$ being defined as the image of the measure by the map $x\mapsto ax$, as above. 

A random variable $Z$ with values in $\mathbb K$ is now called \emph{$\alpha$-stable} if $a^{1/\alpha}Z + b^{1/\alpha}Z'$ is equal in law to $(a+b)^{1/\alpha}$ for every $a,b>0$, where $Z'$ is an independent copy of $Z$ and $\alpha\in\R$. With the theory of Laplace transforms and infinitely divisible distributions on semigroups (the main reference to this subject is \cite{Berg1984}), the authors of \cite{Davydov2008} then show that to every $\alpha$-stable random variable $Z$ there corresponds a \emph{L\'evy measure} $\Lambda$ which is \emph{homogeneous of order $\alpha$}, i.e.\ $\Lambda(aB) = a^\alpha\Lambda(B)$ for any Borel set $B$. Actually, $\Lambda$ is \emph{a priori} only defined on a certain dual of $\mathbb K$, and a considerable part of the work in \cite{Davydov2008} is to give conditions under which $\Lambda$ is supported by $\mathbb K$ itself. These conditions are satisfied for the first example given above, but not for the second, since they require in particular that the cone be pointed.
Moreover, and this is their most important result, under some conditions satisfied by the first example, $Z$ can be expressed as its \emph{LePage series}, i.e.\ the sum over the points of the Poisson process with intensity measure $\Lambda$. 

In order to get to the decomposition  \eqref{eq_definition_dppp}, one must then disintegrate the homogeneous L\'evy measure $\Lambda$ into a radial and an angular component, such that $\Lambda = cr^{-\alpha}\dd r\times \sigma$ for $c>0$ and some measure $\sigma$ on the unit sphere $\mathbb S=\{x\in\mathbb K: \|x\|=1\}$. This is also called the \emph{spectral decomposition} and $\sigma$ is called the \emph{spectral measure}. If $\sigma$ has mass 1, then the LePage series can be written as
\[
Z = \sum_i \xi_i X_i,
\]
where $\xi_1,\xi_2,\ldots$ are the atoms of a Poisson process of intensity $cr^{-\alpha}\dd r$ and $X_1,X_2,\ldots$ are iid according to $\sigma$, independent of the $\xi_i$. This is exactly the decomposition \eqref{eq_definition_dppp}.

\section{A succinct proof of the decomposition (\ref{eq_definition_dppp})}

As mentioned in the introduction, we will give here a short proof of the decomposition \eqref{eq_definition_dppp} and its extension to random measures. We hope that this proof will be more accessible to probabilists who are not familiar with the methods used in \cite{Davydov2008}. Furthermore, the results in \cite{Davydov2008} cannot be directly applied to give the extension of \eqref{eq_definition_dppp} to random measures, such that it may be of interest to give a rigorous proof in that setting.


\subsection{Definitions and notation}

We denote by $\Mcal$ the space of boundedly finite measures on $\R$, i.e.\ measures, which assign finite mass to every bounded Borel set in $\R$, and by $\Ncal$ the subspace of counting measures. It is known (see e.g.\ \cite{DV2003}, p.\ 403ff) that there exists a metric $d$ on $\Mcal$ which induces the vague topology and under which $(\Mcal,d)$ is complete and separable (but \emph{not} locally compact). We further set $\Mstar = \Mcal \backslash\{0\}$, which is an open subset and hence a complete separable metric space as well (\cite{BouGT2}, IX.6.1, Proposition~2), when endowed with the metric $d^*(\mu,\nu) = d(\mu,\nu) + |d(\mu,0)^{-1} - d(\nu,0)^{-1}|$, equivalent to $d$ on $\Mstar$. The spaces $\Ncal$ and $\Nstar = \Ncal\backslash\{0\}$ are closed subsets of $\Mcal$ and $\Mstar$, and therefore complete separable metric spaces as well (\cite{BouGT2}, IX.6.1, Proposition~1).

For every $x\in\R$, we define the translation operator $T_x:\Mcal\to\Mcal$, by $(T_x \mu)(A) = \mu(A-x)$ for every Borel set $A\subset\R$. Furthermore, we define the function $M:\Mstar\to\R$ by 
\[
M(\mu) = \inf\{x\in\R: \mu((x,\infty)) < \min(1,\mu(\R)/2)\}.
\] 
Note that if $\mu\in\Ncal$, then $M(\mu)$ is the position of the right-most atom of $\mu$, i.e.\ $M(\mu) = \sup \supp \mu$. It is easy to show that the maps $(x,\mu)\mapsto T_x\mu$ and $M$ are continuous, hence measurable.

A \emph{random measure} $Z$ on $\R$ is a random variable taking values in $\Mstar$. If $Z$ takes values in $\Ncal$, we also call $Z$ a \emph{point process}. For every non-negative measurable function $f:\R\to\R_+$, we define the \emph{cumulant}
\[K(f) = K_Z(f) = -\log \E\left[\exp(-\langle Z,f\rangle )\right]\in[0,\infty],\]
where $\langle \mu,f\rangle  = \int_\R f(x)\mu(\dd x)$. The cumulant uniquely characterises $Z$ (\cite{DV2003}, p.\ 161).

\begin{theorem}
\label{th_random_measure}
Let $Z$ be a random measure and let $K(f)$ be its cumulant. Then $Z$ is superposable if and only if for every measurable non-negative function $f:\R\to\R_+$,
\begin{equation}
\label{eq_K_superposable}
K(f) = c\int_\R e^{-x}f(x)\,\dd x + \int_\R e^{-x} \int_{\Mstar} [1- \exp(-\langle \mu,f \rangle)] T_x\Delta(\dd \mu)\,\dd x,
\end{equation}
for some constant $c\ge0$ and some measure $\Delta$ on $\Mstar$, such that 
\begin{equation}
\label{eq_condition_Delta}
\int_\R e^x\int_0^\infty (1-e^{-y}) \Delta(\mu(A+x)\in\dd y)\,\dd x < \infty,
\end{equation} for every bounded Borel set $A\subset\R$. Moreover, $\Delta$ can be chosen such that $\Delta(M(\mu)\ne 0) = 0$, and as such, it is unique unless $Z=0$ almost surely.
\end{theorem} 


\begin{corollary}
\label{cor_PP}
A point process $Z$ is superposable if and only if it has the representation \eqref{eq_definition_dppp} for some point process $D$ satisfying 
\begin{equation}
\label{eq_condition_D_PP}
\int^{\infty} \P(D(A+x) > 0) e^x\,\dd x < \infty.
\end{equation}
If $\P(Z\ne 0) > 0$, then there exists a unique pair $(m,D)$, such that $\P(M(D) = m) = 1$.
\end{corollary}

\subsection{Infinitely divisible random measures}
\label{sec_infinite_divisibility}
Our proof is based on the theory of infinitely divisible random measures as exposed in Kallenberg \cite{Kal1983}. A random measure $Z$ is said to be \emph{infinitely divisible} if for every $n\in\N$ there exist iid random measures $Z^{(1)},\ldots,Z^{(n)}$ such that $Z$ is equal in law to $Z^{(1)} + \cdots + Z^{(n)}$. It is said to be infinitely divisible as a point process, if $Z^{(1)}$ can be chosen to be a point process. Note that a (deterministic) counting measure is infinitely divisible as a random measure but not as a point process.

The main result about infinitely divisible random measures is the following (see \cite{Kal1983}, Theorem 6.1 or \cite{DV2008}, Proposition 10.2.IX, however, note the error in the theorem statement of the latter reference: $F_1$ may be infinite as it is defined).

\begin{fact}
\label{fact_idrm}
The random measure $Z$ is infinitely divisible if and only if
\[K(f) = \langle\lambda,f\rangle + \int_{\Mstar} [1- \exp(-\langle \mu,f \rangle)] \Lambda(\dd \mu),\]
where $\lambda\in\Mcal$ and $\Lambda$ is a measure on $\Mstar$ satisfying
\begin{equation}
\label{eq_Lambda}
\int_0^\infty (1-e^{-x}) \Lambda(\mu(A) \in \dd x) < \infty,
\end{equation}
for every bounded Borel set $A\subset\R$.
\end{fact}
The probabilistic interpretation (\cite{Kal1983}, Lemma 6.5) of this fact is that $Z$ is the superposition of the non-random measure $\lambda$ and of the atoms of a Poisson process on $\Mstar$ with intensity $\Lambda$, which is exactly the representation of $Z$ as the LePage series mentioned in Section~\ref{sec:cones}. It has the following analogous result in the case of point processes (\cite{DV2008}, Proposition 10.2.V), where the measure $\Lambda$ is also called the \emph{KLM measure}.

\begin{fact}
\label{fact_idpp}
A point process $Z$ is infinitely divisible as a point process if and only if $\lambda = 0$ and $\Lambda$ is concentrated on $\Nstar$, where $\lambda$ and $\Lambda$ are the measures from Fact \ref{fact_idrm}. Then, \eqref{eq_Lambda} is equivalent to $\Lambda(\mu(A) > 0) < \infty$ for every bounded Borel set $A\subset\R$.
\end{fact}

In particular, the L\'evy/KLM measure of a Poisson process on $\R$ with intensity measure $\nu(\dd x)$ is the image of $\nu$ by the map $x\mapsto \delta_x$.

\subsection{Proof of Theorem~\ref{th_random_measure}}
\label{sec_proof}
We can now prove Theorem~\ref{th_random_measure} and Corollary~\ref{cor_PP}. For the ``if'' part, we note that \eqref{eq_condition_Delta} implies \eqref{eq_Lambda} for the measure $\Lambda = \int e^{-x}T_x\Delta\,\dd x$, such that the process with cumulant given by \eqref{eq_K_superposable} exists. The superposability is readily verified. Further note that for point processes the condition \eqref{eq_condition_D_PP} is equivalent to \eqref{eq_condition_Delta}.

It remains to prove the ``only if'' parts. Let $Z$ be a superposable random measure. Then, for $\alpha,\beta\in\R$, such that $e^\alpha + e^\beta = 1$, we have
\[\begin{split}
K(f) = -\log \E[\exp(-\langle Z,f\rangle )] &= -\log \E[\exp(-\langle T_\alpha Z,f\rangle )] - \log \E[\exp(-\langle T_\beta Z,f\rangle )]\\
&= K(f(\cdot + \alpha)) + K(f(\cdot + \beta)).
\end{split}\]
Setting $\varphi(x) = K(f(\cdot + \log x))$ for $x \in \R_+$ (with $\varphi(0) = 0$) and replacing $f$ by $f(\cdot + \log x)$ in the above equation, we get $\varphi(x) = \varphi(xe^\alpha) + \varphi(xe^\beta)$ for all $x\in\R_+$, or $\varphi(x)+\varphi(y) = \varphi(x+y)$ for all $x,y\in\R_+$. This is the famous Cauchy functional equation and since $\varphi$ is by definition non-negative on $\R_+$, it is known and easy to show \cite{Darboux1880} that $\varphi(x) = \varphi(1) x$ for all $x\in\R_+$. As a consequence, we obtain the following corollary:
\begin{corollary}
\label{cor_exp} $K(f(\cdot+x)) = e^x K(f)$ for all $x\in\R$.
\end{corollary}

Furthermore, it is easy to show that superposability implies infinite divisibility. We then have the following lemma.

\begin{lemma}
\label{lem_density}
Let $\lambda,\Lambda$ be the measures corresponding to $Z$ by Fact~\ref{fact_idrm}. 
\begin{enumerate}[nolistsep]
\item There exists a constant $c\ge0$, such that $\lambda = ce^{-x}\,\dd x$.
\item For every $x\in\R$, we have $\dd T_x\Lambda = e^x\dd\Lambda$ in the sense of Radon-Nikodym derivatives.
\item $\mu(\R_+) = 0$ for $\Lambda$-almost every $\mu$.
\end{enumerate}
\end{lemma}
\begin{proof}
The measures $T_x\lambda$, $T_x\Lambda$ are the measures corresponding to the infinitely divisible random measure $T_xZ$ by Fact \ref{fact_idrm}. But by Corollary \ref{cor_exp}, the measures $e^x\lambda$ and $e^x \Lambda$ correspond to $T_xZ$, as well. Since these measures are unique, we have $T_x\lambda = e^x\lambda$ and $T_x\Lambda = e^x\Lambda$. The second statement follows immediately. For the first statement, note that $c_1 = \lambda([0,1)) < \infty$, since $[0,1)$ is a bounded set. It follows that \[\lambda([0,\infty)) = \sum_{n\ge0} \lambda([n,n+1)) = \sum_{n\ge0} c_1e^{-n} = \frac{c_1e}{e-1} =: c,\]
hence $\lambda([x,\infty)) = ce^{-x}$ for every $x\in\R$. The first statement of the lemma follows. For the third statement, let 
  $I_n = [n,n+1)$ and $I = [0,1)$. By \eqref{eq_Lambda}, we have \[\int_0^1\Lambda(\mu(I) > x) \,\dd x = \int_0^1 x\Lambda(\mu(I)\in \dd x)  <\infty.\] By monotonicity, the first integral is greater than or equal to $x \Lambda(\mu(I) > x)$ for every $x\in[0,1]$, hence $\Lambda(\mu(I) > x) \le C/x$ for some constant $0\le C < \infty$. By the second statement, it follows that
\[
\Lambda(\mu(I_n) > e^{-n/2}) = e^{-n} \Lambda(\mu(I) > e^{-n/2}) \le C e^{-n/2},
\] 
for every $n\in\N$. Hence, $\sum_{n\in\N} \Lambda\left(\mu(I_n) > e^{-n/2}\right) < \infty$. By the Borel-Cantelli lemma,
\[\Lambda\left(\limsup_{n\to\infty} \left\{\mu(I_n) > e^{-n/2}\right\}\right) = 0,\] 
which implies $\Lambda(\mu((0,\infty)) = \infty) = 0$.
\end{proof}

\begin{lemma}
\label{lem_decomp}
The measure $\Lambda$ admits the decomposition $\Lambda = \int e^{-x}T_x\Delta\,\dd x$, where $\Delta$ is a unique measure on $\mathcal M^*$ with $\Delta(M(\mu)\ne 0) = 0$ which satisfies \eqref{eq_condition_Delta}.
\end{lemma}

\begin{proof}
  We follow the proof of Proposition~4.2 in \cite{Rosinski1990}. Set $\Mstar_0 := \{\mu\in\Mstar:M(\mu)=0\}$, which is a closed subset of $\Mstar$, and therefore a complete seperable metric space (\cite{BouGT2}, IX.6.1, Proposition~1). By the continuity of $(x,\mu)\mapsto T_x\mu$, the map $\phi: \Mstar\to \Mstar_0\times \R$ defined by $\phi(\mu) = (T_{-M(\mu)},M(\mu))$ is a Borel isomorphism, i.e.\ it is bijective and $\phi$ and $\phi^{-1}$ are measurable. The translation operator $T_x$ then acts on $\Mstar_0\times \R$ by $T_x(\mu,m) = (\mu,m+x)$. If $A_n = \{\mu\in\Mstar_0:\mu([-2n,2n]) \le 1/n\}$, then $\Lambda(A_n\times[-n,n])<\infty$ for every $n\in\N$ by \eqref{eq_Lambda}. By the theorem on the existence of conditional probability distributions (see e.g.\ \cite{Kallenberg1997}, Theorems~5.3 and 5.4) there exists then a measure $\Delta_0$ on $\Mstar_0$ with $\Delta_0(A_n)<\infty$ for every $n\in\N$ and a measurable kernel $K(\mu,\dd m)$, with $K(\mu,[-n,n])<\infty$ for every $n\in\N$, such that 
  \begin{equation*}
\Lambda(\dd \mu,\dd m) = \int_{\Mstar_0} \Delta_0(\dd\mu)K(\mu,\dd m).
  \end{equation*}
Moreover, we can assume in the above construction that $K(\mu,[0,1]) = K(\mu',[0,1]) = 1$ for every $\mu,\mu'\in\Mstar_0$ and $n\in\N$, and with this normalization, $\Delta_0$ is unique. By Lemma~\ref{lem_density}, we now have $T_xK(\mu,\dd m) = e^x K(\mu,\dd m)$ for every $x\in\R$ and $\mu\in\Mstar_0$. As in the proof of the first statement of Lemma~\ref{lem_density}, we then conclude that $K(\mu,\dd m) = c(\mu)e^{-m}\,\dd m$ for some constant $c(\mu)\ge 0$, and by the above normalization, $c(\mu)\equiv c := e/(e-1)$. Setting $\Delta(\dd m) = c\Delta_0(\dd m)$ then gives
\begin{equation*}
  \Lambda(\dd \mu,\dd m) = \int_{\Mstar_0} \Delta(\dd\mu)e^{-m}\,\dd m,
\end{equation*}
which finishes the proof.
\end{proof}

The ``only if'' part of Theorem~\ref{th_random_measure} now follows from the previous lemmas. If $Z$ is a point process, then Fact~\ref{fact_idpp} implies that $\lambda = 0$ and that $\Lambda$ is concentrated on $\Nstar$, hence $\Delta$ as well. Equation \eqref{eq_condition_Delta} then implies that $\Delta(\mu(A)>0) < \infty$ for any bounded Borel set $A\subset\R$. In particular, this holds for $A=[-1,1]$. But since $\mu\in\Nstar_0$ implies $\mu([-1,1]) > 0$ and since $\Delta$ is concentrated on $\Nstar_0$, it has finite mass. If $\Delta$ has mass zero, then $Z=0$ almost surely. If $\Delta$ has positive mass, set $m=\log \Delta(\Nstar)$. The measure $\Delta' = e^{-m}T_m\Delta$ is then a probability measure and $\Lambda = \int e^{-x}T_x\Delta'\,\dd x$. Furthermore, $Z$ satisfies \eqref{eq_definition_dppp}, where $D$ follows the law $\Delta'$. Uniqueness of the pair $(m,D)$ follows from Lemma~\ref{lem_decomp}. This finishes the proof of Corollary~\ref{cor_PP}.

\subsection{Finiteness of the intensity}
\label{sec_dppp}
If $Z$ is a superposable point process and has finite intensity (i.e.\ $E[Z(A)] < \infty$ for every bounded Borel set $A\subset \R$), then it is easy to show that the intensity is proportional to $e^{-x}\,\dd x$. However, in the process which occurs in the extremal particles of branching Brownian motion or branching random walk, the intensity of the decoration grows with $|x|e^{|x|}$, as $x\to -\infty$ \cite[Section~4.3]{Brunet2011}. The following simple result shows that in these cases, $Z$ does not have finite intensity.
\begin{proposition}
\label{prop_dppp_exp}
Let $Z$ be defined as in \eqref{eq_definition_dppp}. Then $Z$ has finite intensity if and only if $\E[\langle D,e^x\rangle ] < \infty$.
\end{proposition}

\begin{proof}
By Tonelli's theorem,
\[E[Z(A)] = \E\left[\sum_{i\in\N} \E[T_{\xi_i}D(A)\,|\,\xi]\right]= \int_\R \E[D(A-y) e^{-y}]\, \dd y = \E\left[\int_\R D(A-y) e^{-y}\, \dd y\right],\]
for every bounded Borel set $A\subset \R$. Again by Tonelli's theorem we have
\[\int_\R D(A-y) e^{-y}\, \dd y = \int_\R \int_\R \Ind_{A-y}(x)e^{-y}\, \dd y\, D(\dd x) = \langle D,\int_\R \Ind_{A-y}(\cdot)e^{-y}\, \dd y \rangle.\]
For $x\in\R$, $x\in A-y$ implies $y \in [\min A - x,\max A - x]$. Since $e^{-y}$ is decreasing, we therefore have
\[|A|e^{-\max A}e^x \le \int_\R \Ind_{A-y}(x)e^{-y}\, \dd y \le |A|e^{-\min A}e^x,\] where $|A|$ denotes the Lebesgue measure of $A$. We conclude that $E[Z(A)] < \infty$ if and only if $\E[\langle D,e^x\rangle ] < \infty.$
\end{proof}




\bibliographystyle{pascal} 
\bibliography{library,these} 




\cleardoublepage
\chapterstar{Errata}

\paragraph{Chapter~\ref{ch:NBBM}, Section~\ref{sec:doob}.}
The last paragraph contains a wrong statement: Even if $R(x)\equiv R$, the function $h(x)$ is in general not a harmonic function for the stopped process $X^H$ under $P^x$. The truth is that (for general bounded functions $R(x)$) the process 
\begin{equation}
\label{eq:M}
 M_t = \frac{h(X_{t\wedge H})}{h(x)} \exp\left(-\int_0^{t\wedge H} R(X_s)(1-Q(X_s))\,\dd s\right)
\end{equation}
is a mean-one martingale under $P^x$ for every $x\in\mathscr E$. If we define the law $P_h^x$ by a martingale change of measure $\dd P_h^x|_{\F_t} = M_t\,\dd P^x|_{\F_t}$, then the very last statement in that section remains true, in that the law $\overline{P}_h^x$ is obtained from the law $P_h^x$ by killing the process at the time-dependent rate $R(x) Q(x)\Ind_{(t < H)}$.

In order to prove the above, it is helpful to enlarge the space $\mathscr E$: Let $\overline {\mathscr E}$ be a copy of $\mathscr E$ and denote the copy of an element $x\in\mathscr E$ by $\overline x$. We can view $\overline{P}^x$ as a the law of a process on  $\mathscr E\cup \overline {\mathscr E}$ with the generator
\[
 \overline L f(x) = L f(x) + R(x)(f(\overline x)-f(x)) \tand \overline Lf(\overline x)=0 \text{ for all } x\in\mathscr E.
\]
Here, $L$ is the generator of $X$ and $f:\mathscr E\cup \overline {\mathscr E}\to \R$ is any measurable function such that the restriction of $f$ to $\mathscr E$ is in the domain of $L$. We then define the function $\overline h:\mathscr E\cup \overline {\mathscr E}\to[0,1]$ by $\overline h(x) = h(x)$ and $\overline h(\overline x) = h(x)Q(x)$ for all $x\in\mathscr E$. The function $\overline h$ is then a harmonic function for the stopped process $X^H$ under $\overline{P}^x$, which implies that for every $x\in \mathscr E\backslash F$, we have
\begin{equation}
\label{eq:Lh}
 L h(x) = R(x)(1-Q(x))h(x).
\end{equation}
Moreover, the law $\overline{P}_h^x$ is the $\overline h$-transform of $\overline{P}^x$. In particular, if $\overline L_h$ denotes the generator of $X^H$ under $\overline{P}_h^x$, we have for every $x\in \mathscr E\backslash F$ and $f$ as above,
\[
  \overline L_h f(x) = h(x)^{-1}\overline L (\overline h f)(x) = h(x)^{-1} L (hf)(x) - R(x)(1-Q(x))f(x) + R(x)Q(x)(f(\overline x)-f(x)).
\]
Now, \eqref{eq:Lh} implies that the process $M_t$ defined in \eqref{eq:M} is a local martingale\footnote{See Lemma~3.1 in Z.\ Palmowski, T.\ Rolski (2002). A technique for exponential change of measure for Markov processes. \textit{Bernoulli}, 8(6), 767–785.} and even a martingale, because we assumed that $R(x)$ is bounded\footnote{See Proposition~3.2 in the same paper.}. This implies that the operator $L_h$ defined by 
\[
 L_hf = h(x)^{-1} L (hf)(x) - R(x)(1-Q(x))f(x)
\]
is the generator of the process of the process $X^H$ under the law $P_h^x$ defined above. This gives
\[
 \overline L_h f(x) = L_hf(x) + R(x)Q(x)(f(\overline x)-f(x)),
\]
which shows that the law $\overline{P}_h^x$ is indeed obtained from the law $P_h^x$ by killing the process at the rate $R(x) Q(x)\Ind_{(t < H)}$.

\paragraph{Chapter~\ref{ch:NBBM}, Section~\ref{sec:weakly_conditioned}.}
The previous error entails an error at the beginning of the section: Under the law $\Ptilde_f^x$, particles move according to the law obtained from Brownian motion with drift $-\mu_t$ (stopped at $0$ and $a$) through a change of measure by the martingale 
\[
 h(X_{t\wedge H},t\wedge H) \exp\left(-\int_0^{t\wedge H} \beta_0 (1-Q(X_s,s))\,\dd s\right).
\]
As a consequence, for Lemma~\ref{lem:many_to_one_conditioned} to remain true, the definition \eqref{eq:def_e} needs to be changed to read
\[
 e(x,t) = \beta_0 \sum_{k\ge 0} (k+1) (1-h(x,t)^{k-1}) q(k) \le \beta_0 (m_2+m)(1-h(x,t))
\]
and with this estimate all the following ones remain true.

\paragraph{Corrections of some typographical mistakes}
\begin{itemize}
\item p. 126, Figure 2.6: ...the probability that at time $t$ there are less than $N$ particles not descending from $u$ is bounded by $C_\delta e^{-A}(a^{-1} + e^{-A}Y^{(0),\neg u}_{t'})$.
\item p. 127, first paragraph: For $n\ge 0$ and $0\le k \le n-2$, we define $B_{k,n}$ to be the number of particles that turn blue at a time $t\in I_n$, have an ancestor that turned blue at a time $t'\in I_k$, have none that has hit 0 between $t_{k+2}$ and $t_n$, have not hit $a$ between $t'$ and $t_{k+2}$ and did not break out between $t_{k+2}$ and $t$.
\item p. 133, first paragraph of Section~\ref{sec:coupling}: A \emph{selection mechanism} for branching Brownian motion is by definition a stopping line $\mathscr L$, which has the interpretation that if $(u,t)\in\mathscr L$, we think of $u$ being \emph{killed} at the time $t$. The set of particles in the system at time $t$ then consists of all the particles $u\in\mathscr N(t)$ which do not have an ancestor which has been killed at a time $s\le t$, i.e.\ all the particles $u\in\mathscr N(t)$ with $\mathscr L \not\preceq (u,t)$.
\end{itemize}

\paragraph{Chapter~\ref{ch:stable}.} Some typographical mistakes have been corrected in the published version. Reference: [Pascal Maillard, A note on stable point processes occuring in branching Brownian motion, \textit{Electronic
Communications in Probability} \textbf{18}, no. 5, 2013]

\end{document}